%% file: parabolic-surfaces.tex
\documentclass[12pt,twoside,leqno,openany]{amsart}

\pdfoutput=1

\usepackage{amsbsy,amscd,amsfonts,amsmath,amssymb,amsthm,color,
fancybox,fancyhdr,footmisc,graphics,graphicx,ifthen,mathrsfs,
multicol,pdfpages,rotating,times,wasysym,pifont,blkarray}
\usepackage[dvipsnames,svgnames,x11names]{xcolor}

\usepackage[all]{xy}
\usepackage[utf8]{inputenc}
\usepackage[T1]{fontenc}
\usepackage{mathtools}
\newtagform{EngelLie}[\scriptstyle]{$}{$}
\makeatletter\newcommand{\leqnomode}{\tagsleft@true}
\newcommand{\reqnomode}{\tagsleft@false}\makeatother

\input macros.tex


\begin{document}

\setcounter{section}{0}

$\:$

\bigskip\bigskip\bigskip\bigskip\bigskip

\begin{center}

{\large\bf On Differential Invariants}
\label{parabolic-surfaces-differential-invariants}

\medskip

{\large\bf of Parabolic Surfaces}\footnote{1 This work was supported in part by the Polish National Science Centre (NCN) via the grant number 2018/29/B/ST1/0258.}



\bigskip\bigskip

Zhangchi {\sc Chen},
Joël~{\sc Merker},

\end{center}\bigskip

\begin{center}
\begin{minipage}[t]{12.5cm}
\parindent 0.53cm
\footnotesize
\noindent
{\sc Abstract}.
The algebra of differential invariants under ${\sf SA}_3(\R)$ of
generic parabolic surfaces $S^2 \subset \R^3$ with nonvanishing
Pocchiola $4$\textsuperscript{th} invariant $\Waux$ is shown to be
generated, through invariant differentiations, by only one other
invariant, $\Maux$, of order $5$, having $57$
differential monomials. The proof is based on
Fels-Olver's recurrence formulas, pulled back to the
parabolic jet bundles.

\hfill
{\scriptsize
[Message to the busy reader: 
Sections~{\ref{introduction-parabolic-surfaces}}
and~{\ref{presentation-results}}
explain and summarize the contents.]}
\end{minipage}
\end{center}

\Section{\bf Introduction}
\label{introduction-parabolic-surfaces}
\HEAD{{\ref{introduction-parabolic-surfaces}}.~{\sf 
Introduction}
}{
Zhangchi {\sc Chen}, Joël~{\sc Merker}}

In continuous and discrete mathematics, group actions are widespread.
They also arise in various fields of applied science, especially in
classical mechanics
{\cite{Auffray-Kolev-Olive-2017,
Desmorat-Kolev-Olive-2017,
Olive-2017}}.
Invariants may 
be used {\em e.g.} to take into account physical symmetries of a 
body, with the aim of reducing the data size of its analysis.

When available, the associated (algebraic or differential) group
invariants enable one to solve equivalence problems, to classify
geometric objects, to set up canonical forms for them, to know their
symmetries, and to find complete lists of homogeneous models.

This article deals more specifically with {\em differential}
invariants, group actions being in general {\em nonlinear}. It also
attempts to handle {\em extended} explicit expressions.

Our computational requirements are: be algorithmic, be explicit, 
{\em and} be {\sl synthetic}\,\,---\,\,what algorithms usually are not!
Two main goals are in our minds: compute collections of generating
sets of invariants, and understand {\em differential} relations among
such collections.

We will touch neither rewriting procedures in terms of the generating
invariants, nor algorithms for computing inside algebras of
invariants.  In any case, as soon as the number of independent
variables becomes $\geqslant 2$, it is well known that one encounters
a high complexity in symbolic expressions.  Although objects are
inserted in several theoretically satisfactory frameworks, everybody
is often left with frustratingly small achievements while playing on
any computer.

\smallskip

Most problems in the (infinitely wide) Lie-Cartan theory come up with
a given action of a certain $r$-dimensional Lie group $G$ acting on an
$m$-dimensional manifold $M$ with $m \geqslant 2$. Attacking these
questions often
involves studying the induced action of $G$ on submanifolds $S^p
\subset M^m$ of a prescribed dimension $p$, with $1 \leqslant p
\leqslant m-1$.  

While exploring deeper such problems, the occurrence of certain
submanifolds $S^p \subset M^m$ repeats itself as producing certain
(new) sub-submanifolds $S_2^{p_2} \subset S^p \subset M^m$,
sub-sub-submanifolds $S_3^{p_3} \subset S_2^{p_2} \subset S^p \subset
M^m$, and so on, most of the times after replacing $M^m$ with some
appropriate, cut-down and branched, jet {\em sub}space. Then new
(invariant) equations appear, new (invariant) bifurcations are
created. So here $S^p$ and $M^m$ should be thought of as being
sub-objects within some of the branches of a certain root problem
lying at departure.

Since the dimension $r$ of the group is often (much) higher than the
dimension $m$ of the ambient space, in order to somehow `{\sl
des-intricate}' the group action, one prolongs it to the jet bundles
$J_{p,q}^n$ of any order $n \geqslant 0$, where $q := m - p$ is the
codimension. Certain general (but complicated) formulas going back to
Lie ({\cite[Ch.~25]{Lie-Merker-2015}}) show how tangent directions
attached to $S^p \subset M^m$, and also higher order jets as well, 
transfer through diffeomorphisms from 
the base $M$ to the jet bundles $J_{p,q}^n$ 
of any order $n$.  Section~{\ref{graph-transformations-jet-spaces}}
presents these formulas, and thanks to them, the action of $G$ on $M$
lifts as a $G$-action on every jet space $J_{p,q}^n$.

The geometry of submanifolds under Lie transformation groups: their
equivalences, their symmetries, their normal forms; is entirely
governed by what is known as differential invariants.  They are best
visualized inside $J_{p,q}^n$.

A {\sl differential invariant} is a (perhaps locally defined)
real-valued function $\Iaux\colon\, J_{p,q}^n \longrightarrow \R$ that
is invariant under the prolonged group action. Any finite-dimensional
Lie group action admits an infinite number of functionally independent
differential invariants of progressively higher and higher orders.

A universal question is to find a minimal set of generating
differential invariants. Certainly, the minimal number of
generating differential invariants cannot be fixed {\em a priori}: it
strongly depends on the particularities of the group action.  Since
the 19\textsuperscript{th} Century, the question of finite generation
of differential invariants was addressed by several authors, also in
the more general context of (infinite-dimensional) Lie pseudo-groups.
The fundamental {\sl Basis Theorem} states that all the differential
invariants can be generated, from a finite number of low order
invariants, by repeated {\sl invariant differentiations}.

\smallskip

Serendipitously indeed ({\cite{Fels-Olver-1999}}), there always exist
$p = \dim\, S$ linearly independent {\sl invariant differential
operators} $\mathcal{D}_1, \dots, \mathcal{D}_p$ with the property
that each $\mathcal{D}_j$ maps every differential invariant $\Iaux$ to
a differential invariant $\mathcal{D}_j \big( \Iaux \big)$.  The great
value of such invariant derivations is that one can explicitly write
down their action on invariantized jet monomials, and compare the
outcome with higher order invariants. This comparison, which
incorporates appropriate correction terms, is captured by the
celebrated {\sl recurrence formulas}, set up in the widest context
by Fels-Olver in~{\cite{Fels-Olver-1999}}, {\em see} also
Section~{\ref{recurrence-formulae-differential-invariants}}.

Let us repeat that, notwithstanding their power, recent symbolic
implementations are often led to unsurmountable obstacles while
attempting to explicitly compute invariants, 
or even cross-sections to the $G$-action on jet spaces. Our
Section~{\ref{search-resolved-cross-section-parabolic}} 
illustrates this difficulty. But serendipitously again, from the
computational side, a minimal amount of data is necessary to set up
the key recurrence formulas of Fels-Olver. Remarkably, these
formulas can be explicitly determined without knowing the actual
formulas for either the differential invariants, or the invariant
differential operators, or even the moving frame itself (!).

Indeed, the recurrence formulas can be written with only the knowledge
of the infinitesimal generators of the action and the equations of the
cross-section. Therefore, understanding these recurrence formulas is
the `{\sl master key}', according to Olver, that `{\sl unlocks}' the
structure of the algebra of differential invariants, the determination
of generators, and the classification of syzygies. The only required
ingredients are the prolongation formulas for the infinitesimal
generators, or, equivalently, the Lie matrix, along with the
specification of the cross-section normalizations.

\smallskip

Nevertheless, our slogan will be: 

\medskip

\begin{center}
{\em Explicit expressions of invariants are
necessary in exploring classification branches}. 
\end{center}

\medskip

In differential invariant contexts where {\em no} branching is
tracked, explicit expressions of invariants are not crucial.
In~{\cite{Olver-2007}}, Olver showed that the algebra of differential
invariants of a suitably generic hyperbolic or elliptic surface $S
\subset \R^3$ under the equi-affine group action is generated
by a single differential invariant, the third order Pick invariant,
with its invariant derivatives. The proof was based on the
straight equivariant approach to the method of moving frames.  We
believe anyway that classification bifurcations
also exist for hyperbolic
and elliptic surfaces.

\begin{Question}
{\sl 
What about parabolic surfaces?}
\end{Question}

These are (local) surfaces $S^2 \subset \R^3$ of 
the graphed form $\big\{ u = F(x,y) \big\}$ 
whose hessian matrix
$\big( \begin{smallmatrix} F_{xx} & F_{xy} \\ F_{yx} & F_{yy} 
\end{smallmatrix} \big)$ is identically of rank $1$, not $2$
(elliptic or hyperbolic cases). The vanishing of the Hessian
determinant:
\[
0
\,\equiv\,
F_{xx}\,F_{yy}
-
F_{xy}^2
\eqno
{\scriptstyle{(\text{\rm at all points}\,\,(x,y))}},
\]
then creates a differential relation which must be differentiated
again and again to build up the relevant {\sl parabolic jet spaces}
$P\!J_{2,1}^n$, of any order $n \geqslant 2$.  These differential
relations also have strong influence on the recurrence formulas.
Knowing explicitly the Hessian is unavoidable
to study this branch: the category
of parabolic surfaces. Other more complicated
branching invariants will come up in our deeper explorations,
as we will summarize in the next 
Section~{\ref{presentation-results}}.

The present article therefore opens up a new natural context of {\sl
jet spaces with differential relations}, in which the (cut-down)
pulled back recurrence formulas have an entirely invariant
meaning. The principle of passing to (bifurcating) submanifolds will
be illustrated several times on examples.

\smallskip

It is well known ({\em
see}~Section~{\ref{graph-transformations-jet-spaces}}) that the
coefficients of the prolonged infinitesimal generators of any group
action are {\em polynomial functions} of the jet coordinates.  In
particular, if the action of $G$ is transitive on $M$, which is often
the case, and if one normalizes all the order zero coordinates, this
implies that the {\sl Maurer–Cartan invariants}, which appear in the
fundamental recurrence formulas, are also rational functions of the
collection of generating invariants.  As a consequence of the theory,
all the higher order differential invariants are {\em
rational functions} of the generating differential invariants.

Rationality holds true of all the 
structures studied in this article.
The same algebraicity features hold
for a large class of pseudo-group actions: the differential
invariant algebra is intrinsically rational, in the sense that all
recurrence formulas, commutation relations and syzygies, involve
rational functions of the basic differential invariants,
all of order $\geqslant 1$.

Beyond transitivity of the $G$-action, 
in Theorem~{\ref{Thm-translations-transvections-G-1}}, 
we provide a condition on $G$ insuring that all basic 
differential invariants are of order $\geqslant 2$.

\smallskip

Last but not least, following the approach of Fels-Olver ({\em cf.}
the recent~{\cite{Olver-2018}}), we offer in this article another
interpretation of Cartan’s method of $G$-structure reductions, the
application of which goes beyond the plain determination of normal
forms for power series given at the origin.  One of its advantages is
that it includes an explicit approach to finding generating
invariants.  This method will be presented by elaborating on several
examples, before a general theoretical description which will be expressed
in a forthcoming publication.

\medskip\noindent {\bf Acknowledgments.} The realization of this research work in differential invariants has received generous financial support from the scientific grant 2018/29/B/ST1/02583 originating from the Polish National Science Center (NCN).

During spring 2019, in May in Orsay, the authors benefited from countless oral exchanges with Paweł Nurowski (Center For Theoretical Physics), who shared his deep knowledge of Cartan’s method of equivalence by calculating the explicit numerator of the Pick invariant in full details.

In March 2020 in Orsay, the authors learned a lot from Boris Doubrov (Belarusian State University) on geometric structures of homogeneous real hypersurfaces in complex vector spaces.

Grateful thanks are addressed to an anonymous referee for pointing out the classification of developable surfaces.
\Section{\bf Presentation of the Results}
\label{presentation-results}
\HEAD{{\ref{presentation-results}}.~{\sf 
Presentation of the Results}
}{
Zhangchi {\sc Chen}, Joël~{\sc Merker}}

We can now start a precise desciption of our results.
Several aspects are true generally, but we restrict ourselves
to a presentation of the $2$-dimensional case.
Most of our considerations will be of local nature.
We will assume real analyticity throughout.

We consider the special affine group:
\[
\SA_3(\R)
\,=\,
\SL_3(\R)
\ltimes
\R^3,
\]
which consists of invertible linear transformations $(x,y,u)
\longmapsto (s,t,v)$ coupled with translations:
\leqnomode\usetagform{default}
\begin{align}
\label{introduction-a-b-c-d-k-l-m-n-p-q-r-s}
\aligned
s
&
\,=\,
{\sf a}\,x
+
{\sf b}\,y
+
{\sf c}\,u
+
{\sf d},
\\
t
&
\,=\,
{\sf k}\,x
+
{\sf l}\,y
+
{\sf m}\,u
+
{\sf n},
\\
v
&
\,=\,
{\sf p}\,x
+
{\sf q}\,y
+
{\sf r}\,u
+
{\sf s},
\endaligned
\ \ \ \ \ \ \ \ \ \ \ \ \ \ \ \ \ \ \ \
1
\,=\,
\left\vert\!
\begin{array}{ccc}
{\sf a} & {\sf b} & {\sf c}
\\
{\sf k} & {\sf l} & {\sf m}
\\
{\sf p} & {\sf q} & {\sf r}
\end{array}
\!\right\vert,
\end{align}
preserving volume and orientation. We have:
\[
\dim\,
\SA_3(\R)
\,=\,
3\cdot 3-1
+
3
\,=\,
11.
\]
We will always consider special affine transformations not
far from the identity, hence we may view $\SA_3(\R)$ as
a local Lie group.
The full affine group
will be denoted ${\sf A}_3(\R) = \GL_3(\R) \ltimes \R^3$.

In the source space $(x,y,u)$, 
we consider surfaces $S^2 \subset \R_{x,y}^2 \times \R_u^1$ 
graphed as $\big\{ u = F(x,y) \big\}$ with
convergent power series $F \in \R \{x, y\}$, and similarly,
in the target space $(s,t,v)$, we consider 
graphed analytic surfaces
$\big\{ v = G(s,t) \big\}$:
\[
\aligned
u
&
\,=\,
\sum_{j=0}^\infty\,
\sum_{k=0}^\infty\,
F_{j,k}\,
{\textstyle{\frac{x^j}{j!}}}\,
{\textstyle{\frac{y^k}{k!}}},
\\
v
&
\,=\,
\sum_{l=0}^\infty\,
\sum_{m=0}^\infty\,
G_{l,m}\,
{\textstyle{\frac{s^l}{l!}}}\,
{\textstyle{\frac{t^m}{m!}}}.
\endaligned
\]

\begin{Problem}
\label{1-Problem}
{\sl 
Determine when two given surfaces $\big\{ u = F(x,y) \big\}$
and $\big\{ v = G(s,t) \big\}$ are $\SA_3$-equivalent.}
\end{Problem}

When this holds, 
by a special affine transformation, every point $\big(x, y, F(x,y) 
\big)$ is mapped to a point $\big(s, t, G(s,t) \big)$, and
a {\sl fundamental equation} holds in $\R\{ x,y\}$:
\[
{\sf p}\,x
+
{\sf q}\,y
+
{\sf r}\,F(x,y)
+
{\sf s}
\,\,\equiv\,\,
G
\Big(
{\sf a}\,x
+
{\sf b}\,y
+
{\sf c}\,F(x,y)
+
{\sf d},\,\,\,\,
{\sf k}\,x
+
{\sf l}\,y
+
{\sf m}\,F(x,y)
+
{\sf n}
\Big).
\]

\begin{Problem}
\label{2-Problem}
{\sl 
Classify surfaces $\big\{ u = F(x,y) \big\}$ under 
the $\SA_3(\R)$ action, especially, find all
(locally) homogeneous models.}
\end{Problem}

Problem~{\ref{2-Problem}} has been studied by
means of Lie-theoretical
methods. 
The classification of all 
${\sf A}_3(\R)$-homogeneous surfaces with 
identically vanishing Pick
invariant in general affine geometry was obtained 
in~{\cite{Abdalla-Dillen-Vrancken-1997}}.

In~{\cite{Doubrov-Komrakov-Rabinovich-1996}}, 
all locally homogeneous
two-dimensional surfaces in the three-dimensional affine
geometry were described,
with a list of $18$ items.
There exist similar classifications for equiaffine
geometry~{\cite{Guggenheimer-1977,
Jensen-1977}}.

\smallskip

We denote a general element of the special affine group
by $g \in \SA_3(\R)$, and the general transformation as:
\[
s
\,=\,
s\big(g,x,y,u\big),
\ \ \ \ \ \ \ \ \ \ \ \ \
t
\,=\,
t\big(g,x,y,u\big),
\ \ \ \ \ \ \ \ \ \ \ \ \
v
\,=\,
v\big(g,x,y,u\big).
\]

\begin{Definition}
\label{Def-differential-invariant}
A {\sl differential invariant} of order $n$ is a function
of the horizontal coordinates and the partial
derivatives of the graphing function up to order $n$:
\[
\Iaux
\Big(
s,t,\,
\big\{
G_{s^lt^m}
(s,t)
\big\}_{0\leqslant l+m\leqslant n}
\Big)
\,\,\equiv\,\,
\Iaux
\Big(
x,y,\,
\big\{
F_{x^jy^k}(x,y)
\big\}_{0\leqslant j+k\leqslant n}
\Big),
\]
which is unchanged after replacement of $(s,t,v)$ in terms
of $\big(g, x,y,u\big)$, for every $g \in \SA_3(\R)$.
\end{Definition}

\begin{Problem}
\label{3-Problem}
{\sl 
Describe the structure of the algebra of differential invariants
of surfaces under the action of $\SA_3(\R)$.}
\end{Problem}

In this memoir, we will focus on Problems~{\ref{1-Problem}}
and~{\ref{3-Problem}}, in the spirit of~{\cite{Olver-1995, 
Fels-Olver-1999,
Olver-2000, Olver-2007, Olver-2007b, Olver-2018}}.

To a graphed surface $\big\{ u = F(x,y) \big\}$ is associated
its {\sl Hessian matrix:}
\[
\Hessian_F
\,=\,
\left(\!
\begin{array}{cc}
F_{xx} & F_{xy}
\\
F_{yx} & F_{yy}
\end{array}
\!\right).
\]

\begin{Definition}
\label{Def-relative-invariant}
{\cite{Fels-Olver-1998}} A {\sl relative invariant} is a function satisfying:
\[
\Paux
\Big(
s,t,\,
\big\{
G_{s^lt^m}
(s,t)
\big\}_{0\leqslant l+m\leqslant n}
\Big)
\,\,\equiv\,\,
\nonzero
\cdot
\Paux
\Big(
x,y,\,
\big\{
F_{x^jy^k}(x,y)
\big\}_{0\leqslant j+k\leqslant n}
\Big),
\]
with a nowhere vanishing factor, at least when $g \in \SA_3(\R)$ is
not far from the identity.
\end{Definition}

A starting observation 
(Section~{\ref{equiaffine-SA-3-surfaces-S2-R3}})
is that the Hessian {\em determinant} is a relative invariant:
\[
G_{ss}\,
G_{tt}
-
G_{st}^2
\,\,=\,\,
\nonzero
\cdot
\Big(
F_{xx}\,F_{yy}
-
F_{xy}^2
\Big),
\]
even under general affine transformations.
Moreover, Proposition~{\ref{Prp-invariancy-rank-Hessian-matrix}}
shows that the rank of the Hessian matrix 
remains unchanged through any (special) affine transformation.

For the general theory of surfaces, 
this implies an elementary initial branching:
\[
\xymatrix{
&&
\Hessian_F\equiv\text{\scriptsize
$\left(\!\!\begin{array}{cc} 0 & 0 \\ 0 & 0 
\end{array}\!\!\right)$}
\ar[rr]
&&
\{u=0\}
\\
\big\{u=F(x,y)\big\}
\ar[urr]
\ar[rr]
\ar[drr]
&&
{\sf rank}\,\,\Hessian_F\,\equiv\,1
\ar[rr]
&&
\overset{{\sf Root}\,\,{\sf Hypothesis}\,\,
{\sf in}\,\,{\sf this}\,\,{\sf Memoir}}{
\boxed{
F_{xx}\neq 0\equiv F_{xx}F_{yy}-F_{xy}^2}}
\\
&&
{\sf rank}\,\,\Hessian_F\,\equiv\,2
\ar[rr]
&&
{\scriptstyle{{\sf Not}\,\,{\sf treated}\,\,{\sf here}.}}
}
\] 

Geometrically, it is clear that the case where the Hessian
matrix is identically zero:
\[
0
\,\equiv\,
F_{xx}
\,\equiv\,
F_{xy}
\,\equiv\,
F_{yy},
\]
is {\sl flat} in the proper sense, hence there exists a
special affine transformation which maps
any such $\big\{ u = F(x,y) \big\}$
to a reference plane $\{ v = 0 \}$.
This branch is hence trivial.

The rank $2$ case is a wide story
in itself, it conducts to the so-called
{\sl Pick invariant}\footnote{On p.~{\pageref{explicit-Pick}},
the reader will find its explicit expression.
Thanks are adressed to Pawe{\l} Nurowski for his help.}, 
of order $3$, and to further
order $4$ differential invariants, 
{\em cf.}~{\cite{Spivak-1979, Olver-2007}}.

In this memoir, we will study the middle branch only.
After a rotation in the $(x,y)$ space, we can assume that
$F_{xx} (x,y) \neq 0$ is nowhere vanishing (our reasonings are
local). Then our main {\sl root hypothesis} will constantly be:
\[
F_{xx}
\,\,\neq\,\,
0
\,\,\equiv\,\,
F_{xx}\,F_{yy}
-
F_{xy}^2.
\]
Solving:
\[
F_{yy}
\,\equiv\,
\frac{F_{xy}^2}{F_{xx}},
\]
we may differentiate once:
\[
\aligned
F_{xyy}
&
\,=\,
2\,
\frac{F_{xy}\,F_{xxy}}{F_{xx}}
-
\frac{F_{xy}^2\,F_{xxx}}{F_{xx}^2},
\\
F_{yyy}
&
\,=\,
3\,
\frac{F_{xy}^2\,F_{xxy}}{F_{xx}^2}
-
2\,
\frac{F_{xy}^3\,F_{xxx}}{F_{xx}^3},
\endaligned
\]
and so on. 

\begin{center}
\scalebox{1.05}{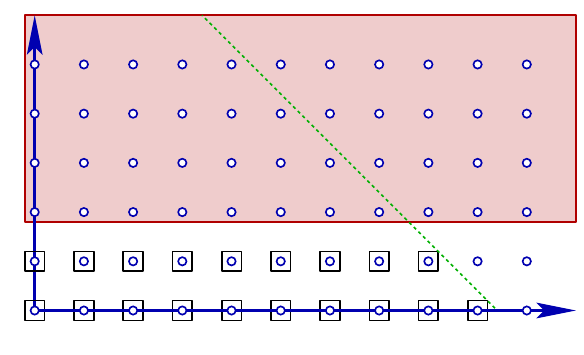}
\end{center}

It is easy to convince oneself 
({\em see}~Section~{\ref{parabolic-jet-relations}})
that every partial derivative $F_{x^j y^k}$
with $k \geqslant 2$
expresses in terms of the partial derivatives:
\[
\big\{
F_{x^{j'}}
\big\}_{j'\leqslant j+k},
\ \ \ \ \
\big\{
F_{x^{j''}y}
\big\}_{j''\leqslant j+k-1}.
\]
This conducts us to introduce the 
{\sl parabolic jet spaces} of any order $n \geqslant 2$:
\[
P\!J_{2,1}^n
\,\,\ni\,\,\,\,
\Big(
x,y,\,\,
\begin{smallmatrix} 
u_y, & \dots & \dots & u_{x^{n-1}y}, &
\\
\\ 
u, & u_x, & \dots, & u_{x^{n-1}}, & u_{x^n}
\end{smallmatrix}
\Big)
\,\,\,\,\in\,\,
\R^{3+2n}.
\]

In effective differential invariant theory, 
for instance in the case of (not necessarily parabolic) 
surfaces,
under any action of a (local) Lie group $G$, 
certain relative invariants are encountered, call them:
\[
\Paux
\,=\,
\Paux
\Big(
x,y,u,\,
\big\{
u_{x^jy^k}
\big\}_{1\leqslant j+k\leqslant n}
\big),
\ \ \ \ \ \ \ \ \ \ \ \
\Qaux,
\ \ \ \ \ \ \ \ \ \ \ \
\Raux,
\ \ \ \ \ \ \ \ \ \ \ \
\dots
\]
According to Definition~{\ref{Def-relative-invariant}},
their zero-sets $\big\{ \Paux = 0 \big\}$, 
$\big\{ \Qaux = 0 \big\}$, \dots, are invariant
under $G$. 
They are responsible for the creation of {\sl branches}
and of further {\sl subbranches:}
\[
\xymatrix{
&&
&&
\Qaux\,\equiv\,0,
\\
&&
\ar[urr]
\Paux\,\equiv\,0
\ar[rr]
&&
\Qaux\,\neq\,0,
\\
\ar[urr]
(\ast)
\ar[drr]
&&
&&
\\
&&
\ar[rr]
\Paux\,\neq\,0
\ar[drr]
&&
\Raux\,\equiv\,0,
\\
&&
&&
\Raux\,\neq\,0.
}
\]

We adopt Lie's principle of thought
({\cite[Chap.~1]{Lie-Merker-2015}),
which admits that either a (relative) differential invariant is
identically zero, or it is assumed to be nowhere zero, after
restriction to an appropriate open subset.  Mixed cases where some
(relative) invariant is nonzero on some nonempty open subset and
vanishes on a nonempty closed subset are excluded from exploration.

Importantly, as soon as some (relative) invariant vanishes identically,
like our Hessian
determinant:
\[
\Haux_F 
\,:=\,
F_{xx}\,F_{yy} 
- 
F_{xy}^2,
\]
one must express {\em all differential consequences} of this
assumption in order to explore properly the concerned branch. When, on
some (sub)branch, there occurs a simultaneous vanishing of two or more
(relative) invariants, one must at first express the differential
consequences under a closed workable form, like setting up a
meaningful Gr\"obner basis for the differential ideal generated.

We can now start to present our results. At first, if we abbreviate:
\[
{\sf root}
\,:=\,
\boxed{
\aligned
0
&
\,\neq\,
F_{xx}
\\
0
&
\,\equiv\,
F_{xx}\,F_{yy}
-
F_{xy}^2
\endaligned}
\,\,=\,\,
\boxed{
F_{xx}\neq 0\equiv\Haux_F}
\]
the branching diagram which summarizes everything is:
\[
\xymatrix{
&&
&&
\Paux\,\equiv\,0
&&
\Caux\,\equiv\,0
&&
\\
&&
\ar[urr]
\Saux\,\equiv\,0
\ar[rr]
&&
\ar[urr]
\Paux\,\neq\,0
\ar[rr]
&&
\Caux\,\neq\,0
&&
\\
\ar[urr]
\ar[rr]
\boxed{\root}
&&
\ar[rr]
\Saux\,\neq\,0
\ar[ddrr]
&&
\ar[rr]
\Waux\,\equiv\,0
\ar[drr]
&&
\Xaux\,\equiv\,0
&&
\\
&&
&&
&&
\ar[rr]
\Xaux\,\neq\,0
\ar[drr]
&&
\Yaux\,\equiv\,0
\\
&&
&&
\ar[rr]
\Waux\,\neq\,0
\ar[drr]
&&
\Maux\,\equiv\,0
&&
\Yaux\,\neq\,0
\\
&&
&&
&&
\Maux\,\neq\,0
&&
}
\]
This tree decomposes in $3$ main branches, 
extracted in three diagrams below, 
just before the statements of 
$3$ associated theorems. 

In the first, top branch, $\Saux$ and $\Paux$ are relative invariants:
\[
\aligned
\Saux
&
\,:=\,
\frac{F_{xx}\,F_{xxy}-F_{xy}\,F_{xxx}}{F_{xx}^2},
\\
\Paux
&
\,:=\,
\frac{1}{3}\,
\frac{-\,5\,F_{xxx}^2+3\,F_{xx}\,F_{xxxx}}{F_{xx}^2},
\endaligned
\]
while $\Caux$ is a differential invariant:
\[
\Caux
\,:=\,
\frac{1}{\sqrt{3}}\,
\frac{9\,F_{xx}^2\,F_{xxxxx}
-45\,F_{xx}\,F_{xxx}\,F_{xxxx}
+40\,F_{xxx}^3}{
\big(
\pm\,3\,F_{xx}\,F_{xxxx}
\mp\,5\,F_{xxx}^2
\big)^{3/2}}.
\]

In the second, middle branch, $\Waux$ is a differential invariant,
but it is assumed to vanish identically, hence it is trivial,
and further, $\Xaux$ and $\Yaux$ are differential invariants:
\[
\Xaux
\,:=\,
\frac{1}{9}\,
\frac{\big(F_{xx}\,F_{xxy}-F_{xy}\,F_{xxx}\big)\,
\big(9\,F_{xx}^2\,F_{xxxxx}
-45\,F_{xx}\,F_{xxx}\,F_{xxxx}
+40\,F_{xxx}^3\big)}{
F_{xx}^6},
\]

\[
\scriptsize
\aligned
\Yaux
&
\,:=\,
\frac{1}{18}\,
\frac{
\left( 
F_{{xxy}}F_{{xx}}-F_{{xxx}}F_{{xy}} 
\right)^{5/3}
}{
{F_{{xx}}}^{10} 
\left( 
9\,F_{{xxxxx}}{F_{{xx}}}^{2}-45\,F_{{xxx}}F_{{xxxx}}F_{{xx}}+40\,{F_{{xxx}}}^{3}
\right)}\,
\bigg\{
\\
&
\bigg\{
11200\,{F_{{xxx}}}^{8}
-
12600\,{F_{{xxx}}}^{3}F_{{xxxxx}}{F_{{xx}}}^{3}F_{{xxxx}}
+
13230\,F_{{xxx}}F_{{xxxxx}}{F_{{xx}}}^{4}{F_{{xxxx}}}^{2}
+
1134\,F_{{xxx}}F_{{xxxxx}}{F_{{xx}}}^{5}F_{{xxxxxx}}
\\
&
-\,
3150\,{F_{{xxx}}}^{2}F_{{xxxx}}{F_{{xx}}}^{4}F_{{xxxxxx}}
-
810\,F_{{xxxxxxx}}{F_{{xx}}}^{5}F_{{xxx}}F_{{xxxx}}
-
33600\,{F_{{xxx}}}^{6}F_{{xxxx}}F_{{xx}}
-
7875\,{F_{{xxx}}}^{2}{F_{{xxxx}}}^{3}{F_{{xx}}}^{3}
\\
&
-\,
756\,{F_{{xxx}}}^{2}{F_{{xxxxx}}}^{2}{F_{{xx}}}^{4}
+
6720\,{F_{{xxx}}}^{5}F_{{xxxxx}}{F_{{xx}}}^{2}
+
31500\,{F_{{xxx}}}^{4}{F_{{xxxx}}}^{2}{F_{{xx}}}^{2}
-
4725\,{F_{{xxxx}}}^{4}{F_{{xx}}}^{4}
\\
&
-\,
189\,{F_{{xxxxxx}}}^{2}{F_{{xx}}}^{6}
+
1890\,{F_{{xxxx}}}^{2}{F_{{xx}}}^{5}F_{{xxxxxx}}
-
2835\,F_{{xxxx}}{F_{{xx}}}^{5}{F_{{xxxxx}}}^{2}
+
162\,F_{{xxxxxxx}}{F_{{xx}}}^{6}F_{{xxxxx}}
\\
&
+
720\,F_{{xxxxxxx}}{F_{{xx}}}^{4}{F_{{xxx}}}^{3}
\bigg\}.
\endaligned
\]

In the third, last, bottom branch, 
$\Waux$ is a nontrivial differential invariant:
\[
\Waux
\,:=\,
\frac{
F_{xx}^2\,F_{xxxy}
-
F_{xx}\,F_{xy}\,F_{xxxx}
+
2\,F_{xy}\,F_{xxx}^2
-
2\,F_{xx}\,F_{xxx}\,F_{xxy}}{
(F_{xx})^2\,\,
\big(
F_{xx}\,F_{xxy}
-
F_{xy}\,F_{xxx}
\big)^{2/3}},
\]
and $\Maux$ is also a differential invariant:

\[
\scriptsize
\aligned
\Maux
&
\,:=\,
\frac{1}{36}\,
\frac{1}{
{F_{{xx}}}^{6}\,\,
\left(-\,F_{{xxy}}F_{{xx}}+F_{{xxx}}F_{{xy}}\right)\,\,
\left(F_{{xx}}F_{{xxxx}}F_{{xy}}
-{F_{{xx}}}^{2}F_{{xxxy}}
+2\,F_{{xxx}}F_{{xxy}}F_{{xx}}
-2\,{F_{{xxx}}}^{2}F_{{xy}}\right)}\,
\bigg\{
\\
&
\bigg\{
\ \ 
270\,{F_{{xx}}}^{6}F_{{xxxxy}}{F_{{xxy}}}^{2}F_{{xxxx}}
-
72\,F_{{xxx}}F_{{xxxxx}}{F_{{xx}}}^{5}{F_{{xxy}}}^{3}
+
820\,F_{{xxx}}{F_{{xx}}}^{3}{F_{{xxxx}}}^{3}{F_{{xy}}}^{3}
-
2195\,{F_{{xxx}}}^{3}{F_{{xx}}}^{2}{F_{{xxxx}}}^{2}{F_{{xy}}}^{3}
\\
&
+\,
2560\,{F_{{xxx}}}^{5}F_{{xx}}F_{{xxxx}}{F_{{xy}}}^{3}
+
2000\,{F_{{xxx}}}^{2}{F_{{xx}}}^{5}{F_{{xxxy}}}^{2}F_{{xxy}}
-
2000\,{F_{{xxx}}}^{3}{F_{{xx}}}^{4}{F_{{xxxy}}}^{2}F_{{xy}}
-
3040\,{F_{{xxx}}}^{3}{F_{{xx}}}^{4}F_{{xxxy}}{F_{{xxy}}}^{2}
\\
&
-
3040\,{F_{{xxx}}}^{5}{F_{{xx}}}^{2}F_{{xxxy}}{F_{{xy}}}^{2}
-
3840\,{F_{{xxx}}}^{5}{F_{{xxy}}}^{2}{F_{{xx}}}^{2}F_{{xy}}
+
3840\,{F_{{xxx}}}^{6}F_{{xxy}}F_{{xx}}{F_{{xy}}}^{2}
-
420\,{F_{{xxxx}}}^{3}{F_{{xx}}}^{4}F_{{xxy}}{F_{{xy}}}^{2}
\\
&
+\,
480\,F_{{xxxx}}{F_{{xx}}}^{4}{F_{{xxy}}}^{3}{F_{{xxx}}}^{2}
-
420\,F_{{xxy}}{F_{{xx}}}^{6}F_{{xxxx}}{F_{{xxxy}}}^{2}
+
192\,{F_{{xxx}}}^{4}{F_{{xx}}}^{2}F_{{xxxxx}}{F_{{xy}}}^{3}
-
120\,{F_{{xxx}}}^{2}{F_{{xx}}}^{5}F_{{xxxxy}}{F_{{xxy}}}^{2}
\\
&
-\,
120\,{F_{{xxx}}}^{4}{F_{{xx}}}^{3}F_{{xxxxy}}{F_{{xy}}}^{2}
+
36\,F_{{xxxxx}}{F_{{xx}}}^{6}{F_{{xxy}}}^{2}F_{{xxxy}}
+
45\,{F_{{xx}}}^{6}{F_{{xxxxy}}}^{2}F_{{xxx}}F_{{xy}}
-
45\,{F_{{xx}}}^{5}{F_{{xxxxx}}}^{2}{F_{{xy}}}^{2}F_{{xxy}}
\\
&
+\,
45\,{F_{{xx}}}^{4}{F_{{xxxxx}}}^{2}{F_{{xy}}}^{3}F_{{xxx}}
-
120\,{F_{{xx}}}^{4}F_{{xxxxx}}{F_{{xy}}}^{3}{F_{{xxxx}}}^{2}
-
120\,{F_{{xx}}}^{6}F_{{xxxxx}}F_{{xy}}{F_{{xxxy}}}^{2}
+
120\,{F_{{xx}}}^{5}F_{{xxxxy}}{F_{{xxxx}}}^{2}{F_{{xy}}}^{2}
\\
&
+\,
1280\,{F_{{xxx}}}^{4}{F_{{xxy}}}^{3}{F_{{xx}}}^{3}
-
400\,F_{{xxx}}{F_{{xx}}}^{6}{F_{{xxxy}}}^{3}
-
45\,{F_{{xx}}}^{7}{F_{{xxxxy}}}^{2}F_{{xxy}}
-
405\,{F_{{xxy}}}^{3}{F_{{xx}}}^{5}{F_{{xxxx}}}^{2}
\\
&
+\,
120\,{F_{{xx}}}^{7}F_{{xxxxy}}{F_{{xxxy}}}^{2}
-
1280\,{F_{{xxx}}}^{7}{F_{{xy}}}^{3}
-
5200\,{F_{{xxx}}}^{2}{F_{{xx}}}^{4}F_{{xxxx}}F_{{xy}}F_{{xxxy}}F_{{xxy}}
+
432\,F_{{xxx}}{F_{{xx}}}^{4}F_{{xxxx}}{F_{{xy}}}^{2}F_{{xxxxx}}F_{{xxy}}
\\
&
-\,
360\,F_{{xxx}}{F_{{xx}}}^{5}F_{{xxxx}}F_{{xy}}F_{{xxxxy}}F_{{xxy}}
+
108\,F_{{xxx}}{F_{{xx}}}^{5}F_{{xxxy}}F_{{xxxxx}}F_{{xy}}F_{{xxy}}
-
2040\,F_{{xxx}}{F_{{xx}}}^{4}{F_{{xxxx}}}^{2}{F_{{xy}}}^{2}F_{{xxxy}}
\\
&
+
1985\,{F_{{xxx}}}^{2}{F_{{xx}}}^{3}{F_{{xxxx}}}^{2}{F_{{xy}}}^{2}F_{{xxy}}
+
1620\,F_{{xxx}}{F_{{xx}}}^{5}F_{{xxxx}}F_{{xy}}{F_{{xxxy}}}^{2}
+
4600\,{F_{{xxx}}}^{3}{F_{{xx}}}^{3}F_{{xxxx}}{F_{{xy}}}^{2}F_{{xxxy}}
\\
&
+
1600\,{F_{{xxx}}}^{3}{F_{{xx}}}^{3}F_{{xxxx}}F_{{xy}}{F_{{xxy}}}^{2}
-
4640\,{F_{{xxx}}}^{4}{F_{{xx}}}^{2}F_{{xxxx}}{F_{{xy}}}^{2}F_{{xxy}}
+
6080\,{F_{{xxx}}}^{4}{F_{{xx}}}^{3}F_{{xxxy}}F_{{xxy}}F_{{xy}}
\\
&
+
840\,{F_{{xxxx}}}^{2}{F_{{xx}}}^{5}F_{{xxy}}F_{{xy}}F_{{xxxy}}
+
615\,{F_{{xxxx}}}^{2}{F_{{xx}}}^{4}{F_{{xxy}}}^{2}F_{{xy}}F_{{xxx}}
+
600\,F_{{xxxx}}{F_{{xx}}}^{5}{F_{{xxy}}}^{2}F_{{xxxy}}F_{{xxx}}
\\
&
+
336\,{F_{{xxx}}}^{2}{F_{{xx}}}^{4}F_{{xxxxx}}F_{{xy}}{F_{{xxy}}}^{2}
-
456\,{F_{{xxx}}}^{3}{F_{{xx}}}^{3}F_{{xxxxx}}{F_{{xy}}}^{2}F_{{xxy}}
-
126\,{F_{{xxx}}}^{2}{F_{{xx}}}^{3}F_{{xxxx}}{F_{{xy}}}^{3}F_{{xxxxx}}
\\
&
+
90\,{F_{{xxx}}}^{2}{F_{{xx}}}^{4}F_{{xxxx}}{F_{{xy}}}^{2}F_{{xxxxy}}
-
144\,{F_{{xxx}}}^{2}{F_{{xx}}}^{4}F_{{xxxy}}F_{{xxxxx}}{F_{{xy}}}^{2}
-
306\,{F_{{xx}}}^{5}F_{{xxxxx}}F_{{xy}}F_{{xxxx}}{F_{{xxy}}}^{2}
\\
&
+
240\,{F_{{xxx}}}^{3}{F_{{xx}}}^{4}F_{{xxxxy}}F_{{xxy}}F_{{xy}}
-
180\,F_{{xxx}}{F_{{xx}}}^{6}F_{{xxxy}}F_{{xxxxy}}F_{{xxy}}
+
180\,{F_{{xxx}}}^{2}{F_{{xx}}}^{5}F_{{xxxy}}F_{{xxxxy}}F_{{xy}}
\\
&
+
90\,{F_{{xx}}}^{6}F_{{xxxxx}}F_{{xy}}F_{{xxxxy}}F_{{xxy}}
-
90\,{F_{{xx}}}^{5}F_{{xxxxx}}{F_{{xy}}}^{2}F_{{xxxxy}}F_{{xxx}}
+
240\,{F_{{xx}}}^{5}F_{{xxxxx}}{F_{{xy}}}^{2}F_{{xxxx}}F_{{xxxy}}
\\
&
\,-
240\,{F_{{xx}}}^{6}F_{{xxxxy}}F_{{xxxx}}F_{{xy}}F_{{xxxy}}
\bigg\}.
\endaligned
\]

It is important to show these invariants, because they not
only cause the branchings, but also, they will constitute
generating collections for the full algebras of differential
invariants.

We may now state our results for the three kinds of branches.
We always start from our root assumption.

\[
\xymatrix{
&&
&&
\Paux\,\equiv\,0
&&
\Caux\,\equiv\,0
&&
\\
&&
\ar[urr]
\Saux\,\equiv\,0
\ar[rr]
&&
\ar[urr]
\Paux\,\neq\,0
\ar[rr]
&&
\Caux\,\neq\,0
&&
\\
\ar[urr]
\boxed{
F_{xx}\neq 0\equiv\Haux_F}
&&
&&
&&
&&
}
\]
The full affine group in two dimensions is 
${\sf A}_2(\R) = \GL_2(\R) \ltimes \R^2$.

\begin{Theorem}
\label{Thm-S-0}
Within the branch $\Saux \equiv 0$:

\smallskip\noindent{\bf (1)}\,
Every surface $S^2 \subset \R^3$ is special affinely equivalent
to the product of a curve in $\R_{x,u}^{1+1}$ and $\R_{y}^1$,
and $\SA_3(\R)$-equivalences amount to ${\sf A}_2(\R)$-equivalences
of such curves;

\smallskip\noindent{\bf (2)}\,
There is a relative invariant $\Paux$ of order $4$;

\smallskip\noindent{\bf (3)}\,
When $\Paux \equiv 0$, the surface is $\SA_3$-equivalent
to $\big\{ u = x^2 \big\}$,
the product of a parabola and $\R_y^1$, and conversely;

\smallskip\noindent{\bf (4)}\,
When $\Paux \neq 0$, the surface is, in a unique way,
$\SA_3$-equivalent to:
\[
u
\,=\,
{\textstyle{\frac{x^2}{2!}}}
\pm
{\textstyle{\frac{x^4}{4!}}}
+
F_{5,0}\,
{\textstyle{\frac{x^5}{5!}}}
+
\sum_{j\geqslant 6}\,
F_{j,0}\,
{\textstyle{\frac{x^j}{j!}}},
\]
and the collection of coefficients $F_{5,0}$, 
$\big\{ F_{j,0} \big\}_{j
\geqslant 6}$ is in one-to-one
correspondence with equivalent classes.
\end{Theorem}

Here:
\[
F_{5,0}
\,=\,
F_{xxxxx}(0)
\,=\,
\text{\rm value of}\,\,
\Caux\,\,
\text{\rm at the origin}.
\]
Infinitely many differential invariants correspond to
these coefficients $F_{j,0}$, as we will soon explain.

\begin{Question}
{\sl 
How to compute explicitly differential invariants?}
\end{Question}

It is clear that $\SA_3(\R)$ contains all translations:
\[
s
\,=\,
x+{\sf d},
\ \ \ \ \ \ \ \ \ \ \ \ \ \ \ \ \ \ \ \
t
\,=\,
y+{\sf n},
\ \ \ \ \ \ \ \ \ \ \ \ \ \ \ \ \ \ \ \
v
\,=\,
u+{\sf s}.
\]
This implies\,\,---\,\,exercise from 
Definition~{\ref{Def-differential-invariant}}, or {\em see}
Theorem~{\ref{Thm-translations-G-0}}\,\,---\,\,that every
differential invariant:
\[
\Iaux
\Big(
\underbrace{
x,y,u}_{\sf absent},\,
\big\{
u_{x^jy^k}
\big\}_{1\leqslant j+k\leqslant n}
\Big),
\]
must depend only on jet derivatives of order $\geqslant 1$.

To compute these invariants $\Iaux$, we start from a power series
at the origin:
\[
u
\,=\,
\sum_{j+k\geqslant 1}\,
F_{j,k}\,
\frac{x^j}{j!}\,
\frac{y^k}{k!},
\]
and we progressively perform (several) `simple', `natural',
special affine transformations in order to annihilate\big/normalize
as much as possible Taylor coefficients $F_{j,k}$.
Rigorous descriptions illustrated by examples
will be provided in 
Sections~{\ref{special-affine-power-series-invariants-curves-R-2}},
{\ref{affine-invariants-curves-R-2}},
{\ref{parabolic-pseudostablization}},
{\ref{relative-invariant-S-first-invariant-W}},
{\ref{branch-W-equiv-0-branch-W-nonzero}},
but here, we only present general
ideas. One main feature of the process is its progressivity.

At the end, we reach a certain `{\sl normal form}':
\[
v
\,=\,
\sum_{l+m\geqslant 1}\,
G_{l,m}\,
\frac{s^l}{l!}\,
\frac{t^m}{m!},
\]
in which several coefficients $G_{l,m}$ are `simplified',
for instance as in Theorem~{\ref{Thm-S-0}} above:
\[
G_{1,0}
\,=\,
G_{0,1}
\,=\,
0,
\ \ \ \ \ \ \ \ \ \ \ \ \ \ \ \ \ \ \ \
G_{2,0}
\,=\,
1,
\ \ \ \ \ \ \ \ \ \ \ \ \ \ \ \ \ \ \ \
G_{3,0}
\,=\,
0,
\ \ \ \ \ \ \ \ \ \ \ \ \ \ \ \ \ \ \ \
\text{\em etc.}
\]
Certainly, 
the full composition of all the progressively normalizing 
maps belongs to $\SA_3(\R)$, hence is of the 
form~({\ref{introduction-a-b-c-d-k-l-m-n-p-q-r-s}})
for some specific constants ${\sf a}$, \dots, ${\sf s}$.
These constants are complicated at the end, but step-by-step
they are simple, only the full composition of normalizing maps
creates complexity.

After the process is pushed at its farthest point, the identity
map of $\SA_3(\R)$ is the only transformation which leaves
untouched the `{\sl normal form}' of the power series
$v = G(s,t)$. 

While normalizing low order Taylor coefficients,
we also keep track (on a computer) of the way how the {\em other}
(higher order) Taylor coefficients are modified.
At the end, we receive formulas:
\[
G_{l,m}
\,=\,
\Pi_{l,m}
\Big(
\big\{
F_{j,k}
\big\}_{1\leqslant j+k\leqslant l+m}
\Big)
\eqno
{\scriptstyle{(l\,+\,m\,\geqslant\,1)}}.
\]
Then granted that:
\[
F_{j,k}
\,=\,
u_{x^jy^k}(0,0),
\]
all the desired
genuine differential invariants are obtained simply by replacing
in these formulas Taylor coefficients by jet coordinates:
\[
\Iaux_{l,m}
\,:=\,
\Pi_{l,m}
\Big(
\big\{
u_{x^jy^k}
(x,y)
\big\}_{1\leqslant j+k\leqslant l+m}
\Big)
\eqno
{\scriptstyle{(l\,+\,m\,\geqslant\,1)}}.
\]
Importantly, the hypothesis that the group contains 
all translations guarantees 
(Theorem~{\ref{Thm-translations-G-0}})
that we obtain the
expressions of {\em all} differential invariants
at {\em every} point $(x,y)$ near the origin.

This process could be explained abstractly in any dimension
(forthcoming).
During normalizations, relative invariants play a crucial role.

\begin{Observation}
Any (relative) invariant $\Paux$:

\smallskip\noindent$\bullet$\,
Either creates a new branch $\Paux \equiv 0$ to be explored
farther;

\smallskip\noindent$\bullet$\,
or is absorbed, when $\Paux \neq 0$, into some constant by normalization.

\end{Observation}

This is, for instance, true of 
$\Saux$, $\Paux$, $\Waux$, $\Xaux$: when they are nonzero,
they will be used to normalize some Taylor coefficients.

\begin{Theorem}
With the assumption $\Saux \neq 0$, there is exactly one
differential invariant of fourth order, $\Waux$.
\end{Theorem}

\smallskip

We can now state our second result for the second, middle branch.
(The third, bottom branch will also assume $\Saux \neq 0$.)
\[
\xymatrix{
\ar[rr]
\boxed{
F_{xx}\neq 0\equiv\Haux_F}
&&
\ar[rr]
\Saux\,\neq\,0
&&
\ar[rr]
\Waux\,\equiv\,0
\ar[drr]
&&
\Xaux\,\equiv\,0
&&
\\
&&
&&
&&
\ar[rr]
\Xaux\,\neq\,0
\ar[drr]
&&
\Yaux\,\equiv\,0
\\
&&
&&
&&
&&
\Yaux\,\neq\,0
}
\]

\begin{Theorem}
\label{Thm-S-nonzero-W-0}
Within the branch $\Saux \neq 0$, $\Waux \equiv 0$:

\smallskip\noindent{\bf (1)}\,
There is a single invariant, $\Xaux$, of order $5$;

\smallskip\noindent{\bf (2)}\,
When $\Xaux \equiv 0$, every surface $S^2 \subset \R^3$
is $\SA_3$-equivalent to the model:
\[
\aligned
u
&
\,=\,
\frac{1}{2}\,
\frac{x^2}{1-y}
\\
&
\,=\,
{\textstyle{\frac{x^2}{2}}}
+
{\textstyle{\frac{x^2\,y}{2}}}
+
{\textstyle{\frac{x^2\,y^2}{2}}}
+
{\textstyle{\frac{x^2\,y^3}{2}}}
+\cdots+
{\textstyle{\frac{x^2\,y^k}{2}}}
+\cdots;
\endaligned
\]

\smallskip\noindent{\bf (3)}\,
When $\Xaux \neq 0$, every surface is $\SA_3$-equivalent to:
\[
\!\!\!\!\!\!\!\!\!\!\!\!\!\!\!
\aligned
u
\,=\,
{\textstyle{\frac{x^2}{2}}}
+
{\textstyle{\frac{x^2\,y}{2}}}
+
{\textstyle{\frac{x^2\,y^2}{2}}}
+
F_{5,0}\,
\frac{x^5}{120}
+
{\textstyle{\frac{x^2\,y^3}{2}}}
+
4\,F_{5,0}\,
\frac{x^5y}{120}
+
{\textstyle{\frac{x^2\,y^4}{2}}}
+
F_{7,0}\,
\frac{x^7}{5\,040}
&
+
20\,F_{5,0}\,
\frac{x^5y^2}{240}
+
{\textstyle{\frac{x^2\,y^5}{2}}}
\,+
\\
&
+
\sum_{j+k\geqslant 8}\,
F_{j,k}\,
x^jy^k,
\endaligned
\]
with:
\[
\aligned
F_{5,0}
\,=\,
\text{\rm value of}\,\,
\Xaux\,\,
\text{\rm at the origin},
\\
F_{7,0}
\,=\,
\text{\rm value of}\,\,
\Yaux\,\,
\text{\rm at the origin};
\endaligned
\]

\smallskip\noindent{\bf (4)}\,
The collection of coefficients $F_{5,0}$,
$F_{7,0}$, $\big\{ F_{j,0} \big\}_{j
\geqslant 8}$ is in one-to-one
correspondence with equivalent classes.
\end{Theorem}

Lastly, we treat the main (thickest) branch:
\[
\xymatrix{
\ar[rr]
\boxed{
F_{xx}\neq 0\equiv\Haux_F}
&&
\Saux\,\neq\,0
\ar[ddrr]
&&
&&
&&
\\
&&
&&
&&
&&
\\
&&
&&
\ar[rr]
\Waux\,\neq\,0
\ar[drr]
&&
\Maux\,\equiv\,0
&&
\\
&&
&&
&&
\Maux\,\neq\,0.
&&
}
\]

\begin{Theorem}
\label{Thm-S-nonzero-W-nonzero}
Within the main branch $S \neq 0$, $\Waux \neq 0$:

\smallskip\noindent{\bf (1)}\,
There is a single differential invariant $\Maux$, of order
$5$, differentiably independent of $\Waux$;

\smallskip\noindent{\bf (2)}\,
Every surface $S^2 \subset \R^3$ is $\SA_3$-equivalent to:
\[
\aligned
u
\,=\,
{\textstyle{\frac{x^2}{2}}}
+
{\textstyle{\frac{x^2\,y}{2}}}
+
F_{3,1}\,
\frac{x^3\,y}{6}
+
{\textstyle{\frac{x^2\,y^2}{2}}}
+
F_{5,0}\,
\frac{x^5}{120}
+
6\,F_{3,1}\,
\frac{x^3\,y^2}{12}
&
+
{\textstyle{\frac{x^2\,y^3}{2}}}
\,+
\\
&
+
\sum_{j+k\geqslant 6}\,
F_{j,k}\,
x^jy^k,
\endaligned
\]
with:
\[
\aligned
F_{3,1}
\,=\,
\text{\rm value of}\,\,
\Waux\,\,
\text{\rm at the origin},
\\
F_{5,0}
\,=\,
\text{\rm value of}\,\,
\Maux\,\,
\text{\rm at the origin};
\endaligned
\]

\smallskip\noindent{\bf (3)}\,
Any other surface $\big\{ v = G(s,t) \big\}$
within the same branch similarly put into the form:
\[
\aligned
v
\,=\,
{\textstyle{\frac{s^2}{2}}}
+
{\textstyle{\frac{s^2\,t}{2}}}
+
G_{3,1}\,
\frac{s^3\,t}{6}
+
{\textstyle{\frac{s^2\,t^2}{2}}}
+
G_{5,0}\,
\frac{s^5}{120}
+
6\,G_{3,1}\,
\frac{s^3\,t^2}{12}
&
+
{\textstyle{\frac{s^2\,t^3}{2}}}
\,+
\\
&
+
\sum_{l+m\geqslant 6}\,
G_{l,m}\,
s^lt^m,
\endaligned
\]
is $\SA_3$-equivalent to $\big\{ u = F(x,y) \big\}$ above
if and only if all (independent) Taylor coefficients 
in the parabolic jet space match:
\reqnomode\usetagform{EngelLie}
\begin{align}
G_{3,1}
\,=\,
F_{3,1},
\ \ \ \ \ \ \ \ \ \ \ \ \ \ \ \ \ \ \ \
G_{5,0}
\,=\,
F_{5,0},
\ \ \ \ \ \ \ \ \ \ \ \ \ \ \ \ \ \ \ \
G_{l,0}
&
\,=\,
F_{l,0}
\tag{(l\,\geqslant\,6),}
\\
G_{l,1}
&
\,=\,
F_{l,1}
\tag{(l\,\geqslant\,5).}
\end{align}
\end{Theorem}

In these three 
Theorems~{\ref{Thm-S-0}},
{\ref{Thm-S-nonzero-W-0}},
{\ref{Thm-S-nonzero-W-nonzero}},
there always exist, according to 
Fels-Olver~{\cite{Fels-Olver-1999}}, 
two invariant differential operators
$\mathcal{D}_1$ and $\mathcal{D}_2$ satisfying:
\[
\mathcal{D}_i
\big(
{\sf differential}\,\,{\sf invariant}
\big)
\,\,=\,\,
{\sf differential}\,\,{\sf invariant}
\eqno
{\scriptstyle{(i\,=\,1,\,2)}},
\]
and they are non-commuting, in general. 
Invariantly derivating an invariant means applying
$\mathcal{D}_1$ and $\mathcal{D}_2$ several times,
in any order.
More explanations
will be given in 
Sections~{\ref{moving-frame-method}},
{\ref{meaning-differential-invariant}},
{\ref{recurrence-formulae-differential-invariants}},
see especially~({\ref{equation-definition-invariant-D-j}}).

\begin{Theorem}
\label{Thm-minimal-generating-set}
The full algebra of differential invariants under the 
action of $\SA_3(\R)$ is generated by:

\smallskip\noindent\text{\ding{192}}\,
In the branch $\Saux \equiv 0$, $\Paux \neq 0$:
\[
\Caux\ \ 
\text{and its invariant derivatives};
\]

\smallskip\noindent\text{\ding{193}}\,
In the branch $\Saux \neq 0$, $\Waux \equiv 0$:
\[
\Xaux,\ \
\Yaux\ \ 
\text{and their invariant derivatives};
\]

\smallskip\noindent\text{\ding{194}}\,
In the branch $\Saux \neq 0$, $\Waux \neq 0$:
\[
\Waux,\ \
\Maux\ \ 
\text{and their invariant derivatives}.
\]
\end{Theorem}

It is well known that parabolic surfaces coincide with developable
surfaces. While normalizing the ${\sf SA}_3(\R)$-moving frames of
parabolic surfaces using Cartan's formalism, 
Guggenheimer obtained degenerate branches of
cylinders and cones in~{\cite[p.~295]{Guggenheimer-1977}}.  His work
can be re-expressed in terms of explicit differential invariants. In
Section~{\ref{relation-developable-surfaces}}, we will prove that a
parabolic surface is a cylinder if and only if $\Saux \equiv 0$; a
cone if and only if $\Saux \neq 0$ and $\Waux \equiv 0$; a tangential
surface (tangents of a space curve) if and only if $\Saux \neq 0$ and
$\Waux \neq 0$.

\Section{\bf Graph Transformations and Jet Spaces}
\label{graph-transformations-jet-spaces}
\HEAD{{\ref{graph-transformations-jet-spaces}}.~{\sf 
Graph Transformations and Jet Spaces}
}{
Zhangchi {\sc Chen}, Joël~{\sc Merker}}

We provide here a reminder on
how graphs and their transformations 
can be lifted to jet spaces of any order. All functions and
geometric objects will be assumed {\em analytic}, or 
{\em smooth}, or even of a
finite but high enough differentiability
class, provided all differentiations necessary
in our processes are permitted.
Over either $\R$ or $\C$, 
we will mostly work in a local framework,
not mentioning restricted open (sub)sets
in which calculations are true. We will deal only
with finite-dimensional Lie groups 
(but see~\cite{Olver-Pohjanpelto-2008, Olver-Pohjanpelto-2009}).

\Subsection{Jet spaces and jet notations}
Consider $p \geqslant 1$ independent variables $(x^1, \dots, x^p)$ and
also $q \geqslant 1$ variables $(u^1, \dots, u^q)$ which are {\sl
dependent} in the sense that they should be components $u^\alpha =
{\tt u}^\alpha (x)$, $1 \leqslant \alpha \leqslant q$, of maps 
${\tt u} \colon \R^p
\longrightarrow \R^q$.  For any $n \geqslant 0$, introduce the
$n$\textsuperscript{th} order jet space of such maps:
\[
J_{x,u}^n
\,:=\,
J^n
\big(
\R_x^p
\longrightarrow
\R_u^q
\big)
\,\equiv\,
J_{p,q}^n.
\]
This $J_{x,u}^n$ is equipped with coordinates:
\[
z^{(n)}
\,:=\,
\big(
x^j,\,u^\alpha,\,
u_{x^J}^\beta
\big),
\]
where $J = (j_1, \dots, j_\lambda)$ is an unordered
multi-index,
with $1 \leqslant
\lambda = \vert J \vert \leqslant
n$, with
$1 \leqslant j_1, \dots, j_\lambda \leqslant p$, and where each
$u_{x^J}^\alpha \equiv u_J^\alpha$ 
is an independent real coordinate corresponding to the
partial derivative:
\[
\frac{\partial {\tt u}^\alpha(x)}{
\partial x^{j_1}\dots\partial x^{j_\lambda}}
\,\,\,\longleftrightarrow\,\,\,
u_{x^J}^\alpha
\,\equiv\,
u_J^\alpha
\eqno
{\scriptstyle{(1\,\leqslant\,\alpha\,\leqslant\,q,\,\,
J\,=\,(j_1,\dots,j_\lambda))}}.
\]

For instance when $p = q = 1$, we will denote jet coordinates
sometimes $u_x$, $u_{xx}$, $u_{xxx}$, $\dots$, sometimes $u_1$, $u_2$,
$u_3$, $\dots$ Thus:
\[
\dim\,
J_{p,q}^n
\,=\,
p
+
q\,
\binom{p+n}{n}.
\]

When $p \geqslant 2$, say with $p = 2$ and
$q = 1$ to simplify, an alternative {\sl
multi-index} notation will sometimes be employed, especially when
working on any current computer algebra system
like {\sc Mathematica}, {\sc Maple}, {\sc Sage}
and others, for instance:
\[
u_{x^1}
\,\equiv\,
u_{1}
\,\equiv\,
u_{[1,0]},
\ \ \ \ \ \ \ \ \ \ \ \ \ \ \ \ \ \ \ \
u_{x^2}
\,\equiv\,
u_{2}
\,\equiv\,
u_{[0,1]},
\ \ \ \ \ \ \ \ \ \ \ \ \ \ \ \ \ \ \ \
u_{x^1x^2x^1}
\,\equiv\,
u_{1,2,1}
\,\equiv\,
u_{[2,1]}.
\]
To translate, we just count the number of times every index
$1 \leqslant i \leqslant p$ appears in the sequence
$j_1, \dots, j_\lambda$: 
\[
\kappa_i
\,:=\,
\Card\,
\big\{
1\leqslant\nu\leqslant\lambda
\colon\,
j_\nu
=
i
\big\},
\]
and we admit the notational coincidences:
\[
u_{x^{j_1}\cdots x^{j_\lambda}}
\,\equiv\,
u_{j_1,\dots,j_\lambda}
\,\equiv\,
u_{[\kappa_1,\dots,\kappa_p]}.
\]
Sometimes even, when employing the second multi-index notation,
we will drop the brackets:
\[
u_{[\kappa_1,\dots,\kappa_p]}
\,\equiv\,
u_{\kappa_1,\dots,\kappa_p},
\]
especially when studying parabolic surfaces $S^2 \subset \R^3$
later on. There will be no risk of confusion.

\Subsection{Prolongations of diffeomorphisms to jet spaces}
Without full proofs but with 
presentations of ideas, let us review the fundamental
prolongation formulas which are anyway rarely used in applications
because of their complexity.  For completeness of the exposition, it
is nevertheless necessary to be a bit specific about that, since in
some cases anyway, symbolic computer softwares happen to be able to
handle such formulas.

Consider a diffeomorphism 
$\phi \colon \R^{p+q} \longrightarrow \R^{p+q}$
between equidimensional spaces:
\[
(x,u)
\,\,\,\longmapsto\,\,\,
\big(X(x,u),\,U(x,u)\big)
\,=\,
(X,U),
\]
the target space $\R_{X,U}^{p+q}$ being equipped with similar
coordinates $X \in \R^p$ and $U \in \R^q$.
In later paragraphs, we will employ the notation
$(y, v)$ instead of $(X,U)$\,\,---\,\,clearer here as it indicates
parallels and formal analogies.

We must also introduce
functional symbols for the first $p$ components and
the last $q$ components of such a diffeomorphism:
\[
X^i
\,=\,
\varphi^i(x,u),
\ \ \ \ \ \ \ \ \ \ \ \ \ \ \ \ \ \ \ \
U^\alpha
\,=\,
\psi^\alpha(x,u)
\eqno
{\scriptstyle{(1\,\leqslant\,i\,\leqslant\,p,\,\,
1\,\leqslant\,\alpha\,\leqslant\,q)}}.
\] 
But at the end, we will come back to the more parallel notation:
\[
X^i
\,=\,
X^i(x,u),
\ \ \ \ \ \ \ \ \ \ \ \ \ \ \ \ \ \ \ \
U^\alpha
\,=\,
U^\alpha(x,u)
\eqno
{\scriptstyle{(1\,\leqslant\,i\,\leqslant\,p,\,\,
1\,\leqslant\,\alpha\,\leqslant\,q)}}.
\] 

The splitting $\R_x^p \times \R_u^q$ is motivated by the fact that we
are interested in submanifolds $\big\{ u = f(x) \big\}$ that are {\em
graphs} of maps $\R^p \ni x \longmapsto f(x) \in \R^q$.
Similar graphs $\big\{ U = F(X) \big\}$ exist in the target space
$\R_X^p \times \R_U^q$.  We are mainly interested in how source graphs
$\big\{ u = f (x) \big\}$ are transformed into target graphs $\big\{ U
= F(X) \big\}$, provided our diffeomorphism: 
\[
(x,u) 
\,\,\,\longmapsto\,\,\,
\big(
\varphi(x,u),\,\psi(x,u) 
\big)
\]
is not too far from the identity map, and especially, does not
"rotate" too much tangent directions to our graph $\big\{ u = f(x)
\big\}$.

More precisely, we assume that on restriction to the graph
$\big\{ u = f(x) \big\}$, the first $p$ coordinates 
$\varphi^1, \dots, \varphi^p$ or the diffeomorphism 
when equated to $p$ target independent horizontal coordinates
$X^1, \dots, X^p$:
\[
X
\,=\,
\varphi
\big(x,f(x)\big)
\,\,\,\longleftrightarrow\,\,\,
\Upsilon(X)
\,=\,
x,
\]
can be solved by means of the implicit function theorem.
Equivalently, we are assuming that the Jacobian determinant
of the map $x \longmapsto \varphi(x, f(x) \big)$ between
equidimensional $\R^p \longrightarrow \R^p$
is nowhere vanishing:
\[
\aligned
0
&
\,\neq\,
\det\,
\bigg(
\frac{\partial}{\partial x^{i}}\,
\Big[
\varphi^{j}
\big(x,\,f(x)\big)
\Big]
\bigg)_{1\leqslant j\leqslant p\,\,{\sf rows}}^{1\leqslant i\leqslant p\,\,
{\sf columns}}
\\
&
\,=\,
\det\,
\bigg(
\frac{\partial\varphi^j}{\partial x^i}
+
\smallsum{1\leqslant\alpha\leqslant q}\,
\frac{\partial f^\alpha}{\partial x^{i}}\,
\frac{\partial\varphi^{j}}{\partial u^\alpha}
\bigg)_{1\leqslant j\leqslant p}^{1\leqslant i\leqslant p},
\endaligned
\]
so that the implicit function theorem really applies to the
$p$ equations $X^{j} = \varphi^{j} \big( x, f(x) \big)$ to solve
$x = (x^1, \dots, x^p)$ in terms of $X = \big( X^1, \dots, X^p\big)$
by means of a certain map $\Upsilon \colon \R^p
\longrightarrow \R^p$ as written above:
\[
x^i
\,=\,
\Upsilon^i
\big(X^1,\dots,X^p\big)
\eqno
{\scriptstyle{(1\,\leqslant\,i\,\leqslant\,p)}}.
\]

By definition, the map $X \longmapsto
\Upsilon(X)$ satisfies $p$ identities:
\[
X^j
\,\equiv\,
\varphi^j
\Big(
\Upsilon(X),\,
f\big(\Upsilon(X)\big)
\Big)
\eqno
{\scriptstyle{(1\,\leqslant\,j\,\leqslant\,q)}},
\]
which can be differentiated with respect to all the $X^\ell$ for
$1 \leqslant \ell \leqslant p$:
\[
\aligned
\delta_\ell^j
&
\,=\,
\sum_{i=1}^p\,
\frac{\partial\varphi^j}{\partial x^i}\,
\frac{\partial\Upsilon^i}{\partial X^\ell}
+
\sum_{\alpha=1}^q\,
\frac{\partial\varphi^j}{\partial u^\alpha}\,
\sum_{i=1}^p\,
\frac{\partial f^\alpha}{\partial x^i}\,
\frac{\partial\Upsilon^i}{\partial X^\ell}
\\
&
\,=\,
\sum_{i=1}^p\,
\bigg(
\underbrace{
\frac{\partial\varphi^j}{\partial x^i}
+
\smallsum{1\leqslant\alpha\leqslant q}
\frac{\partial f^\alpha}{\partial x^i}\,
\frac{\partial\varphi^j}{\partial u^\alpha}}_{
{\sf recognize}\,\,{\sf Jacobian}\,\,{\sf determinant}}
\bigg)\,
\frac{\partial\Upsilon^i}{\partial X^\ell},
\endaligned
\]
and since the determinant is assumed nonzero, this linear system can
be inverted in order to determine the partial derivatives
$\frac{\partial \Upsilon^i}{\partial X^\ell}$.

Before doing this, in order to rewrite this system
in a more compact and conceptual form, 
let us introduce the $p$ 
{\sl total differentiation operators:}
\[
{\sf D}_{x^i}
\,:=\,
\frac{\partial}{\partial x^i}
+
\sum_{\alpha=1}^q\,
u_i^\alpha\,
\frac{\partial}{\partial u^\alpha}
+
\sum_{\#\,J\geqslant 1}\,
\sum_{\alpha=1}^q\,
u_{J,i}^\alpha\,
\frac{\partial}{\partial u_J^\alpha}
\eqno
{\scriptstyle{(1\,\leqslant\,i\,\leqslant\,p)}}.
\]
These ${\sf D}_{x^i}$ are here written as 
{\em infinite} formal operators,
but anyway, they will always act {\em in a finite truncated way}
on every function encountered, for instance:
\[
{\sf D}_{x^i}
\big(
\varphi^j
\big)
\,=\,
\bigg(
\frac{\partial}{\partial x^{i}}
+
\sum_{\alpha=1}^q\,
u_{i}^\alpha\,
\frac{\partial}{\partial u^\alpha}
+
0
\bigg)\,
\big(
\varphi^j
(x,u)
\big).
\]

On restriction (pullback) to a graph $\big(x, f(x) \big)$, 
the independent jet variables $u_i^\alpha
\longleftrightarrow \frac{\partial
f^\alpha}{\partial x^i}(x)$
become of course
partial derivatives of the graphing functions
$f^\alpha(x^1, \dots, x^n)$, and therefore
the above linear system rewrites compactly as:
\[
\delta_\ell^j
\,=\,
\sum_{i=1}^p\,
{\sf D}_{x^i}
\big(\varphi^j\big)\,
\frac{\partial\Upsilon^i}{\partial X^\ell}
\eqno
{\scriptstyle{(1\,\leqslant\,j,\,\ell\,\leqslant\,p)}}.
\]

For instance, in the case where $p = 1$:
\[
1
\,=\,
{\sf D}_x(\varphi)\,
\Upsilon_X,
\]
we can solve:
\[
\Upsilon_X
\,=\,
\frac{1}{{\sf D}_x(\varphi)},
\]
while an unpleasant matrix inversion will be required as soon as
$p \geqslant 2$.

Next, in the simplest case where $p = 1 = q$, 
the graph $\big(x, f(x) \big)$ is transformed into the
target graph:
\[
\aligned
U
&
\,=\,
\psi
\big(x,f(x)\big)
\\
&
\,=\,
\psi
\Big(
\Upsilon(X),\,
f\big(\Upsilon(X)\big)
\Big)
\\
&
\,=:\,
F(X),
\endaligned
\]
whose tangent directions are obtained by plain derivation:
\[
\aligned
F_X
&
\,=\,
\big(
\psi_x
+
f_x\,\psi_u
\big)\,
\Upsilon_X
\\
&
\,=\,
{\sf D}_x(\psi)\,
\Upsilon_X
\\
&
\,=\,
\frac{{\sf D}_x(\psi)}{{\sf D}_x(\varphi)}.
\endaligned
\]
Forgetting the functions $f(x)$ and $F(X)$ which represent
graphs, and replacing their first-order derivatives
$f_x$ and $F_X$ by our independent jet variables $u_x$ and
$U_X$, we have finished to present the basic

\begin{Theorem}
{\rm
{\cite{Engel-Lie-1888, Lie-Merker-2015, Bluman-Kumei-1989,
Olver-1993, Olver-1995, Merker-2008}}}
When $p = 1 = q$, every 
diffeomorphism from $\R_x^1 \times \R_u^1$ to $\R_X^1 \times \R_U^1$:
\[
\phi 
\colon\ \ \
(x,u)
\,\,\,\longmapsto\,\,\, 
\big(X(x,u), U(x,u)\big)
\]
possesses a lift as a diffeomorphism
\[
\phi^{(1)} 
\colon\ \ \
(x,u,u_x)
\,\,\,\longmapsto\,\,\, 
\big(X(x,u), U(x,u),U_X(x,u,u_x)\big)
\]
from the open subset $GJ^1_{x,u}$ (graphed jets) of $J^1_{x,u}$ defined by
\[
0
\,\neq\,
{\sf D}_x 
\big(X\big)
\,=\,
X_x
+
u_x\,X_u,
\]
\begin{enumerate}
\item (automatic) making commutative the diagram:
\[
\xymatrix{
GJ_{x,u}^1
\ar[rr]^{\phi^{(1)}}
\ar[d]
&
&
J_{X,U}^1
\ar[d]
\\
\R_{x,u}^{1+1}
\ar[rr]^{\phi}
&
&
\R_{X,U}^{1+1},
}
\]
\item preserving the contact ideal, i.e. $dU-U_X\,dX=(*)\,\big(du-u_x\,dx\big)$, where $(*)$ is a nowhere zero function,
\end{enumerate}
uniquely determined by the formula:
\[
U_X
\,=\,
\frac{{\sf D}_x(U)}{{\sf D}_x(X)}
\,=\,
\frac{U_x+u_x\,U_u}{X_x+u_x\,X_u}.
\]

For higher jet orders:
\[
U_{XX}
\,=\,
\frac{1}{{\sf D}_x(X)}\,
{\sf D}_x\big(U_X\big),
\ \ \ \ \ \ \ \ \ \ \ \
U_{XXX}
\,=\,
\frac{1}{{\sf D}_x(X)}\,
{\sf D}_x\big(U_{XX}\big),
\ \ \ \ \ \ \ \ \ \ \ \
\cdots\cdots.
\eqno\qed
\]
\end{Theorem}

We do not explain how iterations of the first-order
jet formula 
$U_X = \frac{{\sf D}_x(U)}{{\sf D}_x(X)}$
provide higher order jets formulas.
Once the case $p = q = 1$ has been understood, the technicalities
concerning higher dimensional
cases (involving many indices) become more accessible.

Indeed, coming back to the compact linear
system for the partial
derivatives of the $\Upsilon^i$ left above, 
using the notation $X^j(x,u)$ instead of $\varphi^j(x,u)$,
we see that this system can be rewritten
in matrix form as:
\[
\left(\!
\begin{array}{cccc}
1 & 0 & \cdots & 0
\\
0 & 1 & \cdots & 0
\\
\vdots & \vdots & \ddots & \vdots
\\
0 & 0 & \cdots & 1
\end{array}
\!\right)
\,\,=\,\,
\left(\!
\begin{array}{cccc}
\Upsilon_{X^1}^1 & \Upsilon_{X^1}^2 & \cdots & \Upsilon_{X^1}^p
\\
\Upsilon_{X^2}^1 & \Upsilon_{X^2}^2 & \cdots & \Upsilon_{X^2}^p
\\
\vdots & \vdots & \ddots & \vdots
\\
\Upsilon_{X^p}^1 & \Upsilon_{X^p}^2 & \cdots & \Upsilon_{X^p}^p
\end{array}
\!\right)
\cdot
\left(\!
\begin{array}{cccc}
{\sf D}_{x^1}(X^1) & {\sf D}_{x^1}(X^2) & \cdots & {\sf D}_{x^1}(X^p)
\\
{\sf D}_{x^2}(X^1) & {\sf D}_{x^2}(X^2) & \cdots & {\sf D}_{x^2}(X^p)
\\
\vdots & \vdots & \ddots & \vdots
\\
{\sf D}_{x^p}(X^1) & {\sf D}_{x^p}(X^2) & \cdots & {\sf D}_{x^p}(X^p)
\end{array}
\!\right),
\]
which means that:
\[
\left(\!
\begin{array}{ccc}
\Upsilon_{X^1}^1 & \cdots & \Upsilon_{X^1}^p
\\
\vdots & \ddots & \vdots
\\
\Upsilon_{X^p}^1 & \cdots & \Upsilon_{X^p}^p
\end{array}
\!\right)
\,\,=\,\,
\left(\!
\begin{array}{ccc}
{\sf D}_{x^1}(X^1) & \cdots & {\sf D}_{x^1}(X^p)
\\
\vdots & \ddots & \vdots
\\
{\sf D}_{x^p}(X^1) & \cdots & {\sf D}_{x^p}(X^p)
\end{array}
\!\right)^{\!-1}.
\]

We also remind that the nonvanishing of the Jacobian determinant of
the map $x \longmapsto X\big( x, f(x) \big)$ can be re-expressed after
a transposition as the invertibility of this matrix:
\[
0
\,\neq\,
\det\,
\left(\!
\begin{array}{ccc}
{\sf D}_{x^1}(X^1) & \cdots & {\sf D}_{x^1}(X^p)
\\
\vdots & \ddots & \vdots
\\
{\sf D}_{x^p}(X^1) & \cdots & {\sf D}_{x^p}(X^p)
\end{array}
\!\right)
\bigg\vert_{
u=f(x),\,\,
u_{x^j}=f_{x^j}(x)}.
\]

The geometric meaning of our assumption that graphs are transformed
into graphs is that the composition of three maps: lifting to the
graph; performing the diffeomorphism; projecting horizontally:
\[
\xymatrix{
(x,f(x))
\ar[rr]^{\phi\,\,\,\,\,\,\,\,\,\,\,\,\,\,\,\,\,\,\,\,\,\,\,\,\,\,\,\,\,\,}
&
&
\big(X(x,f(x)),U(x,f(x))\big)
\ar[d]
\\
x
\ar[u]
&
&
X(x,f(x)),
}
\]
is invertible. 

Similarly as we did in the simple case $p = 1 = q$,
we pass the graph $\big(x, f(x) \big)$ to the target space,
and we re-express the result in terms of $X$:
\reqnomode\usetagform{EngelLie}
\begin{align}
U^\alpha
&
\,=\,
\psi^\alpha\big(x,f(x)\big)
\notag
\\
&
\,=\,
\psi^\alpha
\Big(
\Upsilon(X),\,
f\big(\Upsilon(X)\big)
\Big)
\notag
\\
&
\,=:\,
F^\alpha(X)
\tag{(1\,\leqslant\,\alpha\,\leqslant\,q).}
\end{align}
Next, for fixed $\alpha$, 
we differentiate these equations with respect to 
$X^1, \dots, X^p$:
\[
\aligned
U_{X^1}^\alpha
&
\,=\,
\sum_{i=1}^p\,
\frac{\partial\psi^\alpha}{\partial x^i}\,
\frac{\partial\Upsilon^i}{\partial X^1}
+
\sum_{i=1}^p\,
\sum_{\beta=1}^q\,
\frac{\partial\psi^\alpha}{\partial u^\beta}\,
\frac{\partial f^\beta}{\partial x^i}\,
\frac{\partial\Upsilon^i}{\partial X^1},
\\
\cdot
\cdots
&
\cdots\cdots\cdots\cdots\cdots\cdots\cdots\cdots\cdots
\cdots\cdots\cdots\cdots
\\
U_{X^p}^\alpha
&
\,=\,
\sum_{i=1}^p\,
\frac{\partial\psi^\alpha}{\partial x^i}\,
\frac{\partial\Upsilon^i}{\partial X^p}
+
\sum_{i=1}^p\,
\sum_{\beta=1}^q\,
\frac{\partial\psi^\alpha}{\partial u^\beta}\,
\frac{\partial f^\beta}{\partial x^i}\,
\frac{\partial\Upsilon^i}{\partial X^p},
\endaligned
\]
that is to say after reorganization:
\[
\aligned
U_{X^1}^\alpha
&
\,=\,
\sum_{i=1}^p\,
\frac{\partial\Upsilon^i}{\partial X^1}\,
\bigg(
\frac{\partial\psi^\alpha}{\partial x^i}
+
\sum_{\beta=1}^q\,
\frac{\partial f^\beta}{\partial x^i}\,
\frac{\partial\psi^\alpha}{\partial u^\beta}
\bigg),
\\
\cdot
\cdots
&
\cdots\cdots\cdots\cdots\cdots\cdots\cdots\cdots\cdots
\cdots\cdots\cdot
\\
U_{X^p}^\alpha
&
\,=\,
\sum_{i=1}^p\,
\frac{\partial\Upsilon^i}{\partial X^p}\,
\bigg(
\frac{\partial\psi^\alpha}{\partial x^i}
+
\sum_{\beta=1}^q\,
\frac{\partial f^\beta}{\partial x^i}\,
\frac{\partial\psi^\alpha}{\partial u^\beta}
\bigg).
\endaligned
\] 
Coming back to the notation $U(x,u)$ 
instead of $\psi(x,u)$, and recognizing
the actions of the total differentiation operators, we get,
still for any fixed $1 \leqslant \alpha \leqslant q$:
\[
\left(\!
\begin{array}{c}
U_{X^1}^\alpha
\\
\vdots
\\
U_{X^p}^\alpha
\end{array}
\!\right)
\,\,=\,\,
\left(\!
\begin{array}{ccc}
\Upsilon_{X^1}^1 & \cdots & \Upsilon_{X^1}^p
\\
\vdots & \ddots & \vdots 
\\
\Upsilon_{X^p}^1 & \cdots & \Upsilon_{X^p}^p
\end{array}
\!\right)\,\,
\left(\!
\begin{array}{c}
{\sf D}_{x^1}(U^\alpha)
\\
\vdots
\\
{\sf D}_{x^p}(U^\alpha)
\end{array}
\!\right).
\]

We can now state the fundamental general theorem, 
recalling how we denote coordinates:
\[
\big(x^i,u^\alpha,u_{x^j}^\beta\big)
\,\in\,
J_{x,u}^1 
\ \ \ \ \ \ \ \ \ \ \ \ \ \ \ \ \ \ \ \
\text{and}
\ \ \ \ \ \ \ \ \ \ \ \ \ \ \ \ \ \ \ \
\big(X^i, U^\alpha, U_{X^j}^\beta\big)
\,\in\,
J_{X,U}^1,
\]
with trivial projections from the first jet spaces onto the
ground manifolds:
\[
\xymatrix{
\big(x^i,u^\alpha,u_{x^j}^\beta\big)
\ar[d]
&
&
\big(X^k, U^\gamma, U_{X^l}^\delta\big)
\ar[d]
\\
(x^i,u^\alpha),
&
&
\big(X^k, U^\gamma\big).
}
\]

\begin{Theorem}
{\rm
{\cite{Lie-Merker-2015, Bluman-Kumei-1989,
Olver-1993, Olver-1995, Merker-2008}}}
\label{Thm-prolongation-diffeomorphism-phi}
For any $p \geqslant 1$ and any $q \geqslant 1$, every 
diffeomorphism: 
\[
\aligned
\phi\colon\ \ \
\R_{x,u}^{p+q}
&
\,\,\,\longrightarrow\,\,\,
\R_{X,U}^{p+q}
\\
(x^i,u^\alpha)
&
\,\,\,\longmapsto\,\,\,
\big(X^k(x^i,u^\alpha),\,U^\gamma(x^i,u^\alpha)\big),
\endaligned
\]
possesses a lift as a diffeomorphism
\[
\phi^{(1)} 
\colon\ \ \
(x^i,u^{\alpha},u^{\beta}_{x^j})
\,\,\,\longmapsto\,\,\, 
\big(X^k(x^i,u^\alpha),\,U^\gamma(x^i,u^\alpha),U^\delta_{X^l}(x^i,u^\alpha,u^\beta_{x^j})\big)
\]
from the open subset $GJ^1_{x,u}$ of $J^1_{x,u}$ defined by
\[
0
\,\neq\,
\det\,
\left(\!
\begin{array}{ccc}
{\sf D}_{x^1}(X^1) & \cdots & {\sf D}_{x^1}(X^p)
\\
\vdots & \ddots & \vdots
\\
{\sf D}_{x^p}(X^1) & \cdots & {\sf D}_{x^p}(X^p)
\end{array}
\!\right),
\]
\begin{enumerate}
\item (automatic) making commutative the diagram:
\[
\xymatrix{
GJ_{x,u}^1
\ar[rr]^{\phi^{(1)}}
\ar[d]
&
&
J_{X,U}^1
\ar[d]
\\
\R_{x,u}^{p+q}
\ar[rr]_\phi
&
&
\R_{X,U}^{p+q},
}
\]
\item preserving the contact ideal, i.e. 
\[
\left(\!
\begin{array}{ccc}
dU^1-U^1_{X^1}\,dX^1
\\
\vdots
\\
dU^q-U^q_{X^p}\,dX^p
\end{array}
\!\right)
\,\,=\,\,
\left(\!
\begin{array}{ccc}
* & \cdots & *
\\
\vdots & \ddots & \vdots
\\
* & \cdots & *
\end{array}
\!\right)
\cdot
\left(\!
\begin{array}{ccc}
du^1-u^1_{x^1}\,dx^1
\\
\vdots
\\
du^q-u^q_{x^p}\,dx^p
\end{array}
\!\right),
\]
where the transition matrix is everywhere invertible,
\end{enumerate}
uniquely determined by the formula:
\[
\left(\!
\begin{array}{c}
U_{X^1}^\alpha
\\
\vdots
\\
U_{X^p}^\alpha
\end{array}
\!\right)
\,\,=\,\,
\left(\!
\begin{array}{ccc}
{\sf D}_{x^1}(X^1) & \cdots & {\sf D}_{x^1}(X^p)
\\
\vdots & \ddots & \vdots
\\
{\sf D}_{x^p}(X^1) & \cdots & {\sf D}_{x^p}(X^p)
\end{array}
\!\right)^{\!-1}\,\,
\left(\!
\begin{array}{c}
{\sf D}_{x^1}(U^\alpha)
\\
\vdots
\\
{\sf D}_{x^p}(U^\alpha)
\end{array}
\!\right).
\]

For higher jet orders, with $J = j_1, \dots, j_\lambda$,
with $1 \leqslant j_1, \dots, j_\lambda \leqslant p$ arbitrary,
and with $1 \leqslant \alpha \leqslant q$:
\[
\left(\!
\begin{array}{c}
U_{X^JX^1}^\alpha
\\
\vdots
\\
U_{X^JX^p}^\alpha
\end{array}
\!\right)
\,\,=\,\,
\left(\!
\begin{array}{ccc}
{\sf D}_{x^1}(X^1) & \cdots & {\sf D}_{x^1}(X^p)
\\
\vdots & \ddots & \vdots
\\
{\sf D}_{x^p}(X^1) & \cdots & {\sf D}_{x^p}(X^p)
\end{array}
\!\right)^{\!-1}\,\,
\left(\!
\begin{array}{c}
{\sf D}_{x^1}(U_{X^J}^\alpha)
\\
\vdots
\\
{\sf D}_{x^p}(U_{X^J}^\alpha)
\end{array}
\!\right).
\eqno\qed
\]
\end{Theorem}

Again, we do not provide explanations on how iterations
of the first-order jet formula written above
conduct to the higher order jet formulas,
leaving the details\,\,---\,\,which can be found
in~{\cite{Bluman-Kumei-1989}}\,\,---\,\,as an exercise 
to the interested reader.

\smallskip

Introducing the modified total differentiation operators:
\[
\left(\!
\begin{array}{c}
{\sf E}_{x^1}
\\
\vdots
\\
{\sf E}_{x^p}
\end{array}
\!\right)
\,\,:=\,\,
\left(\!
\begin{array}{ccc}
{\sf D}_{x^1}(X^1) & \cdots & {\sf D}_{x^1}(X^p)
\\
\vdots & \ddots & \vdots
\\
{\sf D}_{x^p}(X^1) & \cdots & {\sf D}_{x^p}(X^p)
\end{array}
\!\right)^{\!-1}\,\,
\left(\!
\begin{array}{c}
{\sf D}_{x^1}
\\
\vdots
\\
{\sf D}_{x^p}
\end{array}
\!\right),
\]
the first prolongation of the diffeomorphism $\phi$ rewrites as:
\[
\aligned
U_{X^j}^\alpha
\,=\,
{\sf E}_{x^j}
\big(U^\alpha\big).
\endaligned
\]

\begin{Theorem}
\label{Thm-prolongation-E-j}
For all $1 \leqslant \alpha \leqslant q$ and $1 \leqslant j_1, \dots,
j_\lambda \leqslant p$:
\[
U_{X^{j_1}\cdots X^{j_\lambda}}^\alpha
\,\,=\,\,
{\sf E}_{x^{j_1}}
\big(
\cdots
\big(
{\sf E}_{x^{j_\lambda}}
(U^\alpha)
\big)
\cdots
\big).
\eqno\qed
\]
\end{Theorem}

\Subsection{Prolongations of vector fields to jet spaces}
So we agree that every (local) diffeomorphism $\phi \colon (x,u) 
\longmapsto \big( X(x,u), U(x,u) \big)$ has uniquely
determined lifts $\phi^{(n)}$ to the jet spaces of any order $n \geqslant 1$, even if the explicit formulas for the components of
$\phi^{(n)}$ are very unwieldy,
due to the inversion of the matrix ${\sf D}_{x^i} \big( X^j \big)$,
and due to the well known exponential
symbolic swelling of formulas
through iterated differentiations.
Fortunately, as discovered by Lie, linearization 
of $\phi^{(n)}$ sheds new light, and simplifies the formulas.

Indeed, every vector field on the base:
\[
{\sf v}
\,=\,
\sum_{i=1}^p\,
\xi^i(x,u)\,
\frac{\partial}{\partial x^i}
+
\sum_{\alpha=1}^q\,
\varphi^\alpha(x,u)\,
\frac{\partial}{\partial u^\alpha},
\]
produces by integration a $1$-parameter 
group of diffeomorphisms:
\[
\phi_t(x,u)
\,:=\,
\exp\,\big(t\,{\bf v}\big)(x,u)
\eqno
{\scriptstyle{(t\,\in\,\R\,\,\text{small})}},
\]
satisfying by definition:
\[
\frac{d}{dt}
\bigg\vert_{t=0}\,
\phi_t(x,u)
\,=\,
{\bf v}
\big\vert_{(x,u)}.
\]
Of course, each $\phi_t$ possesses lifts to domains of graphed jets
of any order $n \geqslant 1$:
\[
\xymatrix{
GJ_{x,u}^n
\ar[rr]^{\phi_t^{(n)}}
\ar[d]
&
&
J_{X,U}^n
\ar[d],
\\
\R_{x,u}^{p+q}
\ar[rr]_{\phi_t}
&
&
\R_{X,U}^{p+q},
}
\]
and it comes naturally to mind to differentiate such a
lifted $1$-parameter family $\phi_t^{(n)}$ of diffeomorphisms,
obtaining a uniquely determined vector field upstairs:
\[
{\sf v}^{(n)}
\,:=\,
\frac{d}{dt}
\bigg\vert_{t=0}\,
\phi_t^{(n)}
\big(x^i,u^\alpha,u_{x^J}^\beta\big).
\]

Without recalling the explanations which can be found
in the literature, employing again the $p$ total
differentiation operators:
\[
{\sf D}_{x^i}
\,=\,
\frac{\partial}{\partial x^i}
+
\sum_{\alpha=1}^q\,
u_i^\alpha\,
\frac{\partial}{\partial u^\alpha}
+
\sum_{\#\,J\geqslant 1}\,
\sum_{\alpha=1}^q\,
u_{J,i}^\alpha\,
\frac{\partial}{\partial u_J^\alpha}
\eqno
{\scriptstyle{(1\,\leqslant\,i\,\leqslant\,p)}},
\]
we shall admit the following 
very useful theorem which enables
to compute the coefficients $\varphi_J^\alpha$
of these (infinitely) prolonged
vector field:
\[
{\sf v}^{(\infty)}
\,=\,
{\sf v}
+
\sum_{\# J\geqslant 1}\,
\sum_{\alpha=1}^q\,
\varphi_J^\alpha
\big(
x^i,u^\beta,u_{x^K}^\gamma
\big)\,
\frac{\partial}{\partial u_J^\alpha}.
\]

\begin{Theorem}
\label{Thm-prolongation-vector-fields}
{\rm {\cite{Lie-Merker-2015, Olver-1993, Olver-1995, 
Merker-2008}}}
For any $1 \leqslant \alpha \leqslant q$ and every $J = j_1,
\dots, j_\lambda$ with $1 \leqslant j_\nu \leqslant p$, one has:
\[
\varphi_J^\alpha
\,=\,
{\sf D}_{x^J}
\bigg(
\varphi^\alpha
-
\sum_{1\leqslant i\leqslant p}\,
\xi^i\,
u_i^\alpha
\bigg)
+
\sum_{1\leqslant i\leqslant p}\,
\xi^i\,
u_{J,i}^\alpha.
\eqno\qed
\]
\end{Theorem}

Here of course, we define:
\[
{\sf D}_{x^J}
\,:=\,
{\sf D}_{x^{j_1}}
\cdots
{\sf D}_{x^{j_\lambda}},
\]
where the order does not matter, since these total differentiation
operators commute. It is important to point out that, due to the
obvious relation:
\[
{\sf D}_{x^J}
\big(u_i^\alpha)
\,=\,
u_{J,i}^\alpha,
\]
the last $p$ terms $\xi_i\, u_{J,i}^\alpha$ 
in the formula for $\varphi_J^\alpha$
which incorporate jets of order $\# J +1$ 
do in fact disappear, hence all $\varphi_J^\alpha$ 
are functions of jet variables of order $\leqslant \# J$,
as they must be.

The above "direct formulas" for the coefficients $\varphi_J^\alpha$
are quite convenient when working on a computer.
Sometimes, alternative formulas proceeding by recurrence 
happen to be also useful:
\[
\varphi_{J,j}^\alpha
\,=\,
{\sf D}_{x^j}
\big(
\varphi_J^\alpha
\big)
-
\sum_{1\leqslant i\leqslant p}\,
{\sf D}_{x^j}(\xi^i)\,
u_{J,i}^\alpha.
\]
These recurrence formulas were closed and
synthetized in~{\cite[Chap.~II]{Merker-2008}}
to produce explicit formulas for all the $\varphi_J^\alpha$
by means of multiple Kronecker symbols.

\Section{\bf The Moving Frame Method}
\label{moving-frame-method}
\HEAD{{\ref{moving-frame-method}}.~{\sf 
The Moving Frame Method}
}{
Zhangchi {\sc Chen}, Joël~{\sc Merker}}

We now present the {\sl moving frame method},
developped by Olver and his
collaborators (\cite{Fels-Olver-1999, Olver-2000,
Olver-2001, Olver-2007, Olver-2011, Olver-2018}).

\Subsection{Source and target jet coordinates}
\label{Subsec-jet-notations}
We keep the previous notation
for the coordinates in the source space:
\[
z
\,=\,
(x, u)
\,=\,
\big(
x^1,\dots,x^p,\, 
u^1, \dots,u^q \big), 
\]
but instead of $(X, U)$, we will denote
the coordinates in the target space by:
\[
w
\,=\,
(y,v)
\,=\,
\big(
y^1,\dots,y^p,\,
v^1,\dots,v^q
\big).
\]

For any jet order $n 
\geqslant 0$, jet coordinates will accordingly be denoted by:
\[
\aligned
z^{(n)}
&
\,=\,
\big(
x^i,u^\alpha,\,
u_{x^J}^\beta
\big),
\\
w^{(n)}
&
\,=\,
\big(
y^j,v^\gamma,\,
v_{y^K}^\delta
\big),
\endaligned
\]
with $J = j_1, \dots, j_\lambda$, 
with $1 \leqslant \lambda \leqslant n$,
with $1 \leqslant j_1, \dots, j_\lambda \leqslant p$, 
and similarly,
with $K = k_1, \dots, k_\mu$, 
with $1 \leqslant \mu \leqslant n$,
with $1 \leqslant k_1, \dots, k_\mu \leqslant p$.
Sometimes, we will abbreviate:
\[
m
\,:=\,
p+q.
\]

Later, we will need to renumber the coordinates by increasing
length of derivatives:
\[
\big(
w_1^{(n)},
\dots,
w_\NN^{(n)}
\big)
\,=\,
\big(
y^j,v^\gamma,\,
v_{y^{K_1}}^{\delta_1},
\dots,
v_{y^{K_n}}^{\delta_n}
\big),
\]
where the total number $\NN$ is the dimension of the
jet space $J_{p,q}^n$:
\[
\NN
\,=\,
\NN_{p,q}(n)
\,:=\,
p
+
q\,
\binom{p+n}{n},
\]
where $\# K_1 = 1$, \dots, $\# K_n = n$, and where we choose
any ordering for the subcollection of 
jet coordinates 
$v_{y^{K_h}}^{\delta_h}$ with $1 \leqslant \delta_h \leqslant q$ and
$\# K_h = h$. 

The first $p + q$ coordinates $w_1^{(n)}, \dots, w_{p+q}^{(n)}$
belong in fact to the zero\textsuperscript{th} jet order space,
hence should be written $w_1^{(0)}, \dots, w_{p+q}^{(0)}$.
By convention, for any jet order $n' \geqslant n$,
the trivial projection $J_{p,q}^{n'} \longrightarrow J_{p,q}^n$
allows to identify $w_k^{(n')} \equiv w_k^{(n)}$
as soon as $w_k^{(n)}$ is a coordinate of $J_{p,q}^n$.

\Subsection{Prolongations of groups actions}
Next, assume that a finite-dimensional local Lie group $G$ acts
on $\R_{x,u}^{p+q}$, and denote its dimension by:
\[
r
\,:=\,
\dim\,G
\eqno
{\scriptstyle{(1\,\leqslant\,r\,<\,\infty)}}.
\]

Indeed, 
it is convenient (and theoretically wider) to assume that $G$ is a
{\em local} Lie group, as did Lie in~{\cite{Engel-Lie-1888,
Lie-Merker-2015}}. Therefore, the diffeomorphisms $\phi \colon
\R_{x,u}^{p+q} \longrightarrow \R_{y,v}^{p+q}$ we are considering are
`parametrized' by group elements $g \in G$ in some neighborhood of the
identity element:
\[
e
\,\in\,
G,
\]
and local Lie group axioms about composition
and inversion hold. As we will not require background
reminder about foundational concepts of Lie theory, 
the interested reader
will be referred to~{\cite[pp.~29--32]{Lie-Merker-2015}}.

We will not employ functional symbols $\phi = (\varphi, \psi)$
like in Section~{\ref{graph-transformations-jet-spaces}},
but instead, we will denote diffeomorphisms
`parametrized' by group elements $g \in G$ as Lie did:
\[
y
\,=\,
y(g,x,u),
\ \ \ \ \ \ \ \ \ \ \ \ \ \ \ \ \ \ \ \
v
\,=\,
v(g,x,u).
\]
Sometimes, we will also use the more compact notation $w = g \cdot z$
instead of $w = w(g,x,u)$.

\begin{Definition}
The action of $G$ on $\R_{x,u}^{p+q}$ is 
said to be {\sl free} 
if the isotropy subgroup of every $z = (x,u)$ reduces to the identity:
\[
\{e\}
\,=\,
\big\{
g\in G
\colon\,
g\cdot z
\,=\,
z
\big\}
\,=:\,
G_z
\eqno
{\scriptstyle{(\forall\,z\,\in\,\R^m)}}.
\]
\end{Definition}

In general, the action of $G$ on $\R_{x,u}^{p+q}$ is not free, just
because as soon as $r \geqslant m + 1$, isotropy subgroups $G_z$ are
of dimension $\geqslant r - m \geqslant 1$. As soon as the maps $g
\longmapsto w(g,x,u)$ are not of full rank $r$, freeness is lost.

A good substitute is to assume that the ranks of the maps
$g \longmapsto w(g,x,u)$ are {\em constant}. In this case,
after perhaps passing to smaller open subsets, 
the $G$-orbits of the local Lie group $G$ are of constant
dimension, and their union foliates some open set in $\R_{x,u}^{p+q}$.
In the analytic category, group actions are generically
of local constant rank, in appropriate subdomains which have an 
invariant meaning, hence rank constancy is essentially no
assumption in Lie theory.
This principle of thought is explained in 
Chapter~1 of~{\cite{Lie-Merker-2015}}.

All these observations remain valid in jet spaces of any order $n
\geqslant 0$, because thanks to
Theorem~{\ref{Thm-prolongation-diffeomorphism-phi}},
every group-diffeomorphism $z \longmapsto g \cdot z$
(close the identity mapping)
lifts upstairs as a jet diffeomorphism:
\[
\xymatrix{
J_z^n
\ar[rr]^{g\cdot}
\ar[d]
&
&
J_w^n
\ar[d]
\\
\R_z^m
\ar[rr]_{g\cdot}
&
&
\R_w^m,
}
\ \ \ \ \ \ \ \ \ \ \ \ \ \ \ \ \ \ \ \
\xymatrix{
z^{(n)}
\ar[rr]^{g\cdot}
\ar[d]
&
&
w^{(n)}\big(g,z^{(n)}\big)
\ar[d]
\\
z
\ar[rr]_{g\cdot}
&
&
w(g,z),
}
\]
vertical arrows being (trivial) projections, as before.
The advantage of prolongating is that the ambient dimension
increases polynomially with $n$:
\[
\dim\,
J_{p,q}^n
\,=\,
p
+
q\,
\binom{p+n}{n},
\]
and hence, granted that prolongation commutes with projections
from higher order $n' \geqslant n \geqslant 0$ jet spaces:
\[
\xymatrix{
J_z^{n'}
\ar[rr]^{g\cdot}
\ar[d]
&
&
J_w^{n'}
\ar[d]
\\
J_z^n
\ar[rr]^{g\cdot}
\ar[d]
&
&
J_w^n
\ar[d]
\\
\R_z^m
\ar[rr]_{g\cdot}
&
&
\R_w^m,
}
\] 
it is clear that the (generic) ranks of the maps $g \longmapsto
w^{(n)} \big(g, z^{(n)} \big)$ are increasing with $n$.
Of course, the ranks of these maps $g \longmapsto w^{(n)}
\big(g, z^{(n)} \big)$ 
are considered for $g \sim e$ close to
the identity element, and for $z^{(n)}$ in $J_{p,q}^{(n)}$,
or in some appropriate open subsets of 
$J_{p,q}^{(n)}$\,\,---\,\,more will be said soon about that.

We shall use the abbreviation:
\[
\NN
\,:=\,
\NN_{p,q}(n)
\,:=\,
p
+
q\,
{\textstyle{\binom{p+n}{n}}}.
\]
For most of the existing local Lie group actions (and for all the ones
studied in the present article), there always exists a
minimal finite jet order:
\[
0
\,\leqslant\,
n_\GG
\,<\,
\infty,
\]
such that:
\[
\rank
\Big(
g
\,\,\longmapsto\,\,
w^{(n_\GG)}
\big(
g,\,z^{(n_\GG)}
\big)
\Big)
\,=\,
r
\,=\,
\dim\,
G,
\]
and hence the action of the local Lie group $G$ is free.
A necessary\,\,---\,\,and often sufficient\,\,---\,\,condition 
is that:
\[
p
+
q\,
{\textstyle{\binom{p+n_\GG}{n_\GG}}}
\,=\,
\NN_{p,q}
\big(n_\GG\big)
\,\geqslant\,
r.
\]

\begin{Example}
Let the affine group ${\sf A}_2(\R)$ on the plane 
$\R_{x,u}^2$ have general
transformations $y = {\sf a}\, x + {\sf b}\, u + {\sf c}$ and $v =
{\sf k}\, x + {\sf l}\,u + {\sf m}$ with $0 \neq
\left\vert\! \begin{smallmatrix} {\sf a} & {\sf b} 
\\ {\sf k} & {\sf l} \end{smallmatrix} \!\right\vert$.  
With the modified total differentiation operator:
\[
{\sf E}_x
\,:=\,
\frac{1}{{\sf a}+{\sf b}\,u_x}\,
\bigg[
\frac{\partial}{\partial x}
+
u_x\,
\frac{\partial}{\partial u}
+
\sum_{i=1}^\infty\,
u_{x^{i+1}}\,
\frac{\partial}{\partial u_{x^i}}
\bigg],
\]
an application of Theorem~{\ref{Thm-prolongation-E-j}} 
yields the second prolongation
on graphed curves $\big\{ u = u(x) \big\}$ of such
affine transformations:
\[
v_{yy}
\,=\,
{\sf E}_x
\big(
{\sf E}_x
(v)
\big)
\,=\,
\frac{{\sf a}{\sf l}-{\sf b}{\sf k}}{
({\sf a}+{\sf b}\,u_x)^3}\,\,
u_{xx}.
\]

We readily see that the condition $u_{xx} \neq 0$ is affinely
invariant, hence the second order jet space stratifies ${\sf
A}_2(\R)$-invariantly as:
\[
J_{x,u}^2
\,=\,
\big\{
u_{xx}
=
0
\big\}
\cup
\big\{
u_{xx}
\neq
0
\big\}.
\]

Of course, the flatness condition $u_{xx}(x) \equiv 0$ which means
that the graph $u(x) = \lambda\, x + \mu$ being a straight line
is affinely invariant!
\end{Example}

\Subsection{Foliated $G$-actions and differential invariants}
For general actions of finite-dimensional Lie groups on graphs
in $\R_{x,u}^{p+q}$, a case-by-case study is required to determine
the stratification of jet spaces $J_{x,u}^n$ in appropriate
invariant pieces wherein the lifted action of $G$ has 
constant rank properties.

In later sections, we will conduct a precise stratification
study for parabolic surfaces $S^2 \subset \R^3$.

\begin{Definition}
A lifted action of a local Lie group $G$ to a jet space
$J_{p,q}^n$ is called {\sl foliated} at a point 
$z_0^{(n)} \in J_{p,q}^n$ when all $G$-orbits have constant
dimension equal to a certain integer $s$ with
$0 \leqslant s \leqslant r = \dim\, G$
in a neighborhood of $z_0^{(n)}$,
with maps $g \longmapsto w^{(n)} \big(g, z^{(n)} \big)$
having constant rank $s$ for all $z^{(n)}$ near $z_0^{(n)}$
and all $g$ near $e \in G$.
\end{Definition}

Under this assumption, 
the rank theorem guarantees that 
the collection of $G$-orbits indeed
constitutes a local foliation by $s$-dimensional manifolds 
in a neighborhood of $z_0^{(n)}$.

Most of the time, for any action of a local Lie group $G$ on graphs
$\{ u = u(x) \}$ in $\R_{x,u}^{p+q}$, the induced action of $G$ on
some jet space $J_{x,u}^n$ (even of low order) becomes almost
everywhere {\em foliated}, although not necessarily free.

\begin{center}
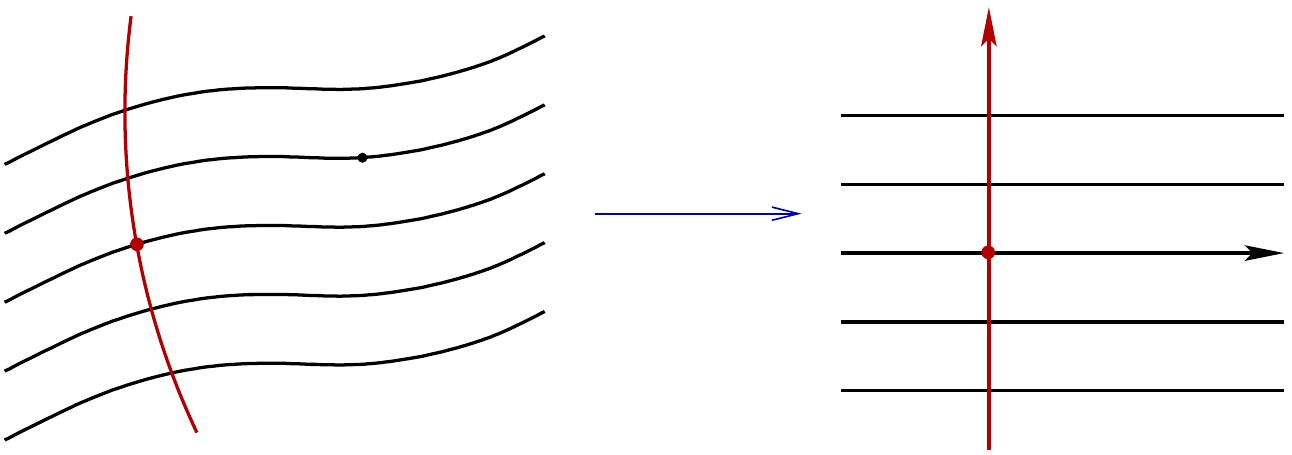
\end{center}

Assuming the action of $G$ is foliated at some $z_0^{(n)}
\in J_{p,q}^n$, let us choose a {\sl transversal} 
$T^n$ to the $G$-orbits at $z_0^{(n)}$,
namely a local submanifold:
\[
z_0^{(n)}
\,\in\,
T^n
\,\subset\,
J_{p,q}^n
\]
of dimension complementary to the dimension of $G$-orbits:
\[
\dim\,T^n
\,=\,
\NN_{p,q}(n)
-
s,
\]
which is also {\sl transversal} to the $G$-orbit 
$G\cdot z_0^{(n)}$ in the sense of transversality theory:
\[
T_{z_0^{(n)}}
T^n
\oplus
T_{z_0^{(n)}}
G\cdot z_0^{(n)}
\,=\,
T_{z_0^{(n)}}
J_{p,q}^n.
\]

After a local straightening diffeomorphism, we can make $G$-orbits
horizontal, directed by $\R^s \times \{0\}$, and the
transversal $T^n$ 
vertical as well, directed by $\{0\} \times \R^{\NN-s}$.
This means that we can get $\NN-s$ functions:
\[
\Iaux_1\big(z^{(n)}\big),\,
\dots\dots,
\Iaux_{\NN-s}\big(z^{(n)}\big),
\]
vanishing at $z_0^{(n)}$, such that the union of $G$-orbits 
near $z_0^{(n)}$ is represented as:
\[
\bigcup_{c_1,\dots,c_{\NN-s}}\,
\big\{
\Iaux_1=c_1,
\dots,
\Iaux_{\NN-s}=c_{\NN-s}
\big\},
\]
with arbitrary constants $c_1, \dots, c_{\NN-s}$ all close to $0$.
In general, producing such functions $\Iaux_1, \dots,
\Iaux_{\NN-s}$ requires an application of the implicit function
theorem. Hence explicitness can be lost.

Certainly, such
$G$-invariants $\Iaux_1, \dots, \Iaux_{\NN-s}$ are
{\sl functionally independent} in the sense that:
\[
\mathcal{R}
\big(
\Iaux_1(z^{(n)}),
\dots,
\Iaux_{\NN-s}(z^{(n)})
\big)
\,\equiv\,
0
\ \ \ \ \ \ \ \ \ \ \ \ \ \ \ \ \ \ \ \
\Longrightarrow
\ \ \ \ \ \ \ \ \ \ \ \ \ \ \ \ \ \ \ \
\mathcal{R}
\,=\,
0,
\]
because this 
property is trivially true after straightening,
and is invariant under diffeomorphisms.
Equivalently, the differentials 
of these invariants are linearly independent:
\[
0
\,\neq\,
d\Iaux_1
\big(z^{(n)}\big)
\wedge\cdots\wedge
d\Iaux_{\NN-s}
\big(z^{(n)}\big)
\eqno
{\scriptstyle{(\forall\,z^{(n)}\,\text{\rm near}\,
z_0^{(n)})}}.
\]

\begin{Definition}
A (local) {\sl differential invariant} $\Iaux \colon J_{p,q}^n 
\longrightarrow \R$ is a function
defined in a neighborhood of some $z_0^{(n)} \in J_{p,q}^n$ 
which is constant on
(local) $G$-orbits:
\[
\Iaux
\big(g\cdot z^{(n)}\big)
\,=\,
\Iaux\big(z^{(n)}\big)
\eqno
{\scriptstyle{(\forall\,g\,\in\,G)}}.
\]
\end{Definition}

Of course, $G$ is a local Lie group as before.  In true applications,
as we will soon {\em e.g.}  in the context of parabolic surfaces
$S^2 \subset \R^3$, the behavior of differential invariants, their
number, their syzygies, will be constant in certain Zariski open
subsets of the jet spaces $J_{p,q}^n$, for every $n \geqslant 0$,
and it is a part of an adequate mathematical work
to determine appropriate invariant stratifications of 
these $J_{p,q}^n$. 

\begin{Theorem}
\label{Thm-invariants-I-1-I-N-s}
{\rm {\cite{Engel-Lie-1888, Olver-1993, Olver-1995, Lie-Merker-2015}}}
Given an action on $\R_x^p \times \R_u^q$ of a local
Lie group $G$, for any jet order $n \geqslant 0$, 
given the induced action of $G$ on 
the space $J_{p,q}^n$ of jets of
maps $x \longmapsto u(x)$, 
with $\NN = \NN_{p,q} (n) = \dim\, J_{p,q}^n$, 
at any point $z_0^{(n)} \in J_{p,q}^n$ 
at which the induced $G$-action is foliated
of rank $s$, 
with $0 \leqslant s \leqslant r$ being the dimension
of leaves, there are exactly $\NN_{p,q}(n) - s$ 
functionally independent differential invariants $\Iaux_1, \dots, 
\Iaux_{\NN-s}$ defined near $z_0^{(n)}$ so that any
other differential invariant $\Iaux$ is a certain,
uniquely defined, function of these:
\[
\Iaux
\,=\,
\mathcal{F}
\big(
\Iaux_1,
\dots,
\Iaux_{\NN-s}
\big).
\eqno\qed
\]
\end{Theorem}

\Subsection{Free actions and moving frames}
Assume now that $n \geqslant 0$ is large enough so that the induced
action of $G$ on $J_{p,q}^n$ is locally free near some $z_0^{(n)} \in
J_{p,q}^{(n)}$, and also foliated. As before, denote by $n_\GG$ the
minimal such $n$. Presently, we will work only with jet orders $n
\geqslant n_\GG$.

Thus, all $G$-orbits on $J_{p,q}^n$ have maximal possible dimension $r
= \dim\, G$. Thanks to the preceding
Theorem~{\ref{Thm-invariants-I-1-I-N-s}}, we know that there are
exactly $\NN_{p,q}(n) - r$ functionally independent differential
invariants.

We now take $n := n_\GG$ minimal possible,
so the number of independent 
differential invariants is $\NN_{p,q} (n_\GG) - r$. 

\begin{center}
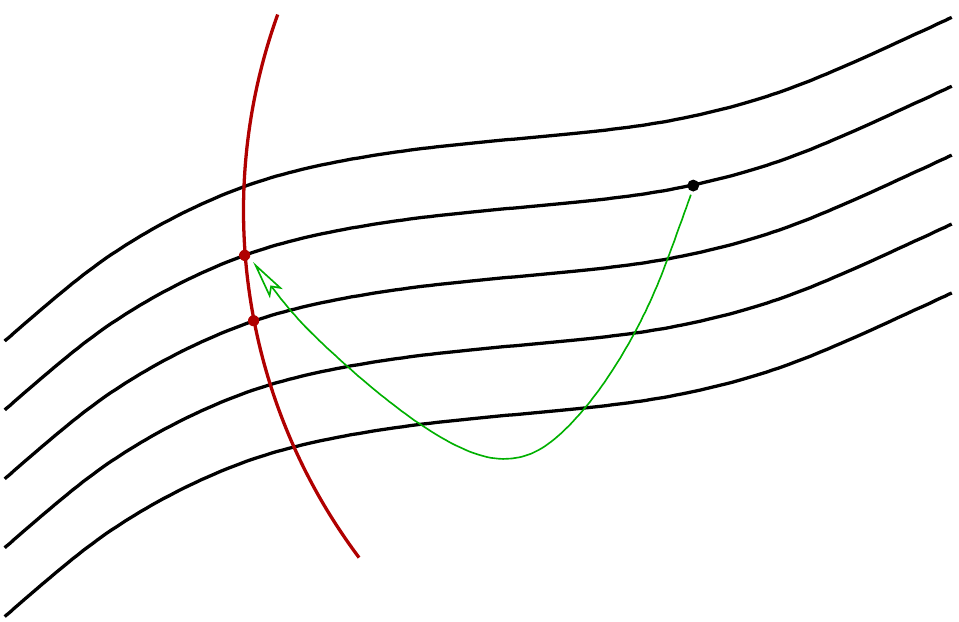
\end{center}

The geometrically evident 
fact that every local $G$-orbit intersects the transversal 
$T^{n_\GG} \subset J_{p,q}^{n_\GG}$ 
in a single point explains the:

\begin{Lemma}
For every $z^{(n_\GG)}$ near $z_0^{(n_\GG)}$, there exists a unique
group element $g \in G$ near $e$ such that:
\[
g
\cdot
z^{(n_\GG)}
\,\in\,
T^{n_\GG}.
\eqno\qed
\]
\end{Lemma}

As in~{\cite{Fels-Olver-1999}}, we will denote this unique group
element by $\rho \big( z^{(n_\GG)} \big)$, so that:
\[
\rho\big(z^{(n_\GG)}\big)
\cdot
z^{(n_\GG)}
\,\in\,
T^{n_\GG}.
\]

We can be more specific about how the map 
$\rho \colon J_{p,q}^{n_\GG} \longmapsto G$ can be
handled in concrete situations. 
Remind from Subsection~{\ref{Subsec-jet-notations}} 
our notation for coordinates on $J_{p,q}^{n_\GG}$:
\[
\big(
w_1^{(n_\GG)},
\dots,
w_{\NN_{p,q}(n_\GG)}^{(n_\GG)}
\big)
\,=\,
\big(
y^j,v^\gamma,\,
v_{y^{K_1}}^{\delta_1},
\dots,
v_{y^{K_{n_\GG}}}^{\delta_{n_\GG}}
\big).
\]
A transversal to the $G$-orbits is usually constructed
after a detailed study of specific features of
the group action. 

\begin{center}
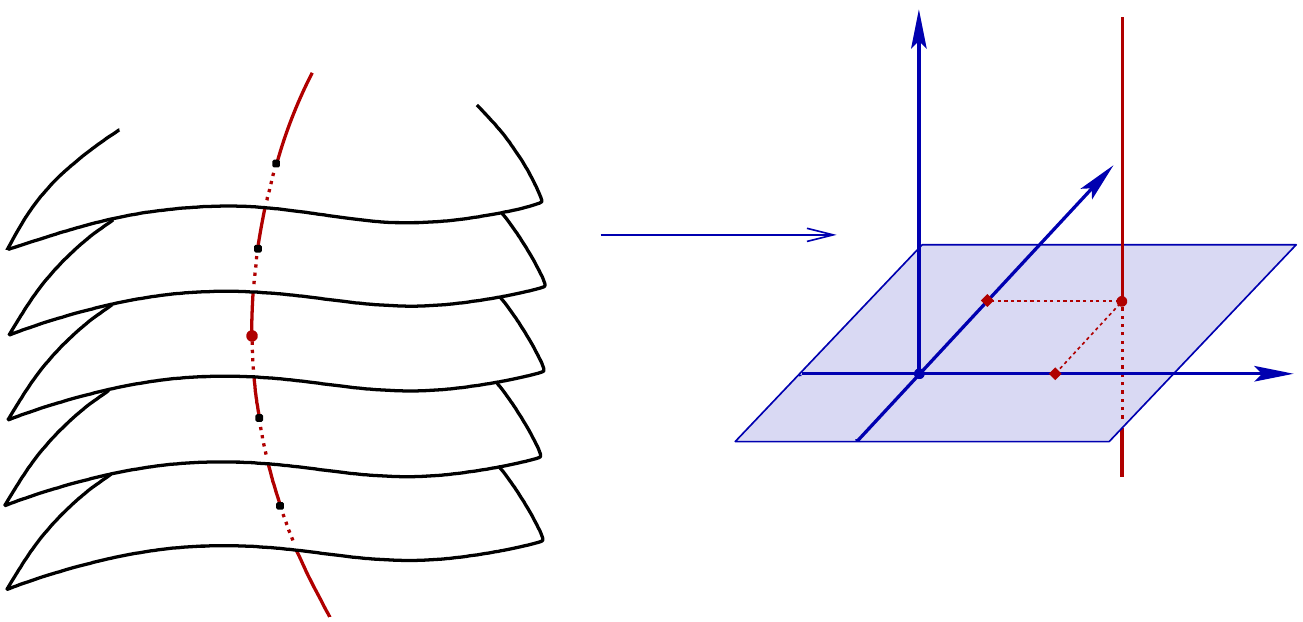
\end{center}

Generally, a number $r = \dim\, G$ of integers:
\[
1
\,\leqslant\,
\nu_1
\,<\,\cdots\,<\,
\nu_h
\,<\,\cdots\,<\,
\nu_r
\,\leqslant\,
\NN,
\]
are then `selected' in some
`natural' way, and constants $c_1, \dots, c_r \in \R$\,\,---\,\,often
equal to $0$ or $1$\,\,---\,\,are also `chosen' so that, viewed in the
target jet space as on the right
hand side of the figure above, the transversal has defining
equations of the shape:
\[
w_{\nu_1}^{(n_1)}
\,=\,
c_1,\,\,
\dots\dots,\,\,
w_{\nu_h}^{(n_h)}
\,=\,
c_h,
\dots\dots,\,\,
w_{\nu_r}^{(n_r)}
\,=\,
c_r,
\]
where for $1 \leqslant h \leqslant r$,
the integer $(n_h)$ is the minimal jet space order to which 
the jet coordinate
$w_{\nu_h}^{(n_\GG)} = w_{\nu_h}^{(n_h)}$ belongs,
and where we have $n_r = n_\GG$, by definition of $n_\GG$. 

Viewed in the source space as on the 
left of the figure above, the equations of the transversal 
$T^{n_\GG} \subset J_{p,q}^{n_\GG}$ 
then read by understanding that the $w_h^{(n_h)}$ are
coordinates of the prolongation of the $G$-action:
\[
w_{\nu_1}^{(n_1)}
\big(g,z^{(n_1)}\big)
\,=\,
c_1,\,\,\,
\dots\dots\dots,\,\,\,
w_{\nu_r}^{(n_r)}
\big(g,z^{(n_r)}\big)
\,=\,
c_r,
\]
and the assumption that $T^{n_\GG}$ is a {\em transversal}
is equivalent to the fact that {\em all the group parameters}
$g = (g_1, \dots, g_r)$ can be solved from these $r$ equations
by means of the implicit function theorem:
\[
g
\,=\,
\rho
\big(
z^{(n_1)},\dots,z^{(n_\GG)}
\big)
\,\equiv\,
\rho
\big(
z^{(n_\GG)}
\big).
\]

For every $n \geqslant n_\GG$, a transversal
$T^n \subset J_{p,q}^n$ can be chosen to have exactly 
the {\em same} equations as a transversal 
$T^{n_\GG} \subset J_{p,q}^{n_\GG}$ in the minimal order $n_\GG$
jet space as written above.
For later applications, {\em e.g.} to parabolic
surfaces $S^2 \subset \R^3$, 
we want to keep memory that the map: 
\[
z^{(n)}
\,\longmapsto\,
\rho\big(z^{(n)}\big)
\,=\,
\rho\big(z^{(n_\GG)}\big)
\]
only depends on jets of this minimal possible order $n_\GG$.
We can now formulate a definition valid for any $n \geqslant n_\GG$.
\[
\xymatrix{
J_{x,u}^n
\ar@/^5pc/[rrrrrrrd]^{\rho}
\ar[rr]_{g\cdot}
\ar[dd]
& &
J_{y,v}^n
\ar@/^1pc/[rrrd]^{\rho}
\ar[dd]
& & &
\\
& & & & &
G
& &
G
\ar[ll]^{\cdot g^{-1}}
\\
\R_{x,u}^{p+q}
\ar[rr]^{g\cdot}
& &
\R_{y,v}^{p+q}
& & &
}
\]

\begin{Definition}
{\rm {\cite{Fels-Olver-1998, Fels-Olver-1999}}}
A (local) {\sl right moving frame} for the action
on $\R_{x,u}^{p+q}$ of a local Lie
group $G$ lifted to an $n$\textsuperscript{th}
jet space $J_{p,q}^n$ is a map:
\[
\rho
\colon\ \ \
J_{p,q}^n
\,\longrightarrow\,
G,
\]
defined near some jet $z_0^{(n)} \in J_{p,q}^n$ 
which satisfies the right equivariancy rule:
\[
\rho
\big(
g\cdot z^{(n)}
\big)
\,=\,
\rho\big(z^{(n)}\big)
\cdot
g^{-1},
\]
for all $g \in G$ near $e$ and all $z^{(n)}$ near
$z_0^{(n)}$.
\end{Definition}

\[
\xymatrix{
z^{(n)}
\ar@/^7pc/[rrrrrrrrdd]^{\rho}
\ar[rr]_{g\cdot}
& &
g\cdot z^{(n)}
\ar@/^1pc/[rrrd]^{\rho}
& & &
\\
& & & & &
\rho\big(g\cdot z^{(n)}\big)
\ar@{=}[d]
& & &
\\
& &
& & &
\rho(z^{(n)})\cdot g^{-1}
& & &
\rho(z^{(n)})
\ar[lll]^{\cdot g^{-1}}
}
\]

\begin{Lemma}
{\rm {\cite[Sec.~4]{Fels-Olver-1999}}}
For any choice of a transveral $T^{n_\GG} \subset J_{p,q}^{n_\GG}$
as above, and for any $n \geqslant n_\GG$, 
the map constructed above by means of the implicit function
theorem:
\[
\rho\big(z^{(n)}\big)
\,=\,
\rho\big(z^{(n_\GG)}\big)
\]
automatically satisfies the right equivariancy rule.\qed
\end{Lemma}

The proof, short, will be skipped here,
as well as the proof of the next elementary

\begin{Theorem}
\label{Thm-existence-moving-frame}
{\rm {\cite[Sec.~4]{Fels-Olver-1999}}}
If the action of $G$ lifted to $J_{p,q}^n$ is foliated at a point
$z_0^{(n)}$, the following three conditions are equivalent:

\smallskip\noindent{\bf (i)}\,
a moving frame exists in a neighborhood of $z_0^{(n)}$;

\smallskip\noindent{\bf (ii)}\, 
$G$ acts locally freely near $z_0^{(n)}$;

\smallskip\noindent{\bf (iii)}\,
$G$-orbits have maximal dimension $r = \dim\, G$ near 
$z_0^{(n)}$.\qed
\end{Theorem}

When studying parabolic surfaces $S^2 \subset \R^3$,
with $G = \SA_3(\R) = \SL_3(\R) \ltimes \R^3$ 
having dimension $11$,
we will see that 
in some appropriately `truncated' $4$\textsuperscript{th} 
jet space 
having {\em also} dimension $11$,
the $G$-orbits have constant
dimension $10 = 11 - 1$, whence 
Theorem~{\ref{Thm-invariants-I-1-I-N-s}}
guarantees that there exists {\em one} differential
invariant of order $4$, which we will call
$\Waux$.
But then because the rank of the $G$-action
degenerates by $1 = 11 - 10$ dimension,
the above 
Theorem~{\ref{Thm-existence-moving-frame}} 
does {\em not} apply. 

\begin{Question}
{\sl 
Is there an appropriate substitute to a moving frame
when $G$-orbits have constant dimension 
$0 \leqslant s \leqslant r-1$ smaller than $\dim\, G$?
}
\end{Question}

One can think of the {\sl isotropy subgroup} $G_{z_0^{(n)}} \subset G$
of a reference point $z_0^{(n)}$, but then, in general,
this subgroup $H := G_{z_0^{(n)}}$ is not distinguished in $G$, 
only left cosets $G \big/ H$ or right cosets
$H \big\backslash G$ can be dealt with.

\smallskip

In presence of a moving frame, 
Theorem~{\ref{Thm-invariants-I-1-I-N-s}} becomes the

\begin{Theorem}
{\rm {\cite[Sec.~4]{Fels-Olver-1999}}}
If $\rho \colon J_{p,q}^n \longrightarrow G$ is a right moving frame,
with $n \geqslant n_\GG$, then in coordinates, all components
of the map:
\[
z^{(n)}
\,\,\longmapsto\,\,
\rho\big(z^{(n_\GG)}\big)
\cdot
z^{(n)}
\]
constitute a complete generating set of differential invariants of order at most
$n$.\qed
\end{Theorem}

Concretely, after having solved the group parameters $g = (g_1, \dots,
g_r)$ from the $r$ equations 
$w_{\nu_h} = c_h$, $1 \leqslant h \leqslant r$, 
of some transversal $T^{n_\GG} \subset
J_{p,q}^{n_\GG}$, one {\em replaces} these solutions $g = \rho \big(
z^{(n_\GG)} \big)$ in {\em all the other} 
coordinate formulas $w_\nu \big( g,
z^{(n)} \big)$ for the prolonged $G$-action,
and one obtains a generating set of 
$\NN_{p,q}(n) - r$ differential invariants of order at most $n$:
\[
I_\nu
\,:=\,
w_\nu
\big(
\rho(z^{n_\GG}),\,
z^{(n)}
\big)
\eqno
{\scriptstyle{(1\,\leqslant\,\nu\,\leqslant\,\NN_{p,q}(n),\,\,\,
\nu\,\neq\,\nu_1,\dots,\nu_r)}},
\]
When $n \geqslant 1 + n_\GG$, this means that {\em all} 
indices:
\[
\NN_{p,q}(n_\GG)
+
1
\,\,\leqslant\,\,
\nu
\,\,\leqslant\,\,
\NN_{p,q}(n)
\]
are concerned by such a replacement.

By definition of how the implicit function theorem
applies to the equations $w_{\nu_h} = c_h$,
if one replaces $g = \rho\big( z^{(n_\GG)} \big)$ in the coordinates
$w_{\nu_1}, \dots, w_{\nu_r}$ which were used to construct the moving
frame, one gets trivial constants:
\[
w_{\nu_1}
\big(
\rho(z^{(n_\GG)}),\,
z^{(n_1)}
\big)
\,\equiv\,
c_1,\,\,
\dots\dots\dots,\,\,
w_{\nu_r}
\big(
\rho(z^{(n_\GG)}),\,
z^{(n_r)}
\big)
\,\equiv\,
c_r.
\]

\begin{Terminology}
{\rm {\cite{Fels-Olver-1999}}}
These $r$ differential quantities 
$w_{\nu_h} \big( \rho(z^{(n_\GG)}), z^{(n_h)} \big) \equiv c_h$
are called {\sl phantom differential invariants}.
\end{Terminology}

These ghost objects will be very useful later.

\Section{\bf What it Really Means to Be a Differential Invariant}
\label{meaning-differential-invariant}
\HEAD{{\ref{meaning-differential-invariant}}.~{\sf What it Really 
Means to Be a Differential Invariant}
}{
Zhangchi {\sc Chen}, Joël~{\sc Merker}}

As before, let $G$ be a 
local Lie group acting on graphs
$\big\{ u = {\tt u} (x) \big\}$ in $\R_x^p \times \R_u^q$,
of finite dimension 
$1 \leqslant r = \dim\, G < \infty$. The action of an
element $g \in G$ lifts to all jet spaces $J_{x,u}^n$ of any
order $n \geqslant 0$, as in the following diagrams:
\[
\xymatrix{
J_{x,u}^n
\ar[rr]^{g\cdot}
\ar[d]
&
&
J_{y,v}^n
\ar[d]
\\
\R_{x,u}^{p+q}
\ar[rr]_{g\cdot}
&
&
\R_{y,v}^{p+q},
}
\ \ \ \ \ \ \ \ \ \ \ \ \ \ \ \ \ \ \ \
\xymatrix{
z^{(n)}
\ar[rr]^{g\cdot}
\ar[d]
&
&
w^{(n)}\big(g,z^{(n)}\big)
\ar[d]
\\
z
\ar[rr]_{g\cdot}
&
&
w(g,z).
}
\]

By definition, a {\sl differential invariant}
of order $n$ is a function
$\Iaux \colon J_{x,u}^n \longrightarrow \R$ satisfying:
\[
\Iaux\big(g\cdot z^{(n)}\big)
\,=\,
\Iaux\big(z^{(n)}\big)
\eqno
{\scriptstyle{(\forall\,g\,\in\,G)}},
\]
a property that can be diagrammatized as follows:
\[
\xymatrix{
&
&
&
&
\R
\\
J_{x,u}^n
\ar[rr]^{g\cdot}
\ar[d]
\ar@/^2pc/[urrrr]^{\Iaux}
&
&
J_{y,v}^n
\ar[d]
\ar[urr]_{\Iaux}
\\
\R_{x,u}^{p+q}
\ar[rr]_{g\cdot}
&
&
\R_{y,v}^{p+q},
}
\ \ \ \ \ \ \ \ \ \ \ \ \ \ \ \ \ \ \ \
\xymatrix{
&
&
&
&
\Iaux\big(z^{(n)}\big)
\ar@{=}[d]
\\
&
&
&
&
\Iaux\big(g\cdot z^{(n)}\big)
\\
z^{(n)}
\ar[rr]_{g\cdot}
\ar@/^2pc/[uurrrr]^{\Iaux}
&
&
g\cdot z^{(n)}
\ar[urr]_{\Iaux}.
}
\]

Viewed as such, a differential invariant is a function
defined on the jet space $J_{x,u}^n$ equipped with {\em independent}
jet coordinates:
\[
z^{(n)}
\,=\,
\big(
x^j,u^\alpha,\,
u_{x^J}^\beta
\big),
\]
but it has another more interesting meaning. 

Indeed, remember that we are considering {\em graphs} 
$\big\{ u = {\tt u}(x) \big\}$ in the source space 
and their {\em transforms}, which are graphs $\big\{ v = 
{\tt v}(y) \big\}$
in the target space. So pulled back to any 
graph $\{ u = {\tt u}(x) \}$,
a differential invariant becomes a function of the coordinates
$x = (x^1, \dots, x^p)$ on the graph:
\[
x
\,\longmapsto\,
\Iaux
\big(
x,\,
{\tt u}^\alpha(x),\,
{\tt u}_{x^J}^\beta(x)
\big).
\]
Hence a differential invariant takes various values at various
points of a graphed manifold $\{ u = {\tt u}(x) \}$.

\begin{Question}
{\sl 
Then what does it really mean, for $\Iaux \big(x^j, 
u^\alpha, u_{x^J}^\beta \big)$, to be a 
{\rm differential invariant}?
}
\end{Question}

Of course, in the target space
$\R_{y,v}^{p+q}$, the {\em same} function
$\Iaux$ of the target arguments must be considered:
\[
\Iaux
\big(
y^k,v^\gamma,\,v_{y^K}^\delta
\big).
\]
But what is the relation with $\Iaux\big( x^j, u^\alpha,
u_{x^J}^\beta \big)$?

Recall that any diffeomorphism $\phi \colon
(x,u) \longmapsto \big( \varphi(x,u),
\psi(x,u) \big)$ not far from the identity induces a
{\sl horizontal diffeomorphism} between the graphing 
horizontal spaces
$\R_x^p$ and $\R_y^p$, simply through
three maps: lifting to the
graph; performing the diffeomorphism; projecting horizontally:
\[
\xymatrix{
(x,{\tt u}(x))
\ar[rr]^{\phi\,\,\,\,\,\,\,\,\,\,\,\,\,\,\,\,\,\,\,\,\,\,\,\,\,\,\,\,\,\,}
&
&
\big(\varphi(x,{\tt u}(x)),\,\psi(x,{\tt u}(x))\big)
\ar[d]
\\
x
\ar[u]
&
&
\varphi(x,{\tt u}(x)).
}
\]
We therefore assume that we have 
a (local) diffeomorphism $\R_x^p \longmapsto
\R_y^p$:
\[
x
\xrightarrow[{\rule[0pt]{60pt}{0pt}}]{}
\varphi\big(x,{\tt u}(x)\big).
\]

\begin{center}
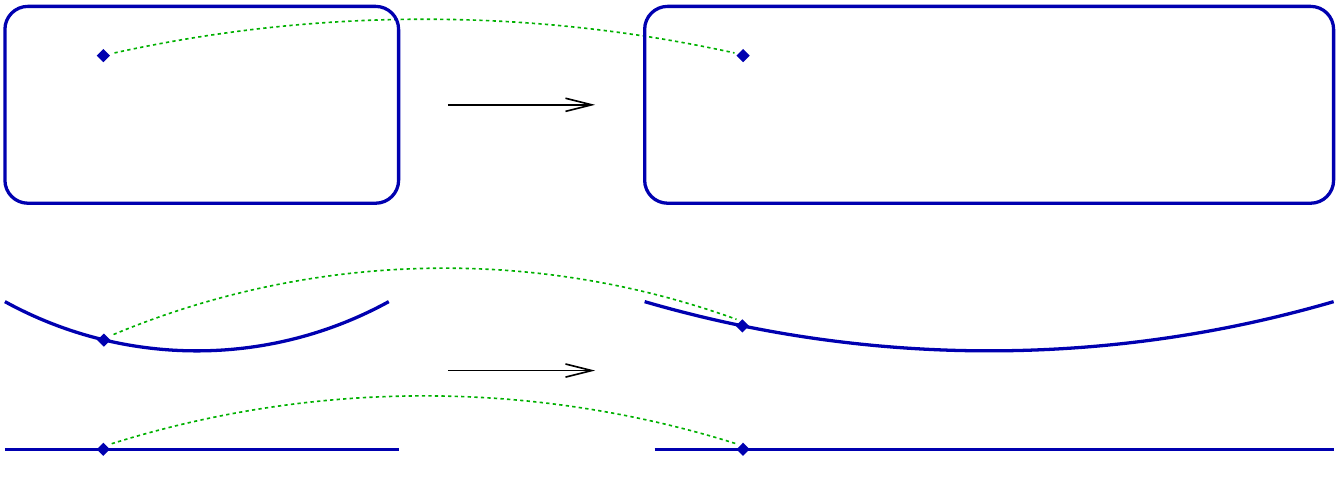
\end{center}

In our case of the action of a local Lie group $G$,
we have a {\em family} of diffeomorphisms parametrized
by $g \in G$, that we can denote as:
\[
(x,u)
\,\,\longmapsto\,\,
\big(
\varphi(g,x,u),\,
\psi(g,x,u)
\big).
\]
On restriction to a graph $\big\{ u = {\sf u}(x) \big\}$, 
this becomes:
\[
\big(x,{\tt u}(x)\big)
\,\,\longmapsto\,\,
\Big(
\varphi\big(g,x,{\tt u}(x)\big),\,
\psi\big(g,x,{\tt u}(x)\big)
\Big).
\]
Furthermore, as we did in 
Section~{\ref{graph-transformations-jet-spaces}},
we must express how a graph $\{ u = {\tt u}(x) \}$ 
is transferred to the target space.

At first, let us abbreviate the family of horizontal 
diffeomorphisms as:
\[
x
\,\,\longmapsto\,\,
\varphi
\big(
g,x,{\tt u}(x)
\big)
\,=:\,
{\tt y}(g,x)
\,=:\,
y.
\]
Also, let us denote the family of inverse diffeomorphisms as:
\[
x
\,:=\,
{\tt x}(g,y)
\,\,\longleftarrow\,\,
y,
\]
by means of certain maps ${\tt x}(g, \cdot)$ 
coming from an application of the
implicit function theorem.  Then the source graph $\{ u = {\sf u}(x)
\}$ becomes a family of graphs over the $\R_y^p$ space:
\[
\aligned
v
&
\,=\,
\psi
\big(
g,\,
{\tt x}(g,y),\,
{\tt u}
\big(
{\tt x}(g,y)
\big)
\big)
\\
&\,
\,=:\,
{\sf v}
(g,y).
\endaligned
\]

In summary, the action of a local Lie group $G$ on $\R^{p+q}$ provides
a family of diffeomorphisms $\phi_g(\cdot) = \big( \varphi (g, \cdot),
\psi (g, \cdot) \big)$ together with:

\smallskip\noindent$\bullet$\,
a $G$-parametrized family of horizontal diffeomorphisms
$x \longmapsto {\tt y} (g,x)$ from $\R_x^p$ to $\R_y^p$;

\smallskip\noindent$\bullet$\,
a $G$-parametrized family of graphs $\big\{ v = {\tt v}(g,y) \big\}$
in the target space; 

\smallskip\noindent
which are both associated 
to any given graph $\{ u = {\tt u} (x) \}$
in the source space.

\smallskip

Hence as in the figure above, for any $g \in G$, when $x \in \R^p$
varies, points $\big(x, {\sf u}(x) \big)$ in a source graph are sent
to points $\big( {\tt y}(g,x),\, {\tt v} \big(g, {\tt y}(g,x) \big)
\big)$ in a target graph so that
the pointwise one-to-one correspondence between these
two graphs reads:
\[
\big(
x,{\tt u}(x)
\big)
\,\,\,\longleftrightarrow\,\,\,
\Big(
{\tt y}(g,x),\,
{\tt v}
\big(g,{\tt y}(g,x)
\big)
\Big).
\]
We can now express what it really means to be a differential 
invariant.

\begin{Theorem}
A differential invariant $\Iaux = \Iaux \big(
x^j, u^\alpha, u_{x^J}^\beta \big)$ takes identical values
at points which correspond one to another
through the action of any group element $g \in G$:
\[
\Iaux
\big(
x,{\tt u}^\alpha(x),{\tt u}_{x^J}^\beta(x)
\big)
\,\,\equiv\,\,
\Iaux
\Big(
{\tt y}(g,x),\,
{\tt v}^\gamma
\big(g,{\tt y}(g,x)\big),\,\,
{\tt v}_{y^K}^\delta
\big(g,{\tt y}(g,x)\big)
\Big)
\eqno
{\scriptstyle{(\forall\,x\,\in\,\R^p)}}.
\eqno\qed
\]
\end{Theorem}

\Section{\bf Solving Cross-Section Equations for Curves
$C^1 \subset \R^2$ Under $\SA_2(\R)$}
\label{solving-cross-section-equations-curves}
\HEAD{{\ref{solving-cross-section-equations-curves}}.~{\sf 
Solving Cross-Section Equations for Curves
$C^1 \subset \R^2$ Under $\SA_2(\R)$}
}{
Zhangchi {\sc Chen}, Joël~{\sc Merker}}

We illustrate this {\sl cross-section} approach 
({\cite{Fels-Olver-1999, Olver-2007, Olver-2018}})
in the case of curves $C^1 \subset \R^2$ under 
the action of the {\sl equiaffine} group $\SA_2(\R)$,
consisting of area-preserving affine transformations
$\R_{x,u}^{1+1} \longrightarrow \R_{y,v}^{1+1}$:
\[
\aligned
y
&
\,=\,
{\sf a}\,x
+
{\sf b}\,u
+
{\sf c},
\\
v
&
\,=\,
{\sf k}\,x
+
{\sf l}\,u
+
{\sf m},
\endaligned
\ \ \ \ \ \ \ \ \ \ \ \ \ \ \ \ \ \ \ \
1
\,=\,
\left\vert\!
\begin{array}{cc}
{\sf a} & {\sf b}
\\
{\sf k} & {\sf l}
\end{array}
\!\right\vert.
\]

According to the preceding sections, any (local) invertible 
diffeomorphism: 
\[
\phi 
\colon
\ \ \
\R_{x,u}^2 
\,\,\,\longrightarrow\,\,\,
\R_{y,v}^2 
\]
with $\det\, \Jac\, \phi \neq 0$ which is not far
from the identity transforms, through an application
of the implicit function theorem, sends every graphed curve
$\big\{ u = F(x) \big\}$ into a similar graphed curve
$\big\{ y = G(v) \big\}$. Then tangents are transferred 
from one curve to the other, 
and higher order jets as well. This means that the 
diffeomorphism $\phi$ lifts as an invertible
transformation between corresponding $n$\textsuperscript{th}
jet spaces:
\[
\xymatrix{
J_{x,u}^n
\ar[rr]^{\phi^{(n)}}
\ar[d]
&
&
J_{y,v}^n
\ar[d]
\\
\R_{x,u}^2
\ar[rr]_g
&
&
\R_{y,v}^2,
}
\] 
with successive:
\[
\aligned
\phi^{(1)}
\colon\ \ \
\big(x,u,u_1\big)
&
\,\,\,\longmapsto\,\,\,
\Big(
y(x,u),\,
v(x,u),\,\,
v_1\big(x,u,u_1\big)
\Big),
\\
\phi^{(2)}
\colon\ \ \
\big(x,u,u_1,u_2\big)
&
\,\,\,\longmapsto\,\,\,
\Big(
y(x,u),\,
v(x,u),\,\,
v_1\big(x,u,u_1\big),\,\,
v_2\big(x,u,u_1,u_2\big)
\Big),
\endaligned
\]
whose first two components are of course those of $\phi$,
with certain uniquely determined functions 
$v_1\big(x,u,u_1\big)$, 
$v_2\big(x,u,u_1,u_2\big)$, {\em etc.},
which use the total differentiation operator:
\[
{\sf D}_x
\,=\,
\frac{\partial}{\partial x}
+
u_1\,\frac{\partial}{\partial u}
+
\sum_{i=1}^\infty\,
u_{i+1}\,
\frac{\partial}{\partial u_i},
\]
namely:
\[
v_1
\,=\,
\frac{{\sf D}_x(v)}{{\sf D}_x(y)},
\ \ \ \ \ \ \ \ \ \ \ \ \ \ \ \ \ \ \ \
v_2
\,=\,
\frac{{\sf D}_x(v_1)}{{\sf D}_x(y)},
\ \ \ \ \ \ \ \ \ \ \ \ \ \ \ \ \ \ \ \
v_3
\,=\,
\frac{{\sf D}_x(v_2)}{{\sf D}_x(y)},
\ \ \ \ \ \ \ \ \ \ \ \ \ \ \ \ \ \ \ \
\cdots\cdots.
\]

However, the expanded formulas are known to rapidly become unwieldy,
especially in higher dimensions. At least here with just
one independent variable $x$, we have up to jet order $2$
({\cite{Lie-1883, Merker-2008}}):
\[
\aligned
v_1
&
\,=\,
\frac{v_x+u_1v_u}{y_x+u_1y_u},
\\
v_2
&
\,=\,
\frac{1}{(y_x+u_1y_u)^3}\,
\bigg\{
u_2
\left\vert\!
\begin{array}{cc}
y_{x} & y_{u}
\\
v_{x} & v_{u}
\end{array}
\!\right\vert
+
\left\vert\!
\begin{array}{cc}
y_{x} & y_{xx}
\\
v_{x} & v_{xx}
\end{array}
\!\right\vert
+
u_1
\left[
2\,
\left\vert\!
\begin{array}{cc}
y_{x} & y_{xu}
\\
v_{x} & v_{xu}
\end{array}
\!\right\vert
-
\left\vert\!
\begin{array}{cc}
y_{xx} & y_{u}
\\
v_{xx} & v_{u}
\end{array}
\!\right\vert\,
\right]
\\
&
\ \ \ \ \ \ \ \ \ \ \ \ \ \ \ \ \ \ \ \ \ \ \ \ \ \ \ \ \ \ \ \ \ \ 
\ \ \ \ \ \ \ \ \ \ \ \ \ \ \ \ \ \,
+
u_1^2
\left[\,
\left\vert\!
\begin{array}{cc}
y_{x} & y_{uu}
\\
v_{x} & v_{uu}
\end{array}
\!\right\vert
-
2\,
\left\vert\!
\begin{array}{cc}
y_{xu} & y_{u}
\\
v_{xu} & v_{u}
\end{array}
\!\right\vert\,
\right]
-
u_1^3
\left\vert\!
\begin{array}{cc}
y_{uu} & y_{u}
\\
v_{uu} & v_{u}
\end{array}
\!\right\vert\
\bigg\}.
\endaligned
\]

When the diffeomorphism is a special affine transformation:
\[
y
\,=\,
{\sf a}\,x
+
{\sf b}\,u
+
{\sf c},
\ \ \ \ \ \ \ \ \ \ \ \ \ \ \ \ \ \ \ \
v
\,=\,
{\sf k}\,x
+
{\textstyle{\frac{1+{\sf b}{\sf k}}{{\sf a}}}}\,
u
+
{\sf m},
\]
a computer yields compact formulas:
\[
\aligned
v_1
&
\,=\,
\frac{{\sf a}{\sf k}
+
(1+{\sf b}{\sf k})\,u_1}{
{\sf a}\,({\sf a}+{\sf b}\,u_1)},
\\
v_2
&
\,=\,
\frac{u_2}{
({\sf a}+{\sf b}\,u_1)^3},
\\
v_3
&
\,=\,
\frac{-3{\sf b}\,u_2^2+{\sf a}\,u_3+{\sf b}\,u_1u_3}{
({\sf a}+{\sf b}\,u_1)^5},
\\
v_4
&
\,=\,
\frac{-10{\sf b}^2\,u_1u_2u_3
-
10{\sf a}{\sf b}\,u_2u_3
+
15{\sf b}^2\,u_2^3
+
2{\sf a}{\sf b}\,
u_1u_4
+
{\sf b}^2\,u_1^2u_4
+
{\sf a}^2\,u_4}{
({\sf a}+{\sf b}\,u_1)^7}.
\endaligned
\]
The natural cross-section ({\em cf.}
Section~{\ref{special-affine-power-series-invariants-curves-R-2}}):
\[
\big\{
y=0,\,\,\,
v=0,\,\,\,
v_1=0,\,\,\,
v_2=1,\,\,\,
v_3=0
\big\},
\]
enables one to solve the $5$ group parameters ${\sf a}$, 
${\sf b}$, ${\sf c}$, ${\sf k}$, ${\sf m}$ as:
\[
\aligned
{\sf a}
&
\,=\,
\frac{3\,u_2^2-u_1u_3}{3\,u_2^{5/3}},
\\
{\sf b}
&
\,=\,
\frac{u_3}{3\,u_2^{5/3}},
\\
{\sf c}
&
\,=\,
\frac{-3\,xu_2^2+xu_1u_3-uu_3}{3\,u_2^{5/3}},
\\
{\sf k}
&
\,=\,
-\,\frac{u_1}{u_2^{1/3}},
\\
{\sf m}
&
\,=\,
\frac{-u+xu_1}{u_2^{1/3}},
\endaligned
\]
which provides the moving frame map $\rho \colon
J_{x,u}^n \longrightarrow \SA_2(\R)$ for any $n \geqslant 3$, 
and then, replacing in the formula for $v_4$, we 
obtain what we will call the {\sl parabolas invariant:}
\[
\Paux
\,:=\,
\frac{1}{3}\,
\frac{-5\,u_3^2+3\,u_2u_4}{u_2^{8/3}},
\]
and which is also known as the {\sl equiaffine curvature}.

\Section{\bf Equiaffine Group $\SA_3(\R)$ and its Action on
Graphed Surfaces $S^2 \subset \R^3$}
\label{equiaffine-SA-3-surfaces-S2-R3}
\HEAD{{\ref{equiaffine-SA-3-surfaces-S2-R3}}.~{\sf 
Equiaffine Group $\SA_3(\R)$ and its Action on
Graphed Surfaces $S^2 \subset \R^3$}
}{
Zhangchi {\sc Chen}, Joël~{\sc Merker}}

Let a source space $\R^3$ be equipped with coordinates $(x,y,u)$,
and let a target space $\R^3$ be equipped with coordinates
$(s,t,v)$. We are interested in how local 
analytic graphed surfaces $\big\{ u = F(x,y) \big\}$
are mapped to local analytic graphed surfaces
$\big\{ v = G(s,t) \big\}$ through simple 
transformations $\R_{x,y,u}^3 \longrightarrow \R_{s,t,v}^3$.

\Subsection{Affine and special affine transformations}
The {\sl affine group:}
\[
{\sf A}_3(\R)
\,:=\,
{\sf GL}_3(\R)
\ltimes
\R^3
\]
consists of invertible linear transformations in ${\sf GL}_3(\R)$,
coupled with translations:
\[
\aligned
s
&
\,=\,
{\sf a}\,x
+
{\sf b}\,y
+
{\sf c}\,u
+
{\sf d},
\notag
\\
t
&
\,=\,
{\sf k}\,x
+
{\sf l}\,y
+
{\sf m}\,u
+
{\sf n},
\\
v
&
\,=\,
{\sf p}\,x
+
{\sf q}\,y
+
{\sf r}\,u
+
{\sf s},
\notag
\endaligned
\ \ \ \ \ \ \ \ \ \ \ \ \ \ \ \ \ \ \ \
0
\,\neq\,
\left\vert\!
\begin{array}{ccc}
{\sf a} & {\sf b} & {\sf c}
\\
{\sf k} & {\sf l} & {\sf m}
\\
{\sf p} & {\sf q} & {\sf r}
\end{array}
\!\right\vert
\,=:\,
\delta.
\]
When the determinant $\delta = 1$,
so that the volume (and the orientation) 
of geometric objects
remains unchanged, the transformation is called {\sl special affine},
or {\sl equiaffine}:
\[
\SA_3(\R)
\,:=\,
\SL_3(\R)
\ltimes
\R^3.
\]
Both groups ${\sf A}_3(\R)$ and $\SA_3(\R)$ 
act on surfaces $S^2 \subset \R^3$ graphed as 
$\big\{ u = F(x,y) \big\}$.

We write $\R \{x, y\}$ to denote the ring of convergent
power series defined in some neighborhood of the origin
$(0, 0) \in \R^2$. Each element $F(x,y) \in \R\{ x, y\}$
admits a power series expansion:
\[
F(x,y)
\,=\,
F_{0,0}
+
F_{1,0}\,
{\textstyle{\frac{x^1}{1!}}}
+
F_{0,1}\,
{\textstyle{\frac{y^1}{1!}}}
+
F_{2,0}\,
{\textstyle{\frac{x^2}{2!}}}
+
F_{1,1}\,
{\textstyle{\frac{x^1y^1}{1!\,1!}}}
+
F_{0,2}\,
{\textstyle{\frac{y^2}{2!}}}
+\cdots,
\]
namely:
\[
F(x,y)
\,=\,
\sum_{k=0}^\infty\,
\sum_{l=0}^\infty\,
F_{k,l}\,
{\textstyle{\frac{x^k}{k!}}}\,
{\textstyle{\frac{y^l}{l!}}}
\eqno
{\scriptstyle{(F_{k,l}\,\in\,\R)}}.
\]
The two monomials $x$ and $y$ generate the maximal ideal
$\langle x, y \rangle \subset \R\{ x,y\}$ formed of
power series with $F(0,0) = 0$.
For any order 
$\oorder \geqslant 0$, 
the quotient:
\[
\R\{x,y\}
\big/
\langle x,y\rangle^{\oorder+1},
\]
is a free $\R$-module of rank $\binom{\oorder+2}{2}$, plainly
generated by all the monomials $x^j y^k$ with $j + k \leqslant
\oorder$.
 
In this article, we will focus on local analytic graphing functions
$F(x,y) \in \R\{ x, y\}$ of this sort, and on similar functions
$G(s,t) \in \R\{ s,t\}$, possibly with $F(0, 0) \neq 0 \neq G(0, 0)$.

\begin{Definition}
Two local graphed surfaces $\big\{ u = F(x,y) \big\}$ and $\big\{ v =
G(s,t) \big\}$ are {\sl (special) affinely equivalent} if there exists
a (special) affine transformation which sends the one to the 
other\,\,---\,\,a concrete criterion follows in a second.
\end{Definition}

We will always consider transformations in ${\sf A}_3(\R)$ 
or in $\SA_3(\R)$ which are not far from the identity, so that
graphs are transformed into graphs. 
Thus, if we let act a (special) affine
transformation $(x,y,u) \longmapsto 
(s,t,v)$ as above, any point $\big( x, y, F(x,y) \big)$ is
sent to the point $\big( s, t, G(s,t) \big)$ and we must have:
\[
{\sf p}\,x
+
{\sf q}\,y
+
{\sf r}\,u
+
{\sf s}
\,\,=\,\,
G
\big(
{\sf a}\,x
+
{\sf b}\,y
+
{\sf c}\,u
+
{\sf d},\,\,
{\sf k}\,x
+
{\sf l}\,y
+
{\sf m}\,u
+
{\sf n}
\big)
\Big\vert_{u=F(x,y)}.
\]
Equivalently, the following {\sl fundamental equation}:
\leqnomode\usetagform{default}
\begin{align}
\label{eqfond-F-x-y-G-s-t}
{\sf p}\,x
+
{\sf q}\,y
+
{\sf r}\,F(x,y)
+
{\sf s}
\,\,=\,\,
G
\big(
{\sf a}\,x
+
{\sf b}\,y
+
{\sf c}\,F(x,y)
+
{\sf d},\,\,
{\sf k}\,x
+
{\sf l}\,y
+
{\sf m}\,F(x,y)
+
{\sf n}
\big),
\end{align}
must hold identically in the domain of convergence of $F(x,y)$. 

\begin{Question}
How to determine equivalence classes of surfaces $S^2 \subset \R^3$
modulo ${\sf A}_3(\R)$ or $\SA_3(\R)$?
\end{Question}

As already presented in a general context, we will see that there
exist rational combinations of derivatives of the graphing function
$F$ which are invariant under ${\sf A}_3(\R)$, or $\SA_3(\R)$, 
called {\sl differential invariants}. We will
realize that if $\big\{ u = F(x,y) \big\}$ is equivalent to $\big\{ v
= G(s,t) \big\}$, differential invariants are `the same', namely
correspond to each other as explained in
Section~{\ref{meaning-differential-invariant}}.  We will also see that
differential invariants determine equivalent classes of (local)
surfaces. All this will be clearer later.

\Subsection{Three types of substitutions}
For now, 
in order to avoid confusion, it is important to differentiate
three types of use of the 
action of a transformation group like
${\sf A}_3(\R)$ or $\SA_3(\R)$.

\smallskip\noindent{\bf (S1)}\,
{\footnotesize\sf Equivalence:} 
Have in hands an explicit 
(special) affine equivalence between two given surfaces
$\big\{ u = F(x,y) \big\}$ and $\big\{ v = G(s,t) \big\}$, 
as above\,\,---\,\,a situation which rarely occurs {\em a priori}.

\smallskip\noindent{\bf (S2)}\,
{\footnotesize\sf Transformation:} 
Start from a given surface $\big\{ u = F(x,y) \big\}$, apply
a (special) affine transformation, and ask what is the equation
of the new surface $\big\{ v = G(s,t) \big\}$\,\,---\,\,the 
answer relying on an application of the implicit function
theorem.

\smallskip\noindent{\bf (S3)}\,
{\footnotesize\sf Normalization:} 
Start from a given surface $\big\{ u = F(x,y) \big\}$
and ask whether there exist appropriate 
special affine transformations with put the target
graphed surface $\big\{ v = G(s,t) \big\}$ into a `simpler',
`normalized', form\,\,---\,\,the core of the problem,
which, after being solved, will
enable to know whether two given surfaces are equivalent
like in {\small\bf (S1)}.

\smallskip

In both situations {\small\bf (S2)} and {\small\bf (S3)},
experience tells us that 
it is more appropriate to consider the {\em inverse} transformation
$(x,y,u) \longleftarrow (s,t,v)$. 
We choose to define such an inverse (special) affine 
transformation using the same letters for group parameters:
\begin{align}
\aligned
\label{x-y-u-A-3-s-t-v}
x
&
\,=\,
{\sf a}\,s
+
{\sf b}\,t
+
{\sf c}\,v
+
{\sf d},
\notag
\\
y
&
\,=\,
{\sf k}\,s
+
{\sf l}\,t
+
{\sf m}\,v
+
{\sf n},
\\
u
&
\,=\,
{\sf p}\,s
+
{\sf q}\,t
+
{\sf r}\,v
+
{\sf s},
\notag
\endaligned
\ \ \ \ \ \ \ \ \ \ \ \ \ \ \ \ \ \ \ \
0
\,\neq\,
\left\vert\!
\begin{array}{ccc}
{\sf a} & {\sf b} & {\sf c}
\\
{\sf k} & {\sf l} & {\sf m}
\\
{\sf p} & {\sf q} & {\sf r}
\end{array}
\!\right\vert
\,=:\,
\delta.
\end{align}

For instance, if $F(x,y)$ is given, we can answer {\small\bf (S2)}
by plugging $(x,y,u)$ in terms of $(s,t,v)$ inside 
$0 = -\, u + F(x,y)$, which gives:
\[
0
\,=\,
-\,
{\sf p}\,s
-
{\sf q}\,t
-
{\sf r}\,v
-
{\sf s}
+
F
\big(
{\sf a}\,s
+
{\sf b}\,t
+
{\sf c}\,v
+
{\sf d},\,\,
{\sf k}\,s
+
{\sf l}\,t
+
{\sf m}\,v
+
{\sf n}
\big),
\]
and then the graphing function $G(s,t)$ of the transformed 
surface is obtained by solving this equation for $v$,
using the implicit function theorem, provided the partial
derivative with respect to $v$ is nonvanishing:
\[
0
\,\neq\,
-\,
{\sf r}
+
{\sf c}\,F_x
+
{\sf m}\,F_y,
\]
a condition which holds automatically when the transformation
is close to the identity, since ${\sf r} \approx 1$, while
${\sf c} \approx 0 \approx {\sf m}$.

Later, we will mainly work in order to perform normalizing 
substitutions of the sort {\small\bf (S3)}, first
in our treatment of curves done in
Sections~{\ref{special-affine-power-series-invariants-curves-R-2}}, 
{\ref{affine-invariants-curves-R-2}}, 
and then more intensively when studdying parabolic surfaces
$S^2 \subset \R^3$, in 
Sections~{\ref{parabolic-pseudostablization}},
{\ref{relative-invariant-S-first-invariant-W}},
{\ref{branch-W-equiv-0-branch-W-nonzero}},
{\ref{recurrence-relations-parabolic-surfaces}}.

\Section{\bf Affine Invariancy of the Ranks of Hessian Matrices}
\label{affine-invariancy-ranks-Hessian-matrices}
\HEAD{{\ref{affine-invariancy-ranks-Hessian-matrices}}.~{\sf 
Affine Invariancy of the Ranks of Hessian Matrices}
}{
Zhangchi {\sc Chen}, Joël~{\sc Merker}}

For now, we come back to the simple situation {\small\bf (S1)}
in which an equivalence is given, not sought for. 
Differentiating the fundamental
equation~({\ref{eqfond-F-x-y-G-s-t}})
with respect to $x$ and to $y$, we get:
\leqnomode\usetagform{default}
\begin{align}
\label{F-x-F-y-G-s-G-t}
\aligned
{\sf p}+{\sf r}\,F_x
&
\,=\,
\big({\sf a}+{\sf c}\,F_x\big)\,
G_s
+
\big({\sf k}+{\sf m}\,F_x\big)\,
G_t,
\\
{\sf q}+{\sf r}\,F_y
&
\,=\,
\big({\sf b}+{\sf c}\,F_y\big)\,
G_s
+
\big({\sf l}+{\sf m}\,F_y\big)\,G_t.
\endaligned
\end{align}
Therefore, for those transformations of $\SA_3(\R)$ satisfying:
\leqnomode\usetagform{default}
\begin{align}
\label{definition-Lambda}
0
\,\neq\,
\Lambda
\,:=\,
\left\vert\!
\begin{array}{cc}
{\sf a}+{\sf c}\,F_x
&
{\sf k}+{\sf m}\,F_x
\\
{\sf b}+{\sf c}\,F_y
&
{\sf l}+{\sf m}\,F_y
\end{array}
\!\right\vert,
\end{align}
a condition 
which certainly holds for transformations close to the identity,
we can solve $G_s$ and $G_t$ in terms of $F_x$, $F_y$:
\[
G_s
\,:=\,
\frac{1}{\Lambda}\,
\left\vert\!
\begin{array}{cc}
{\sf p}+{\sf r}\,F_x
&
{\sf k}+{\sf m}\,F_x
\\
{\sf q}+{\sf r}\,F_y
&
{\sf l}+{\sf m}\,F_y
\end{array}
\!\right\vert
\ \ \ \ \ \ \ \ \ \ \ \ \ \ \ \ \ \ \ \
G_t
\,=\,
\left\vert\!
\begin{array}{cc}
{\sf a}+{\sf c}\,F_x
&
{\sf p}+{\sf r}\,F_x
\\
{\sf b}+{\sf c}\,F_y
&
{\sf q}+{\sf r}\,F_y
\end{array}
\!\right\vert.
\]

Beyond, by keeping differentiating the fundamental
equation~({\ref{eqfond-F-x-y-G-s-t}}),
we claim that we can iteratively
solve every partial derivative $G_{s^l t^m}$ 
in terms of the partial derivatives $\big\{ F_{x^j y^k} \big\}_{
j+k \leqslant l+m}$.

\begin{Lemma}
For every order $\oorder \geqslant 0$, the affine group 
${\sf A}_3(\R)$ and the special affine $\SA_3(\R)$ groups act
on Taylor coefficients with $l + m \leqslant \oorder$:
\[
G_{l,m}
\,=\,
\formula\,
\Big(
\begin{smallmatrix} 
{\sf a}, & {\sf b}, & {\sf c}, & {\sf d},
\\ 
{\sf k}, & {\sf l}, & {\sf m}, & {\sf n}, 
\\ 
{\sf p}, & {\sf q}, & {\sf r}, & {\sf s}, 
\end{smallmatrix}
\ \ \
\big\{
F_{j,k}
\big\}_{j+k\leqslant l+m}
\Big).\eqno\qed
\]
\end{Lemma}

Such (complicated) formulas express the induced action
of ${\sf A}_3(\R)$ or of $\SA_3(\R)$ on the spaces  of Taylor
coefficients:
\[
\R^{\binom{\oorder+2}{2}}
\,\,\cong\,\,
\R\{x,y\}
\big/
\langle x,y\rangle^{\oorder+1}.
\]

In the context of jet spaces, 
Section~{\ref{graph-transformations-jet-spaces}}
already showed 
how to obtain such formulas, provided one restricts 
considerations to the fiber over the origin $(0,0)$ of the jet
space $J_{x,y}^{\oorder}$. Indeed, remember we introduced the
modified total differentiation operators that we can
now write after pull-back to the graph 
$\{ u = F(x,y) \}$ as:
\[
\left(\!
\begin{array}{c}
{\sf E}_x
\\
{\sf E}_y
\end{array}
\!\right)
\,=\,
\left(\!
\begin{array}{cc}
\frac{{\sf l}+{\sf m}F_y}{\Lambda} 
&
\frac{-{\sf k}-{\sf m}F_x}{\Lambda}
\\
\frac{-{\sf b}-{\sf c}F_y}{\Lambda}
&
\frac{{\sf a}+{\sf c}F_x}{\Lambda} 
\end{array}
\!\right)\,
\left(\!
\begin{array}{c}
{\sf D}_x
\\
{\sf D}_y
\end{array}
\!\right),
\]
where ${\sf D}_x$, ${\sf D}_y$ are total differentiation operators,
and then, similarly as in 
Theorem~{\ref{Thm-prolongation-E-j}}\,\,---\,\,but only
over the origin $(0,0)$, not over every point $(x,y)$\,\,---,
we have:
\[
\aligned
G_s
&
\,=\,
{\sf E}_x
\Big(
{\sf p}\,x
+
{\sf q}\,y
+
{\sf r}\,F(x,y)
+
{\sf s}
\Big),
\\
G_t
&
\,=\,
{\sf E}_y
\Big(
{\sf p}\,x
+
{\sf q}\,y
+
{\sf r}\,F(x,y)
+
{\sf s}
\Big).
\endaligned
\]
and beyond:
\[
G_{ss}
\,=\,
{\sf E}_x
\big(
G_s
\big),
\ \ \ \ \ \ \ \ \ \ \ \ \ \ \ \ \ \ \ \
G_{st}
\,=\,
{\sf E}_x
\big(
G_t
\big)
\,=\,
{\sf E}_y
\big(
G_s
\big)
\,=\,
G_{ts},
\ \ \ \ \ \ \ \ \ \ \ \ \ \ \ \ \ \ \ \
G_{tt}
\,=\,
{\sf E}_y
\big(
G_t
\big).
\]

In length, the second-order formulas for Taylor coefficients 
obtained by differentiating~({\ref{F-x-F-y-G-s-G-t}}) 
once more with respect
to $x$ and $y$ 
read:
\[
\footnotesize
\aligned
{\sf r}\,F_{xx}
&
\,=\,
{\sf c}\,F_{xx}\,G_s
+
{\sf m}\,F_{xx}\,G_t
+
\big({\sf a}+{\sf c}\,F_x\big)^2\,
G_{ss}
+
2\,
\big({\sf a}+{\sf c}\,F_x\big)\,
\big({\sf k}+{\sf m}\,F_x\big)\,
G_{st}
+
\big({\sf k}+{\sf m}\,F_x\big)^2\,
G_{tt},
\\
{\sf r}\,F_{xy}
&
\,=\,
{\sf c}\,F_{xy}\,G_s
+
{\sf m}\,F_{xy}\,G_t
\\
&
\ \ \ \ \
+
\big({\sf a}+{\sf c}\,F_x\big)\,
\big({\sf b}+{\sf c}\,F_y\big)\,
G_{ss}
+
\Big[
\big({\sf a}+{\sf c}\,F_x\big)\,
\big({\sf l}+{\sf m}\,F_y\big)
+
\big({\sf b}+{\sf c}\,F_y\big)\,
\big({\sf k}+{\sf m}\,F_x\big)
\Big]\,
G_{st}
\\
&
\ \ \ \ \
+
\big({\sf k}+{\sf m}\,F_x\big)\,
\big({\sf l}+{\sf m}\,F_y\big)\,
G_{tt},
\\
{\sf r}\,F_{yy}
&
\,=\,
{\sf c}\,F_{yy}\,G_s
+
{\sf m}\,F_{yy}\,G_t
+
\big({\sf b}+{\sf c}\,F_y\big)^2\,
G_{ss}
+
2\,
\big({\sf b}+{\sf c}\,F_y\big)\,
\big({\sf l}+{\sf m}\,F_y\big)\,
G_{st}
+
\big({\sf l}+{\sf m}\,F_y\big)^2\,
G_{tt}.
\endaligned
\]

As is known, 
the Hessian matrix of a graphed surface $\big\{ u = F(x,y) \big\}$
is:
\[
\Hessian_F
\,=\,
\left(\!
\begin{array}{cc}
F_{xx} & F_{xy}
\\
F_{yx} & F_{yy}
\end{array}
\!\right).
\]
We will denote its determinant by:
\[
\Haux_F
\,=\,
\det\,
\Hessian_F
\,=\,
F_{xx}\,F_{yy}
-
F_{xy}^2.
\]

When all entries of the Hessian matrix 
are identically $0 \equiv F_{xx} \equiv F_{xy}
\equiv F_{yy}$, the Taylor series of $F$ reduces to its
first order terms $F_{0,0} + 
F_{1,0}\, x + F_{0,1}\, y$, and it is easy to verify the

\begin{Proposition}
For a graphed surface $\big\{ u = F(x,y) \big\}$, 
the following two properties are equivalent:

\smallskip\noindent{\bf (i)}\,
there is a special affine transformation 
$(x,y,u) \longmapsto (s,t,v)$
sending it
to a plane $\{ v = 0\}$;

\smallskip\noindent{\bf (ii)}\,
the Hessian matrix $\Hessian (F) \equiv 
\left(\! \begin{smallmatrix} 0 & 0 
\\ 0 & 0 \end{smallmatrix} \!\right)$
is identically zero.\qed
\end{Proposition}

The preceding section showed formulas 
giving $F_x$, $F_y$, $F_{xx}$, $F_{xy}$, $F_{yy}$ in terms of
$G_s$, $G_t$, $G_{ss}$, $G_{st}$, $G_{tt}$.
Then a direct 
check\,\,---\,\,using a computer!\,\,---\,\,yields 
a formula 
({\cite[5.2]{Merker-2019}}) 
for the transfer of Hessians {\em determinants:}
\[
\Haux_G
\,=\,
\nonzero
\cdot
\Haux_F,
\]
valid for any (not necessarily special)
affine transformation with nonzero quantities:
\[
0
\,\neq\,
\delta
\,=\,
\left\vert\!
\begin{array}{ccc}
{\sf a} & {\sf b} & {\sf c}
\\
{\sf k} & {\sf l} & {\sf m}
\\
{\sf p} & {\sf q} & {\sf r}
\end{array}
\!\right\vert
\ \ \ \ \ \ \ \ \ \ \ \ \ \ \ \ \ \ \ \
\text{and}
\ \ \ \ \ \ \ \ \ \ \ \ \ \ \ \ \ \ \ \
\Lambda
\,=\,
\left\vert\!
\begin{array}{cc}
{\sf a}+{\sf c}\,F_x
&
{\sf k}+{\sf m}\,F_x
\\
{\sf b}+{\sf c}\,F_y
&
{\sf l}+{\sf m}\,F_y
\end{array}
\!\right\vert
\,\neq\,
0.
\]

\begin{Theorem}
\label{Thm-transfer-Hessian-determinants}
One has:
\[
G_{ss}\,G_{tt}
-
G_{st}^2
\,=\,
\frac{\delta^2}{\Lambda^4}\,
\Big(
F_{xx}\,F_{yy}
-
F_{xy}^2
\Big).
\eqno\qed
\]
\end{Theorem}

Consequently, the vanishing or the nonvanishing of the Hessian
determinant is an affinely invariant property.
In fact, this is a consequence of a more informative

\begin{Proposition}
\label{Prp-invariancy-rank-Hessian-matrix}
The rank and the signature of the Hessian matrix of a graphed
surface $\big\{ u = F(x,y) \big\}$ are unchanged after any
affine transformation.
More precisely, under any affine 
transformation, one has:
\[
\left(\!
\begin{array}{cc}
{\sf a}+{\sf c}\,F_x
&
{\sf k}+{\sf m}\,F_x
\\
{\sf b}+{\sf c}\,F_y
&
{\sf l}+{\sf m}\,F_y
\end{array}
\!\right)
\cdot
\left(\!
\begin{array}{cc}
G_{ss} & G_{st}
\\
G_{ts} & G_{tt}
\end{array}
\!\right)
\cdot
\left(\!
\begin{array}{cc}
{\sf a}+{\sf c}\,F_x
&
{\sf k}+{\sf m}\,F_x
\\
{\sf b}+{\sf c}\,F_y
&
{\sf l}+{\sf m}\,F_y
\end{array}
\!\right)^{\!{\tt t}}
\,\,=\,\,
\frac{\delta}{\Lambda}\,
\left(\!
\begin{array}{cc}
F_{xx} & F_{xy}
\\
F_{yx} & F_{yy}
\end{array}
\!\right).
\]
\end{Proposition}

\proof
We already wrote above two formulas giving
$G_s$ and $G_t$ in terms of $F_x$, $F_y$.
Using a computer, starting from the three formulas
for second derivatives written above, 
we may similarly solve 
$G_{ss}$, $G_{st}$, $G_{tt}$ in terms of
$F_{xx}$, $F_{xy}$, $F_{yy}$, $F_x$, $F_y$.
The formulas are quite large. 
Equivalently, with $v := {\sf p}\, x + {\sf q}\, y +
{\sf r}\, F(x,y) + {\sf s}$:
\[
G_{ss}
\,=\,
{\sf E}_x\big({\sf E}_x(v)\big),
\ \ \ \ \ \ \ \ \ \ \ \ \ \ \ \ \ \ \ \
G_{st}
\,=\,
{\sf E}_x\big({\sf E}_y(v)\big),
\ \ \ \ \ \ \ \ \ \ \ \ \ \ \ \ \ \ \ \
G_{tt}
\,=\,
{\sf E}_y\big({\sf E}_y(v)\big).
\]
Still on a computer, we verify that this
matrix identity holds.
\endproof

\begin{Definition}
A point $p = (x_p,y_p)$
on a graphed surface $\big\{ u = F(x,y) \big\}$ is
called:

\smallskip\noindent$\bullet$\,
{\sl flat} if $\Hessian_F(p) = 
\left(\! \begin{smallmatrix} 0 & 0 
\\ 0 & 0 \end{smallmatrix} \!\right)$;

\smallskip\noindent$\bullet$\,
{\sl parabolic} if $\Hessian_F(p)$ has rank $1$; 

\smallskip\noindent$\bullet$\,
{\sl elliptic} if $\Hessian_F(p)$ has rank $2$ and signature
$(2, 0)$ or $(0,2)$;

\smallskip\noindent$\bullet$\,
{\sl hyperbolic} if $\Hessian_F(p)$ has rank $2$ and signature
$(1, 1)$.

\end{Definition}

Of course, these $4$ circumstances are mutually exclusive.  After an
elementary special affine transformation, we can assume that the graph
$u=F(x,y)$ passes through the origin $p = 0$ and that $\Hessian_F(0)$
is as shown below.

\medskip

\begin{center}
\begin{tabular}{r|c|l}
\hline
$F(x,y)$&$\Hessian_F$ at the origin& Type of the origin
\\\hline
${\rm O}(3)$&$\left( \begin{array}{cc}
0 & 0 \\
0 & 0
\end{array} \right)$ & ${\sf flat}$
\\\hline
$\frac{1}{2}x^2+{\rm O}(3)$&$\left( \begin{array}{cc}
1 & 0 \\
0 & 0
\end{array} \right)$ & ${\sf parabolic}$
\\\hline
$\frac{1}{2}(x^2+y^2)+{\rm O}(3)$&$\left( \begin{array}{cc}
1 & 0 \\
0 & 1
\end{array} \right)$ & ${\sf elliptic}$
\\\hline
$\frac{1}{2}(x^2-y^2)+{\rm O}(3)$&$\left( \begin{array}{cc}
1 & 0 \\
0 & -1
\end{array} \right)$ & ${\sf hyperbolic}$
\\\hline
\end{tabular}
\end{center}

\smallskip

The map:
\[
(x,y)
\,\,\longmapsto\,\,
\rank\,
\Hessian_F(x,y)
\]
is lower semicontinuous, namely, if the Hessian matrix
has rank $\tau_p$ with
$0 \leqslant \tau_p \leqslant 2$ at some point
$(x_p, y_p)$, then at all nearby points $q \sim p$, it
has rank $\tau_q \geqslant \tau_p$. If it has rank 2 at some point, then it has rank 2 in some neighborhood of that point.

If we agree to make rank constancy assumptions,
as we will do throughout this article, 
analytic surfaces $S^2 \subset \R^3$ 
can be classified as:

\smallskip\noindent$\bullet$
{\sl everywhere flat}, 
if the Hessian matrix is identically of rank $0$; 

\smallskip\noindent$\bullet$
{\sl everywhere parabolic}, 
if the Hessian matrix is identically of rank $1$;

\smallskip\noindent$\bullet$
{\sl everywhere nondegenerate}, 
if the Hessian matrix is everywhere of rank $2$.

\smallskip

Nevertheless, it remains in general 
two kinds of {\sl mixed types}, 
where a surface $\{ u = F(x,y) \}$ can be:

\smallskip\noindent$\bullet$
flat in a proper closed subset, and parabolic elsewhere,
in some dense open subset;

\smallskip\noindent$\bullet$
flat or parabolic in some proper closed subset, and elliptic
or hyperbolic elsewhere,
in some dense open subset.

\smallskip

In this paper, we will avoid studying mixed types,
because it would engage
towards singularity theory. 
Before we focus our attention on everywhere
parabolic surfaces, let us briefly review known works
about everywhere nondegenerate surfaces $S^2 \subset \R^3$. 

\Section{\bf Everywhere Elliptic or Hyperbolic Surfaces $S^2
\subset \R^3$: A Review}
\label{review-everywhere-elliptic-hyperbolic}
\HEAD{{\ref{review-everywhere-elliptic-hyperbolic}}.~{\sf 
Everywhere Elliptic or Hyperbolic Surfaces $S^2 \subset \R^3$:
A Review}
}{
Zhangchi {\sc Chen}, Joël~{\sc Merker}}

By lower semicontinuity of the rank of a matrix, if the origin is an
elliptic (or hyperbolic) point, then there exists a sufficiently small
neighborhood which is everywhere elliptic (or hyperbolic).

\begin{Theorem}
{\rm {\cite[III, p.~165]{Spivak-1979}}}
Under the action of the equi-affine group $\SA_3(\R)$,
every elliptic surface $S^2 \subset \R^3$ is 
equivalent to:
\[
u
\,=\,
{\textstyle{\frac{1}{2}}}\,
\big(
x^2+y^2
\big)
+
{\textstyle{\frac{\Caux}{6}}}\,
\big(
x^3
-
3\,xy^2
\big)
+
{\rm O}_{x,y}(4),
\]
while every hyperbolic surface is equivalent to one 
and only of the following three:
\[
\aligned
u
&
\,=\,
{\textstyle{\frac{1}{2}}}\,
\big(
x^2-y^2
\big)
+
{\textstyle{\frac{\Caux}{6}}}\,
\big(
x^3
+
3\,xy^2
\big)
+
{\rm O}_{x,y}(4),
\\
u
&
\,=\,
{\textstyle{\frac{1}{2}}}\,
\big(
x^2-y^2
\big)
+
{\textstyle{\frac{\Caux}{6}}}\,
\big(
3\,x^2y
+
y^3
\big)
+
{\rm O}_{x,y}(4),\\
u
&
\,=\,
{\textstyle{\frac{1}{2}}}\,
\big(
x^2-y^2
\big)
+
{\textstyle{\frac{1}{6}}}\,
\big(
x+y
\big)^3
+
{\rm O}_{x,y}(4),
\endaligned
\]
where $\Caux$ is unique up to sign.
\end{Theorem}

This $\Caux$ is a Taylor coefficient at the origin,
but the method of the next Section~{\ref{power-series-method}}
will show, thanks to the fact that the action of
$\SA_3(\R)$ is (trivially) transitive, 
that the computation of $\Caux$ for a
power series at the origin provides the
expression of a corresponding differential invariant
at every point $(x,y)$. 

\begin{Definition}
The quantity:
\[
\Paux
\,:=\,
{\textstyle{\frac{1}{2}}}\,
\Caux^2
\] 
is a 3\textsuperscript{rd}-order equi-affine invariant, called
the {\sl Pick invariant}.
\end{Definition}

Its explicit expression, at an elliptic point, is:
\[
\label{explicit-Pick}
\aligned
\Paux
&
\,=\,
\frac{1}{512}\,
\frac{1}{\big(u_{xx}u_{yy}-u_{xy}^2\big)^{11/2}}\,
\Big(
-\,18\,
u_{xx}\,u_{xxy}\,u_{xy}\,u_{xyy}\,u_{yy}
+
12\,
u_{xxx}\,u_{xy}^2\,u_{xyy}\,u_{yy}
\,+
\\
&
\ \ \ \ \ \ \ \ \ \ \ \ \ \ \ \ \ \ \ \ \ \ \ \ \ \ \ \ \ \ \ \ \ \ 
\ \ \ \ \ \ \ \ \ \ \
+
9\,
u_{xx}^2\,u_{xyy}^2\,u_{yy}
+
9\,
u_{xx}\,u_{xxy}^2\,u_{yy}^2
-
6\,
u_{xxx}\,u_{xxy}\,u_{xy}\,u_{yy}^2
\,-
\\
&
\ \ \ \ \ \ \ \ \ \ \ \ \ \ \ \ \ \ \ \ \ \ \ \ \ \ \ \ \ \ \ \ \ \ 
\ \ \ \ \ \ \ \ \ \ \
-\,6\,
u_{xx}\,u_{xxx}\,u_{xyy}\,u_{yy}^2
+
u_{xxx}^2\,u_{yy}^3
+
12\,
u_{xx}\,u_{xxy}\,u_{xy}^2\,u_{yyy}
\,-
\\
&
\ \ \ \ \ \ \ \ \ \ \ \ \ \ \ \ \ \ \ \ \ \ \ \ \ \ \ \ \ \ \ \ \ \ 
\ \ \ \ \ \ \ \ \ \ \
-\,8\,
u_{xxx}\,u_{xy}^3\,u_{yyy}
-
6\,
u_{xx}^2\,u_{xy}\,u_{xyy}\,u_{yyy}
-
6\,
u_{xx}^2\,u_{xxy}\,u_{yy}\,u_{yyy}
\,+
\\
&
\ \ \ \ \ \ \ \ \ \ \ \ \ \ \ \ \ \ \ \ \ \ \ \ \ \ \ \ \ \ \ \ \ \ 
\ \ \ \ \ \ \ \ \ \ \
+
6\,
u_{xx}\,u_{xxx}\,u_{xy}\,u_{yy}\,u_{yyy}
+
u_{xx}^3\,u_{yyy}^2
\Big)^2.
\endaligned
\]
At a hyperbolic point, we replace the factor $\big(u_{xx}u_{yy}-u_{xy}^2\big)^{-11/2}$ by $\big(u_{xy}^2-u_{xx}u_{yy}\big)^{11/2}$.

When $\Caux$ is nonzero, we may assume $\Caux>0$. Then the only
element in $\SA_3(\R)$ fixing the standard form above is the identity,
hence all the coefficients in the Taylor expansion of ${\rm O}_{x,y}
(4)$ are also differential invariants.

Under some non-degeneracy conditions, Olver proved
in~{\cite{Olver-2007}} that all those higher order differential
invariants can be generated by $\Caux$ and its differentials. Once
$\Caux$ is captured 
in a small neighborhood of the origin,
all the differential
invariants are known there.

\Section{\bf Parabolic Jet Relations}
\label{parabolic-jet-relations}
\HEAD{{\ref{parabolic-jet-relations}}.~{\sf Parabolic Jet Relations}
}{
Zhangchi {\sc Chen}, Joël~{\sc Merker}}

Take a graph $\big\{ u = F(x,y) \big\}$, and assume that
the Hessian matrix of $F$: 
\[
\Hessian_F
\,=\,
\left(\!
\begin{array}{cc}
F_{xx} & F_{xy}
\\
F_{yx} & F_{yy}
\end{array}
\!\right)
\]
has rank $1$ at every point. After an rotation 
in the $(x,y)$-space (if necessary), 
this assumption amounts to:
\[
F_{xx}
\,\neq\,
0
\,\equiv\,
F_{xx}\,F_{yy}
-
F_{xy}^2.
\]
Therefore, we can solve:
\leqnomode\usetagform{default}
\begin{align}
\label{solve-F-yy}
F_{yy}
\,=\,
\frac{F_{xy}^2}{F_{xx}}.
\end{align}

Jet spaces will be equipped with coordinates denoted:
\[
\Big(
x,y,\,u,\,\,
u_{1,0},u_{0,1},\,\,
u_{2,0},u_{1,1},u_{0,2},\,\,
u_{3,0},u_{2,1},u_{1,2},u_{0,3},\,\,
\dots\dots
\Big).
\]
Sometimes, we will fix an order $\oorder \geqslant 0$.
At first, we have to express all the differential consequences
of the resolution~({\ref{solve-F-yy}})
for $F_{yy}$.

\begin{center}
\scalebox{1.25}{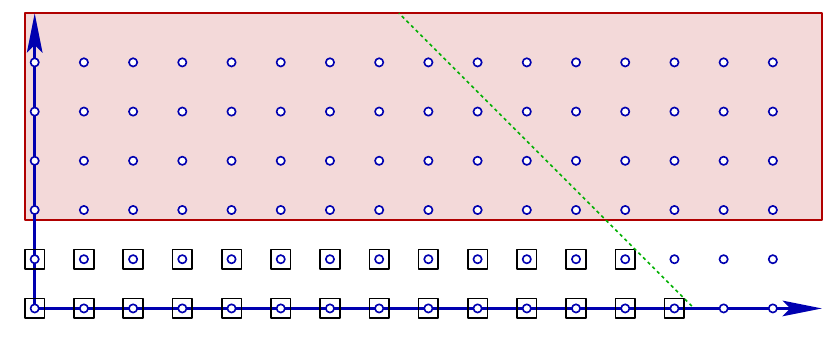}
\end{center}

For jets of order $3$, we obtain by differentiating 
with respect to $x$ and to $y$ and by performing
replacements:
\[
\aligned
F_{xyy}
&
\,=\,
2\,
\frac{F_{xy}\,F_{xxy}}{F_{xx}}
-
\frac{F_{xy}^2\,F_{xxx}}{F_{xx}^2},
\\
F_{yyy}
&
\,=\,
3\,
\frac{F_{xy}^2\,F_{xxy}}{F_{xx}^2}
-
2\,
\frac{F_{xy}^3\,F_{xxx}}{F_{xx}^3},
\endaligned
\]
and also for jets of order $4$:
\[
\footnotesize
\aligned
F_{xxyy}
&
\,=\,
-\,
4\,\frac{F_{xy}\,F_{xxy}\,F_{xxx}}{F_{xx}^2}
+
2\,\frac{F_{xy}\,F_{xxxy}}{F_{xx}}
+
2\,\frac{F_{xxy}^2}{F_{xx}}
+
2\,\frac{F_{xy}^2\,F_{xxx}^2}{F_{xx}^3}
-
\frac{F_{xy}^2\,F_{xxxx}}{F_{xx}^2},
\\
F_{xyyy}
&
\,=\,
-\,
12\,\frac{F_{xy}^2\,F_{xxx}\,F_{xxy}}{F_{xx}^3}
+
6\,\frac{F_{xy}\,F_{xxy}^2}{F_{xx}^2}
+
6\,\frac{F_{xy}^3\,F_{xxx}^2}{F_{xx}^4}
+
3\,\frac{F_{xy}^2\,F_{xxxy}}{F_{xx}^2}
-
2\,\frac{F_{xy}^3\,F_{xxxx}}{F_{xx}^3},
\\
F_{yyyy}
&
\,=\,
-\,
24\,\frac{F_{xy}^3\,F_{xxx}\,F_{xxy}}{F_{xx}^4}
+
12\,\frac{F_{xy}^2\,F_{xxy}^2}{F_{xx}^3}
+
4\,\frac{F_{xy}^3\,F_{xxxy}}{F_{xx}^3}
+
12\,\frac{F_{xy}^4\,F_{xxx}^2}{F_{xx}^5}
-
3\,\frac{F_{xy}^4\,F_{xxxx}}{F_{xx}^4}.
\endaligned
\]
Similar quite longer formulas exist for 
$F_{xxxyy}$, $F_{xxyyy}$, $F_{xyyyy}$, $F_{yyyyy}$,
for the memory of a computer, this is nothing,
even for jets up to order $7$.
It is also easy to prove by induction the

\begin{Observation}
\label{Obs-horizontal-parabolic-jets}
Every partial derivative $F_{x^k y^l}$ 
in the red region with $l \geqslant 2$
expresses rationally in terms of the partial derivatives
in the black region:
\[
\big\{
F_{x^{k'}}
\big\}_{k'\leqslant k+l},
\ \ \ \ \
\big\{
F_{x^{k''}y}
\big\}_{k''\leqslant k+l-1},
\]
with denominators containing only powers $\big( F_{xx} \big)^\ast$.\qed
\end{Observation}

Using the jet notation, this means that we will exclusively work
in the {\em submanifold} of the jet spaces $J_{x,u}^{\oorder}$
defined by:
\[
u_{0,2}
\,=\,
\frac{u_{1,1}^2}{u_{2,0}},
\]
next:
\[
\aligned
u_{1,2}
&
\,=\,
2\,\frac{u_{1,1}\,u_{2,1}}{u_{2,0}}
-
\frac{u_{1,1}^2\,u_{3,0}}{u_{2,0}^2},
\\
u_{0,3}
&
\,=\,
3\,
\frac{u_{1,1}^2\,u_{2,1}}{u_{2,0}^2}
-
2\,\frac{u_{1,1}^3\,u_{3,0}}{u_{2,0}^3},
\endaligned
\]
further:
\[
\aligned
u_{2,2}
&
\,=\,
2\,
\frac{u_{2,1}^2}{u_{2,0}}
-
4\,
\frac{u_{1,1}\,u_{2,1}\,u_{3,0}}{u_{2,0}^2}
+
2\,
\frac{u_{1,1}\,u_{3,1}}{u_{2,0}}
+
2\,
\frac{u_{1,1}^2\,u_{3,0}^2}{u_{2,0}^3}
-
\frac{u_{1,1}^2\,u_{4,0}}{u_{2,0}^2},
\\
u_{1,3}
&
\,=\,
6\,
\frac{u_{1,1}\,u_{2,1}^2}{u_{2,0}^2}
-
12\,
\frac{u_{1,1}^2\,u_{3,0}\,u_{2,1}}{u_{2,0}^3}
+
3\,\frac{u_{1,1}^2\,u_{3,1}}{u_{2,0}^2}
+
6\,
\frac{u_{1,1}^3\,u_{3,0}^2}{u_{2,0}^4}
-
2\,
\frac{u_{1,1}^3\,u_{4,0}}{u_{2,0}^3},
\\
u_{0,4}
&
\,=\,
12\,
\frac{u_{1,1}^2\,u_{2,1}^2}{u_{2,0}^3}
-
24\,
\frac{u_{1,1}^3\,u_{2,1}\,u_{3,0}}{u_{2,0}^4}
+
12\,
\frac{u_{1,1}^4\,u_{3,0}^2}{u_{2,0}^5}
+
4\,
\frac{u_{1,1}^3\,u_{3,1}}{u_{2,0}^3}
-
3\,
\frac{u_{1,1}^4\,u_{4,0}}{u_{2,0}^4},
\endaligned
\]
and so on. Again, this is nothing on a computer.

\Section{\bf In Search of a Resolved Cross-Section for 
Parabolic Surfaces $S^2 \subset \R^3$}
\label{search-resolved-cross-section-parabolic}
\HEAD{{\ref{search-resolved-cross-section-parabolic}}.~{\sf In 
Search of a Resolved Cross-Section for 
Parabolic Surfaces $S^2 \subset \R^3$}
}{
Zhangchi {\sc Chen}, Joël~{\sc Merker}}

Now, if we let a general affine transformation:
\[
\aligned
s
&
\,=\,
{\sf a}\,x
+
{\sf b}\,y
+
{\sf c}\,u
+
{\sf d},
\notag
\\
t
&
\,=\,
{\sf k}\,x
+
{\sf l}\,y
+
{\sf m}\,u
+
{\sf n},
\\
v
&
\,=\,
{\sf p}\,x
+
{\sf q}\,y
+
{\sf r}\,u
+
{\sf s},
\notag
\endaligned
\ \ \ \ \ \ \ \ \ \ \ \ \ \ \ \ \ \ \ \
\delta
\,=\,
\left\vert\!
\begin{array}{ccc}
{\sf a} & {\sf b} & {\sf c}
\\
{\sf k} & {\sf l} & {\sf m}
\\
{\sf p} & {\sf q} & {\sf r}
\end{array}
\!\right\vert
\,\neq\,
0,
\]
act on graphed surfaces, the first prolongation formulas are:
\[
\aligned
v_{1,0}
&
\,=\,
\frac{{\sf l}{\sf p}
-
{\sf k}{\sf q}
+
\big(
{\sf l}{\sf r}
-
{\sf m}{\sf q}
\big)\,
u_{1,0}
+
\big(
{\sf m}{\sf p}
-
{\sf k}{\sf r}
\big)\,
u_{0,1}}{
{\sf a}{\sf l}
-
{\sf b}{\sf k}
+
\big(
{\sf c}{\sf l}
-
{\sf b}{\sf m}
\big)\,
u_{1,0}
+
\big(
{\sf a}{\sf m}
-
{\sf c}{\sf k}
\big)\,
u_{0,1}},
\\
v_{0,1}
&
\,=\,
\frac{{\sf a}{\sf q}
-
{\sf b}{\sf p}
+
\big(
{\sf c}{\sf q}
-
{\sf b}{\sf r}
\big)\,
u_{1,0}
+
\big(
{\sf a}{\sf r}
-
{\sf c}{\sf p}
\big)\,
u_{0,1}}{
{\sf a}{\sf l}
-
{\sf b}{\sf k}
+
\big(
{\sf c}{\sf l}
-
{\sf b}{\sf m}
\big)\,
u_{1,0}
+
\big(
{\sf a}{\sf m}
-
{\sf c}{\sf k}
\big)\,
u_{0,1}},
\endaligned
\]
and we of course recognize from
({\ref{definition-Lambda}}) the denominator:
\[
\Lambda
\,=\,
{\sf a}{\sf l}
-
{\sf b}{\sf k}
+
\big(
{\sf c}{\sf l}
-
{\sf b}{\sf m}
\big)\,
u_{1,0}
+
\big(
{\sf a}{\sf m}
-
{\sf c}{\sf k}
\big)\,
u_{0,1}.
\]

Next, for second-order jets $v_{2,0}$, $v_{1,1}$, $v_{0,2}$,
we consider {\em only} $v_{2,0}$, $v_{1,1}$, because the
jet $v_{0,2} = v_{1,1}^2 \big/ v_{2,0}$ is dependent.
Furthermore, we have to 
take account of the parabolic jet relations 
explained in 
Section~{\ref{parabolic-jet-relations}}.
Using a computer, we obtain:
\[
\aligned
v_{2,0}
&
\,=\,
\frac{\delta}{\Lambda^3\,u_{2,0}}\,\,
\Big(
\underbrace{
{\sf l}\,u_{2,0}
+
{\sf m}\,u_{0,1}\,u_{2,0}
-
{\sf k}\,u_{1,1}
-
{\sf m}\,u_{1,0}\,u_{1,1}}_{
=\,\,\Pi}
\Big)^2,
\\
v_{1,1}
&
\,=\,
\frac{\delta}{\Lambda^3\,u_{2,0}}\,
\Pi\,\,
\big(
{\sf c}\,u_{1,0}\,u_{1,1}
+
{\sf a}\,u_{1,1}
-
{\sf b}\,u_{2,0}
-
{\sf c}\,u_{0,1}\,u_{2,0}
\big),
\endaligned
\]
while the formulas for the two independent third-order jets 
$v_{3,0}$ and $v_{2,1}$ start to become large:
\[
\aligned
v_{3,0}
&
\,=\,
\frac{\delta}{\Lambda^5\,u_{2,0}^3}\,
\Pi^2\,
\Big(
\text{\rm {\bf 48} monomials homogeneous of degree}\,\,3\,\,
\text{\rm in}\,\,
\begin{smallmatrix} 
{\sf a}, & {\sf b}, & {\sf c}, & {\sf d},
\\ 
{\sf k}, & {\sf l}, & {\sf m}, & {\sf n}, 
\\ 
{\sf p}, & {\sf q}, & {\sf r}, & {\sf s} 
\end{smallmatrix}
\Big),
\\
v_{2,1}
&
\,=\,
\frac{\delta}{\Lambda^5\,u_{2,0}^3}\,
\Pi\,
\Big(
\text{\rm {\bf 135} monomials homogeneous of degree}\,\,4\,\,
\text{\rm in}\,\,
\begin{smallmatrix} 
{\sf a}, & {\sf b}, & {\sf c}, & {\sf d},
\\ 
{\sf k}, & {\sf l}, & {\sf m}, & {\sf n}, 
\\ 
{\sf p}, & {\sf q}, & {\sf r}, & {\sf s} 
\end{smallmatrix}
\Big).
\endaligned
\]
and the last two useful formulas start to become huge:
\[
\aligned
v_{4,0}
&
\,\,\ni\,\,
\text{\rm factor in the numerator containing {\bf 720} monomials},
\\
v_{4,1}
&
\,\,\ni\,\,
\text{\rm factor in the numerator containing {\bf 14\,156} monomials}.
\endaligned
\]

We were not able, on a computer, to solve 
for the $11$ group parameters
$\Big(
\begin{smallmatrix} 
{\sf a}, & {\sf b}, & {\sf c}, & {\sf d},
\\ 
{\sf k}, & {\sf l}, & {\sf m}, & {\sf n}, 
\\ 
{\sf p}, & {\sf q}, & {\sf r}, & {\sf s} 
\end{smallmatrix}
\Big)$
from the $11$ natural cross-section equations:
\[
\aligned
s
\,=\,
0,
\ \ \ \ \ \ \ \ \ 
t
\,=\,
0,
\ \ \ \ \ \ \ \ \ 
v
\,=\,
0,
\ \ \ \ \ \ \ \ \ 
v_{1,0}
&
\,=\,
0,
\ \ \ \ \ \ \ \ \ 
v_{0,1}
\,=\,
0,
\\
v_{2,0}
&
\,=\,
1,
\ \ \ \ \ \ \ \ \ 
v_{1,1}
\,=\,
0,
\\
v_{3,0}
&
\,=\,
0,
\ \ \ \ \ \ \ \ \ 
v_{2,1}
\,=\,
1,
\\
v_{4,0}
&
\,=\,
0,
\\
&
\ \ \ \ \ \ \ \ \ \ \ \ \ \ \ \ \ \ \ 
v_{4,1}
\,=\,
0.
\endaligned
\]
Fortunately, an alternative, more progressive, method works,
as we will see in
Sections~{\ref{parabolic-pseudostablization}},
{\ref{relative-invariant-S-first-invariant-W}},
{\ref{branch-W-equiv-0-branch-W-nonzero}},
{\ref{recurrence-relations-parabolic-surfaces}}.

By solving order-by-order the cross-section equations, they 
progressively simplify\,\,---\,\,a lot!\,\,---,
and we never have to deal with huge expresssions.

\Section{\bf The Power Series Method}
\label{power-series-method}
\HEAD{{\ref{power-series-method}}.~{\sf 
The Power Series Method}
}{
Zhangchi {\sc Chen}, Joël~{\sc Merker}}

Consider a local Lie group $G_0$ acting on $\R_{x,u}^{p+q}$,
and assume that all of its elements fix the origin
$z = 0$:
\[
g_0
\cdot
0
\,=\,
0
\eqno
{\scriptstyle{(\forall\,g_0\,\in\,G_0)}}.
\]
As usual, we denote coordinates on the target space 
$\R_{y,v}^{p+q}$ as $w = (y,v)$.
For the moment, we do not necessarily assume
that $G_0$ is the isotropy subgroup of the origin
for the action of a certain larger group $G \supset G_0$.

When working with power series, we will abandon the notation $J = j_1,
\dots, j_\lambda$ used for jet spaces in the preceding sections, with
$1 \leqslant j_1, \dots, j_\lambda \leqslant p$ not recording
repetitions, and instead, we will employ the standard {\sl
multi-index} notation\,\,---\,\,with the same letter(s)\,\,---:
\[
J
\,=\,
\big(j_1,\dots,j_p\big)
\,\in\,
\N^p.
\]
The advantage is that 
we can introduce useful quantities which would be 
otherwise difficult to denote:
\[
J!
\,:=\,
j_1!\,\cdots\,j_p!
\ \ \ \ \ \ \ \ \ \ \ \ \ \ \ \ \ \ \ \
\text{and}
\ \ \ \ \ \ \ \ \ \ \ \ \ \ \ \ \ \ \ \
\vert
J
\vert
\,:=\,
j_1+\cdots+j_p,
\]
and also:
\[
x^J
\,:=\,
x^{j_1}\cdots x_p^{j_p}.
\]

As before, we denote the $G_0$-action by $w = g_0 \cdot z$, or
equivalently:
\[
y
\,=\,
y\big(g_0,x,u),
\ \ \ \ \ \ \ \ \ \ \ \ \ \ \ \ \ \ \ \
v
\,=\,
v\big(g_0,x,u).
\]

Now, consider a (converging) power series mapping 
$\R^p \ni x \longmapsto
F(x) \in \R^q$ in the source space which vanishes at the origin:
\[
u^\alpha
\,=\,
\sum_{\vert J\vert\geqslant1}\,
F_J^\alpha\,
\frac{x^J}{J!}.
\]
Every diffeomorphism $z \longmapsto g_0 \cdot z =: w$ corresponding to
a group parameter $g_0 \in G_0$ close to the identity element then
transforms the graph $\big\{ u = F(x) \big\}$ of this power series
into another graph $\big\{ v = G(g_0, y) \big\}$ depeding on $g_0$,
still passing through the origin $(y,v) = (0,0)$, whose graphing
function also has a power series expansion:
\[
v^\beta
\,:=\,
\sum_{\vert K\vert\geqslant 1}\,
G_K^\beta\,
\frac{y^\beta}{K!},
\]
with coefficients:
\[
G_K^\beta
\,=\,
G_K^\beta
\Big(
g_0,\,
\big\{
F_L^\gamma
\big\}_{1\leqslant\vert L\vert\leqslant\vert K\vert}^{
1\leqslant\gamma\leqslant q}
\Big),
\]
depending on the coefficients of the source power series, and on the
group parameters as well\,\,---\,\,of course.  One easily convinces
oneself that, as is written, the $G_K^\beta$ only depend on power
series coefficients $F_J^\alpha$ of order $\vert J \vert \leqslant
\vert K\vert$.

\begin{Definition}
\label{Def-power-series-invariant}
A {\sl power series invariant} of order $n \geqslant 1$ is a function
of the Taylor coefficients:
\[
\Iaux
\,=\,
\Iaux\,
\Big(
\big\{
F_J^\alpha
\big\}_{1\leqslant\vert J\vert\leqslant n}^{
1\leqslant\alpha\leqslant q}
\Big),
\]
which has unchanged value after the action of 
{\em any} element in our local Lie group:
\[
\Iaux
\left(
\Big\{
G_K^\beta
\big(
g_0,\,
\big\{
F_L^\gamma
\big\}_{1\leqslant\vert L\vert\leqslant\vert K\vert}^{
1\leqslant\gamma\leqslant q}
\big)
\Big\}_{1\leqslant\vert K\vert\leqslant n}^{
1\leqslant\beta\leqslant q}
\right)
\,\,=\,\,
\Iaux\,
\Big(
\big\{
F_J^\alpha
\big\}_{1\leqslant\vert J\vert\leqslant n}^{
1\leqslant\alpha\leqslant q}
\Big).
\eqno
{\scriptstyle{(\forall\,g_0\,\in\,G_0)}}.
\]
\end{Definition}

This concept has a meaning only {\em at the origin!} 
It only concerns the derivatives 
{\em at} $0$ of the graphing functions ${\tt u}^\alpha(x)$:
\[
\big\{
{\tt u}_{x^J}^\alpha
(0)
\big\}_{J\in\N_\ast^p}^{
1\leqslant\alpha\leqslant q}.
\] 
By contrast, the general theory of differential invariants
is able to handle derivatives at {\em all points} in the 
source horizontal space:
\[
\big\{
{\tt u}_{x^J}^\alpha
(x)
\big\}_{J\in\N_\ast^p,\,\,\,
x\,\,
\text{\rm varies in}\,\,\R^p}^{
1\leqslant\alpha\leqslant q}.
\]

So it seems that this notion of {\sl power series invariant} is quite
restrictive! But a bit paradoxically\,\,---\,\,and quite
the contrary!\,\,---, we will 
rapidly realize that the power series invariants do capture all
differential invariants at any point $x \in \R^p$, 
provided only that the larger group $G \supset G_0 = \Iso(G, 0)$ 
contains the ambient translations.

In~{\cite{Olver-2018}}, it is shown how:
\[
\text{\footnotesize\sf
invariants in jet spaces}\,\,\,
\mathbin{\scalebox{1.37}{\ensuremath{\leadsto}}}\,\,\,
\text{\footnotesize\sf
power series invariants}.
\]
Our goal now is to explore the reverse transmission:
\[
\text{\footnotesize\sf
invariants in jet spaces}\,\,\,
\reflectbox{$\mathbin{\scalebox{1.37}{\ensuremath{\leadsto}}}$}\,\,\,
\text{\footnotesize\sf
power series invariants},
\]
which will bring some computational advantages.

\smallskip

Taking a local Lie group $G$ acting on a neighborhood of $0
\in \R_{x,u}^{p+q}$, not necessarily fixing the origin, 
we will make two kinds of assumptions.

\begin{Hypothesis}
The group $G$ contains all translations of the ambient space:
\[
\big(
x^i,u^\alpha
\big)
\,\,\,\longmapsto\,\,\,
\big(
x^i
+
{\sf a}^i,\,\,
u^\alpha
+
{\sf b}^\alpha
\big),
\]
whence $\dim\, G \geqslant p + q$.
\end{Hypothesis}

Equivalently, the Lie algebra $\mathfrak{g} = \Lie(G)$ 
of its action contains
all the unit coordinate infinitesimal generators:
\[
\partial_{x^i},
\ \ \ \ \ \ \ \ \ \ \ \ \ \ \ \ \ \ \ \
\partial_{u^\alpha}.
\]
The proof of the following result is elementary.

\begin{Theorem}
\label{Thm-translations-G-0}
If $G$ contains all ambient translations, 
then for any jet order $n \geqslant 0$, 
all differential invariants of $G$ are
independent of $x^i$, $u^\alpha$:
\[
\Iaux
\,=\,
\Iaux
\Big(
\big\{
u_{x^J}^\beta
\big\}_{1\leqslant\vert J\vert\leqslant n}^{
1\leqslant\beta\leqslant q}
\Big),
\]
and there is a one-to-one correspondence:
\[
\text{\footnotesize\sf
Differential invariants of}\,\,
G\,\,\,
\mathbin{\scalebox{1.37}{\ensuremath{\longleftrightarrow}}}\,\,\,
\text{\footnotesize\sf
Power series invariants of}\,\,
G.
\eqno\qed
\]
\end{Theorem}

More precisely\,\,---\,\,and this is the interesting aspect!\,\,---,
the theorem tells us that 
starting from a power series invariant determined in a way
that will be explained later on:
\[
\Iaux
\,=\,
\Iaux
\Big(
\big\{
F_J^\alpha
\big\}_{2\leqslant\vert J\vert\leqslant n}^{
1\leqslant\alpha\leqslant q}
\Big)
\,=\,
\Iaux
\Big(
\big\{
{\tt u}_{x^J}^\alpha
(0)
\big\}_{2\leqslant\vert J\vert\leqslant n}^{
1\leqslant\alpha\leqslant q}
\Big),
\]
one deduces instantly a {\em differential invariant} by just
replacing the origin with any horizontal $x \in \R^p$:
\[
\Iaux
\,=\,
\Iaux
\Big(
\big\{
{\tt u}_{x^J}^\alpha
(x)
\big\}_{2\leqslant\vert J\vert\leqslant n}^{
1\leqslant\alpha\leqslant q}
\Big),
\]
and in particular, at the origin, one recovers the starting
power series invariant.

\begin{Example}
Consider the action of the special Euclidean group 
$\SE_2(\R) := \SO_2(\R) \ltimes \R^2$ on curves
$\{ u = F(x) \}$ in $\R_{x,u}^2$. Because
this action is transitive, we can restrict our attention 
to curves passing through $(0,0)$
and to rotations:
\[
\aligned
x
&
\,=\,
\cos\,\theta\,\,
y
-
\sin\,\theta\,\,
v,
\\
u
&
\,=\,
\sin\,\theta\,\,
y
+
\cos\,\theta\,\,
v,
\endaligned
\ \ \ \ \ \ \ \ \ \ \ \ \ \ \ \ \ \ \ \
\text{abbreviated as}
\ \ \ \ \ \ \ \ \ \ \ \ \ \ \ \ \ \ \ \
\aligned
x
&
\,=\,
{\sf c}\,y
-
{\sf s}\,v,
\\
u
&
\,=\,
{\sf s}\,y
+
{\sf c}\,v,
\endaligned
\]
with ${\sf s}^2 + {\sf c}^2 = 1$. When we deal with the local Lie group near the identity $(s,c)=(0,1)$, we may assume $c>0$.
Two graphed power series:
\[
u
\,=\,
F(x)
\,=\,
F_1\,
{\textstyle{\frac{x^1}{1!}}}
+
F_2\,
{\textstyle{\frac{x^2}{2!}}}
+
\cdots
\ \ \ \ \ \ \ \ \ \ \ \ \ \ \ \ \ \ \ \
\text{and}
\ \ \ \ \ \ \ \ \ \ \ \ \ \ \ \ \ \ \ \
v
\,=\,
G(y)
\,=\,
G_1\,
{\textstyle{\frac{y^1}{1!}}}
+
G_2\,
{\textstyle{\frac{y^2}{2!}}}
+
\cdots
\]
are mapped one to another if and only if:
\[
0
\,\equiv\,
-\,
{\sf s}\,y
-
{\sf c}\,G(y)
+
F
\big(
{\sf c}\,y
-
{\sf s}\,G(y)
\big)
\eqno
{\scriptstyle{(\text{\rm in}\,\C\{y\})}}.
\]

In this identity, the first and second order terms read:
\[
0
\,\equiv\,
\big(
-\,{\sf s}
-
{\sf c}\,G_1
+
{\sf c}\,F_1
-
{\sf s}\,F_1G_1
\big)\,
{\textstyle{\frac{y^1}{1!}}}
+
\big(
-\,
{\sf c}\,G_2
-
{\sf s}\,F_1G_2
+
{\sf c}^2\,F_2
-
2\,{\sf c}{\sf s}\,F_2G_1
+
{\sf s}^2\,F_2G_1^2
\big)\,
{\textstyle{\frac{y^2}{2!}}}
+
{\rm O}_y(3).
\]
We can make $G_1 := 0$ thanks to the choice of $\theta \in \R$
satisfying:
\[
0
\,=\,
-\,
\sin\,\theta
+
\cos\,\theta\,
F_1
\ \ \ \ \ \ \ \ \ \ \ \ \ \ \ \ \ \ \ \
\Longleftrightarrow
\ \ \ \ \ \ \ \ \ \ \ \ \ \ \ \ \ \ \ \
{\textstyle{\frac{{\sf s}}{{\sf c}}}}
\,=\,
F_1.
\]

Next, we may verify that the transformation of $\SO_2(\R)$ with $c>0$
which stabilizes this normalization of first order terms,
namely sends:
\[
u
\,=\,
0
+
F_2\,
{\textstyle{\frac{x^2}{2!}}}
+
\cdots
\ \ \ \ \ \ \ \ \ \ \ \ \ \ \ \ \ \ \ \
\text{to}
\ \ \ \ \ \ \ \ \ \ \ \ \ \ \ \ \ \ \ \
v
\,=\,
0
+
G_2\,
{\textstyle{\frac{y^2}{2!}}}
+
\cdots,
\]
is the identity.

Looking at the second order term above, with $G_1 = 0$,
we get:
\[
G_2
\,=\,
F_2\,
\frac{{\sf c}}{1+\frac{{\sf s}}{\sf c}\,F_1}
\,=\,
F_2\,
\frac{{\sf c}}{1+F_1^2},
\]
and using $1 = {\sf s}^2 + {\sf c}^2 = {\sf c}^2\, F_1^2 + {\sf c}^2$,
which gives ${\sf c} = \frac{1}{\sqrt{1+F_1^2}}$, 
we conclude that we have obtained a power series invariant:
\[
v_{yy}(0)
\,=\,
u_{xx}(0)\,
\frac{1}{\big(1+u_x(0)^2\big)^{3/2}},
\]
which yields a differential invariant at any point $x$:
\[
u_{xx}(x)\,
\frac{1}{\big(1+u_x(x)^2\big)^{3/2}},
\]
We have thus recovered the Euclidean curvature of curves in the plane
by applying the power series method. This method also works for the full Euclidean group $\SE_2(\R)$, not just the local one. A similar calculation for $c\in[-1,1]$ instead of $c>0$ shows that there are two normal forms: $v=0+G_2\frac{y^2}{2}+\cdots$ and $v=0-G_2\frac{y^2}{2}-\cdots$, for every nonlinear curve. These two forms are equivalent to each other by a semi-circle rotation $(s,c)=(0,1)$. So our power series invariant for the full group of order 2 is $|G_2|=\frac{|F_2|}{(1+F_1^2)^{3/2}}$, which yields a differential invariant at any point $x$:
\[
|u_{xx}(x)|\,
\frac{1}{\big(1+u_x(x)^2\big)^{3/2}}.
\]

Beyond this (too) simple
example, we will see how the method works in many other contexts.
\end{Example}

Furthermore, the interesting reverse transmission:
\[
\text{\footnotesize\sf
invariants in jet spaces}\,\,\,
\reflectbox{$\mathbin{\scalebox{1.37}{\ensuremath{\leadsto}}}$}\,\,\,
\text{\footnotesize\sf
power series invariants},
\]
holds in other situation, 
as is expressed by our second

\begin{Hypothesis}
The group $G$ contains all translations of the ambient space,
{\em and} all {\sl vertical transvections:}
\[
\big(
x^i,u^\alpha
\big)
\,\,\,\longmapsto\,\,\,
\big(
x^i,\,\,
u^\alpha
+
\smallsum{1\leqslant i\leqslant p}\,
{\sf c}_i^\alpha\,
x^i
\big),
\]
whence $\dim\, G \geqslant p + q + p\,q$.
\end{Hypothesis}

Equivalently, the Lie algebra $\mathfrak{g} = \Lie(G)$ 
contains all the infinitesimal generators:
\[
\partial_{x^i},
\ \ \ \ \ \ \ \ \ \ \ \ \ \ \ \ \ \ \ \
\partial_{u^\alpha},
\ \ \ \ \ \ \ \ \ \ \ \ \ \ \ \ \ \ \ \
x^j\,\partial_{u^\beta}.
\]

\begin{Theorem}
\label{Thm-translations-transvections-G-1}
If $G$ contains all ambient translations
{\em and} all vertical transvections, 
then for any jet order $n \geqslant 0$, 
all differential invariants of $G$ are
independent of $x^i$, $u^\alpha$, $u_{x^j}^\beta$:
\[
\Iaux
\,=\,
\Iaux
\Big(
\big\{
u_{x^J}^\gamma
\big\}_{2\leqslant\vert J\vert\leqslant n}^{
1\leqslant\gamma\leqslant q}
\Big),
\]
and there is a one-to-one correspondence:
\[
\text{\footnotesize\sf
Differential invariants of}\,\,
G\,\,\,
\mathbin{\scalebox{1.37}{\ensuremath{\longleftrightarrow}}}\,\,\,
\text{\footnotesize\sf
Power series invariants of}\,\,
G.
\eqno\qed
\]
\end{Theorem}

More precisely, the concerned power series invariants can be 
searched for power series starting only at order $2$:
\[
u^\alpha
\,=\,
\sum_{\vert J\vert\geqslant 2}\,
F_J^\alpha\,
\frac{x^J}{J!}
\ \ \ \ \ \ \ \ \ \ \ \ \ \ \ \ \ \ \ \
\text{and}
\ \ \ \ \ \ \ \ \ \ \ \ \ \ \ \ \ \ \ \
v^\beta
\,=\,
\sum_{\vert K\vert\geqslant 2}\,
G_K^\beta\,
\frac{y^K}{K!}.
\]

\smallskip

At least three aspects of the power series method are attractive.

\smallskip\noindent$\bullet$\,
Prolongations of a $G_0$-action to jet spaces
of any order $n \geqslant 0$ can be performed
`automatically', especially on a computer, for there is no need to
write down the complicated formulas ${\sf E}_k ( v^\alpha )$ of
Theorem~{\ref{Thm-prolongation-E-j}},
they are implicitly performed after replacement of
$y = y(g_0, x, u)$ and $v = v(g_0, x, u)$ 
inside $0 = -\, v + G(y)$.

\smallskip\noindent$\bullet$\,
The search for power series invariant takes place above only $1$
point, the origin $0 \in \R^m$, {\em no differentiations are required
anymore}, just computations with Taylor coefficients, which are
constants, that is to say, the\,\,---\,\,nonlinear!\,\,---\,\,action of
$G_0$ is considered just on a vector space, not on a (jet) bundle.

\smallskip\noindent$\bullet$\,
In presence of differential relations, as will be illustrated in
Sections~{\ref{parabolic-pseudostablization}},
{\ref{relative-invariant-S-first-invariant-W}},
{\ref{branch-W-equiv-0-branch-W-nonzero}},
{\ref{recurrence-relations-parabolic-surfaces}},
such differential relations can also be `automatically'
implemented on Taylor series, and this saves computation time.

\smallskip

Lastly, there is another advantage of the power series method.

\smallskip\noindent$\bullet$\,
A progressive stratification of the normalisations of group parameters
for increasing fixed orders $\vert J \vert = 1, 2, 3, \dots$ conducts
to a high proximity with the Cartan method of equivalence, and its
famous {\sl reductions} and {\sl prolongations} of $G$-structures.

\smallskip

To explain these claims, let us exhibit another elementary example,
that of curves in the plane $\R_{x,u}^{1+1}$ under the special affine
group.

\Section{\bf Special Affine Power Series Invariants
of Curves in $\R^2$}
\label{special-affine-power-series-invariants-curves-R-2}
\HEAD{{\ref{special-affine-power-series-invariants-curves-R-2}}.~{\sf 
Special Affine Power Series Invariants of Curves in $\R^2$}
}{
Zhangchi {\sc Chen}, Joël~{\sc Merker}}

In $\R^2 \ni (x,u)$, we consider a graphed curve passing through the
origin normalized to order $1$:
\[
0
\,=\,
-\,u
+
F(x)
\,=\,
-\,u
+
F_2\,
{\textstyle{\frac{x^2}{2!}}}
+
F_3\,
{\textstyle{\frac{x^3}{3!}}}
+\cdots
\]
which satisfies $F_2 \neq 0$. These {\em formal coefficients} $F_i$
will be re-initialized in later stages of the
process. But at the
beginning, we assign these Taylor coefficients to be
initial {\sl functional jets:}
\[
F_i
\,:=\,
u_i
\,=\,
u_i(x),
\]
corresponding to $i$\textsuperscript{th} derivatives $u_{x^i}(x)$ at
arbitrary points $\big(x, u(x) \big)$ of the curve, as was explained
in Section~{\ref{power-series-method}}. 
All power series computations we
will perform at the origin will therefore have an interesting meaning
at {\em every point} of any curve in $\R^2$: such {\em is} 
the `{\sl power}' of power series!

\Subsection{First loop}
The special affine linear group $\SL_2(\R)$ consists of matrices
$\left(\! \begin{smallmatrix} {\sf a} & {\sf b} \\ {\sf c} & {\sf
d} \end{smallmatrix} \!\right)$ having determinant $1 = {\sf a} {\sf
d} - {\sf b} {\sf c}$.  Without primes or bars, target coordinates
will be denoted $(y,v)$. 

\begin{center}
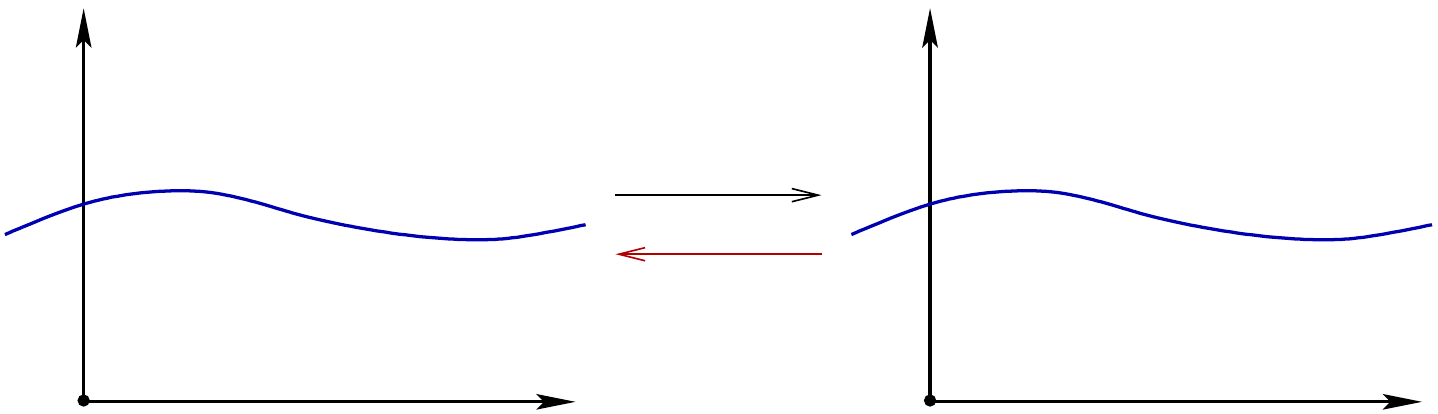
\end{center}

We consider special affine transformations:
\[
\R_{x,u}^2 
\,\,\,\longrightarrow\,\,\,
\R_{y,v}^2
\]
which are not too far from the identity,
so that any analytic graph $\big\{ u = F(x) \big\}$ is sent
to a similar graphed curve $\big\{ v = G(y) \big\}$.

Since the source power series
$F(x) = \sum_{i \geqslant 2}\, F_2\,
\frac{x^i}{i!}$ is given and since we want to simplify 
the target power series $G(y) = \sum_{j\geqslant 2}\, 
G_j\, \frac{y^j}{j!}$, it is more natural to 
work with the {\em inverse} special affine transformation, which
is also represented by means of an $\SL_2(\R)$ matrix as:
\[
x
\,=\,
{\sf a}\,y+{\sf b}\,v,
\ \ \ \ \ \ \ \ \ \ \ \ \ \ \ \ \ \ \ \
u
\,=\,
{\sf c}\,y+{\sf d}\,v,
\] 
Then the graphing function $G(y)$ is uniquely determined,
by a {\sl fundamental equation:}
\[
0
\,\equiv\,
-\,
\big({\sf c}\,y+{\sf d}\,v\big)
+
F\big(
{\sf a}\,y+{\sf b}\,v
\big)
\Big\vert_{\text{\rm replace}\,\,v=G(y)}
\eqno
{\scriptstyle{(\text{\rm in}\,\,\R\{y\})}},
\]
holding identically as a power series of the single horizontal
variable $y$\,\,---\,\,why and how $G$ is thusly determined will be
clear in a while. The game is to use the group parameters 
freedom
${\sf a}$, ${\sf b}$, ${\sf c}$, ${\sf d}$ 
in order to `kill' as much as possible 
coefficients $G_j$.

After an affine transformation, 
we may of course assume 
that our target graph enjoys a similar first-order
normalization $v = {\rm O}_y(2)$, namely:
\[
0
\,=\,
-\,v
+
G(y)
\,=\,
-\,v
+
G_2\,
{\textstyle{\frac{y^2}{2!}}}\,
+
G_2\,
{\textstyle{\frac{y^3}{3!}}}\,
+
\cdots.
\]
Then performing the plain replacement above:
\leqnomode\usetagform{default}
\begin{equation}
\label{eq-fond-F-G}
0
\,\equiv\,
-\,{\sf c}\,y
-
{\sf d}\,G(y)
+
F\big({\sf a}\,y+{\sf b}\,G(y)\big),
\end{equation}
we glean first-order terms which must vanish:
\[
0
\,\equiv\,
-\,{\sf c}\,
y
+
{\rm O}_y(2)
\eqno
{\scriptstyle{(\text{\rm in}\,\,\R\{y\})}}.
\]

\begin{Lemma}
\label{LM-c-zero}
The subgroup of $\SL_2(\R)$ sending
$v = {\rm O}_y(2)$ to $u = {\rm O}_x(2)$ 
is $2$-dimensional and consists of matrices:
\[
G_{\stabsmall}^{(1)}
\colon
\ \ \ \ \ \ \ \ \ \ \ \ \ \ \ \ \ \ \ \
\left(\!
\begin{array}{cc}
{\sf a} & {\sf b}
\\
0 & \frac{1}{{\sf a}}
\end{array}
\!\right)
\eqno
{\scriptstyle{({\sf a}\,\neq\,0)}}.
\qed
\]
\end{Lemma}

Thus, we have computed the subgroup which stabilizes
the current normal form of our power series. In later
stages of the process, higher order jets stabilizer subgroups:
\[
{\sf SL}_2(\R)
\,\supset\,
{\sf G}_{\stabsmall}^{(1)}
\,\supset\,
{\sf G}_{\stabsmall}^{(2)} 
\,\supset\, 
{\sf G}_{\stabsmall}^{(3)}
\,\supset\, 
\,\cdots\, 
\,\supset\,
{\sf G}_{\stabsmall}^{(\tau)}
\,=\,
\{e\}
\]
will naturally appear until final reduction to identity.
A deep proximity exists with Cartan's method
of equivalence.

But before jumping to the second loop of the
`algorithm', we must examine how this reduced subgroup
${\sf G}_{\stabsmall}^{(1)}$ acts on the second order term
$y^2$ in~({\ref{eq-fond-F-G}}), getting:
\[
0
\,\equiv\,
\big(
-\,
{\textstyle{\frac{1}{{\sf a}}}}\,
G_2
+
{\sf a}^2\,F_2
\big)\,
{\textstyle{\frac{y^2}{2!}}}
+
{\rm O}_y(3).
\]
From our assumption
$F_2 \neq 0$, we deduce that $G_2 \neq 0$
is inherited.

Furthermore,
taking ${\sf a} := \frac{1}{F_2^{1/3}}$, 
we can make $G_2 := 1$ by means of the
specific matrix\,\,---\,\,we choose
${\sf b} := 0$ for simplicity\,\,---:
\[
\left(\!
\begin{array}{cc}
\frac{1}{F_2^{1/3}} & 0
\\
0 & F_2^{1/3}
\end{array}
\!\right)
\,\in\,
{\sf G}_{\stabsmall}^{(1)}.
\]

Then with this precise special affine transformation, a
computation of the higher order terms in~({\ref{eq-fond-F-G}}) gives 
us:
\[
0
\,\equiv\,
\big(
-\,F_2^{1/3}G_3
+
{\textstyle{\frac{F_3}{F_2}}}
\big)\,
{\textstyle{\frac{y^3}{3!}}}
+
\big(
-\,F_2^{1/3}G_4
+
{\textstyle{\frac{F_4}{F_2^{4/3}}}}
\big)\,
{\textstyle{\frac{y^4}{4!}}}
+
\big(
-\,F_2^{1/3}G_5
+
{\textstyle{\frac{F_5}{F_2^{5/3}}}}
\big)\,
{\textstyle{\frac{y^5}{5!}}}
+\cdots,
\]
and we obtain as promised unique determinations:
\[
G_i
\,:=\,
\frac{F_i}{F_2^{\frac{1+i}{3}}}
\eqno
{\scriptstyle{(i\,\geqslant\,3)}}.
\]

Finally, remembering that the $F_i$ were formal variables
representing the functional jets $u_{x^i}(x) = u_i(x) = u_i = 
F_i$, we see that the coefficients 
of the transformed curve
$\big\{ v = G(y) \big\}$ have become:
\[
G_2
\,:=\,
1,
\ \ \ \ \ \ \ \ \ \ \ \ \ \ \ \ \ \ \ \
G_i
\,:=\,
\frac{u_i}{u_2^{\frac{1+i}{3}}}
\eqno
{\scriptstyle{(i\,\geqslant\,3)}}.
\] 

\Subsection{Second loop} In order to avoid indices heaviness,
we keep the same notation $\big\{ u = F(x) \big\}$ and
$\big\{ v = G(y) \big\}$, which means that what we now call
$F$ is the $G$ of the end of the previous loop.
So in terms of the {\em initial} 
functional jets $u_i$, we have in fact
{\em re-assigned:}
\[
F_2
\,:=\,
1,
\ \ \ \ \ \ \ \ \ \ \ \ \ \ \ \ \ \ \ \
F_i
\,:=\,
\frac{u_i}{u_2^{\frac{1+i}{3}}}
\eqno
{\scriptstyle{(i\,\geqslant\,3)}}.
\]
Keeping this in memory, we will now work 
{\em formally} with power series coefficients 
$F_i$ and $G_j$, and only at the end of the current
loop will we express the result in terms of these $F_i = 
u_i \big/ u_2^{(1+i)/3}$.

So both our source and target graphed curves
may be assumed to have terms
normalized up to order $2$ included:
\[
u
\,=\,
{\textstyle{\frac{x^2}{2!}}}
+
F_3\,
{\textstyle{\frac{x^3}{3!}}}
+
\cdots
\ \ \ \ \ \ \ \ \ \ \ \ \ \ \ \ \ \ \ \
\text{and}
\ \ \ \ \ \ \ \ \ \ \ \ \ \ \ \ \ \ \ \
v
\,=\,
{\textstyle{\frac{y^2}{2!}}}
+
G_3\,
{\textstyle{\frac{y^3}{3!}}}
+
\cdots.
\]

Remembering that Lemma~{\ref{LM-c-zero}}
already showed that stabilization up to 
order $1$ forces ${\sf c} = 0$, in order to determine
the subgroup ${\sf G}_{\sf stab}^{(2)} \subset \SL_2(\R)$,
we can work within ${\sf G}_{\sf stab}^{(1)}$.

\begin{Lemma}
The subgroup of $\SL_2(\R)$ 
sending $u = \frac{1}{2}\, x^2 + {\rm O}_x(3)$ to
$v = \frac{1}{2}\, y^2 + {\rm O}_y(3)$ is $1$-dimensional and
consists of matrices:
\[
G_{\stabsmall}^{(2)}
\colon
\ \ \ \ \ \ \ \ \ \ \ \ \ \ \ \ \ \ \ \
\left(\!
\begin{array}{cc}
1 & {\sf b}
\\
0 & 1
\end{array}
\!\right).
\]
\end{Lemma}

\proof
Hence with ${\sf c} = 0$, back 
to~({\ref{eq-fond-F-G}}), 
we get ${\sf a}^3 = 1$ from:
\[
0
\,\equiv\,
\big(
-\,
{\textstyle{\frac{1}{{\sf a}}}}
+
{\sf a}^2
\big)\,
{\textstyle{\frac{y^2}{2!}}}
+
{\sf O}_y(3).
\qedhere
\]
\endproof

Thus, as promised, we have determined 
the subgroup $G_{\stabsmall}^{(2)}$ which stabilizes
the current normal form of our power series, and as we will
see soon, the process will stop at stage $3$:
\[
{\sf SL}_2(\R)
\,\supset\,
{\sf G}_{\stabsmall}^{(1)}
\,\supset\,
{\sf G}_{\stabsmall}^{(2)} 
\,\supset\, 
{\sf G}_{\stabsmall}^{(3)}
\,=\,
\{e\}.
\]

For the time being, let us examine how this reduced group
${\sf G}_{\stabsmall}^{(2)}$ acts on the third order term
$y^3$ in~({\ref{eq-fond-F-G}}):
\[
0
\,\equiv\,
\big(
-\,G_3
+
F_3
+
3\,{\sf b}
\big)\,
{\textstyle{\frac{y^3}{3!}}}
+
{\rm O}_y(4).
\]
We can make $G_3 = 0$ with ${\sf b} := - \frac{1}{3}\, F_3$,
by means of the specific matrix:
\[
\left(\!
\begin{array}{cc}
1 & -F_3/3
\\
0 & 1
\end{array}
\!\right).
\]

Then with this precise special affine transformation,
higher order vanishing terms:
\[
0
\,\equiv\,
\big(
-G_4
-
{\textstyle{\frac{5}{3}}}\,
F_3^2
+
F_4
\big)\,
{\textstyle{\frac{y^4}{4!}}}
+
\big(
-\,G_5
-
{\textstyle{\frac{5}{3}}}\,
F_3G_4
+
F_5
+
{\textstyle{\frac{5}{3}}}\,
F_3^2
-
{\textstyle{\frac{10}{3}}}\,
F_3F_4
\big)\,
{\textstyle{\frac{y^5}{5!}}}
+
{\rm O}_y(6)
\]
conduct us to:
\[
G_4
\,:=\,
-\,
{\textstyle{\frac{5}{3}}}\,
F_3^2
+
F_4,
\ \ \ \ \ \ \ \ \ \ \ \ \ \ \ \ \ \ \ \
G_5
\,:=\,
F_5
+
{\textstyle{\frac{40}{9}}}\,
F_3^3
-
5\,F_3F_4.
\]
Coming back to the functional jets, we obtain:
\[
G_2
\,:=\,
1,
\ \ \ \ \ \ \ \ \
G_3
\,:=\,
0,
\ \ \ \ \ \ \ \ \
G_4
\,:=\,
\frac{1}{3}\,
\frac{-5\,u_3^2+3\,u_2u_4}{u_2^{8/3}},
\ \ \ \ \ \ \ \ \
G_5
\,:=\,
\frac{1}{9}\,
\frac{9\,u_2^2u_5
-
45\,u_2u_3u_4
+
40\,u_3^3}{u_2^4}.
\]

\Subsection{Third loop}
We start by re-assigning:
\[
F_2
\,:=\,
1,
\ \ \ \ \ \ \ \ \
F_3
\,:=\,
0,
\ \ \ \ \ \ \ \ \
F_4
\,:=\,
\frac{1}{3}\,
\frac{-5\,u_3^2+3\,u_2u_4}{u_2^{8/3}},
\ \ \ \ \ \ \ \ \
F_5
\,:=\,
\frac{1}{9}\,
\frac{9\,u_2^2u_5
-
45\,u_2u_3u_4
+
40\,u_3^3}{u_2^4}.
\]
We again work with formal $u = \sum\, F_i \, \frac{x^i}{i}$
and $v = \sum\, G_j\, \frac{y^j}{j!}$ assuming the normalizations:
\[
F_2
\,=\,
1,
\ \ \ \ \ \ \ \ \ \ \
F_3
\,=\,
0
\ \ \ \ \ \ \ \ \ \ \ \ \ \ \ \ \ \ \ \
\text{and}
\ \ \ \ \ \ \ \ \ \ \ \ \ \ \ \ \ \ \ \
G_2
\,=\,
1,
\ \ \ \ \ \ \ \ \ \ \
G_3
\,=\,
0.
\]
Naturally, since only $1$ degree of freedom was left at the
previous stage, the condition that the coefficient $0$
of $\frac{x^3}{3!}$ is left unchanged drops the group dimension
by $1$, and it is easy to verify the 

\begin{Lemma}
The subgroup of $\SL_2(\R)$ 
sending 
$v = \frac{y^2}{2} + 0 + {\rm O}_y(4)$ 
to
$u = \frac{x^2}{2} + 0 + {\rm O}_x(4)$ 
is $0$-dimensional and
reduces to the identity:
\[
G_{\stabsmall}^{(3)}
\colon
\ \ \ \ \ \ \ \ \ \ \ \ \ \ \ \ \ \ \ \
\left(\!
\begin{array}{cc}
1 & 0
\\
0 & 1
\end{array}
\!\right).\eqno\qed
\]
\end{Lemma}

The algorithm therefore stops, and a first result, valid only at the
level of power series at the origin, is a corollary of this reduction
to an $\{e\}$-group.

\begin{Theorem}
\label{Thm-all-Taylor-coefficients-invariants}
{\bf (1)}\,
Given a real analytic curve $\big\{ u = F(x) \big\}$ in $\R^2$ 
passing through the origin which satisfies: 
\[
F_{xx}(0)
\,\neq\,
0,
\]
there always exists an ${\sf SL}_2(\R)$ transformation which
puts it into the form
\[
u
\,=\,
{\textstyle{\frac{x^2}{2!}}}
+
0
+
F_4\,
{\textstyle{\frac{x^4}{4!}}}
+
F_5\,
{\textstyle{\frac{x^5}{5!}}}
+
\sum_{i\geqslant 6}\,
F_i\,
{\textstyle{\frac{x^i}{i!}}}.
\]

\noindent{\bf (2)}\,
Any other such real analytic curve $\big\{ v = 
G(y) \big\}$ similarly
put into the form:
\[
v
\,=\,
{\textstyle{\frac{y^2}{2!}}}
+
0
+
G_4\,
{\textstyle{\frac{y^4}{4!}}}
+
G_5\,
{\textstyle{\frac{y^5}{5!}}}
+
\sum_{j\geqslant 6}\,
G_j\,
{\textstyle{\frac{y^j}{j!}}},
\]
is ${\sf SL}_2(\R)$-equivalent to $\big\{ u = F(x) \big\}$ above if
and only if all Taylor coefficients match:
\[
G_4
\,=\,
F_4,
\ \ \ \ \ \ \ \ \ \ \ \ \
G_5
\,=\,
F_5,
\ \ \ \ \ \ \ \ \ \ \ \ \
G_i
\,=\,
F_i
\eqno
{\scriptstyle{(\forall\,i\,\geqslant\,6)}}.
\qed
\]
\end{Theorem}

As explained in Section~{\ref{power-series-method}}, 
the power series coefficients so
obtained $F_4$, $F_5$, are differential invariants
at any point $\big(x, u(x) \big)$
of the curve, and beyond,
the explicit expressions of the next two are:
\[
\aligned
F_6
&
\,:=\,
\frac{1}{9}\,
\frac{9\,u_2^3u_6
-
63\,u_2^2u_3u_5
+
105\,u_2u_3^2u_4
-
35\,u_3^4}{
u_2^{16/3}},
\\
F_7
&
\,:=\,
\frac{1}{9}\,
\frac{9\,u_2^4u_7
-
84\,u_2^3u_3u_6
+
210\,u_2^2u_3^2u_5
-
105\,u_2^2u_3u_4^2
+
210\,u_2u_3^3u_4
-
280\,u_3^5}{
u_2^{20/3}}.
\endaligned
\]

Furthermore, in terms
of the total differentiation operator:
\[
\Daux_x
\,:=\,
\frac{\partial}{\partial x}
+
u_1\,\frac{\partial}{\partial u}
+
\sum_{i=1}^\infty\,
u_{i+1}\,
\frac{\partial}{\partial u_i}
\]
the affine-invariant differentiation operator:
\leqnomode\usetagform{default}
\begin{align}
\label{1-dim-D-x-invariant}
\mathcal{D}_x
\,:=\,
{\textstyle{\frac{1}{u_2^{1/3}}}}\,
\Daux_x
\end{align}
enables to produce higher order invariants, for instance:
\[
\mathcal{D}_x
\Big(
\frac{1}{3}\,
\frac{-5\,u_3^2+3\,u_2u_4}{u_2^{8/3}}
\Big)
\,=\,
\frac{1}{9}\,
\frac{9\,u_2^2u_5-45\,u_2u_3u_4+40\,u_3^3}{u_2^4}.
\]
More will be said in 
Section~{\ref{moving-frame-recurrence-curves-R-2}}.

Next, we present an ancient result ({\cite{Halphen-1878}}).

\begin{Lemma}
\label{LM-vanishing-H-M}
{\bf (1)}\,
A curve $u = u(x)$ with $u_{xx} \neq 0$ is affinely equivalent
to a parabola $v = y^2$ if and only if:
\[
0
\,\equiv\,
\Paux(u)
\,:=\,
\frac{1}{3}\,
\frac{3\,u_{xx}\,u_{xxxx}-5\,u_{xxx}^2}{u_{xx}^{8/3}}.
\]

{\bf (2)}\,
A curve $u = u(x)$ with $u_{xx} \neq 0$ is affinely equivalent
to a nondegenerate conic in the plane if and only if:
\[
0
\,\equiv\,
\Caux(u)
\,:=\,
\frac{1}{9}\,
\frac{9\,u_{xx}^2\,u_{xxxxx}-45\,u_{xx}\,u_{xxx}\,u_{xxxx}
+40\,u_{xxx}^3}{u_{xx}^4}.
\]

{\bf (3)}\,
$\Paux(u) \equiv 0$ implies $\Caux(u) \equiv 0$. 

\end{Lemma}

\proof
{\small\bf (1)}\,
The target being $v = y^2$, we have after an affine transformation
whose linear part 
$\left(\! \begin{smallmatrix} {\sf a} & {\sf b} \\ {\sf k} & {\sf
l} \end{smallmatrix} \right) \sim \left(\! \begin{smallmatrix} 1 &
0 \\ 0 & 1 \end{smallmatrix} \right)$ may be assumed close to
the identity:
\[
{\sf k}\,x
+
{\sf l}\,u
+
{\sf m}
\,=\,
\big(
{\sf a}\,x
+
{\sf b}\,u
+
{\sf c}
\big)^2,
\]
hence we may solve for $u$ thanks to the positivity 
of the discriminant:
\[
u
\,=\,
\frac{
{\sf l}
-
2\,{\sf b}{\sf c}
-
2\,{\sf a}{\sf b}\,x
\pm
\sqrt{
\big(
4\,{\sf b}^2{\sf k}
-
4\,{\sf a}{\sf b}{\sf l}
\big)\,
x
+
{\sf l}^2
-
4\,{\sf b}{\sf c}{\sf l}
+
4\,{\sf b}^2{\sf m}
}
}{2\,{\sf b}^2}.
\]

Thus with different constants the general equation of
parabolas is:
\[
u
\,=\,
{\sf d}\,x
+
{\sf e}
+
\sqrt{2\,{\sf g}\,x+{\sf h}}.
\]
At first, in order 
to eliminate ${\sf d}$ and ${\sf e}$, we just differentiate two times:
\[
u_{xx}
\,=\,
-\,\frac{{\sf g}^2}{
\big(2\,{\sf g}\,x+{\sf h}\big)^{3/2}},
\]
and next, to eliminate the remaining constants, we upside-down:
\[
\frac{1}{u_{xx}^{2/3}}
\,=\,
\frac{2}{{\sf g}^{1/3}}\,
x
+
\frac{{\sf h}}{{\sf g}^{4/3}},
\]
and we again differentiate twice:
\[
\Big(
\frac{1}{u_{xx}^{2/3}}
\Big)_{xx}
\,=\,
-\frac{2}{9}\,
\frac{3\,u_{xx}\,u_{xxxx}-5\,u_{xxx}^2}{u_{xx}^{8/3}}
\,\equiv\,
0.
\]
Conversely, from $\Paux(u) \equiv 0$, using $u_{xx} \neq 0$, 
one reconstitutes by integration (exercise)
the general equation $u = {\sf d}\, x + {\sf e} + 
\sqrt{ 2\, {\sf g}\, x + {\sf h}}$.

\smallskip{\small\bf (2)}\,
Now, the target is a general conic in the $\R_{y,v}^2$-plane, 
hence an $x^2$ monomial must be present under the square root:

\[
u
\,=\,
{\sf d}\,x
+
{\sf e}
+
\sqrt{{\sf f}\,x^2+2\,{\sf g}\,x+{\sf h}}.
\]
Quite similarly:
\[
u_{xx}
\,=\,
\,\frac{{\sf f}\,{\sf h}-{\sf g}^2}{
\big({\sf f}\,x^2+2\,{\sf g}\,x+{\sf h}\big)^{3/2}},
\]
whence:
\[
\frac{1}{u_{xx}^{2/3}}
\,=\,
\frac{{\sf f}}{({\sf f}\,{\sf h}-{\sf g}^2)^{2/3}}\,
x^2
+
\frac{2\,{\sf g}}{({\sf f}\,{\sf h}-{\sf g}^2)^{2/3}}\,
x
+
\frac{{\sf h}}{({\sf f}\,{\sf h}-{\sf g}^2)^{2/3}},
\]
and lastly we have to differentiate {\em three times} to get rid
of all remaining constants:
\[
\Big(
\frac{1}{u_{xx}^{2/3}}
\Big)_{xxx}
\,=\,
-\,
\frac{2}{27}\,
\frac{9\,u_{xx}^2\,u_{xxxxx}-45\,u_{xx}\,u_{xxx}\,u_{xxxx}
+40\,u_{xxx}^3}{u_{xx}^{11/3}}
\,\equiv\,
0.
\]
The converse is also left as an exercise.

\smallskip{\small\bf (3)}\,
follows from a direct differentiation of
$0 \equiv 3\, u_{xx}(x)\, u_{xxxx}(x) - 5\, u_{xxx}(x)^2$, 
and will now be explained in a more theoretical framework.
\endproof

\Section{\bf Recurrence formulas for Differential Invariants}
\label{recurrence-formulae-differential-invariants}
\HEAD{{\ref{recurrence-formulae-differential-invariants}}.~{\sf 
Recurrence formulas for Differential Invariants}
}{
Zhangchi {\sc Chen}, Joël~{\sc Merker}}

\Subsection{Prolongations of infinitesimal generators}
As before, let $G$ be a local Lie group of finite dimension
$1 \leqslant r < \infty$ acting on graphs in $\R_{x,u}^{p+q}$.
Abbreviate $z = (x,u)$ and $z^{(n)} = \big(x, u^{(n)} \big)$.
Choose any basis $e_1, \dots, e_r$ for the Lie algebra
$\mathfrak{g} := \Lie (G)$ and introduce the 
$r$ infinitesimal generators of the $G$-action 
on the considered space:
\[
{\bf v}_\sigma
\,:=\,
\frac{d}{dt}\,
\Big\vert_{t=0}
\exp\,
\big(
t\,e_\sigma
\big)
\cdot
z
\eqno
{\scriptstyle{(1\,\leqslant\,\sigma\,\leqslant\,r)}}.
\]
We thus get $r$ vector fields forming a Lie algebra, which we
will write as:
\[
{\bf v}_\sigma
\,=\,
\sum_{i=1}^p\,
\xi_\sigma^i(x,u)\,
\frac{\partial}{\partial x^i}
+
\sum_{\alpha=1}^q\,
\varphi_\sigma^\alpha\,
\frac{\partial}{\partial u^\alpha}
\eqno
{\scriptstyle{(1\,\leqslant\,\sigma\,\leqslant\,r)}}.
\]

For any jet order $\oorder \geqslant 0$, and even up to infinity,
the prolongations of these fields:
\[
{\sf v}_\sigma^{(\infty)}
\,=\,
{\bf v}_\sigma
+
\sum_{\#\,J\geqslant 1}\,
\sum_{\alpha=1}^q\,
{\varphi_{\sigma}}_J^\alpha\,
\big(
x,u^{(J)}
\big)\,
\frac{\partial}{\partial u_J^\alpha},
\]
have coefficients computed by applying 
Theorem~{\ref{Thm-prolongation-vector-fields}}:
\[
{\varphi_{\sigma}}_J^\alpha
\,=\,
{\sf D}_{x^J}
\Big(
\varphi_\sigma^\alpha
-
\sum_{1\leqslant i\leqslant p}\,
\xi_\sigma^i\,
u_i^\alpha
\Big)
+
\sum_{1\leqslant i\leqslant p}\,
\xi_\sigma^i\,
u_{J,i}^\alpha,
\]
where, by a slight abuse of notation,
we have denoted $u^{(J)}$ instead of $u^{(\# J)}$.

\smallskip

Next, we recall that the modified total differentiation
operators are:
\[
\left(\!
\begin{array}{c}
{\sf E}_{x^1}
\\
\vdots
\\
{\sf E}_{x^p}
\end{array}
\!\right)
\,\,:=\,\,
\left(\!
\begin{array}{ccc}
{\sf D}_{x^1}(y^1) & \cdots & {\sf D}_{x^1}(y^p)
\\
\vdots & \ddots & \vdots
\\
{\sf D}_{x^p}(y^1) & \cdots & {\sf D}_{x^p}(y^p)
\end{array}
\!\right)^{\!-1}\,\,
\left(\!
\begin{array}{c}
{\sf D}_{x^1}
\\
\vdots
\\
{\sf D}_{x^p}
\end{array}
\!\right),
\]
and after perfoming abstractly the matrix inversion, 
this means that there exist coefficient-functions such that:
\leqnomode\usetagform{default}
\begin{align}
\label{equation-definition-invariant-D-j}
{\sf E}_{x^i}
\,=\,
\sum_{1\leqslant j\leqslant p}\,
Z_i^j
\big(
g,\,x,u^{(1)}
\big)\,
{\sf D}_{x^j}.
\end{align}
Following~{\cite{Fels-Olver-1999}}, we may {\sl invariantize} these
${\sf E}_{x^i}$ by replacing the group parameters $g = (g_1, \dots,
g_r)$ by their values $g := \rho \big( z^{(n_\GG)} \big)$ solved from
the cross-section equations, and in this way, we produce $p$ {\sl
invariant horizontal differential operators:}
\[
\mathcal{D}_i
\,:=\,
\inv
\big(
{\sf E}_{x^i}
\big)
\,:=\,
\sum_{1\leqslant j\leqslant p}\,
Z_i^j
\big(
\rho(z^{(n_\GG)}),\,
x,u^{(1)}
\big)
\big)\,
{\sf D}_{x^j}.
\]

\begin{Observation}
Such invariant differential operators $\mathcal{D}_1, \dots, 
\mathcal{D}_p$ have coefficients which depend on jets
of order $n_\GG$ ,where $n_\GG$ is the minimal jet order
for which the action of $G$ on $J_{x,u}^n$ becomes
locally foliated (free) of rank equal to $r = \dim\, G$.\qed
\end{Observation}

Knowing that a moving frame map $\rho \colon J_{p,q}^{n_\GG}
\longrightarrow G$ is often difficult to construct explicitly,
these $\mathcal{D}_i$ are {\em almost never} explicit!

\begin{Definition}
The {\sl invariantization} of an arbitrary function $\mathcal{F} =
\mathcal{F} \big( z^{(n)} \big)$ defined on 
any $n$\textsuperscript{th} order jet space
$J_{x,u}^n$ with $n \geqslant 0$
is defined by
replacing firstly $z^{(n)}$ in its argument by
the target value $w^{(n)} = w^{(n)}
\big(g, z^{(n)} \big)$ and secondly the group parameters $g = \big(
g_1, \dots, g_r \big)$ by their expressions solved from 
the (not canonical) cross-section equations:
\[
\Faux
\,:=\,
\inv
\big(
\mathcal{F}\big(z^{(n)}\big)
\big)
\,:=\,
\mathcal{F}
\Big(
w^{(n)}
\big(
\rho(z^{(n_\GG)}),\,\,
z^{(n)}
\big)
\Big).
\]
\end{Definition}

Of course, $\inv ( \mathcal{F} )$ is always
a differential invariant.
In particular, one can invariantize all jet monomials:
\[
\aligned
\Jaux^i
&
\,:=\,
\inv
\big(x^i\big)
\,=\,
y^i\,
\big(
\rho(z^{(n_\GG)}),
x,u
\big),
\\
\Iaux_K^\alpha
&
\,:=\,
\inv\,
\big(
u_K^\alpha
\big)
\,=\,
v_K^\alpha
\Big(
\rho(z^{(n_\GG)}),x,u^{(K)}
\Big),
\endaligned
\]
where, as before, $1 \leqslant i \leqslant p$, $1 \leqslant \alpha
\leqslant q$, $K = k_1, \dots, k_\lambda$, $1 \leqslant k_1, \dots,
k_\lambda \leqslant p$.
We will call them {\sl monomial} differential invariants.

\begin{Convention}
When $K = \emptyset$, we agree that $\Iaux_\emptyset^\alpha = 
\Iaux^\alpha$ is the invariantization of 
the un-differentiated $u^\alpha$.
\end{Convention}

Most of the times\,\,---\,\,at least for {\em all} group actions
considered in this article\,\,---, the action of $G$ is 
transitive on $\R_{x,u}^{p+q}$, hence the cross-section
equations can be assumed to contain the $p+q$ equations:
\[
y^i
\,=\,
0,
\ \ \ \ \ \ \ \ \ \ \ \ \ \ \ \ \ \ \ \
u^\alpha
\,=\,
0,
\]
and then:
\[
\Jaux^i
\,=\,
0,
\ \ \ \ \ \ \ \ \ \ \ \ \ \ \ \ \ \ \ \
\Iaux^\alpha
\,=\,
0.
\]

We will abbreviate:
\[
\Iaux^{(n)}
\,:=\,
\Big(
\Jaux^i,\,
\big\{
\Iaux_K^\alpha
\big\}_{0\leqslant \# K\leqslant n}^{1\leqslant\alpha\leqslant q}
\Big),
\]
and we can at last 
present the fundamental {\sl recurrence formulas}.

\begin{Theorem}
\label{Thm-recurrence-formulae}
{\rm {\cite{Fels-Olver-1999}}}
{\bf [Recurrence relations]}
The invariant derivatives of all the monomial invariants
$\Iaux_K^\alpha$ with respect to $\mathcal{D}_1, \dots, \mathcal{D}_p$
express as differential invariants of order $\# K + 1$ raised by one
unit:
\[
\mathcal{D}_j
\Iaux_K^\alpha
\,=\,
\Iaux_{K,j}^\alpha
+
\sum_{1\leqslant\sigma\leqslant r}\,
{\varphi_\sigma}_K^\alpha
\big(
\Iaux^{(K)}
\big)
\cdot
\Kaux_j^\sigma
\big(
\Iaux^{(n_\GG+1)}
\big),
\]
with the peculiar case of:
\[
\mathcal{D}_j
\Jaux^i
\,=\,
\delta_j^i
+
\sum_{1\leqslant\sigma\leqslant r}\,
\xi_\sigma^i
\big(
\Iaux^{(0)}
\big)\cdot
\Kaux_j^\sigma
\big(
\Iaux^{(n_\GG+1)}
\big),
\]
plus some correction remainders which incorporate the invariantizations
${\varphi_\sigma}_K^\alpha \big( \Iaux^{(K)} \big)$ of the
coefficients ${\varphi_\sigma}_K^\alpha \big( z^{(K)} \big)$ of the
prolonged vector fields ${\bf v}_1^{(\infty)}, \dots,
{\bf v}_r^{(\infty)}$, as well as
certain special differential invariants $\Kaux_j^\sigma$
of order $\leqslant n_\GG + 1$\,\,---\,\,called
of `{\sl Maurer-Cartan}' type.

Moreover, for each fixed $1 \leqslant j \leqslant p$, the linear
system of $r$ equations in the $r$ unknows $\Kaux_j^1, \dots,
\Kaux_j^r$ constituted of the $r$ recurrence relations
written only for the $r$ phantom invariants:
\[
\inv\big(w_{\nu_1}^{(n_1)}\big)
\,=\,
c_1,\,\,\,
\dots\dots,\,\,\,
\inv\big( w_{\nu_r}^{(n_r)}\big)
\,=\,
c_r,
\] 
has left-hand sides all equal to $0 = \mathcal{D}_j(c_1) = \cdots =
\mathcal{D}_j(c_r)$, is of Cramér type, and can be solved uniquely for
$\Kaux_j^1, \dots, \Kaux_j^r$ in terms of 
the other present differential invariants, all of
order at most:
\[
1
+
\underset{1\leqslant h\leqslant r}{\max}\,
n_h
\,=\,
1
+
n_\GG.
\eqno\qed
\]
\end{Theorem}

The last sentence justifies why these Maurer-Cartan invariants can be
written $\Kaux_j^\sigma = \Kaux_j^\sigma \big( \Iaux^{(n_\GG+1)}
\big)$.  They appear in a wider theoretical context
({\cite{Fels-Olver-1999}}), but in all applications, the most
efficient way to determine these $\Kaux_j^\sigma$ is to set up and
solve the mentioned $p$ Cramér systems.

\smallskip

This fundamental theorem has (at least) three outstanding deep 
qualities.

\smallskip\noindent$\bullet$\,
To provide {\em all} possible
{\em differential} relations between differential invariants
({\cite{Fels-Olver-1999}}).

\smallskip\noindent$\bullet$\,
To be applicable without the need of computing explicitly 
{\em any}\,\,---\,\,even a single one which is not 
phantom!\,\,---\,\,differential invariant, a task which reveals
itself desperately beyond the reach
of any modern symbolic software.

\smallskip\noindent$\bullet$\,
To jump at a synthetic level at which calculations amount
to plain linear algebra in low dimension, hence are very accessible.

\Subsection{Groups containing translations}
Now, assume that our local Lie group $G$ contains
{\em all translations:}
\[
\big(x^i,\,u^\alpha\big)
\,\,\,\longmapsto\,\,\,
\big(
x^i+{\sf a}^i,\,\,
u^\alpha+{\sf b}^\alpha
\big).
\]
This assumption is justified throughout the article, because we will
mainly deal with the standard affine or special affine groups of
transformations acting on $\R^{p+q}$.

We take a reference point $p_0 \in \R^{p+q}$, at which we center the
coordinates $(x,u)$, so that $p_0 = (0, 0)$ is the origin.  With $z =
(x,u)$ and $w = (y,v)$, we recall that we denote the $G$-action from
$\R^m$ to $\R^m$ as $z \longmapsto w(g,z)$, where $m = p + q$.

So the initial graphs $\big\{ u^\alpha = F^\alpha (x) \big\}$ we are
considering satisfy $F^\alpha(0) = 0$, while the target
$G$-transformed graphs $\big\{ v^\beta = G^\beta(g, y) \big\}$ do not
necessarily pass through the origin\,\,---\,\,but our goal is to
reduce ourselves to having $G^\beta(g, 0) = 0$ too, at the price of
{\em reducing} $G$. Like for the $G$-structures of the Cartan method
of equivalence, this {\em principle of progressive reduction} will be
practically very powerful.

To this aim, we may introduce the {\sl isotropy subgroup} of the
origin:
\[
G_0
\,:=\,
\big\{
g\in G
\colon\,
w(g,0)
=
0
\big\}.
\] 
It {\em is} a local Lie group too, just because
fixing a point is preserved by composition.
Since our group $G$
which contains all translations is trivially transitive,
we have:
\[
\dim\,
G_0
\,=\,
r
-
(p+q)
\,=:\,
r_0.
\] 

After a reordering, we may assume that a basis for $\mathfrak{g}
= \Lie(G)$ as introduced before is organized so that:
\[
\aligned
\mathfrak{g}
&
\,=\,
\Span\,
\big\{
e_1,\dots,e_m,\,
e_{m+1},\dots,e_r
\big\},
\\
\Lie(G_0)
\,=:\,
\mathfrak{g}_0
&
\,=\,
\Span\,
\big\{
\ \ \ \ \ \ \ \ \ \ \ \ \ \ \ \ \ \ \ \
e_{m+1},\dots,e_r
\big\}.
\endaligned
\]
If as before, 
we also introduce the infinitesimal generators of the action of $G$ on
$\R^m$:
\[
{\bf v}_\sigma
\,:=\,
\frac{d}{dt}\,
\Big\vert_{t=0}
\exp\,
\big(
t\,e_\sigma
\big)
\cdot
z
\eqno
{\scriptstyle{(1\,\leqslant\,\sigma\,\leqslant\,r)}},
\]
then as is well known since Lie
({\cite{Engel-Lie-1888, Lie-Merker-2015}}), 
all the ones associated to $G_0$ do vanish at the origin:
\[
0
\,=\,
{\bf v}_{m+1}
\big\vert_0
\,=\,\cdots\,=\,
{\bf v}_r
\big\vert_0.
\]

Since $G$ contains all translations, it is clear that we may 
assume that the first $p+q$ infinitesimal generators are plainly:
\[
{\bf v}_i
\,=\,
\partial_{x^i},
\ \ \ \ \ \ \ \ \ \ \ \ \ \ \ \ \ \ \ \
{\bf v}_{p+\alpha}
\,=\,
\partial_{u^\alpha}
\eqno
{\scriptstyle{(1\,\leqslant\,i\,\leqslant\,p,\,\,
1\,\leqslant\,\alpha\,\leqslant\,q)}}.
\]
Furthermore, on the jet space $J_{x,u}^{n_\GG}$
of the appropriate order $n_\GG$,
we may assume that the cross-section equations to the
lifted action of $G$ on $J_{x,u}^{n_\GG}$
contain the $(p+q)$ equations:
\[
y^i\big(g,x,u\big)
\,=\,
0,
\ \ \ \ \ \ \ \ \ \ \ \ \ \ \ \ \ \ \ \
v^\alpha\big(g,x,u\big)
\,=\,
0,
\]
and that these equations are the first $m$ among the $r$
cross-section equations we wrote abstractly as:
\[
0
\,=\,
w_{\nu_1=1}^{(n_1)}
\big(z,g\big)
\,=\,\cdots\,=\,
w_{\nu_m=m}^{(n_m)}
\big(z,g\big).
\]
Therefore, among the phantom differential invariants, there are
the $p+q$:
\[
\Jaux^i
\,=\,
\inv
\big(
y^i
\big),
\ \ \ \ \ \ \ \ \ \ \ \ \ \ \ \ \ \ \ \
\Iaux_\emptyset^\alpha
\,=\,
\inv
\big(
v^\alpha
\big).
\]

A quick inspection of the prolongation 
formulas of Theorem~{\ref{Thm-prolongation-vector-fields}}
convinces of an

\begin{Observation} 
The prolongations
of the translational vector fields are all trivial:
\[
{\bf v}_i^{(\infty)}
\,=\,
{\bf v}_i
+
0,
\ \ \ \ \ \ \ \ \ \ \ \ \ \ \ \ \ \ \ \
{\bf v}_{p+\alpha}^{(\infty)}
\,=\,
{\bf v}_{p+\alpha}
+
0.
\eqno\qed
\]
\end{Observation}

Then in the fundamental recurrence formulas of  
Theorem~{\ref{Thm-recurrence-formulae}},
we decide to set aside the $p+q$ ones concerning these
first phantom invariants, namely:
\[
\aligned
\mathcal{D}_j
\Jaux^i
&
\,=\,
\delta_j^i
+
1
\cdot
\Kaux_j^i
+
\sum_{m+1\leqslant\sigma\leqslant r}\,
\xi_\sigma^i
\big(
\Iaux^{(0)}
\big)
\cdot
\Kaux_j^\sigma
\big(
\Iaux^{(n_\GG+1)}
\big),
\\
\mathcal{D}_j
\Iaux_{\emptyset}^\alpha
&
\,=\,
\Iaux_j^\alpha
+
1
\cdot
\Kaux_j^{p+\alpha}
+
\sum_{m+1\leqslant\sigma\leqslant r}\,
\varphi_\sigma^\alpha
\big(
\Iaux^{(0)}
\big)
\cdot
\Kaux_j^\sigma
\big(
\Iaux^{(n_\GG+1)}
\big),
\endaligned
\]
and to conserve only the recurrence relations 
that are of order $\# K \geqslant 1$:
\[
\aligned
\mathcal{D}_j
\Iaux_K^\alpha
&
\,=\,
\Iaux_{K,j}^\alpha
+
\sum_{1\leqslant\sigma\leqslant r}\,
{\varphi_\sigma}_K^\alpha
\big(
\Iaux^{(K)}
\big)
\cdot
\Kaux_j^\sigma
\big(
\Iaux^{(n_\GG+1)}
\big)
\\
&
\,=\,
\Iaux_{K,j}^\alpha
+
0
+
\sum_{m+1\leqslant\sigma\leqslant r}\,
{\varphi_\sigma}_K^\alpha
\big(
\Iaux^{(K)}
\big)
\cdot
\Kaux_j^\sigma
\big(
\Iaux^{(n_\GG+1)}
\big)
\endaligned
\]
but then our observation that ${\varphi_\sigma}_K^\alpha = 0$ 
or all $1 \leqslant \sigma \leqslant m$ 
shows that $m$ terms drop in all the correction sums.

If we now write only the $r-m$ recurrence formulas
which concern the remaining $r-m$ 
phantom differential invariants:
\[
\inv\big(
w_{\nu_{m+1}}^{(n_{m+1})}
(z,g)\big),\,\,
\dots\dots,\,\,
\inv\big(
w_{\nu_r}^{(n_r)}
(z,g)
\big),
\]
we receive for every fixed $j$ a system of $r-m$ linear equations in
the $r-m$ Maurer-Cartan invariants $\Kaux_j^{m+1}, \dots, \Kaux_j^r$.
This system is a subsystem of the full $r \times r$ system of
Theorem~{\ref{Thm-recurrence-formulae}} for $\Kaux_j^1, \dots,
\Kaux_j^{m},\, \Kaux_j^{m+1}, \dots, \Kaux_j^r$, and by
triangularity (by blocks), the $(r-m) \times (r-m)$ subsystem is 
{\em also} of Cramér type, hence can be solved uniquely for
$\Kaux_j^{m+1}, \dots, \Kaux_j^r$ in terms of $\Iaux^{(n_\GG+1)}$.

\Subsection{Groups containing translations and transvections}
\label{Subsection-G-1-MC-invariants}
Assume now that $G$ contains not only translations, but also
all vertical transvections:
\[
v^\alpha
\,=\,
u^\alpha
+
{\sf q}_1^\alpha\,x^1
+\cdots+
{\sf q}_p^\alpha\,x^p
\eqno
{\scriptstyle{(1\,\leqslant\,\alpha\,\leqslant\,q)}}.
\]
Thus, $\dim\, G = r \geqslant p + q + p\,q$.
The $p\, q$ infinitesimal generators are:
\[
x^i\,
\partial_{u^\alpha}
\eqno
{\scriptstyle{(1\,\leqslant\,i\,\leqslant\,p,\,\,
1\,\leqslant\,\alpha\,\leqslant\,q)}}.
\]

An application of
Theorem~{\ref{Thm-prolongation-vector-fields}}
shows that their infinite prolongations are truncated after
first order:
\[
\big(
x^i\,\partial_{u^\alpha}
\big)^{(\infty)}
\,=\,
x^i\,\partial_{u^\alpha}
+
\partial_{u_i^\alpha}
+
0.
\]
We list the first $p + q + p\,q$ infinitesimal generators as:
\[
{\bf v}_i
\,:=\,
\partial_{x^i},
\ \ \ \ \ \ \ \ \ \ \ \ \ \ \ \ \ \ \ \
{\bf v}_{p+\alpha}
\,:=\,
\partial_{u^\alpha},
\ \ \ \ \ \ \ \ \ \ \ \ \ \ \ \ \ \ \ \
\big\{
{\bf v}_{p+q+1},\dots,
{\bf v}_{p+q+pq}
\big\}
\,=\,
\big\{
x^i\,\partial_{u^\alpha}
\big\}.
\]
Hence we have:
\[
\aligned
\Big(
1\leqslant\sigma\leqslant p+q,
\ \ \ \ \
\# K\geqslant 1
\Big)
\ \ \ \ \ \ \ \ \ 
&
\Longrightarrow
\ \ \ \ \ \ \ \ \ 
{\varphi_\sigma}_K^\alpha
\,=\,
0,
\\
\Big(
p+q+1\leqslant\sigma\leqslant p+q+p\,q,
\ \ \ \ \
\# K\geqslant 2
\Big)
\ \ \ \ \ \ \ \ \ 
&
\Longrightarrow
\ \ \ \ \ \ \ \ \ 
{\varphi_\sigma}_K^\alpha
\,=\,
0.
\endaligned
\]

Among the cross-section equations, we may take:
\[
y^i
\,=\,
0,
\ \ \ \ \ \ \ \ \ \ \ \ \ \ \ \ \ \ \ \
v^\alpha
\,=\,
0,
\ \ \ \ \ \ \ \ \ \ \ \ \ \ \ \ \ \ \ \
v_i^\alpha
\,=\,
0.
\]
This conducts to $p + q + p\, q$ phantom invariants:
\[
\Jaux^i
\,=\,
0,
\ \ \ \ \ \ \ \ \ \ \ \ \ \ \ \ \ \ \ \
\Iaux_\emptyset^\alpha
\,=\,
0,
\ \ \ \ \ \ \ \ \ \ \ \ \ \ \ \ \ \ \ \
\Iaux_i^\alpha
\,=\,
0.
\]

Then in the fundamental recurrence formulas of  
Theorem~{\ref{Thm-recurrence-formulae}},
we decide to set aside the $p+q+p\,q$ ones concerning these
first phantom invariants, 
and the interesting recurrence relations 
are, for any $\# K \geqslant 2$:
\[
\aligned
\mathcal{D}_j
\Iaux_K^\alpha
&
\,=\,
\Iaux_{K,j}^\alpha
+
\sum_{1\leqslant\sigma\leqslant r}\,
{\varphi_\sigma}_K^\alpha
\big(
\Iaux^{(K)}
\big)
\cdot
\Kaux_j^\sigma
\big(
\Iaux^{(n_\GG+1)}
\big)
\\
&
\,=\,
\Iaux_{K,j}^\alpha
+
0
+
\sum_{m_1+1\leqslant\sigma\leqslant r}\,
{\varphi_\sigma}_K^\alpha
\big(
\Iaux^{(K)}
\big)
\cdot
\Kaux_j^\sigma
\big(
\Iaux^{(n_\GG+1)}
\big),
\endaligned
\]
where we abbreviate $m_1 := p + q + p\, q$. 

It remains
$r - m_1$ infinitesimal generators ${\bf v}_{m_1+1}, \dots, {\bf 
v}_r$. 
If we now write only the $r-m_1$ recurrence formulas
which concern the remaining $r-m_1$ 
phantom differential invariants:
\[
\inv\big(
w_{\nu_{m_1+1}}^{(n_{m_1+1})}
(z,g)\big),\,\,
\dots\dots,\,\,
\inv\big(
w_{\nu_r}^{(n_r)}
(z,g)
\big),
\]
we receive for every fixed $j$ a system of $r-m_1$ linear equations in
the $r-m_1$ Maurer-Cartan invariants $\Kaux_j^{m_1+1}, \dots, 
\Kaux_j^r$.
This system is a subsystem of the full $r \times r$ system of
Theorem~{\ref{Thm-recurrence-formulae}} for $\Kaux_j^1, \dots,
\Kaux_j^{m_1},\, \Kaux_j^{m_1+1}, \dots, \Kaux_j^r$, and by
triangularity (by blocks), the $(r-m_1) \times (r-m_1)$ subsystem is 
{\em also} of Cramér type, hence can be solved uniquely for
$\Kaux_j^{m_1+1}, \dots, \Kaux_j^r$ in terms of $\Iaux^{(n_\GG+1)}$.
This observation will be useful later on, in 
Section~{\ref{recurrence-relations-parabolic-surfaces}}.

\Section{\bf Moving Frame Method and Recurrence Relations for
Curves}
\label{moving-frame-recurrence-curves-R-2}
\HEAD{{\ref{moving-frame-recurrence-curves-R-2}}.~{\sf Moving Frame 
Method and Recurrence Relations for Curves}
}{
Zhangchi {\sc Chen}, Joël~{\sc Merker}}

In general, the structure of the full algebra of differential
invariants can be understood thanks to the {\sl recurrence
formulas} obtained by Fels and Olver 
in~{\cite{Fels-Olver-1999}}.

For any integer $n \geqslant 0$, 
let $J_{x,u}^n$ be the space of $n$\textsuperscript{th}
order jets of functions $\R \ni x \longmapsto u(x) \in \R$,
equipped with coordinates:
\[
\big(x,u,u_1,\dots,u_n\big),
\] 
where the $u_i$ stand for $\frac{\partial^i u}{\partial x^i}(x)$,
as abstract independent variables.

According to~\cite[Ch.~23]{Lie-Merker-2015}
and to~{\cite{Olver-1995, Merker-2008},
the $n$\textsuperscript{th} prolongation of a general vector field
on the $(x,u)$ space:
\[
{\bf v}
\,=\,
\xi(x,u)\,\frac{\partial}{\partial x}
+
\eta(x,u)\,\frac{\partial}{\partial u}
\]
is a vector field on $J_{x,u}^n$:
\[
{\bf v}^{(n)}
\,:=\,
{\bf v}
+
\sum_{1\leqslant i\leqslant n}\,
\Phi^i\big(x,u,u_1,\dots,u_i\big)\,
\frac{\partial}{\partial u_i},
\]
which has extended coefficients uniquely defined as:
\[
\Phi^i
\,:=\,
\underbrace{
{\sf D}_x
\big(\cdots
\big(
{\sf D}_x}_{i\,\,\text{\sf times}}
\big(
\eta
-
\xi\,u_1
\big)\big)
\cdots
\big)
+
\xi\,u_{i+1}
\eqno
{\scriptstyle{(1\,\leqslant\,i\,\leqslant\,n)}}.
\]
In~{\cite{Merker-2008}}, 
one finds {\em closed} explicit formulas for all these coefficients
$\Phi^i$ in any jet order $n \geqslant 1$, even for an arbitrary
number of independent variables $\big(x^1, \dots, x^p\big)$ and for an
arbitrary number of dependent variables $\big( u^1, \dots, u^q \big)$.
However, in small dimensions and for small jet order, it is
better and almost straightforward to capture these $\Phi^i$
using a computer.

In the present context, 
the group $\SL_2(\R)$ has $3$ natural infinitesimal generators:
\[
{\bf v}_1
\,:=\,
x\,\partial_x
-
u\,\partial_u,
\ \ \ \ \ \ \ \ \ \ \ \ \ \ \ \ 
{\bf v}_2
\,:=\,
u\,\partial_x,
\ \ \ \ \ \ \ \ \ \ \ \ \ \ \ \ 
{\bf v}_3
\,:=\,
x\,\partial_u,
\]
whose prolongations to the $n$\textsuperscript{th}
order jet space are, as a computer tells us:
\[
\aligned
{\bf v}_1^{(n)}
&
\,=\,
x\,\partial_x
-
u\,\partial_u
-
2\,u_1\,\partial_{u_1}
-
3\,u_2\,\partial_{u_2}
-
4\,u_3\,\partial_{u_3}
-\cdots-
(n+1)\,u_n\,\partial_{u_n},
\\
{\bf v}_2^{(n)}
&
\,=\,
u\,\partial_x
-
u_1^2\,\partial_{u_1}
-
3\,u_1\,u_2\,\partial_{u_2}
-
\big(4\,u_1u_3+3\,u_2^2\big)\,
\partial_{u_3}
-
\big(
5\,u_1u_4+10\,u_2u_3
\big)\,
\partial_{u_4}
\,-
\\
&
\ \ \ \ \ \ \ \ \ \ \ \ \
-\,
\big(6\,u_1u_5+15\,u_2u_4+10\,u_3^2\big)\,
\partial_{u_5}
-
\big(7\,u_1u_6+21\,u_2u_5+35\,u_3u_4\big)\,
\partial_{u_6}
\,-
\\
&
\ \ \ \ \ \ \ \ \ \ \ \ \
-\,
\big(8\,u_1u_7+28\,u_2u_6+56\,u_3u_5+35\,u_4^2\big)\,
\partial_{u_7}
-
\cdots,
\\
{\bf v}_3^{(n)}
&
\,=\,
x\,\partial_u
+
\partial_{u_1},
\endaligned
\]
For any jet order $n \geqslant 1$, we would write:
\[
{\bf v}_\kappa^{(n)}
\,=\,
{\sf v}_\kappa
+
\sum_{1\leqslant i\leqslant n}\,
\Phi_\kappa^i\,
\partial_{u_i}
\eqno
{\scriptstyle{(1\,\leqslant\,\kappa\,\leqslant\,3)}}.
\] 
and reading the first terms of
${\bf v}_1^{(n)}$, ${\bf v}_3^{(n)}$, it is visible that:
\[
\Phi_3^{(1)}
\,=\,
1,
\ \ \ \ \ \ 
\Phi_3^{(2)}
\,=\,
\cdots
\,=\,
\Phi_3^{(n)}
\,=\,
0,
\ \ \ \ \ \ \ \ \ \ \ \ \ \ \ \ \ \ \ \
\Phi_1^{(k)}
\,=\,
-\,(k+1)\,
u_k
\eqno
{\scriptstyle{(1\,\leqslant\,k\,\leqslant\,n)}},
\]
while it can be proved (exercise) that when $n =: 2n$ is even:
\[
\Phi_2^{2n}
\,=\,
-\,\sum_{1\leqslant k\leqslant n}\,
{\textstyle{\binom{2n+1}{k}}}\,
u_k\,u_{2n+1-k},
\]
and when $n =: 2n+1$ is odd:
\[
\Phi_2^{2n+1}
\,=\,
-\,\sum_{1\leqslant k\leqslant n}\,
{\textstyle{\binom{2n+2}{k}}}\,
u_k\,u_{2n+2-k}
-
{\textstyle{\frac{1}{2}}}\,
{\textstyle{\binom{2n+2}{n+1}}}\,
\big(
u_{n+1}
\big)^2.
\]

Fels and Olver in~{\cite{Fels-Olver-1999}} introduced an
{\sl invariantization operator} `$\inv$' which transforms
every function on the jet space into a differential invariant.
After having centered our surface at the origin by 
means of a translation, 
and after having set to 
constant three other Taylor coefficients
as in Theorem~{\ref{Thm-all-Taylor-coefficients-invariants}}, 
we come to:
\[
\aligned
\inv(x)
&
\,=\,
0,
\ \ \ \ \ \ \ \ \ \ \ \ \ \ \ \ \ \ \ \
\inv(u)
\,=\,
0,
\\
\Iaux_1
\,:=\,
\inv\big(u_x\big)
&
\,=\,
0,
\ \ \ \ \ \ \ \ \ \ \ \ \ \ \ \ \ \ \ \
\Iaux_2
\,:=\,
\inv\big(u_{xx}\big)
\,=\,
1,
\ \ \ \ \ \ \ \ \ \ \ \ \ \ \ \ \ \ \ \
\Iaux_3
\,:=\,
\inv\big(u_{xxx}\big)
\,=\,
0.
\endaligned
\]
Generally, every pure jet monomial produces a differential invariant:
\[
\Iaux_k
\,:=\,
\inv\big(u_{x^k}\big)
\eqno
{\scriptstyle{(k\,\geqslant 0)}},
\]
whose expression is most of the time lengthy and complex, 
while the $5$ constantified
invariants $\inv(x)$, $\inv(u)$,
$\Iaux_1$, $\Iaux_2$, $\Iaux_3$ are
called {\sl phantom invariants}.

The recurrence formulas obtained by Fels and Olver 
in~{\cite{Fels-Olver-1999}}
can be set up {\em without knowing}
the (huge) explicit expressions of differential invariants and
are of the form:
\begin{equation}
\label{1-D-recurrence-relations}
\Iaux_{k+1}
\,=\,
\mathcal{D}_x\big(\Iaux_k\big)
-
\sum_{1\leqslant\kappa\leqslant 3}\,
\inv\big(
\Phi_\kappa^k
\big)\,
\Raux^\kappa
\ \ \ \ \ \ \ \ \ \ \ \ \ \ \ \ \ \ \ \
{\scriptstyle{(k\,\geqslant\,0)}},
\end{equation}
where the invariant differentiation
operator $\mathcal{D}_x$ already introduced above
in~({\ref{1-dim-D-x-invariant}}) 
is a nonzero multiple of the total differentiation operator
${\sf D}_x$,
where all invariantized coefficients of the
prolonged vector fields are simply:
\[
\inv
\Big(
\Phi_\kappa^k
\big(
x,u,u_1,u_2,u_3,u_4,\dots,u_k
\big)
\Big)
\,=\,
\Phi_\kappa^k
\big(
0,0,0,1,0,\Iaux_4,\dots,\Iaux_k
\big)
\]
due to the fundamental commutation:
\[
{\sf invariantization}
\circ
{\sf function}
\,=\,
{\sf function}
\circ
{\sf invariantization},
\]
and where the Maurer-Cartan 
coefficients $\Raux^1$, $\Raux^2$, $\Raux^3$\,\,---\,\,which are also
differential invariants\,\,---\,\,can be determined 
by solving the linear system composed of the three 
recurrence relations applied to the three phantom invariants,
technically as follows. 

One writes the $3 \times 3$
matrix of the coefficients of the three generators
${\bf v}_1$, ${\bf v}_2$, ${\bf v}_3$ of $\mathfrak{sl}_2(\R)$ with
respect to the three basic fields $\partial_{u_1}$, $\partial_{u_2}$,
$\partial_{u_3}$ corresponding to the three invariants
$\inv(u_1)$, $\inv(u_2)$, $\inv(u_3)$ which
were phantom:
\[
\begin{blockarray}{cccc}
{} & \partial_{u_1} & \partial_{u_2} & \partial_{u_3} 
\\
\begin{block}{c(ccc)}
{\bf v}_1 & -2u_1 & -3u_2 & -4u_3 
\\
{\bf v}_2 & -u_1^2 & -3u_1u_2 & -4u_1u_3-3u_2^2 
\\
{\bf v}_3 & 1 & 0 & 0
\\ 
\end{block}
\end{blockarray}
\ \ \ \ \ \ \ \ \ \ \ \ \ \ \ \ \ \ \ \
\xrightarrow[{\rule[0pt]{50pt}{0pt}}]{{\sf invariantize}}
\ \ \ \ \ \ \ \ \ \ \ \ \ \ \ \ \ \ \ \
\left(\!
\begin{array}{ccc}
0 & -3 & 0
\\
0 & 0 & -3
\\
1 & 0 & 0
\end{array}
\!\right),
\]
one invariantizes this matrix, one transposes it, and one 
gets the
phantom recurrence relations~({\ref{1-D-recurrence-relations}}) 
for $k = 1, 2, 3$:
\[
\left(\!
\begin{array}{c}
\Iaux_2
\\
\Iaux_3
\\
\Iaux_4
\end{array}
\!\right)
\,=\,
\left(\!
\begin{array}{c}
\mathcal{D}_x(\Iaux_1)
\\
\mathcal{D}_x(\Iaux_2)
\\
\mathcal{D}_x(\Iaux_3)
\end{array}
\!\right)
-
\left(\!
\begin{array}{ccc}
0 & 0 & 1
\\
-3 & 0 & 0 
\\
0 & -3 & 0
\end{array}
\!\right)\,
\left(\!
\begin{array}{c}
\Raux^1
\\
\Raux^2
\\
\Raux^3
\end{array}
\!\right),
\]
which, if we abbreviate the parabolas invariant as:
\[
\Paux
\,:=\,
\inv\big(u_{xxxx}\big)
\]
become:
\[
\left(\!
\begin{array}{c}
1
\\
0
\\
\Paux
\end{array}
\!\right)
\,=\,
\left(\!
\begin{array}{c}
0
\\
0
\\
0
\end{array}
\!\right)
-
\left(\!
\begin{array}{ccc}
0 & 0 & 1
\\
-3 & 0 & 0 
\\
0 & -3 & 0
\end{array}
\!\right)\,
\left(\!
\begin{array}{c}
\Raux^1
\\
\Raux^2
\\
\Raux^3
\end{array}
\!\right),
\]
and this simple linear system has as the unique solution:
\[
\Raux^1
\,:=\,
0
\ \ \ \ \ \ \ \ \ \ \ \ \ \ \ \ \ \ \ \
\Raux^2
\,:=\,
{\textstyle{\frac{1}{3}}}\,
\Paux,
\ \ \ \ \ \ \ \ \ \ \ \ \ \ \ \ \ \ \ \
\Raux^3
\,:=\,
\,-1.
\]

Once these three Maurer-Cartan invariants have been determined,
they can be plugged in {\em all} the other recurrence
formulas~({\ref{1-D-recurrence-relations}}). 
In particular for $k = 4, 5$, determining 
first the coefficients of $\partial_{u_4}$, $\partial_{u_5}$:
\[
\begin{blockarray}{ccc}
& \partial_{u_4} & \partial_{u_5}
\\
\begin{block}{c(cc)}
{\bf v}_1 & -5u_4 & -6u_5
\\
{\bf v}_2 & -5u_1u_4-10u_2u_3 & -6u_1u_5-15u_2u_4-10u_3^2
\\
{\bf v}_3 & 0 & 0
\\ 
\end{block}
\end{blockarray}
\ \ \ \ \ \ \ \ \ \ \ \ \ \ \ \ \ \ \ \
\xrightarrow[{\rule[0pt]{50pt}{0pt}}]{{\sf invariantize}}
\ \ \ \ \ \ \ \ \ \ \ \ \ \ \ \ \ \ \ \
\left(\!
\begin{array}{cc}
-5\,\Paux & -6\,\Iaux_5
\\
0 & -15\,\Paux
\\
0 & 0
\end{array}
\!\right),
\]
we reach an invariant 
which we abbreviate as:
\[
\aligned
\Maux
\,:=\,
\Iaux_5
&
\,=\,
\mathcal{D}_x\big(\Paux\big)
-
\big(-5\,\Paux\big)
\cdot
0
-
0
\cdot
{\textstyle{\frac{1}{3}}}\,
\Paux
-
0
\cdot
(-1)
\\
&
\,=\,
\mathcal{D}_x\big(\Paux\big),
\\
\Iaux_6
&
\,=\,
\mathcal{D}_x\big(\Iaux_5\big)
-
\big(-6\,\Maux\big)
\cdot
0
-
\big(-15\,\Paux\big)
\cdot
{\textstyle{\frac{1}{3}}}\,
\Paux
-
0
\cdot
(-1)
\\
&
\,=\,
\mathcal{D}_x\big(\mathcal{D}_x(\Paux)\big)
+
5\,\Paux^2.
\endaligned
\]
These first two formulas suggest that the parabolas invariant $\Paux$
together with all its invariants derivatives $\mathcal{D}_x^\nu \Paux$
generate the algebra of differential invariants, a result established
by Olver in~{\cite{Olver-2018}}. Let us explain this a bit more.

Indeed, if we abbreviate the relevant 
invariantizations of the coefficients of ${\bf v}_2^{(n)}$ as:
\[
\Lambda_k
\,:=\,
\inv
\big(
\Phi_2^k
\big)
\eqno
{\scriptstyle{(k\,\geqslant 4)}},
\]
then we have:
\[
\footnotesize
\aligned
\Lambda_4
&
\,=\,
0,
\ \ \ \ & \ \ \ \
\Iaux_5
&
\,=\,
\mathcal{D}_x\Iaux_4
+
{\textstyle{\frac{\Paux}{3}}}\,
\Lambda_4
\,=\,
\mathcal{D}_x\Paux,
\\
\Lambda_5
&
\,=\,
15\,\Iaux_4,
\ \ \ \ & \ \ \ \
\Iaux_6
&
\,=\,
\mathcal{D}_x\Iaux_5
+
{\textstyle{\frac{\Paux}{3}}}\,
\Lambda_5
\,=\,
\mathcal{D}_x^2\Paux
+
5\,\Paux^2,
\\
\Lambda_6
&
\,=\,
21\,\Iaux_5,
\ \ \ \ & \ \ \ \
\Iaux_7
&
\,=\,
\mathcal{D}_x\Iaux_6
+
{\textstyle{\frac{\Paux}{3}}}\,
\Lambda_6
\,=\,
\mathcal{D}_x^3\Paux
+
17\,\mathcal{D}_x\Paux\,\Paux,
\\
\Lambda_7
&
\,=\,
28\,\Iaux_6
+
35\,\Iaux_4^2,
\ \ \ \ & \ \ \ \
\Iaux_8
&
\,=\,
\mathcal{D}_x^4\Paux
+
17\,\big(\mathcal{D}_x\Paux\big)^2
+
{\textstyle{\frac{79}{3}}}\,
\mathcal{D}_x^2\Paux\,\Paux
+
{\textstyle{\frac{175}{3}}}\,
\Paux^3,
\\
\Lambda_8
&
\,=\,
36\,\Iaux_7
+
126\,\Iaux_4\Iaux_5,
\ \ \ \ & \ \ \ \
\Iaux_9
&
\,=\,
\mathcal{D}_x^5\Paux
+
{\textstyle{\frac{181}{3}}}\,
\mathcal{D}_x^2\Paux\,
\mathcal{D}_x\Paux
+
{\textstyle{\frac{115}{3}}}\,
\mathcal{D}_x^3\Paux\,\Paux
+
421\,
\mathcal{D}_x\Paux\,
\Paux^2,
\endaligned
\]
Although it is clear that $\big\{ \mathcal{D}_x^\nu \Paux \big\}_{\nu
\in \N}$ generates the full algebra of differential invariants
of curves under $\SA_2(\R)$, it is an open question
({\cite{Olver-2018}})
to get {\em closed} explicit formulas for all the $\Iaux_k$.

\smallskip

In summary, what we have done in our basic
Section~{\ref{special-affine-power-series-invariants-curves-R-2}} is
to decompose into successive steps\,\,---\,\,stratified by increasing
jet orders viewed in power series centered at the
origin\,\,---\,\,all the elimination computations which are required
to {\em normalize} plane curves under special affine transformations,
while the famous moving frame method proceeds with a choice of a
suitable cross-section in a sufficiently high order
jet space so as to
perform all the elimination computations in one stroke.

Of course, the two processes are mathematically equivalent, but the
interest to decompose computations into steps is to be able to grasp
some explicit expressions of at least a few differential invariants,
and in dimension $2$ for parabolic surfaces $S^2 \subset \R^3$, we
will soon see that an explicit knowledge of certain `bifurcating'
invariants is necessary in order to study the full ramified tree of
equivalent classes.

\Section{\bf Affine Differential Invariants
of Curves in $\R^2$}
\label{affine-invariants-curves-R-2}
\HEAD{{\ref{affine-invariants-curves-R-2}}.~{\sf Affine 
Differential Invariants of Curves in $\R^2$}
}{
Zhangchi {\sc Chen}, Joël~{\sc Merker}}

But before we pass to dimension $2$, let us show how the story changes
much when considering the full affine group ${\sf A}_2(\R)$, with $6 >
5$ parameters.  In fact, such a study is not artificial, for it will
be required to complete our understanding of a certain (thin) branch
in our later study of parabolic surfaces $S^2 \subset \R^2$ modulo the
special affine group $\SA_3(\R)$.

As all symbolic computation softwares 
do handle polynomials very efficiently, 
we will employ the power series method,
and we will proceed as explained in the beginning
of Section~{\ref{special-affine-power-series-invariants-curves-R-2}}.

\Subsection{First loop}
Leaving aside the condition ${\sf a}{\sf d} - {\sf b}{\sf c} = 1$,
we therefore start with general $\GL_2(\R)$ matrices:
\[
\left(\!
\begin{array}{cc}
{\sf a} & {\sf b}
\\
{\sf c} & {\sf d}
\end{array}
\!\right)
\eqno
{\scriptstyle{({\sf a}{\sf d}\,-\,{\sf b}{\sf c}\,\neq\,0)}}.
\]
With $x = {\sf a}\, y + {\sf b}\, v$ and $u = {\sf c}\, y + 
{\sf d}\, v$, 
the fundamental equation writes again:
\leqnomode\usetagform{default}
\begin{equation}
\label{eq-fond-F-G-GL-2}
0
\,\equiv\,
-\,{\sf c}\,y
-
{\sf d}\,G(y)
+
F\big({\sf a}\,y+{\sf b}\,G(y)\big).
\end{equation}
As for $\SL_2(\R)$, we may assume that first-order terms are 
normalized to vanish:
\[
u
\,=\,
{\rm O}_x(2)
\ \ \ \ \ \ \ \ \ \ \ \ \ \ \ \ \ \ \ \
\text{and}
\ \ \ \ \ \ \ \ \ \ \ \ \ \ \ \ \ \ \ \
v
\,=\,
{\rm O}_y(2).
\]

\begin{Lemma}
\label{LM-c-zero-GL-2}
The subgroup of $\GL_2(\R)$ sending
$v = {\rm O}_y(2)$ 
to 
$u = {\rm O}_x(2)$ 
is $3$-dimensional and consists of matrices:
\[
G_{\stabsmall}^{(1)}
\colon
\ \ \ \ \ \ \ \ \ \ \ \ \ \ \ \ \ \ \ \
\left(\!
\begin{array}{cc}
{\sf a} & {\sf b}
\\
0 & {\sf d}
\end{array}
\!\right)
\eqno
{\scriptstyle{({\sf a}{\sf d}\,\neq\,0)}}.
\qed
\]
\end{Lemma}

Next, assuming of course still that $F_2 \neq 0 \neq G_2$ in:
\[
u
\,=\,
F_2\,
{\textstyle{\frac{x^2}{2!}}}
+
{\rm O}_x(3)
\ \ \ \ \ \ \ \ \ \ \ \ \ \ \ \ \ \ \ \
\text{and}
\ \ \ \ \ \ \ \ \ \ \ \ \ \ \ \ \ \ \ \
v
\,=\,
G_2\,
{\textstyle{\frac{y^2}{2!}}}\,
+
{\rm O}_y(3),
\] 
we examine how this reduced group $G_{\stabsmall}^{(1)}$ acts
on the second order term $y^2$ of the
fundamental equation~({\ref{eq-fond-F-G-GL-2}}):
\[
0
\,\equiv\,
\big(
-\,{\sf d}\,G_2
+
{\sf a}^2\,F_2
\big)\,
{\textstyle{\frac{y^2}{2!}}}\,
+
{\rm O}_y(3).
\]
Consequently, we can make $G_2 := 1$ 
for the curve in the target space $\R_{y,v}^2$ 
by choosing ${\sf d} := {\sf a}^2\, 
F_2$, and with the simplest choice ${\sf a} := 1$, 
we apply the invertible matrix:
\[
\left(\!
\begin{array}{cc}
1 & 0
\\
0 & F_2
\end{array}
\!\right)
\,\,\in\,\,
G_{\stabsmall}^{(1)}.
\]
Calculations of further terms in 
the fundamental equation~({\ref{eq-fond-F-G-GL-2}}) 
provide:
\[
G_3
\,=\,
\frac{F_3}{F_2},
\ \ \ \ \ \ \ \ \ \ \ \ \ \ \ \ \ \ \ \
G_4
\,=\,
\frac{F_4}{F_2},
\ \ \ \ \ \ \ \ \ \ \ \ \ \ \ \ \ \ \ \
G_5
\,=\,
\frac{F_5}{F_2},
\]
and if we come back to initial functional jets, we get new values:
\[
G_2
\,:=\,
1,
\ \ \ \ \ \ \ \ \ \ \ \ \ \ \ \ \ \ \ \
G_3
\,:=\,
\frac{u_3}{u_2},
\ \ \ \ \ \ \ \ \ \ \ \ \ \ \ \ \ \ \ \
G_4
\,:=\,
\frac{u_4}{u_2},
\ \ \ \ \ \ \ \ \ \ \ \ \ \ \ \ \ \ \ \
G_5
\,:=\,
\frac{u_5}{u_2}.
\]

\Subsection{Second loop} 
We restart from equations normalized up to order $2$:
\[
u
\,=\,
{\textstyle{\frac{x^2}{2!}}}
+
F_3\,
{\textstyle{\frac{x^3}{3!}}}
+\cdots,
\ \ \ \ \ \ \ \ \ \ \ \ \ \ \ \ \ \ \ \
v
\,=\,
{\textstyle{\frac{y^2}{2!}}}
+
G_3\,
{\textstyle{\frac{x^3}{3!}}}
+\cdots.
\]
To determine the subgroup 
${\sf G}_{\stabsmall}^{(2)} \subset {\sf G}_{\stabsmall}^{(1)}$
which stabilizes these normalizations, we compute
the first term in 
the fundamental equation~({\ref{eq-fond-F-G-GL-2}}):
\[
0
\,\equiv\,
\big(-{\sf d}+{\sf a}^2\big)\,
{\textstyle{\frac{y^2}{2!}}}
+
{\rm O}_y(3).
\]

\begin{Lemma}
The subgroup of $\GL_2(\R)$ sending
$v = \frac{y^2}{2!} + {\rm O}_y(3)$ 
to 
$u = \frac{x^2}{2!} + {\rm O}_x(3)$ 
is $2$-dimensional and consists of matrices:
\[
G_{\stabsmall}^{(2)}
\colon
\ \ \ \ \ \ \ \ \ \ \ \ \ \ \ \ \ \ \ \
\left(\!
\begin{array}{cc}
{\sf a} & {\sf b}
\\
0 & {\sf a}^2
\end{array}
\!\right)
\eqno
{\scriptstyle{({\sf a}\,\neq\,0)}}.
\qed
\]
\end{Lemma}

Next, we examine how this reduced group ${\sf G}_{\stabsmall}^{(2)}$
acts on the third order term $y^3$:
\[
0
\,\equiv\,
\big(
-\,{\sf a}^2\,G_3
+
3\,{\sf a}{\sf b}
+
{\sf a}^3\,F_3
\big)\,
{\textstyle{\frac{y^3}{3!}}}
+
{\rm O}_y(4).
\]
Consequently, we can make $G_3 := 0$ by choosing 
${\sf b} := - \frac{1}{3}\,
{\sf a}^2\,F_3$, and with the simple matrix:
\[
\left(\!
\begin{array}{cc}
1 & -\frac{1}{3}F_3
\\
0 & 1
\end{array}
\!\right),
\]
we transform $\big\{ u = F(x) \big\}$ into:
\[
v
\,=\,
{\textstyle{\frac{y^2}{2!}}}
+
0
+
\big(
F_4
-
{\textstyle{\frac{5}{3}}}\,
F_3^2
\big)\,
{\textstyle{\frac{y^4}{4!}}}
+
\big(
F_5
-
5\,F_3F_4
+
{\textstyle{\frac{40}{9}}}\,
F_3^3
\big)\,
{\textstyle{\frac{y^5}{5!}}}
+
{\rm O}_y(6).
\] 

Lastly, coming back to the initial functional jets, 
we get new values:
\[
\aligned
G_2
\,:=\,
1,
\ \ \ \ 
G_3
\,:=\,
0,
\ \ \ \
G_4
&
\,:=\,
F_4-{\textstyle{\frac{5}{3}}}\,F_3^2
\,=\,
{\textstyle{\frac{1}{3}}}\,
\frac{-\,5\,u_3^2+3\,u_2u_4}{u_2^2},
\\
G_5
&
\,:=\,
F_5-5\,F_3\,F_4+{\textstyle{\frac{40}{9}}\,F_3^3}
\,=\,
{\textstyle{\frac{1}{9}}}\,
\frac{9\,u_2^2u_5-45\,u_2u_3u_4+40\,u_3^3}{u_2^3}.
\endaligned
\]

\Subsection{Third loop} 
We restart from equations normalized up to order $3$:
\[
u
\,=\,
{\textstyle{\frac{x^2}{2!}}}
+
0
+
F_4\,
{\textstyle{\frac{x^4}{4!}}}
+\cdots,
\ \ \ \ \ \ \ \ \ \ \ \ \ \ \ \ \ \ \ \
v
\,=\,
{\textstyle{\frac{y^2}{2!}}}
+
0
+
G_4\,
{\textstyle{\frac{y^4}{4!}}}
+\cdots.
\]
To determine the subgroup 
${\sf G}_{\stabsmall}^{(3)} \subset {\sf G}_{\stabsmall}^{(2)}$
which stabilizes these normalizations, we compute
the first term in 
the fundamental equation~({\ref{eq-fond-F-G-GL-2}}):
\[
0
\,\equiv\,
\big(3\,{\sf a}{\sf b}\big)\,
{\textstyle{\frac{y^3}{3!}}}
+
{\rm O}_y(4).
\]

\begin{Lemma}
The subgroup of $\GL_2(\R)$ sending
$v = \frac{y^2}{2!} + 0 + {\rm O}_y(4)$ 
to 
$u = \frac{x^2}{2!} + 0 + {\rm O}_x(4)$ 
is $1$-dimensional and consists of matrices:
\[
G_{\stabsmall}^{(3)}
\colon
\ \ \ \ \ \ \ \ \ \ \ \ \ \ \ \ \ \ \ \
\left(\!
\begin{array}{cc}
{\sf a} & 0
\\
0 & {\sf a}^2
\end{array}
\!\right)
\eqno
{\scriptstyle{({\sf a}\,\neq\,0)}}.
\qed
\]
\end{Lemma}

Next, we examine how this reduced group ${\sf G}_{\stabsmall}^{(3)}$
acts on the third order term $y^4$:
\[
0
\,\equiv\,
\big(
-\,{\sf a}^2\,G_4
+
{\sf a}^4\,F_4
\big)\,
{\textstyle{\frac{y^4}{4!}}}
+
{\rm O}_y(5).
\]

At this stage, because ${\sf a} \neq 0$, 
a bifurcation necessarily occurs, according to whether
$F_4 = 0$ or $F_4 \neq 0$. In fact, coming back
to functional jets:
\[
F_4
\,=\,
{\textstyle{\frac{1}{3}}}\,
\frac{-\,5\,u_3^2+3\,u_2u_4}{u_2^2},
\]
the branch $F_4 \equiv 0$ must be interpreted as an identical 
vanishing in the $4$\textsuperscript{th} 
order jet space, and we have the elementary

\begin{Lemma}
\label{LM-parabola-Halphen-R}
A plane curve $\big\{ u = F(x) \big\}$ with $F_{xx} \neq 0$
is affinely equivalent to the model parabola $\big\{ v = y^2 \big\}$
if and only if:
\[
0
\,\equiv\,
{\textstyle{\frac{1}{3}}}\,
\frac{
-\,5\,F_{xxx}^2
+
3\,F_{xx}\,F_{xxxx}}{F_{xx}^{2}}.
\]
\end{Lemma}

\proof
A possible proof has already been left
as an exercise in Lemma~{\ref{LM-vanishing-H-M}},
but it is enlightening to provide another proof which requires
no integration, only basic differential algebra,
{\em cf.}~{\cite[Thm.~4.1]{Merker-2019}}.

The ordinary differential equation $0 \equiv -\,5 F_{xxx}^2 + 
3\, F_{xx} F_{xxxx}$ is invariant
not only under
the $\SA_2(\R)$ action, but also under general affine transformations,
for one can prove directly
({\cite[p.~26]{Merker-2019}})
that
if a graph $\big\{ u = F(x) \big\}$ is transformed into
a graph
$\big\{ v = G(y) \big\}$ through:
\[
x
\,=\,
{\sf a}\,y
+
{\sf b}\,v
+
{\sf c},
\ \ \ \ \ \ \ \ \ \ \ \ \ \ \ \ \ \ \ \
u
\,=\,
{\sf k}\,y
+
{\sf l}\,v
+
{\sf m}
\eqno
{\scriptstyle{({\sf a}{\sf l}\,-\,{\sf b}{\sf k}\,\neq\,0)}},
\]
then:
\[
\aligned
F_{xx}
&
\,=\,
\frac{{\sf a}{\sf l}-{\sf b}{\sf k}}{
\big({\sf a}+{\sf b}\,G_y\big)^3}\,
G_{yy},
\\
-\,
5\,F_{xxx}^2
+
3\,F_{xx}\,F_{xxxx}
&
\,=\,
\frac{({\sf a}{\sf l}-{\sf b}{\sf k})^2}{
\big({\sf a}+{\sf b}\,G_y\big)^8}\,
\Big[
-\,
5\,G_{yyy}^2
+
3\,G_{yy}\,G_{yyyy}
\Big].
\endaligned
\]

Therefore, after some preliminary affine transformation
as done in 
Section~{\ref{special-affine-power-series-invariants-curves-R-2}},
we can assume that the equation of our curve is already
normalized as:
\[
v
\,=\,
G(y)
\,=\,
y^2
+
0
+
{\rm O}_y(4).
\]
Again, we have $0 \equiv -5\, G_{yyy}^2 + 3\, G_{yy} G_{yyyy}$
and $G_{yy} \neq 0$ as well.
To conclude, it would then suffice to show that the remainder
${\rm O}_y(4) \equiv 0$ vanishes identically\,\,---\,\,very easy!

Indeed, we can solve:
\[
G_{yyyy}
\,=\,
{\textstyle{\frac{5}{3}}}\,
\frac{G_{yyy}^2}{G_{yy}},
\]
and from $G_{yyy}(0) = 0$, we deduce $G_{yyyy}(0) = 0$.
Let us abbreviate this relation as:
\[
G_{yyyy}
\,=\,
\mathcal{R}\,
G_{yyy},
\]
where $\mathcal{R} = \mathcal{R}(y)$ denotes an unspecified function.
Then we may differentiate it and replace:
\[
G_{yyyyy}
\,=\,
\mathcal{R}_y\,G_{yyy}
+
\mathcal{R}\,G_{yyyy}
\,=\,
\mathcal{R}\,G_{yyy},
\]
to deduce similarly $G_{yyyyy}(0) = 0$.
Obviously, induction yields $G_{y^k} = \mathcal{R}\, G_{yyy}$
for any $k \geqslant 4$, 
whence $G_{y^k}(0) = 0$, 
and  we conclude ${\rm O}_y(4) \equiv 0$,
since we assume analyticity.
\endproof

The case $F_4 = 0$ being now fully understood, we can assume 
that $F_4 \neq 0$, and then, two symmetric subcases occur:
\[
F_4
\,>\,
0
\ \ \ \ \ \ \ \ \ \ \ \ \ \ \ \ \ \ \ \
\text{and}
\ \ \ \ \ \ \ \ \ \ \ \ \ \ \ \ \ \ \ \
F_4
\,<\,
0.
\]
Assuming $F_4 > 0$, we can then make $G_4 := 1$ by means of the
simple matrix:
\[
\left(\!
\begin{array}{cc}
F_4^{-1/2} & 0
\\
0 & F_4^{-1}
\end{array}
\!\right).
\]
Looking at next terms in the fundamental equation, we get:
\[
G_5
\,:=\,
\frac{F_5}{F_4^{3/2}},
\ \ \ \ \ \ \ \ \ \ \ \ \ \ \ \ \ \ \ \
G_6
\,:=\,
\frac{F_6}{F_4^2},
\ \ \ \ \ \ \ \ \ \ \ \ \ \ \ \ \ \ \ \
G_7
\,:=\,
\frac{F_7}{F_4^{5/2}}.
\]
Coming back to functional jets:
\[
\aligned
F_5
&
\,:=\,
\frac{1}{\sqrt{3}}\,
\frac{9\,u_2^2u_5-45\,u_2u_3u_4+40\,u_3^3}{
\big(3\,u_2u_4-5\,u_3^2\big)^{3/2}},
\\
F_6
&
\,:=\,
\frac{9\,u_2^3u_6
-
63\,u_2^2u_3u_5
+
105\,u_2u_3^2u_4
-
35\,u_3^4}{
\big(3\,u_2u_4-5\,u_3^2\big)^2}.
\endaligned
\]

\Subsection{Fourth loop} 
The subgroup stabilizing equations normalized up to order $4$:

\[
u
\,=\,
{\textstyle{\frac{x^2}{2!}}}
+
0
+
{\textstyle{\frac{x^4}{4!}}}
+\cdots,
\ \ \ \ \ \ \ \ \ \ \ \ \ \ \ \ \ \ \ \
v
\,=\,
{\textstyle{\frac{y^2}{2!}}}
+
0
+
{\textstyle{\frac{y^4}{4!}}}
+\cdots.
\]
is easily seen to reduce to the identity.
The algorithm therefore stops, and a first result, valid only at the
level of power series at the origin, is a corollary of this reduction
to an $\{e\}$-group.

\begin{Theorem}
\label{Thm-GL-2-coefficients-invariants}
{\bf (1)}\,
Given a real analytic curve $\big\{ u = F(x) \big\}$ in $\R^2$ 
passing through the origin which satisfies: 
\[
F_{xx}(0)
\,\neq\,
0
\ \ \ \ \ \ \ \ \ \ \ \ \ \ \ \ \ \ \ \
\text{and}
\ \ \ \ \ \ \ \ \ \ \ \ \ \ \ \ \ \ \ \
\pm\,
\big(
F_{xx}F_{xxxx}
-
{\textstyle{\frac{5}{3}}}\,
F_{xxx}^2
\big)
\,>\,
0,
\]
there always exist a ${\sf GL}_2(\R)$ transformation which
puts it into the form
\[
u
\,=\,
{\textstyle{\frac{x^2}{2!}}}
+
0
\pm
{\textstyle{\frac{x^4}{4!}}}
+
F_5\,
{\textstyle{\frac{x^5}{5!}}}
+
\sum_{i\geqslant 6}\,
F_i\,
{\textstyle{\frac{x^i}{i!}}}.
\]

\noindent{\bf (2)}\,
Any other such real analytic curve $\big\{ v = 
G(y) \big\}$ similarly
put into the form:
\[
v
\,=\,
{\textstyle{\frac{y^2}{2!}}}
+
0
\pm
{\textstyle{\frac{y^4}{4!}}}
+
G_5\,
{\textstyle{\frac{y^5}{5!}}}
+
\sum_{j\geqslant 6}\,
G_j\,
{\textstyle{\frac{y^j}{j!}}}.
\]
is ${\sf GL}_2(\R)$-equivalent to $\big\{ u = F(x) \big\}$ above if
and only if all Taylor coefficients match:
\[
G_5
\,=\,
F_5,
\ \ \ \ \ \ \ \ \ \ \ \ \
G_i
\,=\,
F_i
\eqno
{\scriptstyle{(\forall\,i\,\geqslant\,6)}}.
\qed
\]
\end{Theorem}

\Subsection{Recurrence relations}
Similarly as in Section~{\ref{moving-frame-recurrence-curves-R-2}},
let us set up the corresponding recurrence relations.
By Theorem~{\ref{Thm-GL-2-coefficients-invariants}},
we come to:
\[
\!\!\!\!\!\!\!\!\!\!\!\!\!\!\!
\aligned
\inv(x)
&
\,=\,
0,
\ \ \ \ & \ \ \ \ 
\inv(u)
&
\,=\,
0,
\ \ \ \ & \ \ \ \
\Iaux_1
\,:=\,
\inv\big(u_x\big)
&
\,=\,
0,
\\
\Iaux_2
\,:=\,
\inv\big(u_{xx}\big)
&
\,=\,
1,
\ \ \ \ & \ \ \ \ 
\Iaux_3
\,:=\,
\inv\big(u_{xxx}\big)
&
\,=\,
0,
\ \ \ \ & \ \ \ \ 
\Iaux_4
\,:=\,
\inv\big(u_{xxxx}\big)
&
\,=\,
\pm 1.
\endaligned
\]
We take the $4$ infinitesimal generators of the the action
of ${\sf GL}_2(\R)\subset{\sf A}_2(\R)$ on $\R_{x,u}^2$:
\[
{\bf v}_1
\,:=\,
x\,\partial_x,
\ \ \ \ \ \ \ \ \ 
{\bf v}_2
\,:=\,
u\,\partial_u,
\ \ \ \ \ \ \ \ \ 
{\bf v}_3
\,:=\,
u\,\partial_x,
\ \ \ \ \ \ \ \ \ 
{\bf v}_4
\,:=\,
x\,\partial_u,
\]
the prolongations of which we write as:
\[
v_\kappa^{(\infty)}
\,=\,
v_\kappa
+
\sum_{k\geqslant 0}\,
\Phi_\kappa^k
\big(
x,u,u_1,\dots,u_k
\big)\,
\frac{\partial}{\partial u_k}
\eqno
{\scriptstyle{(1\,\leqslant\,\kappa\,\leqslant\,4)}}.
\]

Generally, every pure jet monomial produces a differential invariant:
\[
\Iaux_k
\,:=\,
\inv\big(u_{x^k}\big)
\eqno
{\scriptstyle{(k\,\geqslant 0)}},
\]
and the recurrence formulas are of the form:
\[
\Iaux_{k+1}
\,=\,
\mathcal{D}_x\big(\Iaux_k\big)
-
\sum_{1\leqslant\kappa\leqslant 4}\,
\inv\big(
\Phi_\kappa^k
\big)\,
\Raux^\kappa
\ \ \ \ \ \ \ \ \ \ \ \ \ \ \ \ \ \ \ \
{\scriptstyle{(k\,\geqslant\,0)}},
\]
for some uniquely defined invariant differentiation
operator $\mathcal{D}_x$ which  
is a nonzero multiple of the total differentiation operator
${\sf D}_x$, and
where, as before, all invariantized coefficients of the
prolonged vector fields are simply:
\[
\inv
\Big(
\Phi_\kappa^k
\big(
x,u,u_1,u_2,u_3,u_4,u_5,\dots,u_k
\big)
\Big)
\,=\,
\Phi_\kappa^k
\big(
0,0,0,1,0,\pm 1,\Iaux_5,\dots,\Iaux_k
\big),
\]
and where the Maurer-Cartan 
invariants $\Raux^1$, $\Raux^2$, $\Raux^3$, 
$\Raux^4$ can be determined as follows. 

One writes the $4 \times 4$
matrix of the coefficients of the four generators
${\bf v}_1$, ${\bf v}_2$, ${\bf v}_3$, ${\bf v}_4$ of $\mathfrak{gl}_2(\R)$ with
respect to the four basic fields $\partial_{u_1}$, $\partial_{u_2}$,
$\partial_{u_3}$,
$\partial_{u_4}$ corresponding to the four invariants
$\inv(u_1)$, $\inv(u_2)$, $\inv(u_3)$, $\inv(u_4)$ which
were phantom:
\[
\!\!\!\!\!\!\!\!\!\!\!\!\!\!\!
\begin{blockarray}{ccccc}
{} & \partial_{u_1} & \partial_{u_2} & \partial_{u_3} & \partial_{u_4} 
\\
\begin{block}{c(cccc)}
{\bf v}_1 & -u_1 & -2u_2 & -3u_3 & -4u_4
\\
{\bf v}_2 & u_1 & u_2 & u_3 & u_4
\\
{\bf v}_3 & -u_1^2 & -3u_1u_2 & -4u_1u_3-3u_2^2
& -5u_1u_4-10u_2u_3
\\
{\bf v}_4 & 1 & 0 & 0 & 0
\\ 
\end{block}
\end{blockarray}
\ \
\xrightarrow[{\rule[0pt]{50pt}{0pt}}]{{\sf invariantize}}
\ \
\left(\!
\begin{array}{cccc}
0 & -2 & 0 & \mp 4
\\
0 & 1 & 0 & \pm 1
\\
0 & 0 & -3 & 0
\\
1 & 0 & 0 & 0
\end{array}
\!\right),
\]
one invariantizes this matrix, one transposes it, and one 
gets the
phantom recurrence relations 
for $k = 1, 2, 3, 4$:
\[
\left(\!
\begin{array}{c}
\Iaux_2
\\
\Iaux_3
\\
\Iaux_4
\\
\Iaux_5
\end{array}
\!\right)
\,=\,
\left(\!
\begin{array}{c}
\mathcal{D}_x(\Iaux_1)
\\
\mathcal{D}_x(\Iaux_2)
\\
\mathcal{D}_x(\Iaux_3)
\\
\mathcal{D}_x(\Iaux_4)
\end{array}
\!\right)
-
\left(\!
\begin{array}{cccc}
0 & 0 & 0 & 1
\\
-2 & 1 & 0 & 0
\\
0 & 0 & -3 & 0
\\
\mp 4 & \pm 1 & 0 & 0
\end{array}
\!\right)\,
\left(\!
\begin{array}{c}
\Raux^1
\\
\Raux^2
\\
\Raux^3
\\
\Raux^4
\end{array}
\!\right),
\]
which become:
\[
\left(\!
\begin{array}{c}
1
\\
0
\\
\pm 1
\\
\Iaux_5
\end{array}
\!\right)
\,=\,
\left(\!
\begin{array}{c}
0
\\
0
\\
0
\\
0
\end{array}
\!\right)
-
\left(\!
\begin{array}{cccc}
0 & 0 & 0 & 1
\\
-2 & 1 & 0 & 0
\\
0 & 0 & -3 & 0
\\
\mp 4 & \pm 1 & 0 & 0
\end{array}
\!\right)\,
\left(\!
\begin{array}{c}
\Raux^1
\\
\Raux^2
\\
\Raux^3
\\
\Raux^4
\end{array}
\!\right),
\]
and this simple linear system has as the unique solution:
\[
\Raux^1
\,:=\,
\pm\,
{\textstyle{\frac{1}{2}}}\,
\Iaux_5,
\ \ \ \ \ \ \ \ \ \ \ \ \ \ \ \ \ \ \ \
\Raux^2
\,:=\,
\pm\,
\Iaux_5,
\ \ \ \ \ \ \ \ \ \ \ \ \ \ \ \ \ \ \ \
\Raux^3
\,:=\,
\pm\,
{\textstyle{\frac{1}{3}}},
\ \ \ \ \ \ \ \ \ \ \ \ \ \ \ \ \ \ \ \
\Raux^4
\,:=\,
-\,1.
\]

Once these four Maurer-Cartan invariants have been determined,
they can be plugged in:
\[
\begin{blockarray}{ccc}
& \partial_{u_5} & \partial_{u_6}
\\
\begin{block}{c(cc)}
{\bf v}_1 & -5u_5 & -6u_6
\\
{\bf v}_2 & u_5 & u_6
\\
{\bf v}_3 & -6u_1u_5-15u_2u_4-10u_3^2 &
-7u_1u_6-21u_2u_5-35u_3u_4
\\
{\bf v}_4 & 0 & 0
\\ 
\end{block}
\end{blockarray}
\ \
\xrightarrow[{\rule[0pt]{50pt}{0pt}}]{{\sf invariantize}}
\ \
\left(\!
\begin{array}{cc}
-5\,\Iaux_5 & -6\,\Iaux_6
\\
\Iaux_5 & \Iaux_6
\\
\mp 15 & -21\Iaux_5
\\
0 & 0
\end{array}
\!\right).
\]
We deduce:
\[
\aligned
\Iaux_6
&
\,=\,
\mathcal{D}_x(\Iaux_5)
\pm
{\textstyle{\frac{3}{2}}}
\,
\Iaux_5^2
+
5,
\\
\Iaux_7
&
\,=\,
\mathcal{D}_x(\Iaux_6)
\pm
2\,\Iaux_5\Iaux_6
\pm
7\,\Iaux_5,
\endaligned
\]
and it is easy to verify that the algebra of differential
invariants
is generated by $\Iaux_5$ and its invariant derivatives 
$\mathcal{D}_x^\nu(\Iaux_5)$ of any order $\nu \geqslant 1$.

\Section{\bf Parabolic Surfaces $S^2 \subset \R^3$: 
Invariant $\Waux$ of Order $4$}
\label{parabolic-pseudostablization}
\HEAD{{\ref{parabolic-pseudostablization}}.~{\sf Parabolic Surfaces 
$S^2 \subset \R^3$: Invariant $\Waux$ of Order $4$
}
}{
Zhangchi {\sc Chen}, Joël~{\sc Merker}}

Recall from Section~{\ref{parabolic-jet-relations}} that
we are interested in parabolic surfaces $S^2 \subset \R^3$
graphed as $\big\{ u = F(x,y) \big\}$,
whose Hessian matrix is degenerate, of rank $1 < 2$,
whence after a rotation in the $(x,y)$-space, we may assume:
\[
F_{xx}
\,\neq\,0
\,\equiv\,
F_{xx}F_{yy}
-
F_{xy}^2.
\]

In the same Section~{\ref{parabolic-jet-relations}}, we also showed
that the relation $F_{yy} = \frac{F_{xy}^2}{F_{xx}}$ has infinitely
many differential consequences, namely, for every $(j,k)$ with $k
\geqslant 2$, there exists a certain universal rational expressions
$R_{j,k}$ such that:
\[
F_{x^jy^k}
\,=\,
R_{j,k}
\Big(
\big\{
F_{j',0}
\big\}_{0\leqslant j'\leqslant j+k},\,\,
\big\{
F_{j'',1}
\big\}_{0\leqslant j''\leqslant j+k-1}
\Big),
\]
the denominators of these $R_{j,k}$ being
certain powers $\big( F_{xx} \big)^\ast$.
Inside the order $n$ jet space $J_{x,u}^n$, this means that:
\leqnomode\usetagform{default}
\begin{align}
u_{j,k}
\,=\,
R_{j,k}
\Big(
\big\{
u_{j',0}
\big\}_{0\leqslant j'\leqslant j+k},\,\,
\big\{
u_{j'',1}
\big\}_{0\leqslant j''\leqslant j+k-1}
\Big).
\end{align}

Theorem~{\ref{Thm-transfer-Hessian-determinants}} 
already showed that the assumption that the Hessian determinant
is identically zero is stable under (special) affine
transformations.
It is then interesting to have a confirmation of this 
basic fact
by means of prolongations of vector fields.

We start by selecting an appropriate basis for the $11$-dimensional
Lie algebra of infinitesimal generators of the action of
$\SA_3(\R)$ on the space $\R_{x,y,u}^3$, 
and we choose the same generators as in~{\cite{Olver-2007}}:
\[
\aligned
{\bf v_1}
&
\,:=\,
x\,\partial_x
-
u\,\partial_u, 
\ \ \ \ \ \ \ \ \ \ \
{\bf v}_2
\,:=\,
y\,\partial_y
-
u\,\partial_u,
\\
{\bf v}_3
\,:=\,
y\,\partial_x,
\ \ \ \ \ \ \
{\bf v}_4
\,:=\,
u\,\partial_x,
&
\ \ \ \ \ \ \
{\bf v}_5
\,:=\,
x\,\partial_y,
\ \ \ \ \ \ \
{\bf v}_6
\,:=\,
u\,\partial_y,
\ \ \ \ \ \ \
{\bf v}_7
\,:=\,
x\,\partial_u,
\ \ \ \ \ \ \
{\bf v}_8
\,:=\,
y\,\partial_u,
\\
{\bf w}_1
&
\,:=\,
\partial_x,
\ \ \ \ \ \ \ \ \ \ \
{\bf w}_2
\,:=\,
\partial_y,
\ \ \ \ \ \ \ \ \ \ \
{\bf w}_3
\,:=\,
\partial_u.
\endaligned
\]

Among these $11$ vector fields, $5$ are essentially useless,
namely the three translational ${\bf w}_1$, ${\bf w}_2$, ${\bf w}_3$,
and the two transvectional ${\bf v}_7$, ${\bf v}_8$,
as was already explained and understood by
Theorem~{\ref{Thm-translations-transvections-G-1}}.
A confirmation of this fact is provided by a computation of the
prolongations of these $5$ vector fields to any order:
\leqnomode\usetagform{default}
\begin{align}
\label{triangular-by-blocks-5x5}
\begin{blockarray}{ccccccc}
& \partial_x & \partial_y & \partial_u & \partial_{u_{1,0}} &
\partial_{u_{0,1}} & 
\big\{\partial_{u_{k,l}}\big\}_{2\leqslant k+l}
\\
\begin{block}{c(cccccc)}
{\bf v}_7 & 0 & 0 & x & 1 & 0 & 0
\\
{\bf v}_8 & 0 & 0 & y & 0 & 1 & 0
\\
{\bf w}_1 & 1 & 0 & 0 & 0 & 0 & 0
\\
{\bf w}_2 & 0 & 1 & 0 & 0 & 0 & 0
\\
{\bf w}_3 & 0 & 0 & 1 & 0 & 0 & 0
\\ 
\end{block}
\end{blockarray}\,\,,
\end{align}
with absolutely zero coefficients in front of all $\partial_{u_{j,k}}$
with $j + k \geqslant 2$.
Up to a permutation of rows and columns, we thus have a matrix
which is
triangular-by-blocks $5 = 2 + 3$, essentially the $5 \times 5$ identity
matrix, followed by a zero matrix. 

Consequently, if, for any fixed jet order $n \geqslant 2$, 
above any point $(x,y)$, 
we want to determine:
\[
\rank\,
\Span\,
\Big(
{\bf v}_1^{(n)},\,\,
{\bf v}_2^{(n)},\,\,
{\bf v}_3^{(n)},\,\,
{\bf v}_4^{(n)},\,\,
{\bf v}_5^{(n)},\,\,
{\bf v}_6^{(n)},\,\,
{\bf v}_7^{(n)},\,\,
{\bf v}_8^{(n)},\,\,
w_1^{(n)},\,\,
w_2^{(n)},\,\,
w_3^{(n)}
\Big),
\]
it suffices to determine the rank of only
$6$ prolonged vector fields ${\bf v}_1^{(n)}, \dots,
{\bf v}_6^{(n)}$, with first-order parts truncated\big/dropped:
\[
\begin{blockarray}{ccccccccc}
& \partial_{u_{2,0}} & \partial_{u_{1,1}} & 
\partial_{u_{3,0}} & \partial_{u_{2,1}} & 
\partial_{u_{1,2}} & \partial_{u_{0,3}} & 
\big\{\partial_{u_{k,l}}\big\}_{4\leqslant k+l}
\\
\begin{block}{c(cccccccc)}
{\bf v}_1 & \ast & \ast & \ast & \ast & \ast & \ast & \cdots
\\
{\bf v}_2 & \ast & \ast & \ast & \ast & \ast & \ast & \cdots
\\
{\bf v}_3 & \ast & \ast & \ast & \ast & \ast & \ast & \cdots
\\
{\bf v}_4 & \ast & \ast & \ast & \ast & \ast & \ast & \cdots
\\
{\bf v}_5 & \ast & \ast & \ast & \ast & \ast & \ast & \cdots
\\
{\bf v}_6 & \ast & \ast & \ast & \ast & \ast & \ast & \cdots
\\ 
\end{block}
\end{blockarray}\,\,.
\]

Let us denote the submanifold 
$\big\{ u_{j,k} = R_{j,k} \big\}$
of {\sl parabolic jets} as:
\[
P\!J_{2,1}^n
\,\subset\,
J_{2,1}^n
\eqno
{\scriptstyle{(n\,\geqslant\,0)}},
\]
with:
\[
\dim\,
P\!J_{2,1}^n
\,=\,
3+2\,n
\,<\,
2+{\textstyle{\binom{2+n}{n}}}
\,=\,
\dim\,
J_{2,1}^n,
\eqno
\]
the inequality being strict as soon as $n \geqslant 2$, {\em e.g.}
when $n = 4$:
\[
11
\,<\,
17,
\]
with the strange coincidence that this dimension $11$ of 
$P\!J_{2,1}^4$ is equal to $\dim\, \SA_3(\R)$!
A nice surprise is waiting for us about that!

\begin{center}
\scalebox{1.5}{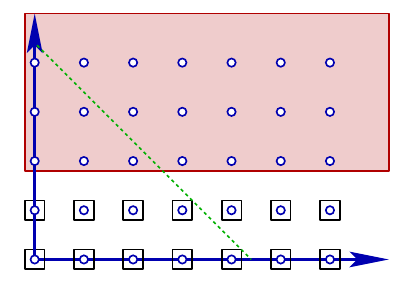}
\end{center}

Also, let us introduce the projection from the full jet space
onto the parabolic jet space:
\[
\aligned
{\sf pj}\,
\Big(
x,y,\,
\big\{
u_{j,k}
\big\}_{j+k\leqslant n}
\Big)
&
\,\,:=\,\,
\Big(
x,y,\,\,
\begin{smallmatrix} 
u_{0,1}, & \dots & \dots & u_{n-1,1}, &
\\
\\ 
u, & u_{1,0}, & \dots, & u_{n-1,0}, & u_{n,0}
\end{smallmatrix}
\Big)
\\
&
\,\,\,\,\in\,\,
\R^{3+2n}.
\endaligned
\]

A general vector field:
\[
{\bf v}
\,=\,
\xi(x,y,u)\,
\frac{\partial}{\partial x}
+
\eta(x,y,u)\,
\frac{\partial}{\partial y}
+
\varphi(x,y,u)\,
\frac{\partial}{\partial u},
\]
prolongs infinitely as:
\[
{\bf v}^\infty
\,=\,
{\bf v}
+
\sum_{1\leqslant j+k}\,
\Phi^{j,k}
\big(
x,y,\,
u^{(j+k)}
\big)\,
\frac{\partial}{\partial u_{j,k}},
\]
with coefficients $\Phi^{j,k}$ given by mean of the two 
(infinite, commuting) total differentiation operators:
\[
\aligned
{\sf D}_x
&
\,:=\,
\frac{\partial}{\partial x}
+
u_{1,0}\,
\frac{\partial}{\partial u}
+
\sum_{1\leqslant j+k}\,
u_{j+1,k}\,
\frac{\partial}{\partial u_{j,k}},
\\
{\sf D}_y
&
\,:=\,
\frac{\partial}{\partial y}
+
u_{0,1}\,
\frac{\partial}{\partial u}
+
\sum_{1\leqslant j+k}\,
u_{j,k+1}\,
\frac{\partial}{\partial u_{j,k}},
\endaligned
\]
by the formulas of 
Theorem~{\ref{Thm-prolongation-vector-fields}}:
\[
\Phi^{j,k}
\,:=\,
{\sf D}_x^j\,
{\sf D}_y^k\,
\Big(
\varphi
-
\xi\,
u_{1,0}
-
\eta\,
u_{0,1}
\Big)
+
\xi\,u_{j+1,k}
+
\eta\,u_{j,k+1},
\]
are functions on $J^{j+k}_{2,1}$.

\begin{Assertion}
For any jet order $n \geqslant 2$, all the $11$ 
prolonged vector field on $J_{2,1}^n$:
\[
{\bf v}_1^{(n)},\ \
{\bf v}_2^{(n)},\ \
{\bf v}_3^{(n)},\ \
{\bf v}_4^{(n)},\ \
{\bf v}_5^{(n)},\ \
{\bf v}_6^{(n)},\ \
\ \ \ \ \ 
{\bf v}_7^{(n)},\ \
{\bf v}_8^{(n)},\ \
\ \ \ \ \ 
{\bf w}_1^{(n)},\ \
{\bf w}_2^{(n)},\ \
{\bf w}_3^{(n)},
\]
are tangent to the (smooth, graphed) 
submanifold $P\!J_{2,1}^n \subset J_{2,1}^n$
of parabolic jets.
\end{Assertion}

\proof
By differentiating one-parameter subgroups, 
this is a consequence of
Theorem~{\ref{Thm-transfer-Hessian-determinants}}.
\endproof

Next, the push-forward to the (horizontal) space $\R^{3+2n}$ 
of parabolic jets of such a prolonged vector field is obtained
by just `killing' all the $\partial \big/ \partial u_{j,k}$ having
$k \geqslant 2$:
\[
{\sf pj}_\ast
\big(
{\bf v}^{(\infty)}
\big)
\,=\,
{\bf v}
+
\sum_{1\leqslant j'}\,
{\sf pj}_\ast
\big(\Phi^{j',0}\big)\,
\frac{\partial}{\partial u_{j',0}}
+
\sum_{0\leqslant j''}\,
{\sf pj}_\ast
\big(\Phi^{j'',1}\big)\,
\frac{\partial}{\partial u_{j'',1}},
\]
and when we restrict ourselves to a finite order $n \geqslant 0$, 
this becomes:
\[
{\sf pj}_\ast
\big(
{\bf v}^{(n)}
\big)
\,=\,
{\bf v}
+
\sum_{1\leqslant j'\leqslant n}\,
{\sf pj}_\ast
\big(\Phi^{j',0}\big)\,
\frac{\partial}{\partial u_{j',0}}
+
\sum_{0\leqslant j''\leqslant n-1}\,
{\sf pj}_\ast
\big(\Phi^{j'',1}\big)\,
\frac{\partial}{\partial u_{j'',1}},
\]
where ${\sf pj}_\ast\big(\Phi^{j'}\big)$ and ${\sf pj}_\ast\big(\Phi^{j'',1}\big)$ are push-forwards of $\Phi^{j',0}$ and $\Phi^{j'',1}$ to $P\!J^{j'}_{2,1}$ and $P\!J^{j''+1}_{2,1}$.

To summarize, the push-forward process consists of

\smallskip\noindent$\bullet$\,
dropping all $\partial \big/ \partial u_{j,k}$ with $k \geqslant 2$;

\smallskip\noindent$\bullet$\,
in the kept coefficients:
\[
\Phi^{j',0}
\,=\,
\Phi^{j',0}
\big(
x,y,\,
u^{(j')}
\big),
\ \ \ \ \ \ \ \ \ \ \ \ \ \ \ \ \ \ \ \
\Phi^{j'',1}
\,=\,
\Phi^{j'',1}
\big(
x,y,\,
u^{(j''+1)}
\big),
\]
replacing all jet coordinates $u_{j''', k'''}$ with $k''' \geqslant 2$ 
by $R_{j''', k'''}$.

For confirmation, for confidence, for coherence, 
the following observation may (for fun) 
be checked on whatever computer 
for many low jet orders $n \geqslant 2$.

\begin{Observation}
{\bf [Exercise]}
In order 
to push forward the $11$ generators of the action of $\SA_3(\R)$
via ${\sf pj}_\ast$,
it suffices to know parabolic jets relations only for
$k = 2$, and for any $j \geqslant 0$:
\[
u_{j,2}
\,=\,
R_{j,2}
\Big(
\big\{
u_{j',0}
\big\}_{2\leqslant j'\leqslant j+2},\,\,
\big\{
u_{j'',1}
\big\}_{1\leqslant j''\leqslant j+1}
\Big).
\eqno\qed
\]
\end{Observation}

\begin{Observation}
{\bf [Exercise]}
The process of replacing $u_{j,2}$ by $R_{j,2}$ 
occurs only when one calculates ${\sf pj}_\ast\big({\bf v}_5^{(n)}\big)$ and ${\sf pj}_\ast\big({\bf v}_6^{(n)}\big)$ for $n\geqslant 2$.
\end{Observation}

Now, although 
Theorem~{\ref{Thm-transfer-Hessian-determinants}} already gave a
full explanation, we want to verify in another way that
the hypothesis that the Hessian is identically zero 
is stable under (special) affine transformations.
So we fix the jet order $n = 2$, whence parabolic
jets are defined by a single equation:
\[
P\!J_{2,1}^2
\,=\,
\big\{
u_{2,0}\,
u_{0,2}
-
u_{1,1}^2
\big\},
\]
and equivalently, because we are working on the space $\big\{ u_{2, 0}
\neq 0 \big\}$, we can solve:
\[
u_{0,2}
\,=\,
\frac{u_{1,1}^2}{u_{2,0}}.
\]
Clearly, ${\bf v}_7^{(2)}$, ${\bf v}_8^{(2)}$,
${\bf w}_1^{(2)}$, ${\bf w}_2^{(2)}$, ${\bf w}_3^{(3)}$
are (trivially) tangent to this submanifold $P\!J_{2,1}^2 \subset 
J_{2,1}^2$, since they incorporate only derivatives $\partial_x$,
$\partial_y$, $\partial_u$, $\partial_{u_{1,0}}$, 
$\partial_{u_{0,1}}$, nothing of higher order $\partial_{u_{j,k}}$
with $j + k \geqslant 2$.

Next, we must verify that the $6$ remaining vector
fields ${\bf v}_1^{(2)}$, ${\bf v}_2^{(2)}$, ${\bf v}_3^{(2)}$,
${\bf v}_4^{(2)}$, ${\bf v}_5^{(2)}$, ${\bf v}_6^{(2)}$ are
also tangent to $P\!J_{2,1}^n$.
We write their coefficients matrix as:
\[
\begin{blockarray}{cccc}
& \partial_{u_{2,0}} & \partial_{u_{1,1}} & 
\partial_{u_{0,2}}
\\
\begin{block}{c(ccc)}
{\bf v}_1 & -3u_{2,0} & -2u_{1,1} & -u_{0,2}
\\
{\bf v}_2 & -u_{2,0} & -2u_{1,1} & -3u_{0,2}
\\
{\bf v}_3 & 0 & -u_{2,0} & -2u_{1,1}
\\
{\bf v}_4 & -3u_{2,0}u_{1,0} & -2u_{1,1}u_{1,0}-u_{2,0}u_{0,1} & 
-2u_{1,1}u_{0,1}-u_{0,2}u_{1,0}
\\
{\bf v}_5 & -2u_{1,1} & -u_{0,2} & 0
\\
{\bf v}_6 & -2u_{1,1}u_{1,0}-u_{2,0}u_{0,1} & 
-2u_{1,1}u_{0,1}-u_{0,2}u_{1,0} & -3u_{0,2}u_{0,1}
\\ 
\end{block}
\end{blockarray}\,\,.
\]
Also, let us abbreviate:
\[
\Haux_u
\,:=\,
u_{2,0}\,
u_{0,2}
-
u_{1,1}^2.
\]
The tangency in question is verified on a computer,
or else with a pen, as follows.

\begin{Assertion}
The quotients:
\[
{\bf v}_\sigma
\big(
\Haux_u
\big)
\Big/
\Haux_u
\]
have values, for $\sigma = 1, 2, 3, 4, 5, 6$:
\[
-\,4,
\ \ \ \ \ \ \ \ \ \
-\,4,
\ \ \ \ \ \ \ \ \ \
0,
\ \ \ \ \ \ \ \ \ \
-\,4\,u_{1,0},
\ \ \ \ \ \ \ \ \ \
0,
\ \ \ \ \ \ \ \ \ \
-\,4\,u_{0,1}.
\eqno\qed
\]
\end{Assertion}

Thus, $\Haux_u = 0$ implies ${\bf v}_\sigma \big( \Haux_u \big) = 0$,
which is tangency!
One also observes that there is a $2 \times 2$ minor
of the above matrix which is nowhere vanishing:
\[
\left\vert\!
\begin{array}{cc}
-u_{2,0} & -2u_{1,1}
\\
0 & -u_{2,0}
\end{array}
\!\right\vert
\,=\,
u_{2,0}^2.
\]

\begin{Question}
{\sl Are there differential invariants of order $2$?}
\end{Question}

{\em No!} Because at every point of $\big\{ u_{2,0} \neq 0 = 
\Haux_u \big\}$, we realize that:
\[
\aligned
\dim\,
P\!J_{2,1}^2
\,=\,
7
\,\,=\,\,
\rank\,\Span\,
\Big(
&
{\sf pj}_\ast\big({\bf v}_2^{(2)}\big),\ \ 
{\sf pj}_\ast\big({\bf v}_3^{(2)}\big),\ \ 
{\sf pj}_\ast\big({\bf v}_7^{(2)}\big),\ \ 
{\sf pj}_\ast\big({\bf v}_8^{(2)}\big),
\\
&
{\sf pj}_\ast\big({\bf w}_1^{(2)}\big),\ \ 
{\sf pj}_\ast\big({\bf w}_2^{(2)}\big),\ \ 
{\sf pj}_\ast\big({\bf w}_3^{(2)}\big)
\Big),
\endaligned
\]
by looking at the triangular-by-blocks
$7 = 2 + 2 + 3$ matrix:
\[
\begin{blockarray}{ccccccccc}
&
\partial_x & \partial_y & \partial_u &
\partial_{u_{1,0}} & \partial_{u_{0,1}} &
\partial_{u_{2,0}} & \partial_{u_{1,1}}
\\
\begin{block}{c(cccccccc)}
{\bf v}_2 & 0 & y & -u & -u_{1,0} & -2u_{0,1} & -u_{2,0} & -2u_{1,1}
\\
{\bf v}_3 & y & 0 & 0 & 0 & -u_{1,0} & 0 & -u_{2,0}
\\
{\bf v}_7 & 0 & 0 & x & 1 & 0 & 0 & 0
\\
{\bf v}_8 & 0 & 0 & y & 0 & 1 & 0 & 0
\\
{\bf w}_1 & 1 & 0 & 0 & 0 & 0 & 0 & 0
\\
{\bf w}_2 & 0 & 1 & 0 & 0 & 0 & 0 & 0
\\
{\bf w}_3 & 0 & 0 & 1 & 0 & 0 & 0 & 0
\\ 
\end{block}
\end{blockarray}\,\,,
\]
whose determinant equals $u_{2,0}^2 \neq 0$.
Consequently, the set we are working on:
\[
\big\{
u_{2,0}
\,\neq\,
0
\,=\,
u_{2,0}u_{0,2}
-
u_{1,1}^2
\big\}
\]
contains only 7-dimension orbits, in fact only a unique orbit, for the prolonged action of $\SA_3(\R)$.
So any transversal to the orbits is zero-dimensional,
and there are {\em no} differential invariants of order $2$.
In particular, the Hessian $\Haux_u$ is {\em not} 
a differential invariant. Taking the notations of
Section~{\ref{moving-frame-method}}, we introduce a

\begin{Terminology}
For the action of a local Lie group $G$ on graphs 
$\{ u = {\tt u}(x) \}$ in $\R_z^{p+q}$,
a function $\Paux = \Paux \big( z^{(n)} \big)$
on the $n$\textsuperscript{th} order jet space
$J_z^n$ is called a {\sl relative invariant} if:
\[
\Paux
\Big(
w^{(n)}
\big(
g,
z^{(n)}
\big)
\Big)
\,=\,
\nonzero
\cdot
\Paux
\big(
z^{(n)}
\big)
\eqno
{\scriptstyle{(\forall\,g\,\in\,G)}}.
\]
\end{Terminology}

Immediately, the zero-set of a relative invariant {\em is}
invariant:
\[
\Big\{
\Paux
\big(
w^{(n)}
\big(g,z^{(n)}\big)
\big)
=
0
\Big\}
\,\,=\,\,
\big\{
\Paux(z^{(n)})
=
0
\big\}.
\]
Precisely what we already observed of the Hessian!

\smallskip

Next, we examine what occurs in the $3$\textsuperscript{rd} order
jet space. Recall that we are working only in the domain:
\[
\big\{
u_{2,0}
\,\neq\,
0
\,=\,
u_{2,0}\,u_{0,2}
-
u_{1,1}^2
\big\}.
\]

The submanifold
$P\!J_{2,1}^3 \subset J_{2,1}^3$ 
of parabolic $3$\textsuperscript{rd} order jets is now
defined by three equations (the last one being in fact not used):
\[
\aligned
u_{0,2}
\,:=\,
\frac{u_{1,1}^2}{u_{2,0}}
\ \ \ \ \ \ \ \ \ \ \ \ \ \ \ \ \ \ \ \
\ \ \ \ \ \ \ \ \ \ \ \ \ \ \ \ \ \ \ \
u_{1,2}
&
\,=\,
2\,\frac{u_{1,1}\,u_{2,1}}{u_{2,0}}
-
\frac{u_{1,1}^2\,u_{3,0}}{u_{2,0}^2},
\\
u_{0,3}
&
\,=\,
3\,
\frac{u_{1,1}^2\,u_{2,1}}{u_{2,0}^2}
-
2\,\frac{u_{1,1}^3\,u_{3,0}}{u_{2,0}^3}.
\endaligned
\]

Beyond 
Theorem~{\ref{Thm-transfer-Hessian-determinants}}, 
again with affine transformations not far from the identity, 
we have the

\begin{Proposition}
{\rm {\cite[5.8]{Merker-2019}}}
With this hypothesis, in terms of the nowhere vanishing quantity:
\[
\Upsilon
\,:=\,
\big(
{\sf l}
+
{\sf m}\,F_y
\big)\,
F_{xx}
-
\big(
{\sf k}
+
{\sf m}\,F_x
\big)\,
F_{xy},
\]
it holds:
\[
\frac{G_{ss}\,G_{sst}-G_{st}\,G_{sss}}{G_{ss}^2}
\,\,=\,\,
\frac{F_{xx}}{\Upsilon}\,
\bigg(
\frac{F_{xx}\,F_{xxy}
-
F_{xy}\,F_{xxx}}{
F_{xx}^2}
\bigg).
\eqno\qed
\]
\end{Proposition}

This implies that, in the submanifold $P\!J_{2,1}^3$ of parabolic
jets of order $3$, the zero-set:
\[
\big\{
u_{2,0}\,u_{2,1}
-
u_{1,1}\,u_{3,0}
=
0
\big\},
\]
is invariant under the prolongation of the $\SA_3$-action,
namely, for every $g \in \SA_3(\R)$ not far from the identity:
\[
g^{(3)}
\Big(
\big\{
u_{2,0}\,u_{2,1}
-
u_{1,1}\,u_{3,0}
=
0
\big\}
\Big)
\,\,\subset\,\,
\big\{
v_{2,0}\,v_{2,1}
-
v_{1,1}\,v_{3,0}
=
0
\big\}.
\]

We now want to see the same property from the vector fields
point of view. Let us abbreviate:
\leqnomode\usetagform{default}
\begin{align}
\label{invariant-S-u}
\Saux_u
\,:=\,
\frac{
u_{2,0}\,u_{2,1}
-
u_{1,1}\,u_{3,0}}{
u_{2,0}^2}.
\end{align}
The tangencies to $\big\{ \Saux_u = 0 \big\}$ 
follow from an

\begin{Assertion}
The quotients:
\[
{\sf pj}_\ast\big({\bf v}_\sigma\big)
\big(
\Saux_u
\big)
\Big/
\Saux_u
\]
have values, for $\sigma = 1, 2, 3, 4, 5, 6$:
\[
0,
\ \ \ \ \ \ \ \ \ \
-\,1,
\ \ \ \ \ \ \ \ \ \
0,
\ \ \ \ \ \ \ \ \ \
0,
\ \ \ \ \ \ \ \ \ \
\frac{u_{1,1}}{u_{2,0}},
\ \ \ \ \ \ \ \ \ \
\frac{u_{1,1}\,u_{1,0}-u_{2,0}\,u_{0,1}}{u_{2,0}}.
\eqno\qed
\]
\end{Assertion}

\begin{Question}
{\sl Are there differential invariants of order $3$?}
\end{Question}

Again: {\em No!} Because at every point of $\big\{ u_{2,0} \neq 0 = 
\Haux_u \big\}$, we realize that:
\[
\aligned
\dim\,
P\!J_{2,1}^3
\,=\,
9
\,\,=\,\,
\rank\,\Span\,
\Big(
&
{\sf pj}_\ast\big({\bf v}_1^{(3)}\big),\ \ 
{\sf pj}_\ast\big({\bf v}_2^{(3)}\big),\ \ 
{\sf pj}_\ast\big({\bf v}_3^{(3)}\big),\ \ 
{\sf pj}_\ast\big({\bf v}_4^{(3)}\big),
\\
&
{\sf pj}_\ast\big({\bf v}_7^{(3)}\big),\ \ 
{\sf pj}_\ast\big({\bf v}_8^{(3)}\big),\ \ 
{\sf pj}_\ast\big({\bf w}_1^{(3)}\big),\ \ 
{\sf pj}_\ast\big({\bf w}_2^{(3)}\big),\ \ 
{\sf pj}_\ast\big({\bf w}_3^{(3)}\big)
\Big),
\endaligned
\]
by looking at the triangular-by-blocks
$9 = 4 + 2 + 3$ matrix:
\[
\footnotesize
\aligned
\begin{blockarray}{ccccccccccc}
&
\partial_x & \partial_y & \partial_u &
\partial_{u_{1,0}} & \partial_{u_{0,1}} &
\partial_{u_{2,0}} & \partial_{u_{1,1}} &
\partial_{u_{3,0}} & \partial_{u_{2,1}}
\\
\begin{block}{c(cccccccccc)}
{\bf v}_1 & x & 0 & -u & -2u_{1,0} & -u_{0,1} & -3u_{2,0} & -2u_{1,1}
& -4u_{3,0} & -3u_{2,1}
\\
{\bf v}_2 & 0 & y & -u & -u_{1,0} & -2u_{0,1} & -u_{2,0} & -2u_{1,1}
& -u_{3,0} & -2u_{2,1}
\\
{\bf v}_3 & y & 0 & 0 & 0 & -u_{1,0} & 0 & -u_{2,0} 
& 0 & -u_{3,0}
\\
{\bf v}_4 & 0 & 0 & u & -u_{1,0}^2 & -u_{1,0}u_{0,1} &
-3u_{2,0}u_{1,0} & 
\begin{smallmatrix} 
-2u_{1,1}u_{1,0} \\
-u_{2,0}u_{0,1} 
\end{smallmatrix}
&
\begin{smallmatrix} 
-4u_{3,0}u_{1,0} \\
-3u_{2,0}^2 
\end{smallmatrix} 
&
\begin{smallmatrix} 
-3u_{2,1}u_{1,0} \\
-3u_{1,1}u_{2,0}-u_{3,0}u_{0,1} 
\end{smallmatrix}
\\
{\bf v}_7 & 0 & 0 & x & 1 & 0 & 0 & 0 & 0 & 0
\\
{\bf v}_8 & 0 & 0 & y & 0 & 1 & 0 & 0 & 0 & 0
\\
{\bf w}_1 & 1 & 0 & 0 & 0 & 0 & 0 & 0 & 0 & 0
\\
{\bf w}_2 & 0 & 1 & 0 & 0 & 0 & 0 & 0 & 0 & 0
\\
{\bf w}_3 & 0 & 0 & 1 & 0 & 0 & 0 & 0 & 0 & 0
\\ 
\end{block}
\end{blockarray}\,\,,
\endaligned
\]
and by realizing that the upper-right $4 \times 4$ block has
a not identically zero determinant:
\[
\left\vert\!
\begin{array}{cccc}
-3u_{2,0} & -2u_{1,1} & -4u_{3,0} & -3u_{2,1}
\\
-u_{2,0} & -2u_{1,1} & -u_{3,0} & -2u_{2,1}
\\
0 & -u_{2,0} & 0 & -u_{3,0}
\\
-3u_{2,0}u_{1,0} & 
\begin{smallmatrix} 
-2u_{1,1}u_{1,0} \\
-u_{2,0}u_{0,1} 
\end{smallmatrix}
&
\begin{smallmatrix} 
-4u_{3,0}u_{1,0} \\
-3u_{2,0}^2 
\end{smallmatrix} 
&
\begin{smallmatrix} 
-3u_{2,1}u_{1,0} \\
-3u_{1,1}u_{2,0}-u_{3,0}u_{0,1} 
\end{smallmatrix}
\end{array}
\!\right\vert
\,\,=\,\,
9\,
u_{2,0}^3\,
\big(
\underbrace{
u_{2,0}\,u_{2,1}
-
u_{1,1}\,u_{3,0}}_{
{\sf recognize}\,\,
\Saux_u}
\big).
\]

Our domain $\big\{ u_{2,0} \neq 0 = u_{2,0}u_{0,2} - 
u_{1,1}^2 \big\}$ stratifies as:
\[
\big\{
\Saux_u
\equiv
0
\big\}
\medcup
\big\{
\Saux_u
\neq
0
\big\}.
\]
So on the dense open subset where $\Saux_u \neq 0$, the above
$9$ vector fields span the tangent space to the $9$-dimensional
parabolic jet space $P\!J_{2,1}^3$, whence any transversal 
to the orbits of the (prolonged) $\SA_3$-action is zero-dimensional,
and there are {\em no} differential invariants at a generic point.

Be careful! The closed subset $\big\{ \Saux_u = 0 \big\}$ 
is $\SA_2$-invariant. In it, the rank in question degenerates,
since the above determinant has value $0$. But since
$u_{2,0} \neq 0$, the equation $\Saux_u = 0$ can be solved for:
\[
u_{2,1}
\,=\,
\frac{u_{1,1}}{u_{2,0}}\,
u_{3,0}.
\]
According to Lie's general principle of thought, 
in the study of graphs $\big\{ u = F(x,y) \big\}$,
there is bifurcation branching:
\[
\xymatrix{
& & & &
\text{\footnotesize\sf Identical degeneracy}\,\,
\Saux_F\,\equiv\,0, 
\\
\Saux_F
\ar[urrr]
\ar[drrr]
& & & &
\\
& & & &
\text{\footnotesize\sf Nowhere vanishing}\,\,
\Saux_F\,\neq\,0,
}
\]
namely one studies either graphs for which $F_{xxy} \equiv 
\frac{F_{xy}}{F_{xx}}\, F_{xxx}$, or graphes for which
$F_{xx}\, F_{xxy} - F_{xy}\, F_{xxx} \neq 0$ at {\em every}
point $(x,y)$. In other words, we disregard mixed types
for which $\Saux_F(x,y) \not\equiv 0$,
while $\big\{ \Saux_F(x,y) = 0 \big\} \neq\emptyset$.

Our main study, in the next 
Sections~{\ref{relative-invariant-S-first-invariant-W}},
{\ref{branch-W-equiv-0-branch-W-nonzero}},
{\ref{recurrence-relations-parabolic-surfaces}},
will concern the branch $\Saux_F \neq 0$. For now, let us summarize
how the branch $\Saux_F \equiv 0$ can be easily terminated.

So our assumptions are:
\leqnomode\usetagform{default}
\begin{align}
\label{asssumptions-branch-H-0-S-0}
F_{xx}
\,\neq\,
0,
\ \ \ \ \ \ \ \ \ \ \ \ \ \ \ \ \ \ \ \
F_{yy}
\,\equiv\,
\frac{F_{xy}^2}{F_{xx}},
\ \ \ \ \ \ \ \ \ \ \ \ \ \ \ \ \ \ \ \
F_{xxy}
\,\equiv\,
\frac{F_{xy}}{F_{xx}}\,
F_{xxx}.
\end{align}
This leads to a {\em new} (smooth) submanifold of the 
$n$\textsuperscript{th} order jet space:
\[
C\!P\!J_{2,1}^n
\,\,\subset\,\,
P\!J_{2,1}^n
\,\,\subset\,\,
J_{2,1}^n,
\]
defined by the equations:
\[
u_{2,0}
\,\neq\,
0,
\ \ \ \ \ \ \ \ \ \ \ \ \ \ \ \ \ \ \ \
u_{0,2}
\,=\,
\frac{u_{1,1}^2}{u_{2,0}},
\ \ \ \ \ \ \ \ \ \ \ \ \ \ \ \ \ \ \ \
u_{2,1}
\,=\,
\frac{u_{1,1}}{u_{2,0}}\,
u_{3,0},
\]
together with their differential consequences (exercise).

\begin{center}
\scalebox{1.5}{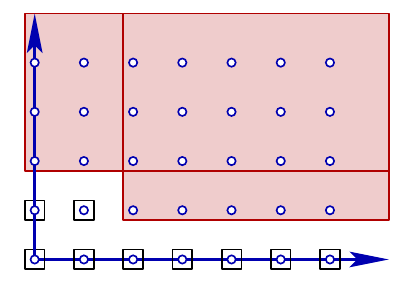}
\end{center}

When $n \geqslant 2$, this submanifold is equipped
with horizontal coordinates:
\[
C\!P\!J_{2,1}^n
\,\,=\,\,
\R^{5+n}
\,\,\ni\,\,
\Big(
x,y,\,\,
\begin{smallmatrix} 
u_{0,1}, & u_{1,1},
\\
\\ 
u, & u_{1,0}, & u_{2,0}, & u_{3,0}, & u_{4,0}, & u_{5,0},
\,\,\dots, & u_{n,0}
\end{smallmatrix}
\Big).
\]

After an elementary affine transformation, we can assume
that a graphed surface satisfying these assumptions
starts at order $2$ terms:
\[
u
\,=\,
F_{2,0}\,
{\textstyle{\frac{x^2}{2!}}}
+
F_{1,1}\,
xy
+
F_{0,2}\,
{\textstyle{\frac{y^2}{2!}}}
+
{\rm O}_{x,y}(3).
\]
Since $F_{0,2} = F_{1,1}^2 \big/ F_{2,0}$, we can write:
\[
u
\,=\,
\frac{F_{2,0}}{2}\,
\Big(
x
+
\frac{F_{1,1}}{F_{2,0}}\,
y
\Big)^2
+
{\rm O}_{x,y}(3),
\]
and making the special affine change of coordinates:
\[
x'
\,:=\,
x
+
{\textstyle{\frac{F_{1,1}}{F_{2,0}}}}\,
y,
\ \ \ \ \ \ \ \ \ \ \ \ \ \ \ \ \ \ \ \
y'
\,:=\,
F_{2,0}\,
y,
\ \ \ \ \ \ \ \ \ \ \ \ \ \ \ \ \ \ \ \
u'
\,:=\,
{\textstyle{\frac{1}{F_{2,0}}}}\,
u,
\]
we come, dropping the primes, to a normalized form:
\[
u
\,=\,
{\textstyle{\frac{x^2}{2!}}}
+
{\rm O}_{x,y}(3).
\]

Importantly, our assumptions~({\ref{asssumptions-branch-H-0-S-0}})
are (special) affinely invariants, and then, in the new system of
coordinates, they continue to hold. 
As in the proof of Lemma~{\ref{LM-parabola-Halphen-R}},
let us write the two {\sc pde}s in question as:
\[
F_{yy}
\,=\,
\mathcal{R}\,
F_{xy},
\ \ \ \ \ \ \ \ \ \ \ \ \ \ \ \ \ \ \ \
F_{xxy}
\,=\,
\mathcal{R}\,
F_{xy},
\]
where $\mathcal{R} = \mathcal{R}(x,y)$ denotes an unspecified 
function.

\begin{Assertion}
In the normalized form, the remainder ${\rm O}_{x,y}(3) = {\rm O}_x(3)$ depends only 
on $x$, not on $y$.
\end{Assertion}

\proof
At first, we see $F_{xxy}(0) = 0$ since $F_{xy}(0)=0$. Differentiating with respect to
$x$, we get:
\[
F_{xxxy}
\,=\,
\mathcal{R}_x\,
F_{xy}
+
\mathcal{R}\,
F_{xxy}
\,=\,
\mathcal{R}\,
F_{xy},
\]
whence similarly $F_{xxxy}(0) = 0$. An easy induction yields
$F_{x^k y} = \mathcal{R}\, F_{xy}$ for every
$k \geqslant 2$, whence $F_{x^ky}(0) = 0$.

Next, we see $F_{yy}(0) = 0$. Differentiating and inducting,
we get $F_{x^kyy} = \mathcal{R}\, F_{xy}$ for every $k \geqslant 0$,
whence $F_{x^kyy}(0) = 0$. By another induction, $F_{x^k y^l} =
\mathcal{R}\, F_{xy}$ for all $k \geqslant 0$, $l \geqslant 2$, whence
$F_{x^k y^l} (0) = 0$, and this concludes 
because we assume analyticity.
\endproof

So surfaces $S^2 \subset \R^3$ satisfying
assumptions~({\ref{asssumptions-branch-H-0-S-0}}) are of the form:
\[
u
\,=\,
f(x)
\,=\,
{\textstyle{\frac{x^2}{2!}}}
+
{\rm O}_x(3),
\]
they are 
{\em products} of curves in $\R_{x,u}^2$ and the straight $\R_y^1$,
which is a degenerate situation.

From any special affine equivalence in $\SA_3(\R)$ we can extract
a $2$-dimensional equivalence:
\[
\!\!\!\!\!\!\!\!\!\!\!\!\!\!\!
\aligned
s
&
\,=\,
{\sf a}\,x
+
{\sf b}\,y
+
{\sf c}\,u
+
{\sf d},
\\
t
&
\,=\,
{\sf k}\,x
+
{\sf l}\,y
+
{\sf m}\,u
+
{\sf n},
\\
v
&
\,=\,
{\sf p}\,x
+
{\sf q}\,y
+
{\sf r}\,u
+
{\sf s}
\endaligned
\ \ \ \ \ \ \
\xrightarrow[{\rule[0pt]{50pt}{0pt}}]{{\sf specialize}}
\ \ \ \ \ \ \
\aligned
s
&
\,=\,
{\sf a}\,x
\ \ \ \ \ \ \ \ \
+
{\sf c}\,u
+
{\sf d},
\\
t
&
\,=\,
\ \ \ \ \ \ \ \ 
{\sf l}\,y
\ \ \ \ \ \ \ \ \ \ \ \,
+
{\sf n},
\\
v
&
\,=\,
{\sf p}\,x
\ \ \ \ \ \ \ \ \
+
{\sf r}\,u
+
{\sf s}
\endaligned
\ \ \ \ \ \ \
\xrightarrow[{\rule[0pt]{50pt}{0pt}}]{{\sf extract}}
\ \ \ \ \ \ \
\aligned
s
&
\,=\,
{\sf a}\,x
+
{\sf c}\,u
+
{\sf d},
\\
v
&
\,=\,
{\sf p}\,x
+
{\sf r}\,u
+
{\sf s},
\endaligned
\]
but then, since the transformation in the $y$-space has zero 
influence, we can adjust the central parameter:
\[
{\sf l}
\,:=\,
\frac{1}{{\sf a}\,{\sf r}-{\sf c}\,{\sf p}},
\]
in order to guarantee that the determinant has value $1$,
while the determinant of the $2$-dimensional transformation:
\[
\left\vert\!
\begin{array}{cc}
{\sf a} & {\sf c}
\\
{\sf p} & {\sf r}
\end{array}
\!\right\vert
\,\neq\,
0
\]
must only be nonzero, not necessarily equal to $1$. In other words,
our curves $\big\{ u = f(x) \big\}$ with $f_{xx} \neq 0$ are
considered modulo {\em affine} transformations of the plane
$\R_{x,u}^2$, not special affine ones.

We have already studied plane curves modulo
${\sf A}_2(\R)$ in 
Section~{\ref{affine-invariants-curves-R-2}}.
Hence we only need to adapt (slightly) 
these results in our present context.

\medskip\noindent$\bullet$\,
There is a relative invariant:
\[
\Paux
\,:=\,
-\,5\,
u_{3,0}^2
+
3\,u_{2,0}\,u_{4,0},
\]
which creates a branching. 

\medskip\noindent$\bullet$\,
One has $\Paux \equiv 0$ if and
only if the surface $\big\{ u = F(x,y) \big\}$ is 
special affinely equivalent to the product 
$\{ v = s^2 \}$ of a parabola
with $\R_t^1$. 

\medskip\noindent$\bullet$\,
When $\Paux \neq 0$, as in 
Theorem~{\ref{Thm-GL-2-coefficients-invariants}},
the surface is $\SA_3$-equivalent to:
\[
u
\,=\,
{\textstyle{\frac{x^2}{2!}}}
+
0
\pm
{\textstyle{\frac{x^4}{4!}}}
+
F_5\,
{\textstyle{\frac{x^5}{5!}}}
+
\sum_{i\geqslant 6}\,
F_i\,
{\textstyle{\frac{x^i}{i!}}},
\]
where $F_5, F_6, F_7, \dots$ are diffferential invariants, 
for instance in the case of $+\, \frac{x^4}{4!}$:
\[
F_5
\,=\,
\Iaux_5
\,:=\,
\frac{1}{\sqrt{3}}\,
\frac{9\,u_{2,0}^2\,u_{5,0}
-
45\,u_{2,0}\,u_{3,0}\,u_{4,0}
+
40\,u_{3,0}^3}{
\big(
3\,u_{2,0}\,u_{4,0}
-
5\,u_{3,0}^3
\big)^{3/2}}.
\]
All the other differential invariants $\Iaux_6, \Iaux_7, \dots$
express in terms of $\Iaux_5$ and its invariant
derivatives $\mathcal{D}_x^\nu (\Iaux_5)$.

\medskip\noindent$\bullet$\,
The degenerate case where $\Iaux_5 \equiv 0$ means,
according to Lemma~{\ref{LM-vanishing-H-M}}, 
that there exist constants such that:
\[
u
\,=\,
{\sf d}\,x
+
{\sf e}
+
\sqrt{{\sf f}\,x^2+2{\sf g}\,x+{\sf h}}.
\]
There are unique such constants so that the right-hand side 
has fourth-order terms normalized as above:
\[
\begin{aligned}
u
&
\,=\,
3
-
\sqrt{9-3\,x^2},
&
\ \ \ \ \ \ \ \ \ \ \ \ \ \ \ \ \ \ \ \ \ \ \ \ \
u
&
\,=\,
-\,3
+
\sqrt{9+3\,x^2}
\\
&
\,=\,
{\textstyle{\frac{x^2}{2!}}}
+
0
+
{\textstyle{\frac{x^4}{4!}}}
+
0
+
{\rm O}_x(6),
&
\ \ \ \ \ \ \ \ \
&
\,=\,
{\textstyle{\frac{x^2}{2!}}}
+
0
-
{\textstyle{\frac{x^4}{4!}}}
+
0
+
{\rm O}_x(6).
\end{aligned}
\]

All this study justifies that we can from now on and until the
end of the paper assume that at every point $(x,y)$:
\[
0
\,\neq\,
\Saux_F
\,=\,
F_{xx}\,F_{xxy}
-
F_{xy}\,F_{xxx}.
\]

\begin{Question}
{\sl 
Are there differential invariants of order $4$?}
\end{Question}

{\em Yes}, at last!
Remind that we already pointed out the coincidence:
\[
\dim\,
P\!J_{2,1}^4
\,\,=\,\,
11
\,\,=\,\,
\dim\,\SA_3(\R).
\]
One could expect that the $11$ vector fields prolonged
to $J_{2,1}^4$ and tangent to 
the submanifold $P\!J_{2,1}^4$ of parabolic jets:
\[
\aligned
&
{\sf pj}_\ast\big({\bf v}_1^{(4)}\big),\ \ 
{\sf pj}_\ast\big({\bf v}_2^{(4)}\big),\ \ 
{\sf pj}_\ast\big({\bf v}_3^{(4)}\big),\ \ 
{\sf pj}_\ast\big({\bf v}_4^{(4)}\big),\ \ 
{\sf pj}_\ast\big({\bf v}_5^{(4)}\big),\ \ 
{\sf pj}_\ast\big({\bf v}_6^{(4)}\big),\ \ 
\\
&
{\sf pj}_\ast\big({\bf v}_7^{(4)}\big),\ \ 
{\sf pj}_\ast\big({\bf v}_8^{(4)}\big),\ \ 
{\sf pj}_\ast\big({\bf w}_1^{(4)}\big),\ \ 
{\sf pj}_\ast\big({\bf w}_2^{(4)}\big),\ \ 
{\sf pj}_\ast\big({\bf w}_3^{(4)}\big),
\endaligned
\]
are linearly independent at a generic point, hence span
the tangent space to $P\!J_{2,1}^4$ almost everywhere,
but this is not the case.

Indeed,  
thanks to our earlier observation~({\ref{triangular-by-blocks-5x5}})
that the last $5$ vector fields have identically zero 
coefficients in front of all the $\partial \big/ \partial u_{j,k}$ 
with $j + k \geqslant 2$,
it suffices to examine the rank of the first $6$ vector fields
in the space of jets of order $\geqslant 2$, namely,
it suffices to determine the rank at a generic point of the
following matrix of coefficients:
\[
\footnotesize
\aligned
\begin{blockarray}{ccccccc}
&
\partial_{u_{2,0}} & \partial_{u_{1,1}} &
\partial_{u_{3,0}} & \partial_{u_{2,1}} &
\partial_{u_{4,0}} & \partial_{u_{3,1}}
\\
\begin{block}{c(cccccc)}
{\bf v}_1 & -3u_{2,0} & -2u_{1,1}
& -4u_{3,0} & -3u_{2,1} & -5u_{4,0} & -4u_{3,1}
\\
{\bf v}_2 & -u_{2,0} & -2u_{1,1}
& -u_{3,0} & -2u_{2,1} & -u_{4,0} & -2u_{3,1}
\\
{\bf v}_3 & 0 & -u_{2,0} 
& 0 & -u_{3,0} & 0 & -u_{4,0}
\\
{\bf v}_4 & 
-3u_{2,0}u_{1,0} & 
\begin{smallmatrix} 
-2u_{1,1}u_{1,0} \\
-u_{2,0}u_{0,1} 
\end{smallmatrix}
&
\begin{smallmatrix} 
-4u_{3,0}u_{1,0} \\
-3u_{2,0}^2 
\end{smallmatrix} 
&
\begin{smallmatrix} 
-3u_{2,1}u_{1,0} \\
-3u_{1,1}u_{2,0}-u_{3,0}u_{0,1} 
\end{smallmatrix}
&
\begin{smallmatrix} 
-5u_{4,0}u_{1,0} \\
-10u_{2,0}u_{3,0} 
\end{smallmatrix}
&
\begin{smallmatrix} 
-4u_{3,1}u_{1,0}-4u_{1,1}u_{3,0} \\
-6u_{2,0}u_{2,1}-u_{4,0}u_{0,1} 
\end{smallmatrix}
\\
{\bf v}_5 & -2u_{1,1} & \Phi_5^{1,1}
& -3u_{2,1} & \Phi_5^{2,1} & -4u_{3,1} & \Phi_5^{3,1}
\\
{\bf v}_6 & 
\begin{smallmatrix}
-2u_{1,1}u_{1,0} \\
-u_{2,0}u_{0,1}  
\end{smallmatrix}
& \Phi_6^{1,1}
& 
\begin{smallmatrix}
-3u_{2,1}u_{1,0} \\
-3u_{1,1}u_{2,0}-u_{3,0}u_{0,1}
\end{smallmatrix}
& \Phi_6^{2,1} 
& 
\begin{smallmatrix}
-4u_{3,1}u_{1,0}-4u_{1,1}u_{3,0} \\
-6u_{2,0}u_{2,1}-u_{4,0}u_{0,1}  
\end{smallmatrix}
& \Phi_6^{3,1}
\\ 
\end{block}
\end{blockarray}\,\,,
\endaligned
\]
in which precisely $6$ entries 
require take account of the parabolic jet
relations shown at the end of 
Section~{\ref{parabolic-jet-relations}}:
\[
\aligned
\Phi_5^{1,1}
&
\,=\,
-\,
\frac{u_{1,1}^2}{u_{2,0}},
\\
\Phi_5^{2,1}
&
\,=\,
-\,4\,
\frac{u_{1,1}\,u_{2,1}}{u_{2,0}}
+
2\,
\frac{u_{1,1}^2\,u_{3,0}}{u_{2,0}^2},
\\
\Phi_5^{3,1}
&
\,=\,
-\,6\,
\frac{u_{1,1}\,u_{3,1}}{u_{2,0}}
+
12\,
\frac{u_{3,0}\,u_{1,1}\,u_{2,1}}{u_{2,0}^2}
-
6\,
\frac{u_{3,0}^2\,u_{1,1}^2}{u_{2,0}^3}
-
6\,
\frac{u_{2,1}^2}{u_{2,0}}
+
3\,
\frac{u_{4,0}\,u_{1,1}^2}{u_{2,0}^2},
\endaligned
\]
and where:
\[
\aligned
\Phi_6^{1,1}
&
\,=\,
-\,
\frac{u_{1,0}\,u_{1,1}^2}{u_{2,0}}
-
2\,
u_{1,1}\,u_{0,1},
\\
\Phi_6^{2,1}
&
\,=\,
-\,4\,
\frac{u_{1,0}\,u_{1,1}\,u_{2,1}}{u_{2,0}}
+
2\,
\frac{u_{1,0}\,u_{1,1}^2\,u_{3,0}}{u_{2,0}^2}
-
3\,u_{1,1}^2
-
2\,u_{2,1}\,u_{0,1},
\\
\Phi_6^{3,1}
&
\,=\,
-\,6\,
\frac{u_{1,0}\,u_{1,1}\,u_{3,1}}{u_{2,0}}
+
12\,
\frac{u_{1,0}\,u_{3,0}\,u_{1,1}\,u_{2,1}}{u_{2,0}^2}
-
6\,
\frac{u_{1,0}\,u_{3,0}^2\,u_{1,1}^2}{u_{2,0}^3}
-
6\,
\frac{u_{1,0}\,u_{2,1}^2}{u_{2,0}}
\,+
\\
&
\ \ \ \ \
+
3\,
\frac{u_{1,0}\,u_{4,0}\,u_{1,1}^2}{u_{2,0}^2}
-
12\,
u_{1,1}\,u_{2,1}
+
2\,
\frac{u_{3,0}\,u_{1,1}^2}{u_{2,0}}
-
2\,u_{3,1}\,u_{0,1}.
\endaligned
\]

\begin{Assertion}
{\bf [On a computer]}
The determinant of the above $6 \times 6$ matrix vanishes 
identically.\qed
\end{Assertion}

\begin{Question}
{\sl 
Then what is the dimension of the generic orbits?}
\end{Question}

\begin{Assertion}
{\bf [On a computer]}
The above $6 \times 6$ matrix has rank $5$ at every point in the domain
of $P\!J_{2,1}^4$ defined by:
\[
u_{2,0}
\,\neq\,0,
\ \ \ \ \ \ \ \ \ \ \ \ \ \ \ \ \ \ \ \
u_{2,0}\,u_{2,1}
-
u_{1,1}\,u_{3,0}
\,\neq\,
0.
\]
\end{Assertion}

\proof 
Denote by $M^{i,j}$ the $5 \times 5$
submatrix of the above $6 \times 6$ matrix
for which the $i$\textsuperscript{th} 
row and the $j$\textsuperscript{th} column are deleted. 
A computer gives:
\[
\aligned
\det\,M^{4,6}\,=\,
&
-18\,u_{2,0}\,
\big(
u_{2,0}\,u_{2,1}
-
u_{1,1}\,u_{3,0}
\big)
\\
&
\big(
3\,u_{2,0}^2\,u_{2,1}^2
-
u_{1,1}\,u_{2,0}\,u_{2,1}\,u_{3,0}
-
2\,u_{1,1}^2\,u_{3,0}^2
-
2\,u_{1,1}\,u_{2,0}^2\,u_{3,1}
+
2\,u_{1,1}^2\,u_{2,0}\,u_{4,0}
\big),
\\
\det\,M^{5,6}\,=\,
&
-18\,u_{2,0}\,
\big(
u_{2,0}\,u_{2,1}
-
u_{1,1}\,u_{3,0}
\big)
\\
&
\big(
3\,u_{2,0}^3\,u_{2,1}
-
3\,u_{1,1}\,u_{2,0}^2\,u_{3,0}
-
5\,u_{1,0}\,u_{2,0}\,u_{2,1}\,u_{3,0}
+
5\,u_{1,0}\,u_{1,1}\,u_{3,0}^2
\\
&
+
2\,u_{1,0}\,u_{2,0}^2\,u_{3,1}
-
2\,u_{1,1}^2\,u_{2,0}\,u_{4,0}
\big),
\\
\det\,M^{6,6}\,=\,
&
-18\,u_{2,0}\,
\big(
u_{2,0}\,u_{2,1}
-
u_{1,1}\,u_{3,0}
\big)
\\
&
\big(
-5\,u_{2,0}\,u_{2,1}\,u_{3,0}
+
5\,u_{1,1}\,u_{3,0}^2
+
2\,u_{2,0}^2\,u_{3,1}
-
2\,u_{1,1}\,u_{2,0}\,u_{4,0}
\big).
\endaligned
\]

Suppose these three determinants are simultaneously $0$. 
Then in our domain $u_{2,0}\,\neq\,0 \neq   
u_{2,0}\,u_{2,1}-u_{1,1}\,u_{3,0}$ we have:
\[
\aligned
0\,=\,
&
3\,u_{2,0}^2\,u_{2,1}^2
-
u_{1,1}\,u_{2,0}\,u_{2,1}\,u_{3,0}
-
2\,u_{1,1}^2\,u_{3,0}^2
-
2\,u_{1,1}\,u_{2,0}^2\,u_{3,1}
+
2\,u_{1,1}^2\,u_{2,0}\,u_{4,0},
\\
0\,=\,
&
3\,u_{2,0}^3\,u_{2,1}
-
3\,u_{1,1}\,u_{2,0}^2\,u_{3,0}
-
5\,u_{1,0}\,u_{2,0}\,u_{2,1}\,u_{3,0}
+
5\,u_{1,0}\,u_{1,1}\,u_{3,0}^2
+
2\,u_{1,0}\,u_{2,0}^2\,u_{3,1}
-
2\,u_{1,1}^2\,u_{2,0}\,u_{4,0},
\\
0\,=\,
&
-5\,u_{2,0}\,u_{2,1}\,u_{3,0}
+
5\,u_{1,1}\,u_{3,0}^2
+
2\,u_{2,0}^2\,u_{3,1}
-
2\,u_{1,1}\,u_{2,0}\,u_{4,0}.
\endaligned
\]
Abbreviating $A := u_{2,0}\,u_{2,1}-u_{1,1}\,u_{3,0}$ 
and $B := u_{2,0}\,u_{3,1}-u_{1,1}\,u_{4,0}$, this is:
\[
\aligned
0\,=\,
&
\big(
3\,u_{2,0}\,u_{2,1}
+
2\,u_{1,1}\,u_{3,0}
\big)\,A
-
2\,u_{1,1}\,u_{2,0}\,B,
\\
0\,=\,
&
\big(
3\,u_{2,0}^2
-
5\,u_{1,0}\,u_{3,0}
\big)\,A
+
2\,u_{1,0}\,u_{2,0}\,B,
\\
0\,=\,
&
-5\,u_{3,0}\,A
+
2\,u_{2,0}\,B.
\endaligned
\]
One can eliminate $u_{2,0}\,B$ and get:
\[
0\,=\,3\,A^2,
\ \ \ \ \ \ \ \ \ \ \ \
0\,=\,3\,u_{2,0}^2\,A.
\]
These equations clearly have 
no solution in our domain. 
So the three minors cannot degenerate simultaneously, 
although each one can degenerate somewhere.
\endproof

Consequently, there is one and only $1 = 6 - 5$ invariant. 
On a computer, one can solve 
the system of $11$ {\sc pde}s:
\[
0
\,\equiv\,
{\sf pj}_\ast\big({\bf v}_\sigma^{(4)}\big)
\big(\Waux\big)
\,\equiv\,
{\sf pj}_\ast\big({\bf w}_\tau^{(4)}\big)
\big(\Waux\big)
\eqno
{\scriptstyle{(1\,\leqslant\,\sigma\,\leqslant\,8,\,\,
1\,\leqslant\,\tau\,\leqslant\,3)}},
\]
for an unknown function $\Waux$ of the $11$ jet variables:
\[
\Big(
x,y,\,\,
\begin{smallmatrix} 
u_{0,1}, & u_{1,1}, & u_{2,1}, & u_{3,1},
\\
\\ 
u, & u_{1,0}, & u_{2,0}, & u_{3,0}, & u_{4,0}
\end{smallmatrix}
\Big),
\]
and obtain its explicit expression.

\begin{Proposition}
\label{Prp-explicit-W-via-vector-fields}
{\bf [On a computer]}
This invariant is:
\[
\Waux_F
\,:=\,
\frac{
F_{xx}^2\,F_{xxxy}
-
F_{xx}\,F_{xy}\,F_{xxxx}
+
2\,F_{xy}\,F_{xxx}^2
-
2\,F_{xx}\,F_{xxx}\,F_{xxy}}{
(F_{xx})^2\,\,
\big(
F_{xx}\,F_{xxy}
-
F_{xy}\,F_{xxx}
\big)^{2/3}}.
\eqno\qed
\]
\end{Proposition}

Because computers rarely succeed as soon as either the number
of variables or the degree increases,
it would be better to have a method for computing $\Waux$ which
would also apply 
to determine the explicit expressions of some
higher order invariants.
To this aim are devoted the next two
Sections~{\ref{relative-invariant-S-first-invariant-W}},
{\ref{branch-W-equiv-0-branch-W-nonzero}}.

\Section{\bf Relative Invariant $\Saux$ and First Invariant $\Waux$}
\label{relative-invariant-S-first-invariant-W}
\HEAD{{\ref{relative-invariant-S-first-invariant-W}}.~{\sf
 Relative Invariant 
$\Saux$ and First Invariant $\Waux$}
}{
Zhangchi {\sc Chen}, Joël~{\sc Merker}}

As in 
Section~{\ref{special-affine-power-series-invariants-curves-R-2}}
for the case of curves,
we apply the power series method 
presented in 
Section~{\ref{power-series-method}} to study equivalences
of analytic
parabolic surfaces $S^2 \subset \R^3$ 
centered at the origin:
\[
u
\,=\,
\sum_{j+k\geqslant 2}\,
F_{j,k}\,
\frac{x^j}{j!}\,
\frac{y^k}{k!},
\]
under the action of the group $\SL_3(\R)$ of special {\em linear}
transformations. Because the special {\em affine} group $\SA_3(\R)$
contains all translations, Theorem~{\ref{Thm-translations-G-0}}
already justified that the computation at the origin of the {\em power
series} $\SL_3$-invariants, in the sense of
Definition~{\ref{Def-power-series-invariant}}, is sufficient to know
the {\em differential} 
$\SA_3$-invariants at every nearby point $(x,y)$.

Furthermore, because $\SL_3(\R)$ also contains all vertical
transvections, Theorem~{\ref{Thm-translations-transvections-G-1}}
already justified that we may assume the first order terms $F_{0,0} +
F_{1,0}\, x + F_{0,1}\, y = 0$ are already normalized to zero, as was
written above.

We can therefore apply the progressive cross-section method of
Section~{\ref{special-affine-power-series-invariants-curves-R-2}}.
The formal coefficients $F_{j,k}$ will be re-initialized
in later stages of the process. At the beginning, we assign
them to be {\sl functional jets:}
\[
F_{j,k}
\,:=\,
u_{j,k},
\]
having in mind $u_{x^j y^k}(x,y)$ in order to receive 
true differential
invariants\,\,---\,\,but 
all computations will be performed at the origin.
We consider special affine transformations:
\[
\R_{x,y,u}^3
\,\,\,\longrightarrow\,\,\,
\R_{s,t,v}^3,
\]
which are not too far from the identity, so that any analytic graph
$\big\{ u = F(x,y) \big\}$ is sent to a similar graphed surface
$\big\{ v = G(s,t) \big\}$.

Since the source power series
$F(x,y) = \sum\, F_{j,k}\,
\frac{x^j\, y^k}{j!\, k!}$ is given and since we want to simplify 
the target power series $G(s,t) = \sum\,  
G_{l,m}\, \frac{s^l\, t^m}{l!\, m!}$, it is more natural to 
work with the {\em inverse} special affine transformation:
\[
\aligned
x
\,=\,
{\sf a}\,s+{\sf b}\,t+{\sf c}\,v,
\\
y
\,=\,
{\sf k}\,s+{\sf l}\,t+{\sf m}\,v,
\\
u
\,=\,
{\sf p}\,s+{\sf q}\,t+{\sf r}\,v.
\endaligned
\] 
Then the graphing function $G(s,t)$ is uniquely determined,
by a {\sl fundamental equation:}
\[
0
\,\equiv\,
-\,
{\sf p}\,s-{\sf q}\,t-{\sf r}\,v
+
F\big(
{\sf a}\,s+{\sf b}\,t+{\sf c}\,v,\,\,
{\sf k}\,s+{\sf l}\,t+{\sf m}\,v
\big)
\Big\vert_{\text{\rm replace}\,\,v\,=\,G(s,t)},
\]
which holds identically as a power series in the two horizontal
variables $(s,t)$.

The goal is to use the group parameters 
freedom
${\sf a}$, ${\sf b}$, ${\sf c}$,
${\sf k}$, ${\sf l}$, ${\sf m}$,
${\sf p}$, ${\sf q}$, ${\sf r}$,
in order to normalize as much as possible 
coefficients $G_{l,m}$.

\Subsection{First loop}
After an affine transformation, we may of course assume that our
target graph enjoys a similar first-order normalization:
\[
v
\,=\,
\sum_{l+m\geqslant 2}\,
G_{l,m}\,
\frac{s^l}{l!}\,
\frac{t^m}{m!}.
\]
Then we perform the replacement:
\leqnomode\usetagform{default}
\begin{equation}
\label{eq-fond-F-G-surfaces}
0
\,\equiv\,
-\,
{\sf p}\,s-{\sf q}\,t-{\sf r}\,G(s,t)
+
F\big(
{\sf a}\,s+{\sf b}\,t+{\sf c}\,G(s,t),\,\,
{\sf k}\,s+{\sf l}\,t+{\sf m}\,G(s,t)
\big),
\end{equation}
and we glean first-order terms which must vanish:
\[
0
\,\equiv\,
-\,{\sf p}\,s
-
{\sf q}\,t
+
{\rm O}_{s,t}(2).
\]

\begin{Lemma}
The subgroup of $\SL_3(\R)$ sending
$v = {\rm O}_{s,t}(2)$
to  
$u = {\rm O}_{x,y}(2)$ 
is $6$-dimensional and consists of matrices:
\[
G_{\stabsmall}^{(1)}
\colon
\ \ \ \ \ \ \ \ \ \ \ \ \ \ \ \ \ \ \ \
\left(\!
\begin{array}{ccc}
{\sf a} & {\sf b} & {\sf c}
\\
{\sf k} & {\sf l} & {\sf m}
\\
0 & 0 & {\sf r}
\end{array}
\!\right),
\]
with $1 = \big( {\sf a}{\sf l} - {\sf b}{\sf k} \big)\,
{\sf r}$.\qed
\end{Lemma}

Next, second order terms in~({\ref{eq-fond-F-G-surfaces}}) are:
\[
\aligned
0
&
\,\equiv\,
\Big(
{\sf k}^2\,
\frac{F_{1,1}^2}{F_{2,0}}
+
{\sf a}^2\,F_{2,0}
+
2\,{\sf a}{\sf k}\,
F_{1,1}-
{\sf r}\,G_{2,0}
\Big)\,
{\textstyle{\frac{s^2}{2}}}
\,+
\\
&
\ \ \ \ \
+
\Big(
{\sf a}{\sf b}\,F_{2,0}
+
{\sf b}{\sf k}\,F_{1,1}
+
{\sf a}{\sf l}\,F_{1,1}
-
{\sf r}\,G_{1,1}
+
{\sf k}{\sf l}\,
\frac{F_{1,1}^2}{F_{2,0}}
\Big)\,
st
+
(\ast)\,
t^2
+
{\rm O}_{s,t}(3),
\endaligned
\]
where $(\ast)$ is unimportant.
We can normalize $G_{2,0} := 1$ and $G_{1,1} := 0$ with the choice
of the simple transformation:
\[
\left(\!
\begin{array}{ccc}
\frac{1}{F_{2,0}^{1/3}} & -\frac{F_{1,1}}{F_{2,0}} & 0
\\
0 & 1 & 0
\\
0 & 0 & F_{2,0}^{1/3}
\end{array}
\!\right)
\,\,\,\in\,\,
G_{\stabsmall}^{(1)},
\]
and the parabolic jet relation $u_{0,2} = \frac{u_{1,1}^2}{u_{2,0}}$
satisfied by $G(s,t)$ gives:
\[
G_{0,2}
\,=\,
0.
\]

After these normalizations,
third order terms become:
\[
0
\,\equiv\,
\Big(
-\,F_{2,0}^{1/3}\,G_{3,0}
+
\frac{F_{3,0}}{F_{2,0}}
\Big)\,
{\textstyle{\frac{s^3}{6}}}
+
\Big(
-\,F_{2,0}^{1/3}\,G_{2,1}
-
\frac{F_{3,0}\,F_{1,1}}{F_{2,0}^{5/3}}
+
\frac{F_{2,1}}{F_{2,0}^{2/3}}
\Big)\,
{\textstyle{\frac{s^2t}{2}}}
+
(\ast)\,st^2
+
(\ast)\,t^3
+
{\rm O}_{s,t}(4).
\]
We solve for the $G_{l,m}$ and we come back to the initial
functional jets:
\[
\aligned
G_{3,0}
&
\,=\,
\frac{F_{3,0}}{F_{2,0}^{4/3}}
\,=\,
\frac{u_{3,0}}{u_{2,0}^{4/3}},
&
\ \ \ \ \ \ \ \ \ \ \ \ \ \ \ \ \ \ \ \
G_{2,1}
&
\,=\,
\frac{F_{2,0}\,F_{2,1}-F_{1,1}\,F_{3,0}}{F_{2,0}^2}
\,=\,
\frac{u_{2,0}\,u_{2,1}-u_{1,1}\,u_{3,0}}{u_{2,0}^2}
\\
G_{4,0}
&
\,=\,
\frac{F_{4,0}}{F_{2,0}^{5/3}}
\,=\,
\frac{u_{4,0}}{u_{2,0}^{5/3}},
&
\ \ \ \ \ \ \ \ \ \ \ \ \ \ \ \ \ \ \ \
G_{3,1}
&
\,=\,
\frac{F_{2,0}\,F_{3,1}-F_{1,1}\,F_{4,0}}{F_{2,0}^{7/3}}
\,=\,
\frac{u_{2,0}\,u_{3,1}-u_{1,1}\,u_{4,0}}{u_{2,0}^{7/3}}.
\endaligned
\]

\Subsection{Second loop}
We restart with two formal power series normalized up to order $2$:
\[
u
\,=\,
\frac{x^2}{2!}
+
\sum_{j+k\geqslant 3}\,
F_{j,k}\,
\frac{x^j}{j!}\,
\frac{y^k}{k!}
\ \ \ \ \ \ \ \ \ \ \ \ \ \ \ \ \ \ \ \
\text{and}
\ \ \ \ \ \ \ \ \ \ \ \ \ \ \ \ \ \ \ \
v
\,=\,
\frac{s^2}{2!}
+
\sum_{l+m\geqslant 3}\,
G_{l,m}\,
\frac{s^l}{l!}\,
\frac{t^m}{m!},
\]
and with the previous stabilizer subgroup
$G_{\stabsmall}^{(1)}$. The fundamental
equation~({\ref{eq-fond-F-G-surfaces}}) is:
\[
\aligned
0
&
\,\equiv\,
-\,{\sf r}\,G(s,t)
+
F\big(
{\sf a}\,s+{\sf b}\,t+{\sf c}\,G(s,t),\,\,
{\sf k}\,s+{\sf l}\,t+{\sf m}\,G(s,t)
\big)
\\
&
\,\equiv\,
\big(
-\,{\sf r}+{\sf a}^2
\big)\,
{\textstyle{\frac{s^2}{2}}}
+
\big({\sf a}{\sf b}\big)\,
st
+
(\ast)\,
{\textstyle{\frac{t^2}{2}}}
+
{\rm O}_{s,t}(3).
\endaligned
\]
Thus ${\sf r} = {\sf a}^2$, then ${\sf b} = 0$, and ${\sf 1} = 
{\sf a} {\sf l} {\sf r}$, so ${\sf l} = {\sf a}^{-3}$.

\begin{Lemma}
The subgroup $G_{\stabsmall}^{(2)} \subset G_{\stabsmall}^{(1)}$
which sends $v = \frac{1}{2}\, s^2 + \cdots$ to $u = \frac{1}{2}\,
x^2 + \cdots$ consists of matrices:
\[
G_{\stabsmall}^{(2)}
\colon
\ \ \ \ \ \ \ \ \ \ \ \ \ \ \ \ \ \ \ \
\left(\!
\begin{array}{ccc}
{\sf a} & 0 & {\sf c}
\\
{\sf k} & {\sf a}^{-3} & {\sf m}
\\
0 & 0 & {\sf a}^2
\end{array}
\!\right).
\eqno\qed
\]
\end{Lemma}

Next, third order terms are, with no stars present:
\[
\aligned
0
&
\,\equiv\,
\Big(
-\,{\sf a}^2\,G_{3,0}
+
3\,{\sf a}^2{\sf k}\,F_{2,1}
+
3\,{\sf a}{\sf c}
+
{\sf a}^3\,F_{3,0}
\Big)\,
{\textstyle{\frac{s^3}{6}}}
\,+
\\
&
\ \ \ \ \
+
\Big(
-\,{\sf a}^2G_{2,1}
+
\frac{1}{{\sf a}}\,
F_{2,1}
\Big)\,
{\textstyle{\frac{s^2t}{2}}}
+
{\rm O}_{s,t}(4).
\endaligned
\]
Solving:
\[
G_{2,1}
\,=\,
\frac{1}{{\sf a}^3}\,
F_{2,1},
\]
we deduce from ${\sf a} \neq 0$ an

\begin{Observation}
After normalization of second order terms,
the properties $F_{2,1} = 0$ and $F_{2,1} \neq 0$ are
$\SL_3(\R)$-invariant.\qed
\end{Observation}

Coming back to the initial functional jets, remind 
we have obtained just above in terms of functional jets:
\[
G_{2,1}
\,=\,
\frac{
u_{2,0}\,u_{2,1}
-
u_{1,1}\,u_{3,0}}{
u_{2,0}^2},
\]
hence we recognize
the relative invariant $\Saux_u$ already
shown in~({\ref{invariant-S-u}}). 
After our normalization, $\Saux_u$ becomes simply:
\[
\Saux_u
\big\vert_{
u_{2,0}=1
\atop
u_{1,1}=0}
\,=\,
\frac{
u_{2,0}\,u_{2,1}
-
u_{1,1}\,u_{3,0}}{
u_{2,0}^2}
\bigg\vert_{
u_{2,0}=1
\atop
u_{1,1}=0}
\,=\,
u_{2,1},
\]
and this explains why we obtained 
the relation $G_{2,1} = F_{2,1} \big/ {\sf a}^3$ between
the plain monomials
$G_{2,1}$ and $F_{2,1}$.

\smallskip

Recall that the branch $\Saux_u \equiv 0$ has already been
studied fully in
Section~{\ref{parabolic-pseudostablization}}.
So we can assume $\Saux_u(x,y) \neq 0$ at every point $(x,y)$,
and in the present context
of power series invariants, this means that we can
assume $F_{2,1} \neq 0 \neq G_{2,1}$.

Then we can normalize $G_{2,1} := 1$ and $G_{3,0} := 0$
by means of the simple transformation:
\[
\left(\!
\begin{array}{ccc}
F_{2,1}^{1/3} & 0 & 0
\\
-\frac{1}{3}\frac{F_{3,0}}{F_{2,1}^{2/3}} & 
\frac{1}{F_{2,1}} & 0
\\
0 & 0 & F_{2,1}^{2/3}
\end{array}
\!\right)
\,\,\,\in\,\,
G_{\stabsmall}^{(2)}.
\]

Looking at terms of order $4$ and solving, we obtain:
\[
\aligned
G_{4,0}
&
\,=\,
-\,
\frac{4}{3}\,
\frac{F_{3,1}\,F_{3,0}}{F_{2,1}^{1/3}}
+\,
\frac{4}{3}\,
\frac{F_{2,1}\,F_{3,0}^2}{F_{2,1}^{1/3}}
+\,
F_{2,1}^{2/3}\,F_{4,0},
\\
G_{3,1}
&
\,=\,
\frac{F_{3,1}}{F_{2,1}^{2/3}}
-\,
2\,F_{2,1}^{1/3}\,F_{3,0},
\endaligned
\]
hence coming back to functional jets:
\[
\aligned
G_{4,0}
&
\,=\,
\frac{1}{3}\,
\frac{4\,u_{3,0}\,u_{2,0}\,u_{4,0}\,u_{1,1}
-4\,u_{3,0}\,u_{2,0}^2\,u_{3,1}
+4\,u_{3,0}^2\,u_{2,1}\,u_{2,0}
-4\,u_{3,0}^3\,u_{1,1}
+3\,u_{4,0}\,u_{2,0}^2\,u_{2,1}}{
u_{2,0}^4\,\,
\big(
u_{2,0}\,u_{2,1}
-u_{1,1}\,u_{3,0}
\big)^{1/3}}
\\
G_{3,1}
&
\,=\,
\frac{-\,u_{2,0}\,u_{4,0}\,u_{1,1}
+u_{2,0}^2\,u_{3,1}
-2\,u_{3,0}\,u_{2,1}\,u_{2,0}
+2\,u_{3,0}^2\,u_{1,1}}{
u_{2,0}^2\,\,
\big(
u_{2,0}\,u_{2,1}
-u_{1,1}\,u_{3,0}
\big)^{2/3}}.
\endaligned
\]

Before continuing, we observe that from the parabolic jet 
relations shown in Section~{\ref{parabolic-jet-relations}},
it comes thanks to $G_{1,1} = 0$:
\[
G_{1,2}
\,=\,
0,
\ \ \ \ \ \ \ \ \ \ \ \ \ \ \ \ \ \ \ \
G_{0,3}
\,=\,
0.
\]

\Subsection{Third loop}
We restart with two formal power series normalized up to order $3$:
\[
u
\,=\,
\frac{x^2}{2}
+
\frac{x^2\,y}{2}
+
\sum_{j+k\geqslant 4}\,
F_{j,k}\,
\frac{x^j}{j!}\,
\frac{y^k}{k!}
\ \ \ \ \ \ \ \ \ \ \ \ \ \ \ \ \ \ \ \
\text{and}
\ \ \ \ \ \ \ \ \ \ \ \ \ \ \ \ \ \ \ \
v
\,=\,
\frac{s^2}{2}
+
\frac{s^2\,t}{2}
+
\sum_{l+m\geqslant 4}\,
G_{l,m}\,
\frac{s^l}{l!}\,
\frac{t^m}{m!},
\]
and with the previous stabilizer subgroup
$G_{\stabsmall}^{(2)}$. The fundamental
equation~({\ref{eq-fond-F-G-surfaces}}) is:
\[
\aligned
0
&
\,\equiv\,
-\,{\sf a}^2\,G(s,t)
+
F\big(
{\sf a}\,s+{\sf c}\,G(s,t),\,\,
{\sf k}\,s+{\sf a}^{-3}\,t+{\sf m}\,G(s,t)
\big)
\\
&
\,\equiv\,
\big(
3\,{\sf a}^2{\sf k}
+
3\,{\sf a}{\sf c}
\big)\,
{\textstyle{\frac{s^3}{6}}}
+
\big(
-\,{\sf a}^2
+
{\textstyle{\frac{1}{{\sf a}}}}
\big)\,
{\textstyle{\frac{s^2t}{2}}}
+
{\rm O}_{s,t}(4).
\endaligned
\]
Thus ${\sf a} = 1$ and ${\sf k} = -\, {\sf c}$.

\begin{Lemma}
The subgroup $G_{\stabsmall}^{(3)} \subset G_{\stabsmall}^{(2)}$
which sends $v = \frac{1}{2}\, s^2 + \frac{s^2 t}{2}
+ \cdots$ to $u = \frac{1}{2}\,
x^2 + \frac{x^2y}{2} + \cdots$ consists of matrices:
\[
G_{\stabsmall}^{(3)}
\colon
\ \ \ \ \ \ \ \ \ \ \ \ \ \ \ \ \ \ \ \
\left(\!
\begin{array}{ccc}
1 & 0 & {\sf c}
\\
-{\sf c} & 1 & {\sf m}
\\
0 & 0 & 1
\end{array}
\!\right).
\eqno\qed
\]
\end{Lemma}

Next, fourth order terms are:
\[
\aligned
0
&
\,\equiv\,
\Big(
-\,G_{4,0}
+
3\,{\sf c}^2
+
6\,{\sf m}
-
4\,{\sf c}\,F_{3,1}
+
F_{4,0}
\Big)\,
{\textstyle{\frac{s^4}{24}}}
\,+
\\
&
\ \ \ \ \
+
\Big(
-\,G_{3,1}
+
F_{3,1}
\Big)\,
{\textstyle{\frac{s^3t}{6}}}
+
{\rm O}_{s,t}(5).
\endaligned
\]

\begin{Observation}
There exists a (single)
$4$\textsuperscript{th} order invariant $G_{3,1} = F_{3,1}$.\qed
\end{Observation}

Its expression in terms of functional jets was already finalized
at the end of the previous loop, and we attribute a name to 
this invariant:
\[
\Waux
\,:=\,
\frac{
u_{2,0}^2\,u_{3,1}
-
u_{2,0}\,u_{4,0}\,u_{1,1}
+
2\,u_{3,0}^2\,u_{1,1}
-
2\,u_{3,0}\,u_{2,1}\,u_{2,0}
}{
u_{2,0}^2\,\,
\big(
u_{2,0}\,u_{2,1}
-u_{1,1}\,u_{3,0}
\big)^{2/3}}.
\]
This confirms what has already been presented with prolonged
vector fields in Section~{\ref{parabolic-pseudostablization}},
{\em cf.} Proposition~{\ref{Prp-explicit-W-via-vector-fields}}. 

\smallskip

Furthermore, we can make $G_{4,0} := 0$ with the 
simple transformation:
\[
\left(\!
\begin{array}{ccc}
1 & 0 & 0
\\
0 & 1 & -\frac{1}{6}F_{4,0}
\\
0 & 0 & 1
\end{array}
\!\right)
\,\,\,\in\,\,
G_{\stabsmall}^{(3)}.
\]

Looking at terms of order $5$ and solving, we obtain:
\[
\aligned
G_{5,0}
&
\,=\,
-\,
\frac{5}{3}\,
F_{4,0}\,F_{3,1}
+
F_{5,0},
\\
G_{4,1}
&
\,=\,
F_{4,1}
-
3\,F_{4,0},
\endaligned
\]
and coming back to functional jets:
\[
\footnotesize
\aligned
G_{5,0}
&
\,=\,
\frac{1}{9}\,
\frac{1}{u_{2,0}^6\,
\big(u_{2,0}\,u_{2,1}-u_{1,1}\,u_{3,0}\big)}\,\,
\bigg\{
\ \ \
5\,
u_{3,0}\,u_{2,0}^2\,u_{4,0}^2\,u_{1,1}^2
-
25\,
u_{3,0}\,u_{2,0}^3\,u_{4,0}\,u_{1,1}\,u_{3,1}
\\
&
\ \ \ \ \ \ \ \ \ \ \ \ \ \ \ \ \ \ \ \ \ \ \ \ \ \ \ \ \ \ \ \ \ \ 
\ \ \ \ \ \ \ \ \ \ \ \ \ \ \ \ \ \ \ \ \ \ \ \
+
30\,u_{3,0}^3\,u_{2,0}\,u_{4,0}\,u_{1,1}^2
+
20\,u_{3,0}\,u_{2,0}^4\,u_{3,1}^2
\\
&
\ \ \ \ \ \ \ \ \ \ \ \ \ \ \ \ \ \ \ \ \ \ \ \ \ \ \ \ \ \ \ \ \ \ 
\ \ \ \ \ \ \ \ \ \ \ \ \ \ \ \ \ \ \ \ \ \ \ \
-\,40\,
u_{3,0}^3\,u_{2,1}^2\,u_{2,0}^2
+
80\,
u_{3,0}^4\,u_{2,1}\,u_{2,0}\,u_{1,1}
\\
&
\ \ \ \ \ \ \ \ \ \ \ \ \ \ \ \ \ \ \ \ \ \ \ \ \ \ \ \ \ \ \ \ \ \ 
\ \ \ \ \ \ \ \ \ \ \ \ \ \ \ \ \ \ \ \ \ \ \ \
-\,40\,
u_{3,0}^5\,u_{1,1}^2
+
15\,
u_{4,0}^2\,u_{2,0}^3\,u_{2,1}\,u_{1,1}
\\
&
\ \ \ \ \ \ \ \ \ \ \ \ \ \ \ \ \ \ \ \ \ \ \ \ \ \ \ \ \ \ \ \ \ \ 
\ \ \ \ \ \ \ \ \ \ \ \ \ \ \ \ \ \ \ \ \ \ \ \
-\,15\,u_{4,0}\,u_{2,0}^4\,u_{2,1}\,u_{3,1}
+
30\,u_{3,0}\,u_{2,0}^3\,u_{4,0}\,u_{2,1}^2
\\
&
\ \ \ \ \ \ \ \ \ \ \ \ \ \ \ \ \ \ \ \ \ \ \ \ \ \ \ \ \ \ \ \ \ \ 
\ \ \ \ \ \ \ \ \ \ \ \ \ \ \ \ \ \ \ \ \ \ \ \
+
9\,
u_{5,0}\,u_{2,0}^4\,u_{2,1}^2
-
3\,u_{3,0}\,u_{2,0}^3\,u_{5,0}\,u_{1,1}\,u_{2,1}
\\
&
\ \ \ \ \ \ \ \ \ \ \ \ \ \ \ \ \ \ \ \ \ \ \ \ \ \ \ \ \ \ \ \ \ \ 
\ \ \ \ \ \ \ \ \ \ \ \ \ \ \ \ \ \ \ \ \ \ \ \
-\,6\,
u_{3,0}^2\,u_{2,0}^2\,u_{5,0}\,u_{1,1}^2
-
15\,
u_{3,0}\,u_{2,0}^4\,u_{4,1}\,u_{2,1}
\\
&
\ \ \ \ \ \ \ \ \ \ \ \ \ \ \ \ \ \ \ \ \ \ \ \ \ \ \ \ \ \ \ \ \ \ 
\ \ \ \ \ \ \ \ \ \ \ \ \ \ \ \ \ \ \ \ \ \ \ \
+
15\,
u_{3,0}^2\,u_{2,0}^3\,u_{4,1}\,u_{1,1}
-
60\,
u_{3,0}^2\,u_{2,0}^2\,u_{4,0}\,u_{1,1}\,u_{2,1}
\ \ \
\bigg\},
\\
G_{4,1}
&
\,=\,
\frac{8\,u_{3,0}^2\,u_{2,1}\,u_{2,0}
-
8\,u_{3,0}^3\,u_{1,1}
+
7\,u_{3,0}\,u_{2,0}\,u_{4,0}\,u_{1,1}
-
4\,u_{3,0}\,u_{2,0}^2\,u_{3,1}
-
u_{2,0}^2\,u_{5,0}\,u_{1,1}
+
u_{2,0}^3\,u_{4,1}
-
3\,u_{4,0}\,u_{2,0}^2\,u_{2,1}}{
u_{2,0}^4\,\,
\big(
u_{2,0}\,u_{2,1}
-
u_{1,1}\,u_{3,0}
\big)^{1/3}}.
\endaligned
\]

According to our general principles, we must take account of the
branching that the invariant $\Waux$ causes:
\[
\xymatrix{
& & & &
\text{\footnotesize\sf Identical degeneracy}\,\,
\Waux_F\,\equiv\,0, 
\\
\Waux_F
\ar[urrr]
\ar[drrr]
& & & &
\\
& & & &
\text{\footnotesize\sf Nowhere vanishing}\,\,
\Waux_F\,\neq\,0.
}
\]
As usual, it is then natural to study first $\Waux_F \equiv 0$, but
before doing this in the next 
Section~{\ref{branch-W-equiv-0-branch-W-nonzero}}, we
prefer to push further a bit our power series algorithm, since it will
definitely show that such a branching lies inside the problem.

Before continuing, we observe that from the parabolic jet 
relations shown in Section~{\ref{parabolic-jet-relations}},
it comes:
\[
G_{2,2}
\,=\,
2,
\ \ \ \ \ \ \ \ \ \ \ \ \ \ \ \ \ \ \ \
G_{1,3}
\,=\,
0,
\ \ \ \ \ \ \ \ \ \ \ \ \ \ \ \ \ \ \ \
G_{0,4}
\,=\,
0.
\]

\Subsection{Fourth loop}
Therefore, we restart with two formal power series normalized 
up to order $4$:
\[
\aligned
u
&
\,=\,
\frac{x^2}{2}
+
\frac{x^2\,y}{2}
+
\Waux\,
\frac{x^3\,y}{6}
+
2\,
\frac{x^2\,y^2}{4}
+
\sum_{j+k\geqslant 4}\,
F_{j,k}\,
\frac{x^j}{j!}\,
\frac{y^k}{k!},
\\
v
&
\,=\,
\frac{s^2}{2}
+
\frac{s^2\,t}{2}
+
\Waux\,
\frac{s^3t}{6}
+
2\,
\frac{s^2\,t^2}{4}
+
\sum_{l+m\geqslant 4}\,
G_{l,m}\,
\frac{s^l}{l!}\,
\frac{t^m}{m!}.
\endaligned
\]
We repeat that for them to be $\SA_3$-equivalent,
their $(3,1)$ coefficients {\em must} be equal:
\[
F_{3,1}
\,=\,
\Waux
\,=\,
G_{3,1}.
\]
We also restart with the previous stabilizer subgroup
$G_{\stabsmall}^{(3)}$. The fundamental
equation~({\ref{eq-fond-F-G-surfaces}}) is:
\[
\aligned
0
&
\,\equiv\,
-\,G(s,t)
+
F\big(
s+{\sf c}\,G(s,t),\,\,
-\,{\sf c}\,s+t+{\sf m}\,G(s,t)
\big)
\\
&
\,\equiv\,
\big(
6\,{\sf m}
+
3\,{\sf c}^2
-
4\,{\sf c}\,F_{3,1}
\big)\,
{\textstyle{\frac{s^4}{24}}}
+
{\rm O}_{s,t}(5).
\endaligned
\]

\begin{Lemma}
\label{LM-G-4-stab-vers-W-0}
The subgroup $G_{\stabsmall}^{(4)} \subset G_{\stabsmall}^{(3)}$
is only $1$-dimensional and consists of matrices:
\[
G_{\stabsmall}^{(4)}
\colon
\ \ \ \ \ \ \ \ \ \ \ \ \ \ \ \ \ \ \ \
\left(\!
\begin{array}{ccc}
1 & 0 & {\sf c}
\\
-{\sf c} & 1 & 
\frac{2}{3}\,{\sf c}\,\Waux
-
\frac{1}{2}\,{\sf c}^2
\\
0 & 0 & 1
\end{array}
\!\right).
\eqno\qed
\]
\end{Lemma}

Next, fifth order terms are:
\[
\aligned
0
&
\,\equiv\,
\Big(
5\,{\sf c}^2\,\Waux
+
\frac{20}{3}\,
{\sf c}\,\Waux^2
-
G_{5,0}
-
5\,{\sf c}\,F_{4,1}
+
F_{5,0}
\Big)\,
{\textstyle{\frac{s^5}{120}}}
\,+
\\
&
\ \ \ \ \
+
\Big(
-\,G_{4,1}
-
2\,{\sf c}\,\Waux
+
F_{4,1}
\Big)\,
{\textstyle{\frac{s^4t}{24}}}
+
{\rm O}_{s,t}(6).
\endaligned
\]
Only when $\Waux \neq 0$, we can normalize $G_{4,1} := 0$ with ${\sf c}
:= \frac{1}{2}\, \frac{F_{4,1}}{\Waux}$, with 
the simple transformation:
\[
\left(\!
\begin{array}{ccc}
1 & 0 & \frac{1}{2}\,\frac{F_{4,1}}{\Waux}
\\
-\frac{1}{2}\frac{F_{4,1}}{\Waux} & 1 &
\frac{1}{3}\,F_{4,1}-\frac{1}{8}\,\frac{F_{4,1}^2}{\Waux^2}
\\
0 & 0 & 1
\end{array}
\!\right)
\,\,\,\in\,\,
G_{\stabsmall}^{(4)}.
\]

The reason we normalize $G_{4,1}$ not $G_{5,0}$ when $\Waux\neq0$ is that the range of $G_{5,0}=5\,c^2\,W+\dots$ is the range of a quadratic function of $c$. It has either a maximum or a minimum. We cannot always normalize it to 0 or any fixed real number. In other words, every orbit crosses the hyperplane $\{G_{4,1}=0\}$ exactly once, but some do not touch $\{G_{5,0}=0\}$.

This confirms our claim 
that the invariant $\Waux$ causes a bifurcation, and we will
now explore the two branches.

\Section{\bf Branch $\Waux \equiv 0$ and Branch $\Waux \neq 0$}
\label{branch-W-equiv-0-branch-W-nonzero}
\HEAD{{\ref{branch-W-equiv-0-branch-W-nonzero}}.~{\sf Branch 
$\Waux \equiv 0$ and Branch $\Waux \neq 0$}
}{
Zhangchi {\sc Chen}, Joël~{\sc Merker}}

Let us treat the branch $\Waux \equiv 0$. We thus have two
differential relations:
\[
\aligned
0
&
\,\equiv\,
F_{xx}\,F_{yy}
-
F_{xy}^2,
\\
0
&
\,\equiv\,
F_{xx}^2\,F_{xxxy}
-
F_{xx}\,F_{xy}\,F_{xxxx}
+
2\,
F_{xy}\,F_{xxx}^2
-
2\,
F_{xx}\,F_{xxx}\,F_{xxy},
\endaligned
\]
which, thanks to the assumption $F_{xx} \neq 0$, can be solved as:
\[
\aligned
F_{yy}
&
\,=\,
\frac{F_{xy}^2}{F_{xx}},
\\
F_{xxxy}
&
\,=\,
\frac{F_{xy}\,F_{xxxx}}{F_{xx}}
-
2\,
\frac{F_{xy}\,F_{xxx}^2}{F_{xx}^2}
+
2\,
\frac{F_{xxx}\,F_{xxy}}{F_{xx}}.
\endaligned
\]
These two relations have differential consequences as well.

\begin{center}
\scalebox{1.5}{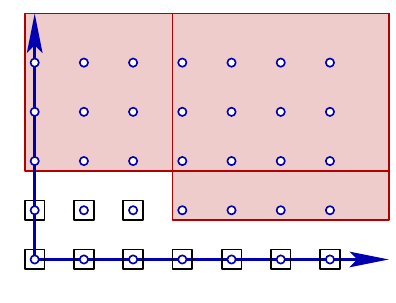}
\end{center}

\begin{Observation}
Every partial derivative $F_{x^k y^l}$ 
in the red region
with $k + l \leqslant \oorder$ and $\oorder \geqslant 2$ 
expresses rationally in terms of the partial derivatives
in the black region:
\[
\big\{
F_{x^{k'}}
\big\}_{k'\leqslant\oorder},
\ \ \ \ \
\big\{
F_y,\ 
F_{xy},\ 
F_{xxy}
\big\},
\]
with denominators containing only powers $\big( F_{xx} \big)^\ast$.\qed
\end{Observation}

Using the jet notation, this means that we will exclusively work
in the {\em submanifold} of the jet space $J_{x,u}^{\oorder}$
defined by up to order $4$ by:
\[
\footnotesize
\aligned
F_{0,2}
&
\,=\,
\frac{F_{1,1}^2}{F_{2,0}},
\\
F_{0,3}
&
\,=\,
3\,\frac{F_{1,1}^2\,F_{2,1}}{F_{2,0}^2}
-
2\,
\frac{F_{1,1}^3\,F_{3,0}}{F_{2,0}^3},
\\
F_{0,4}
&
\,=\,
12\,
\frac{F_{1,1}^2\,F_{2,1}^2}{F_{2,0}^3}
+
\frac{F_{1,1}^4\,F_{4,0}}{F_{2,0}^4}
-
16\,
\frac{F_{1,1}^3\,F_{3,0}\,F_{2,1}}{F_{2,0}^4}
+
4\,
\frac{F_{1,1}^4\,F_{3,0}^2}{F_{2,0}^5},
\\
F_{1,2}
&
\,=\,
-\,
\frac{F_{1,1}^2\,F_{3,0}}{F_{2,0}^2}
+
2\,
\frac{F_{1,1}\,F_{2,1}}{F_{2,0}},
\\
F_{1,3}
&
\,=\,
6\,
\frac{F_{1,1}\,F_{2,1}^2}{F_{2,0}^2}
+
\frac{F_{1,1}^3\,F_{4,0}}{F_{2,0}^3}
-
6\,
\frac{F_{1,1}^2\,F_{3,0}\,F_{2,1}}{F_{2,0}^3},
\\
F_{2,2}
&
\,=\,
\frac{F_{4,0}\,F_{1,1}^2}{F_{2,0}^2}
-
2\,
\frac{F_{1,1}^2\,F_{3,0}^2}{F_{2,0}^3}
+
2\,
\frac{F_{2,1}^2}{F_{2,0}},
\\
F_{3,1}
&
\,=\,
\frac{F_{1,1}\,F_{4,0}}{F_{2,0}}
-
2\,
\frac{F_{1,1}\,F_{3,0}^2}{F_{2,0}^2}
+
2\,
\frac{F_{3,0}\,F_{2,1}}{F_{2,0}}.
\endaligned
\]
Observe that in contrast to the general branch
$\Haux \equiv 0 \neq \Waux$, 
the jet coordinates $F_{3,1}$, $F_{4,1}$, \dots, 
are {\em dependent}.
Higher differential relations will also be needed, for instance:
\[
\footnotesize
\aligned
F_{0,5}
&
\,=\,
60\,
\frac{F_{1,1}^2\,F_{2,1}^3}{F_{2,0}^4}
+
\frac{F_{1,1}^5\,F_{5,0}}{F_{2,0}^5}
-
120\,
\frac{F_{1,1}^3\,F_{2,1}^2\,F_{3,0}}{F_{2,0}^5}
\,+
\\
&
\ \ \ \ \
+
15\,
\frac{F_{1,1}^4\,F_{4,0}\,F_{2,1}}{F_{2,0}^5}
-
15\,
\frac{F_{1,1}^5\,F_{3,0}\,F_{4,0}}{F_{2,0}^6}
+
60\,
\frac{F_{1,1}^4\,F_{2,1}\,F_{3,0}^2}{F_{2,0}^6},
\\
F_{1,4}
&
\,=\,
-\,12\,
\frac{F_{1,1}^4\,F_{3,0}\,F_{4,0}}{F_{2,0}^5}
+
12\,
\frac{F_{1,1}^4\,F_{3,0}^3}{F_{2,0}^6}
+
12\,
\frac{F_{1,1}^3\,F_{2,1}\,F_{4,0}}{F_{2,0}^4}
\,+
\\
&
\ \ \ \ \
+
\frac{F_{1,1}^4\,F_{5,0}}{F_{2,0}^4}
+
24\,
\frac{F_{1,1}\,F_{2,1}^3}{F_{2,0}^3}
-
36\,
\frac{F_{1,1}^2\,F_{3,0}\,F_{2,1}^2}{F_{2,0}^4},
\\
F_{2,3}
&
\,=\,
-\,9\,
\frac{F_{3,0}\,F_{1,1}^3\,F_{4,0}}{F_{2,0}^4}
-
18\,
\frac{F_{2,1}\,F_{1,1}^2\,F_{3,0}^2}{F_{2,0}^4}
+
12\,
\frac{F_{3,0}^3\,F_{1,1}^3}{F_{2,0}^5}
\,+
\\
&
\ \ \ \ \
+
9\,
\frac{F_{2,1}\,F_{1,1}^2\,F_{4,0}}{F_{2,0}^3}
+
\frac{F_{1,1}^3\,F_{5,0}}{F_{2,0}^3}
+
6\,
\frac{F_{2,1}^3}{F_{2,0}^2},
\\
F_{3,2}
&
\,=\,
6\,
\frac{F_{3,0}^3\,F_{1,1}^2}{F_{2,0}^4}
-
6\,
\frac{F_{3,0}\,F_{1,1}^2\,F_{4,0}}{F_{2,0}^3}
+
6\,
\frac{F_{2,1}^2\,F_{3,0}}{F_{2,0}^2}
\,+
\\
&
\ \ \ \ \
+
6\,
\frac{F_{1,1}\,F_{4,0}\,F_{2,1}}{F_{2,0}^2}
-
12\,
\frac{F_{3,0}^2\,F_{1,1}\,F_{2,1}}{F_{2,0}^3}
+
\frac{F_{5,0}\,F_{1,1}^2}{F_{2,0}^2},
\\
F_{4,1}
&
\,=\,
3\,
\frac{F_{4,0}\,F_{2,1}}{F_{2,0}}
-
3\,
\frac{F_{3,0}\,F_{1,1}\,F_{4,0}}{F_{2,0}^2}
+
\frac{F_{1,1}\,F_{5,0}}{F_{2,0}}.
\endaligned
\]

In the branch $\Haux \equiv 0 \equiv \Waux$ the first three
loops are the same as in the preceding
Section~{\ref{relative-invariant-S-first-invariant-W}}.
So we do not repeat the constructions.

\Subsection{Fourth loop}
Assuming $F_{3,1} = G_{3,1} = \Waux = 0$, we restart from:
\[
\aligned
u
&
\,=\,
{\textstyle{\frac{x^2}{2}}}\,
+
{\textstyle{\frac{x^2y}{2}}}
+
0
+
{\textstyle{\frac{x^2y^2}{2}}}
+
\sum_{j+k\geqslant 5}\,
F_{j,k}\,
{\textstyle{\frac{x^j}{j!}}}\,
{\textstyle{\frac{y^k}{k!}}},
\\
v
&
\,=\,
{\textstyle{\frac{s^2}{2}}}\,
+
{\textstyle{\frac{s^2t}{2}}}
+
0
+
{\textstyle{\frac{s^2t^2}{2}}}
+
\sum_{l+m\geqslant 5}\,
G_{l,m}\,
{\textstyle{\frac{s^l}{l!}}}\,
{\textstyle{\frac{t^m}{m!}}}.
\endaligned
\]
When computing the coefficient $F_{5,0}$ in terms of the initial
functional jets, we must take account of the differential
relations written above, and we obtain:
\leqnomode\usetagform{default}
\begin{align}
\label{invariant-F-5-0}
F_{5,0}
\,=\,
\frac{1}{9}\,
\frac{\big(u_{2,0}\,u_{2,1}-u_{1,1}\,u_{3,0}\big)\,
\big(9\,u_{2,0}^2\,u_{5,0}
-45\,u_{2,0}\,u_{3,0}\,u_{4,0}
+40\,u_{3,0}^3\big)}{
u_{2,0}^6}.
\end{align}
We also see that:
\[
F_{4,1}
\,=\,
0,
\ \ \ \ \ \ \ \ \
F_{3,2}
\,=\,
0,
\ \ \ \ \ \ \ \ \
F_{2,3}
\,=\,
6,
\ \ \ \ \ \ \ \ \
F_{1,4}
\,=\,
0,
\ \ \ \ \ \ \ \ \
F_{0,5}
\,=\,
0.
\]

We also restart with the group $G_{\stabsmall}^{(4)}$ of
Lemma~{\ref{LM-G-4-stab-vers-W-0}}
but in which we set $\Waux := 0$:
\[
G_{\stabsmall}^{(4)}
\colon
\ \ \ \ \ \ \ \ \ \ \ \ \ \ \ \ \ \ \ \
\left(\!
\begin{array}{ccc}
1 & 0 & {\sf c}
\\
-{\sf c} & 1 & 
-
\frac{1}{2}\,{\sf c}^2
\\
0 & 0 & 1
\end{array}
\!\right).
\]
This group (stabilizes) sends $v = \frac{s^2}{2} + \frac{s^2t}{2} +
\frac{s^2t^2}{2} + {\rm O}_{s,t}(5)$ to
$u = \frac{x^2}{2} + \frac{x^2y}{2} + \frac{x^2y^2}{2} + 
{\rm O}_{x,y}(5)$.

Next, we look at $5$\textsuperscript{th} order terms, and we realize
that:
\[
G_{5,0}
\,=\,
F_{5,0},
\]
hence the expression~({\ref{invariant-F-5-0}}) of $F_{5,0}$ is 
a differential invariant. We will call it:
\[
\Xaux
\,:=\,
\frac{1}{9}\,
\frac{\big(u_{xx}\,u_{xxy}-u_{xy}\,u_{xxx}\big)\,
\big(9\,u_{xx}^2\,u_{xxxxx}
-45\,u_{xx}\,u_{xxx}\,u_{xxxx}
+40\,u_{xxx}^3\big)}{
u_{xx}^6}.
\]
We observe that in the current branch $u_{xx} \neq 0 \neq 
\Saux$, we have:
\[
\Xaux
\,\equiv\,
0
\ \ \ \ \ \ \ \ \ \ \ \ \ \ \ \ \ \ \ \
\Longleftrightarrow
\ \ \ \ \ \ \ \ \ \ \ \ \ \ \ \ \ \ \ \
9\,u_{xx}^2\,u_{xxxxx}
-45\,u_{xx}\,u_{xxx}\,u_{xxxx}
+40\,u_{xxx}^3
\,\equiv\,
0.
\]

The following result has been proved for the larger, full
affine group ${\sf A}_3(\R)$
in~{\cite[Thm.~4.1]{Merker-2019}}.
We verify that it also holds for $\SA_3(\R)$.

\begin{Theorem}
\label{Thm-flat-cone}
For a local real analytic graph $\big\{ u = F(x,y) \big\}$ 
at the origin of $\R_{x,y}^2 \times \R_u^1$ which satisfies:
\[
F_{xx}
\,\neq\,
0,
\ \ \ \ \ \ \ \ \ \ \ \ \ \ \ \ \ \ \ \
F_{xx}\,F_{yy}
-
F_{xy}^2
\,\equiv\,
0,
\ \ \ \ \ \ \ \ \ \ \ \ \ \ \ \ \ \ \ \
F_{xx}\,F_{xxy}
-
F_{xy}\,F_{xxx}
\,\neq\,
0,
\]
the following two conditions are equivalent:

\smallskip\noindent{\bf (i)}\,
it is special affinely equivalent to:
\[
v
\,=\,
\frac{1}{2}\,
\frac{s^2}{1-t};
\]

\smallskip\noindent{\bf (ii)}\,
$\Waux(F) \equiv 0 \equiv \Xaux(F)$, that is to say:
\[
\aligned
0
&
\,\equiv\,
F_{xx}^2\,F_{xxxy}
-
F_{xx}\,F_{xy}\,F_{xxxx}
+
2\,
F_{xy}\,F_{xxx}^2
-
2\,
F_{xx}\,F_{xxx}\,F_{xxy},
\\
0
&
\,\equiv\,
9\,F_{xx}^2\,F_{xxxxx}
-
45\,F_{xx}\,F_{xxx}\,F_{xxxx}
+
40\,F_{xxx}^3.
\endaligned
\]
\end{Theorem}

\proof
Only {\small\bf (ii)} $\Longrightarrow$ {\small\bf (i)}
matters.
We have seen that after some $\SA_3(\R)$ transformation,
we have:
\[
u
\,=\,
F(x,y)
\,=\,
{\textstyle{\frac{x^2}{2}}}
+
{\textstyle{\frac{x^2\,y}{2}}}
+
{\textstyle{\frac{x^2\,y^2}{2}}}
+
{\rm O}_{x,y}(5).
\]
Because $\Waux$ and $\Maux$ {\em are} differential
invariants, the conditions $\Waux \equiv 0$ and
$\Maux \equiv 0$ still hold. This graphing function $F$
satisfies:
\[
\aligned
0
&
\,=\,
F_{xy}(0),
&
\ \ \ \ \ \ \ \ \ \ \ 
0
&
\,=\,
F_{yy}(0),
&
\ \ \ \ \ \ \ \ \ \ \
\\
0
&
\,=\,
F_{xxx}(0),
&
\ \ \ \ \ \ \ \ \ \ \
1
&
\,=\,
F_{xxy}(0),
&
\ \ \ \ \ \ \ \ \ \ \
0
&
\,=\,
F_{xyy}(0)
\,=\,
F_{yyy}(0),
\\
0
&
\,=\,
F_{xxxx}(0),
&
\ \ \ \ \ \ \ \ \ \ \
0
&
\,=\,
F_{xxxy}(0),
&
\ \ \ \ \ \ \ \ \ \ \
2
&
\,=\,
F_{xxyy}(0),
\ \ \ \ \ \ \ \ \ \ \
0
\,=\,
F_{xyyy}(0)
\,=\,
F_{yyyy}(0).
\endaligned
\]

The condition $\Maux(x,y) \equiv 0$, valid at every point,
can be solved as:
\[
\aligned
F_{xxxxx}
&
\,=\,
\Big(
5\,
\frac{F_{xxxx}}{F_{xx}}
-
\frac{40}{9}\,
\frac{F_{xxx}^2}{F_{xx}^2}
\Big)\,
F_{xxx}
\\
&
\,=\,
\mathcal{R}
\cdot
F_{xxx}.
\endaligned
\]

We claim that $F_{x^j} = \mathcal{R}\, F_{xxx} + \mathcal{R}\,
F_{xxxx}$ for all $j \geqslant 4$. This is true for
$j = 4, 5$, and the induction is:
\[
\aligned
F_{x^{j+1}}
&
\,=\,
\mathcal{R}_x\,
F_{xxx}
+
\mathcal{R}\,
F_{xxxx}
+
\mathcal{R}_x\,
F_{xxxx}
+
\mathcal{R}\,
F_{xxxxx}
\\
&
\,=\,
\mathcal{R}\,
F_{xxx}
+
\mathcal{R}\,
F_{xxxx}.
\endaligned
\]
Since $F_{xxx}(0) = 0 = F_{xxxx}(0)$, we get $F_{x^j}(0) = 0$
for all $j \geqslant 3$.

Next, we solve from $\Waux (x,y) \equiv 0$:
\[
\aligned
F_{xxxy}
&
\,=\,
\Big(
-\,2\,
\frac{F_{xy}\,F_{xxx}}{F_{xx}^2}
+
2\,\frac{F_{xxy}}{F_{xx}}
\Big)\,
F_{xxx}
+
\Big(
\frac{F_{xy}}{F_{xx}}
\Big)\,
F_{xxxx}
\\
&
\,=\,
\mathcal{R}\,
F_{xxx}
+
\mathcal{R}\,
F_{xxxx}.
\endaligned
\]
The same argument gives $F_{x^j y}(0) = 0$ for all $j \geqslant 3$.

We claim that $F_{x^jy^k}(0) = 0$ for
all $j \geqslant 3$ and all $k \geqslant 2$.
Indeed, by induction from $F_{x^j y^{k-1}} = \mathcal{R}\,
F_{xxx} + \mathcal{R}\, F_{xxxx}$, we get:
\[
\aligned
F_{x^jy^k}
&
\,=\,
\mathcal{R}_y\,
F_{xxx}
+
\mathcal{R}\,F_{xxxy}
+
\mathcal{R}_y\,
F_{xxxx}
+
\mathcal{R}\,
F_{xxxxy}
\\
&
\,=\,
\mathcal{R}\,
F_{xxx}
+
\mathcal{R}\,
F_{xxxx}.
\endaligned
\]

\begin{center}
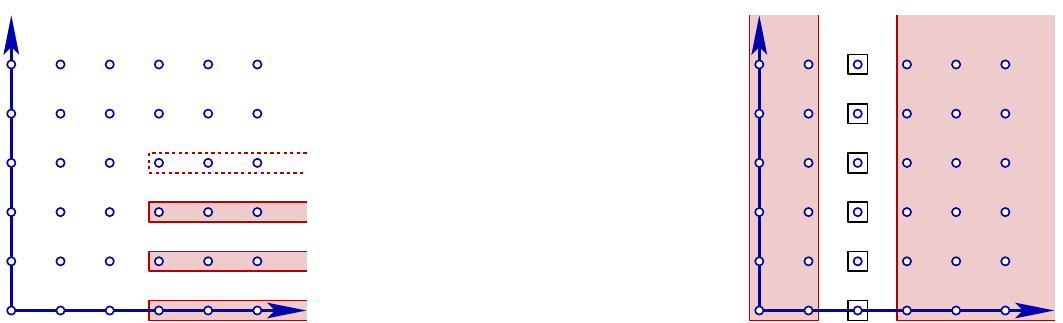
\end{center}

\noindent
This offers $F_{x^jy^k}(0) = 0$ for all $j \geqslant 3$, all
all $k \geqslant 2$, and hence $F$ reduces to:
\[
F(x,y)
\,=\,
F_0(y)
+
x\,F_1(y)
+
x^2\,F_2(y),
\]
with $F_2(0) = \frac{1}{2}$, and due to the normalization above,
we have $F_0(y) = {\rm O}_y(5)$ and $F_1(y) = {\rm O}_y(4)$.

We can now take account of our constant hypothesis that the
Hessian determinant vanishes:
\[
0
\,\equiv\,
\big(
2\,F_2(y)
\big)\,
\Big(
F_{0,yy}
+
x\,F_{1,yy}
+
x^2\,F_{2,yy}
\Big)
-
\Big(
F_{1,y}
+
2x\,F_{2,y}
\Big)^2.
\]
The coefficients of $x^0$, $x^1$ yield\,\,---\,\,remind 
$F_2(0) \neq 0$\,\,---:
\[
\aligned
F_{0,yy}
&
\,=\,
\frac{F_{1,y}^2}{2\,F_2}
\,=\,
\mathcal{R}\,
F_{1,y},
\\
F_{1,yy}
&
\,=\,
\frac{4\,F_{1,y}\,F_{2,y}}{2\,F_2}
\,=\,
\mathcal{R}\,
F_{1,y}.
\endaligned
\]
From $F_{1,y}(0) = 0$, we get $F_{1,yy}(0) = 0$. By iteration:
\reqnomode\usetagform{EngelLie}
\begin{align}
F_{1,yyy}
&
\,=\,
\mathcal{R}_y\,
F_{1,y}+
\mathcal{R}\,F_{1,yy}
\,=\,
\mathcal{R}\,F_{1,y},
\notag
\\
F_{1,y^k}
&
\,=\,
\mathcal{R}\,F_{1,y}
\tag{(k\,\geqslant\,3),}
\end{align}
so $F_{1,y^k}(0) = 0$, whence $F_1(y) \equiv 0$. It comes 
$F_{0, yy}(y) \equiv 0$, and lastly, $F_0(y) \equiv 0$.
In sum:
\[
F(x,y)
\,=\,
x^2\,F_2(y)
\,=:\,
x^2\,T(y),
\]
with $T(0) = \frac{1}{2} = T_y(0)$. Back to the Hessian:
\[
0
\,\equiv\,
2\,T\cdot
x^2\,T_{yy}
-
\big(
2\,x\,T_y
\big)^2,
\]
we solve:
\[
\aligned
T_{yy}
\,=\,
2!\,\frac{(T_y)^2}{T}
\ \ \ \ \ \ \ \
&
\Longrightarrow
\ \ \ \ \ \ \ \ 
T_{yyy}
\,=\,
2!\,
\frac{2\,T_y\,T_{yy}}{T}
-
2!\,\frac{(T_y)^2\,T_y}{T^2}
\,=\,
3!\,
\frac{(T_y)^3}{T^2}
\\
&
\Longrightarrow
\ \ \ \ \ \ \ \ 
T_{y^k}
\,=\,
k!\,
\frac{(T_y)^k}{T^{k-1}},
\endaligned
\]
whence $T_{y^k}(0) = k!\, \frac{1}{2}$, and thus
$T(y) = \frac{1}{2} + \frac{1}{2}\, y + 
\frac{1}{2}\, y^2 + \cdots + \frac{1}{2}\, y^k + \cdots$,
and finally, after having performed only 
special affine transformations:
\[
u
\,=\,
\frac{1}{2}\,
\frac{x^2}{1-y}.
\qedhere
\]
\endproof

Consequently, we can assume that $\Xaux \neq 0$.

\Subsection{Fifth loop}
We restart from:
\[
\aligned
u
&
\,=\,
{\textstyle{\frac{x^2}{2}}}\,
+
{\textstyle{\frac{x^2y}{2}}}
+
0
+
{\textstyle{\frac{x^2y^2}{2}}}
+
\Xaux\,
{\textstyle{\frac{x^5}{5!}}}\,
+
{\textstyle{\frac{x^2y^3}{2}}}
+
\sum_{j+k\geqslant 6}\,
F_{j,k}\,
{\textstyle{\frac{x^j}{j!}}}\,
{\textstyle{\frac{y^k}{k!}}},
\\
v
&
\,=\,
{\textstyle{\frac{s^2}{2}}}\,
+
{\textstyle{\frac{s^2t}{2}}}
+
0
+
{\textstyle{\frac{s^2t^2}{2}}}
+
\Xaux\,
{\textstyle{\frac{s^5}{5!}}}\,
+
{\textstyle{\frac{s^2t^3}{2}}}
+
\sum_{j+k\geqslant 6}\,
G_{l,m}\,
{\textstyle{\frac{s^l}{l!}}}\,
{\textstyle{\frac{t^m}{m!}}},
\endaligned
\]
with a common coefficient $\Xaux \neq 0$. 
We examine how the $1$-dimensional subgroup 
$G_{\stabsmall}^{(4)}$ acts on 
the single independent $6$\textsuperscript{th} order
coefficient: 
\[
\aligned
0
&
\,\equiv\,
-\,G(s,t)
+
F\big({\sf s}+{\sf c}\,G(s,t),\,\,
-\,{\sf c}\,s+t-
{\textstyle{\frac{1}{2}}}\,
{\sf c}^2\,G(s,t)
\big)
\\
&
\,\equiv\,
\Big(
-3\,
{\sf c}\,\Xaux
+
F_{6,0}
-
G_{6,0}
\Big)\,
{\textstyle{\frac{s^6}{6!}}}
+
{\rm O}_{s,t}(7),
\endaligned
\]
Thus, we can normalize $G_{6,0} := 0$ by means of the simple
transformation:
\[
\left(\!
\begin{array}{ccc}
1 & 0 & \frac{1}{3}\,\frac{F_{6,0}}{\Xaux}
\\
-\frac{1}{3}\frac{F_{6,0}}{\Xaux} & 1 &
-\frac{1}{18}\,\frac{F_{6,0}^2}{\Xaux^2}
\\
0 & 0 & 1
\end{array}
\!\right)
\,\,\,\in\,\,
G_{\stabsmall}^{(4)}.
\]
Then in terms of functional jets, there is an invariant
of order $7$, that we call $\Yaux$:

\begin{center}
\includegraphics[scale=0.85]{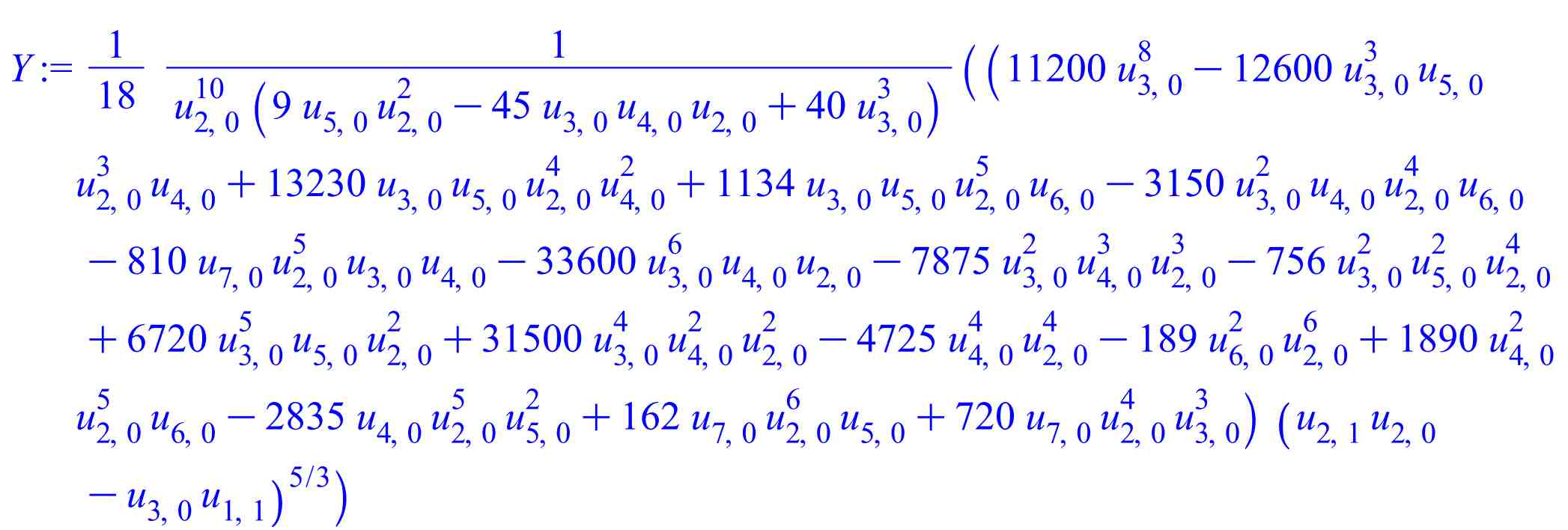}
\end{center}

The synthesis for this branch is made in 
Section~{\ref{introduction-parabolic-surfaces}},
especially in Theorem~{\ref{Thm-S-nonzero-W-0}}.

\medskip

Next, we treat the branch $\Waux \neq 0$.
At the end of Section~{\ref{relative-invariant-S-first-invariant-W}},
we have seen that when $\Waux \neq 0$, we can normalize
$G_{4,1} = 0$. Solving the other coefficient there, we get:
\[
G_{5,0}
\,=\,
-\,\frac{5}{4}\,
\frac{F_{4,1}^2}{F_{3,1}}
+
\frac{10}{3}\,
F_{3,1}\,F_{4,1}
+
F_{5,0}.
\]
Since the last group parameter ${\sf c}$ has been consumed,
it is not necessary to do a fifth loop, and we conclude
that:

\begin{Observation}
In the branch $\Waux \neq 0$, there is a (single) 
$5$\textsuperscript{th} order $\SA_3$-invariant.\qed
\end{Observation}

In terms of functional jets, its explicit expression 
incorporates {\bf 57} monomials in the numerator, and we
call it $\Maux$.

\begin{center}
\includegraphics[scale=0.75]{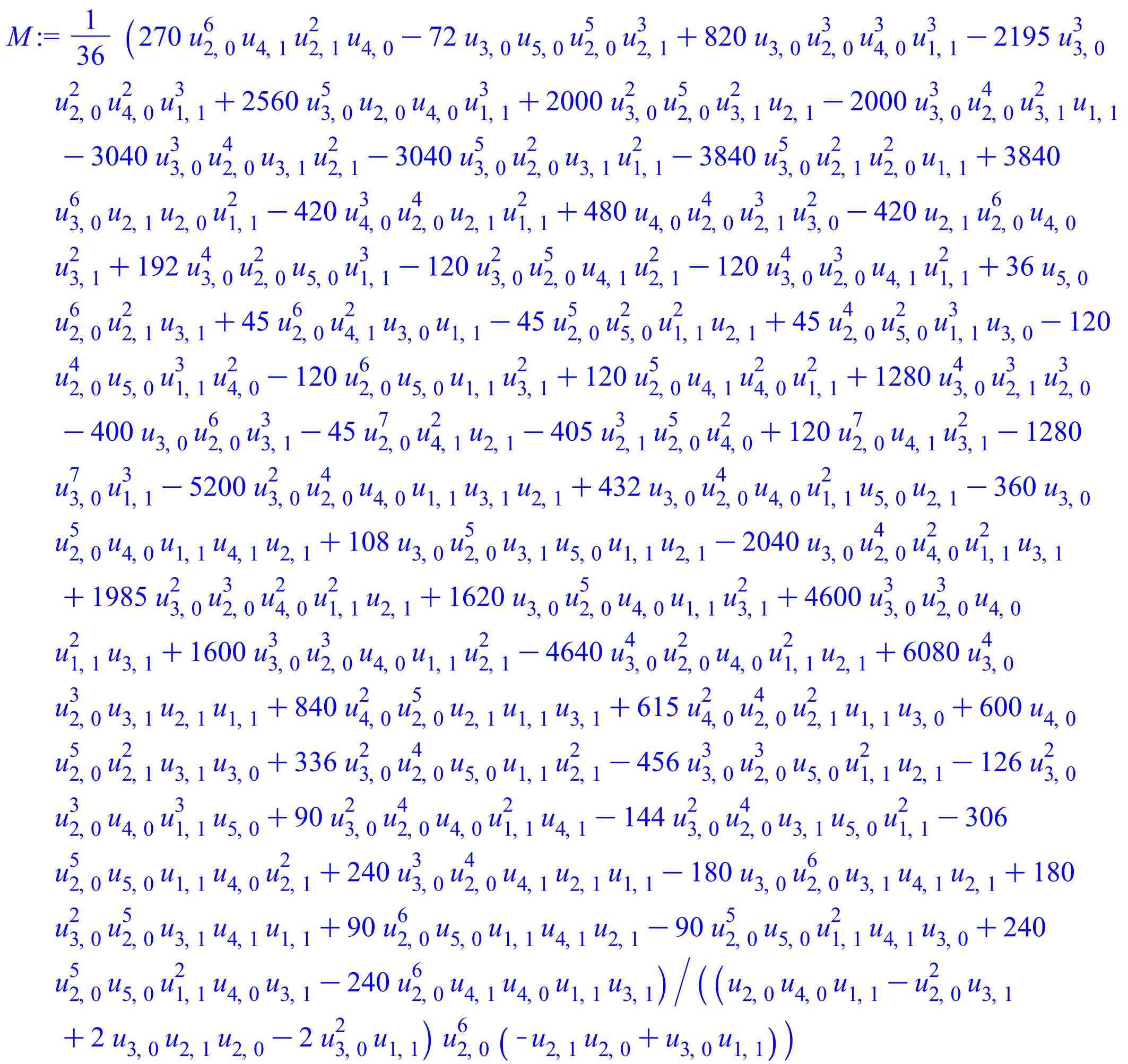}
\end{center}

The synthesis for this branch is made in 
Section~{\ref{introduction-parabolic-surfaces}},
especially in Theorem~{\ref{Thm-S-nonzero-W-nonzero}}.

\Section{\bf Recurrence Relations for Parabolic Surfaces}
\label{recurrence-relations-parabolic-surfaces}
\HEAD{{\ref{recurrence-relations-parabolic-surfaces}}.~{\sf 
Recurrence Relations for Parabolic Surfaces}
}{
Zhangchi {\sc Chen}, Joël~{\sc Merker}}

Since the branch $\Saux \equiv 0$ has already been done in 
Section~{\ref{affine-invariants-curves-R-2}},
we assume from now on that $\Saux \neq 0$, still with 
$F_{xx} \neq 0 \equiv \Haux$. 

\Subsection{Branch $\Waux \neq 0$}
\label{Subsection-W-nonzero}
By what precedes, the normalization at the origin takes the form:
\[
\aligned
u
&
\,=\,
F(x,y)
\\
&
\,=\,
{\textstyle{\frac{x^2}{2}}}
+
{\textstyle{\frac{x^2y}{2}}}
+
\Waux\,
{\textstyle{\frac{x^3y}{6}}}
+
{\textstyle{\frac{x^2y^2}{2}}}
+
\Maux\,
{\textstyle{\frac{x^5}{5!}}}
+
6\,\Waux\,
{\textstyle{\frac{x^3y^2}{3!\,2!}}}
+
{\textstyle{\frac{x^2y^3}{2}}}
+
{\rm O}_{x,y}(6).
\endaligned
\]
Later, we will need from parabolic jet relations:
\[
F_{4,2}
\,=\,
6\,
\big(
F_{3,1}
\big)^2
\,=\,
6\,\Waux^2.
\]

All invariants are:
\[
\Iaux_{j,k}
\,:=\,
\inv\,
\big(
u_{x^jy^k}
\big)
\eqno
{\scriptstyle{(j\,+\,k\,\geqslant\,2)}}.
\]
Of course, $\SA_3(\R)$ contains all translations and all vertical
transvections, hence 
Theorem~{\ref{Thm-translations-transvections-G-1}} 
applies.

Moreover, by Section~{\ref{Subsection-G-1-MC-invariants}},
we can reduce ourselves to the $6$ vector fields:

\[
\aligned
{\bf v_1}
&
\,:=\,
x\,\partial_x
-
u\,\partial_u, 
\ \ \ \ \ \ \ \ \ \ \
{\bf v}_2
\,:=\,
y\,\partial_y
-
u\,\partial_u,
\\
{\bf v}_3
\,:=\,
y\,\partial_x,
&
\ \ \ \ \ \ \
{\bf v}_4
\,:=\,
u\,\partial_x,
\ \ \ \ \ \ \
{\bf v}_5
\,:=\,
x\,\partial_y,
\ \ \ \ \ \ \
{\bf v}_6
\,:=\,
u\,\partial_y,
\endaligned
\]
and we have the recurrence formul{ae}:
\[
\aligned
\mathcal{D}_1
\Iaux_{j,k}
&
\,=\,
\Iaux_{j+1,k}
+
\sum_{1\leqslant\sigma\leqslant 6}\,
\Phi_\sigma^{j,k}
\big(
\Iaux^{(j+k)}
\big)
\cdot
\Kaux_1^\sigma,
\\
\mathcal{D}_2
\Iaux_{j,k}
&
\,=\,
\Iaux_{j,k+1}
+
\sum_{1\leqslant\sigma\leqslant 6}\,
\Phi_\sigma^{j,k}
\big(
\Iaux^{(j+k)}
\big)
\cdot
\Kaux_2^\sigma.
\endaligned
\]
According to what precedes, the $6$ relevant phantom invariants
are:
\[
\Iaux_{2,0}
\,=\,
1,
\ \ \ \ \ \ \ \ \ \ \
\Iaux_{1,1}
\,=\,
0,
\ \ \ \ \ \ \ \ \ \ \
\Iaux_{3,0}
\,=\,
0,
\ \ \ \ \ \ \ \ \ \ \
\Iaux_{2,1}
\,=\,
1,
\ \ \ \ \ \ \ \ \ \ \
\Iaux_{4,0}
\,=\,
0,
\ \ \ \ \ \ \ \ \ \ \
\Iaux_{4,1}
\,=\,
0.
\]
We denote specially:
\[
\Waux
\,:=\,
\Iaux_{3,1}
\ \ \ \ \ \ \ \ \ \ \ \ \ \ \ \ \ \ \ \
\text{and}
\ \ \ \ \ \ \ \ \ \ \ \ \ \ \ \ \ \ \ \
\Maux
\,:=\,
\Iaux_{5,0},
\]
and we apply the first collection of recurrence formulas:
\[
\left(\!
\begin{array}{c}
0
\\
0
\\
0
\\
0
\\
0
\\
0
\end{array}
\!\right)
\,=\,
\left(\!
\begin{array}{c}
\mathcal{D}_1\Iaux_{2,0}
\\
\mathcal{D}_1\Iaux_{1,1}
\\
\mathcal{D}_1\Iaux_{3,0}
\\
\mathcal{D}_1\Iaux_{2,1}
\\
\mathcal{D}_1\Iaux_{4,0}
\\
\mathcal{D}_1\Iaux_{4,1}
\end{array}
\!\right)
\,\,=\,\,
\left(\!
\begin{array}{c}
0
\\
1
\\
0
\\
\Waux
\\
\Maux
\\
\Iaux_{5,1}
\end{array}
\!\right)
+
\left(\!
\begin{array}{cccccc}
-3 & -1 & 0 & 0 & 0 & 0 
\\
0 & 0 & -1 & 0 & 0 & 0
\\
0 & 0 & 0 & -3 & -3 & 0
\\
-3 & -2 & 0 & 0 & 0 & 0
\\
0 & 0 & 0 & 0 & -4\,\Waux & -6
\\
0 & 0 & -\Maux & -10\,\Waux & -24\,\Waux & -18
\end{array}
\!\right)\,
\left(\!
\begin{array}{c}
\Kaux_1^1
\\
\Kaux_1^2
\\
\Kaux_1^3
\\
\Kaux_1^4
\\
\Kaux_1^5
\\
\Kaux_1^6
\end{array}
\!\right).
\]
This Cramér system has the unique solution:
\[
\aligned
&
\Kaux_1^1
\,=\,
-\,\frac{1}{3}\,
\Waux,
\ \ \ \ \ \ \ \ \ \
& &
\Kaux_1^2
\,=\,
\Waux,
\ \ \ \ \ \ \ \ \ \
& &
\Kaux_1^3
\,=\,
1,
\\
&
\Kaux_1^4
\,=\,
2\,\frac{\Maux}{\Waux}
-
\frac{1}{2}\,
\frac{\Iaux_{5,1}}{\Waux},
\ \ \ \ \ \ \ \ \ \
& &
\Kaux_1^5
\,=\,
-2\,\frac{\Maux}{\Waux}
+
\frac{1}{2}\,
\frac{\Iaux_{5,1}}{\Waux},
\ \ \ \ \ \ \ \ \ \
& &
\Kaux_1^6
\,=\,
\frac{3}{2}\,
\Maux
-
\frac{1}{3}\,
\Iaux_{5,1}.
\endaligned
\]

The second collection is:
\[
\left(\!
\begin{array}{c}
0
\\
0
\\
0
\\
0
\\
0
\\
0
\end{array}
\!\right)
\,=\,
\left(\!
\begin{array}{c}
\mathcal{D}_2\Iaux_{2,0}
\\
\mathcal{D}_2\Iaux_{1,1}
\\
\mathcal{D}_2\Iaux_{3,0}
\\
\mathcal{D}_2\Iaux_{2,1}
\\
\mathcal{D}_2\Iaux_{4,0}
\\
\mathcal{D}_2\Iaux_{4,1}
\end{array}
\!\right)
\,\,=\,\,
\left(\!
\begin{array}{c}
1
\\
0
\\
\Waux
\\
2
\\
0
\\
6\,\Waux^2
\end{array}
\!\right)
+
\left(\!
\begin{array}{cccccc}
-3 & -1 & 0 & 0 & 0 & 0 
\\
0 & 0 & -1 & 0 & 0 & 0
\\
0 & 0 & 0 & -3 & -3 & 0
\\
-3 & -2 & 0 & 0 & 0 & 0
\\
0 & 0 & 0 & 0 & -4\,\Waux & -6
\\
0 & 0 & -\Maux & -10\,\Waux & -24\,\Waux & -18
\end{array}
\!\right)\,
\left(\!
\begin{array}{c}
\Kaux_2^1
\\
\Kaux_2^2
\\
\Kaux_2^3
\\
\Kaux_2^4
\\
\Kaux_2^5
\\
\Kaux_2^6
\end{array}
\!\right),
\]
and has the unique solution:
\[
\aligned
\Kaux_2^1
\,=\,
0,
\ \ \ \ \ \ \ \ \ \
\Kaux_2^2
\,=\,
1,
\ \ \ \ \ \ \ \ \ \
\Kaux_2^3
\,=\,
0,
\ \ \ \ \ \ \ \ \ \
\Kaux_2^4
&
\,=\,
-\,\Waux,
\ \ \ \ \ \ \ \ \ \
\Kaux_2^5
&
\,=\,
\frac{4}{3}\
\Waux,
\ \ \ \ \ \ \ \ \ \
\Kaux_2^6
\,=\,
-\,\frac{8}{9}\,
\Waux^2.
\endaligned
\]

For the three non-phantom invariants:
\[
\Iaux_{3,1}
\,=:\,
\Waux,
\ \ \ \ \ \ \ \ \ \ \ \ \ \ \ \ \ \ \ \
\Iaux_{5,0}
\,=:\,
\Maux,
\ \ \ \ \ \ \ \ \ \ \ \ \ \ \ \ \ \ \ \
\Iaux_{6,0},
\]
the recurrence formul{ae} are:
\[
\left(
\begin{array}{c}
\mathcal{D}_1\Waux
\\
\mathcal{D}_1\Maux
\\
\mathcal{D}_1\Iaux_{6,0}
\end{array}
\!\right)
\,\,=\,\,
\left(\!
\begin{array}{c}
0
\\
\Iaux_{6,0}
\\
\Iaux_{7,0}
\end{array}
\!\right)
+
\left(\!
\begin{array}{cccccc}
-4\Waux & -2\Waux & 0 & -6 & -6 & 0
\\
-6\Maux & -\Maux & 0 & 0 & 0 & -10\Waux
\\
-7\Iaux_{6,0} & -\Iaux_{6,0} & 0 & -21\Maux & -6\Iaux_{5,1} & 0
\end{array}
\!\right)\,
\left(\!
\begin{array}{c}
\Kaux_1^1
\\
\Kaux_1^2
\\
\Kaux_1^3
\\
\Kaux_1^4
\\
\Kaux_1^5
\\
\Kaux_1^6
\end{array}
\!\right),
\]
and:
\[
\left(
\begin{array}{c}
\mathcal{D}_2\Waux
\\
\mathcal{D}_2\Maux
\\
\mathcal{D}_2\Iaux_{6,0}
\end{array}
\!\right)
\,\,=\,\,
\left(\!
\begin{array}{c}
6\Waux
\\
\Iaux_{5,1}
\\
\Iaux_{6,1}
\end{array}
\!\right)
+
\left(\!
\begin{array}{cccccc}
-4\Waux & -2\Waux & 0 & -6 & -6 & 0
\\
-6\Maux & -\Maux & 0 & 0 & 0 & -10\Waux
\\
-7\Iaux_{6,0} & -\Iaux_{6,0} & 0 & -21\Maux & -6\Iaux_{5,1} & 0
\end{array}
\!\right)\,
\left(\!
\begin{array}{c}
\Kaux_2^1
\\
\Kaux_2^2
\\
\Kaux_2^3
\\
\Kaux_2^4
\\
\Kaux_2^5
\\
\Kaux_2^6
\end{array}
\!\right).
\]
So we obtain:
\[
\aligned
\mathcal{D}_1\Waux
&
\,=\,
-\,\frac{2\,\Waux^2}{3},
&
\ \ \ \ \ \ \ \ \ \ \ \ \ \ \ \ \ \ \ \
\mathcal{D}_2\Waux
&
\,=\,
2\,\Waux,
\\
\mathcal{D}_1\Maux
&
\,=\,
\Iaux_{6,0}
-
14\,\Maux\,\Waux
+
\frac{10}{3}\,
\Iaux_{5,1}\,\Waux,
&
\ \ \ \ \ \ \ \ \ \ \ \ \ \ \ \ \ \ \ \
\mathcal{D}_2\Maux
&
\,=\,
\Iaux_{5,1}
-
\Maux
+
\frac{80}{9}\,
\Waux^3,
\endaligned
\]
and:
\[
\aligned
\mathcal{D}_1\Iaux_{6,0}
&
\,=\,
\Iaux_{7,0}
-
\frac{3}{2\,\Waux}\,
\big(
7\,\Maux
-
2\,\Iaux_{5,1}
\big)\,
\big(
4\,\Maux
-
\Iaux_{5,1}
\big)
+
\frac{4}{3}\,\Waux\,\Iaux_{6,0}
\\
\mathcal{D}_2\Iaux_{6,0}
&
\,=\,
\Iaux_{6,1}
-
\Iaux_{6,0}
+
21\,\Waux\,\Maux
-
8\,\Waux\,\Iaux_{5,1}.
\endaligned
\]
Looking at these equations, we see that we can solve the 
$6$\textsuperscript{th} order invariants $\Iaux_{6,0}$
and $\Iaux_{5,1}$ in terms of $\Waux$, $\Maux$,
$\mathcal{D}_1 \Maux$, $\mathcal{D}_2 \Maux$.
We can also solve $\Iaux_{7,0}$, $\Iaux_{6,1}$ 
in terms of $\Waux$, $\Maux$ and their invariant derivatives.
An elementary induction yields:

\begin{Proposition}
Within the branch $\Saux \neq 0$, $\Waux \neq 0$, 
all invariants are generated by $\Waux$, $\Maux$
and their invariant derivatives.

Moreover, $\Maux$ cannot be obtained from $\Waux$ and
its invariant derivatives.
\end{Proposition}

\proof
Indeed, both $\mathcal{D}_1 \Waux = -\frac{2}{3}\, \Waux^2$ and
$\mathcal{D}_2 \Waux = 2\, \Waux$ do not raise the jet order.
\endproof

\Subsection{Branch $\Waux \equiv 0$}
\label{Subsection-W-zero}
Thus, we assume $\Saux \neq 0$ and $\Waux \equiv 0$.
Since the case $\Xaux \equiv 0$ has already been covered by
Theorem~{\ref{Thm-flat-cone}}, we may also assume $\Xaux \neq 0$.

According to 
Theorem~{\ref{Thm-S-nonzero-W-0}}, the normal form is:
\[
\!\!\!\!\!\!\!\!\!\!\!\!\!\!\!
\aligned
u
\,=\,
{\textstyle{\frac{x^2}{2}}}
+
{\textstyle{\frac{x^2\,y}{2}}}
+
{\textstyle{\frac{x^2\,y^2}{2}}}
+
F_{5,0}\,
\frac{x^5}{120}
+
{\textstyle{\frac{x^2\,y^3}{2}}}
+
4\,F_{5,0}\,
\frac{x^5y}{120}
+
{\textstyle{\frac{x^2\,y^4}{2}}}
+
F_{7,0}\,
\frac{x^7}{5\,040}
&
+
20\,F_{5,0}\,
\frac{x^5y^2}{240}
+
{\textstyle{\frac{x^2\,y^5}{2}}}
+
\\
&
+
\sum_{j+k\geqslant 8}\,
F_{j,k}\,
x^jy^k,
\endaligned
\]
and this means that we have $6$ phantom invariants:
\[
\Iaux_{2,0}
\,=\,
1,
\ \ \ \ \ \ \ \ \ \ \
\Iaux_{1,1}
\,=\,
0,
\ \ \ \ \ \ \ \ \ \ \
\Iaux_{3,0}
\,=\,
0,
\ \ \ \ \ \ \ \ \ \ \
\Iaux_{2,1}
\,=\,
1,
\ \ \ \ \ \ \ \ \ \ \
\Iaux_{4,0}
\,=\,
0,
\ \ \ \ \ \ \ \ \ \ \
\Iaux_{6,0}
\,=\,
0.
\]
We denote specially:
\[
\Xaux
\,:=\,
\Iaux_{5,0}
\ \ \ \ \ \ \ \ \ \ \ \ \ \ \ \ \ \ \ \
\text{and}
\ \ \ \ \ \ \ \ \ \ \ \ \ \ \ \ \ \ \ \
\Yaux
\,:=\,
\Iaux_{7,0},
\]
and we apply the first collection of recurrence formulas:
\[
\left(\!
\begin{array}{c}
0
\\
0
\\
0
\\
0
\\
0
\\
0
\end{array}
\!\right)
\,=\,
\left(\!
\begin{array}{c}
\mathcal{D}_1\Iaux_{2,0}
\\
\mathcal{D}_1\Iaux_{1,1}
\\
\mathcal{D}_1\Iaux_{3,0}
\\
\mathcal{D}_1\Iaux_{2,1}
\\
\mathcal{D}_1\Iaux_{4,0}
\\
\mathcal{D}_1\Iaux_{6,0}
\end{array}
\!\right)
\,\,=\,\,
\left(\!
\begin{array}{c}
0
\\
1
\\
0
\\
0
\\
\Xaux
\\
\Yaux
\end{array}
\!\right)
+
\left(\!
\begin{array}{cccccc}
-3 & -1 & 0 & 0 & 0 & 0 
\\
0 & 0 & -1 & 0 & 0 & 0
\\
0 & 0 & 0 & -3 & -3 & 0
\\
-3 & -2 & 0 & 0 & 0 & 0
\\
0 & 0 & 0 & 0 & 0 & -6
\\
0 & 0 & 0 & -21\,\Xaux & -24\,\Xaux & 0
\end{array}
\!\right)\,
\left(\!
\begin{array}{c}
\Kaux_1^1
\\
\Kaux_1^2
\\
\Kaux_1^3
\\
\Kaux_1^4
\\
\Kaux_1^5
\\
\Kaux_1^6
\end{array}
\!\right).
\]
This Cramér system has the unique solution:
\[
\Kaux_1^1
\,=\,
0,
\ \ \ \ \ \ \ \ \ \
\Kaux_1^2
\,=\,
0,
\ \ \ \ \ \ \ \ \ \
\Kaux_1^3
\,=\,
1,
\ \ \ \ \ \ \ \ \ \
\Kaux_1^4
\,=\,
-\,\frac{\Yaux}{3\,\Xaux},
\ \ \ \ \ \ \ \ \ \
\Kaux_1^5
\,=\,
\frac{\Yaux}{3\,\Xaux},
\ \ \ \ \ \ \ \ \ \
\Kaux_1^6
\,=\,
\frac{1}{6}\,
\Xaux.
\]

The second collection is:
\[
\left(\!
\begin{array}{c}
0
\\
0
\\
0
\\
0
\\
0
\\
0
\end{array}
\!\right)
\,=\,
\left(\!
\begin{array}{c}
\mathcal{D}_2\Iaux_{2,0}
\\
\mathcal{D}_2\Iaux_{1,1}
\\
\mathcal{D}_2\Iaux_{3,0}
\\
\mathcal{D}_2\Iaux_{2,1}
\\
\mathcal{D}_2\Iaux_{4,0}
\\
\mathcal{D}_2\Iaux_{6,0}
\end{array}
\!\right)
\,\,=\,\,
\left(\!
\begin{array}{c}
1
\\
0
\\
0
\\
2
\\
0
\\
0
\end{array}
\!\right)
+
\left(\!
\begin{array}{cccccc}
-3 & -1 & 0 & 0 & 0 & 0 
\\
0 & 0 & -1 & 0 & 0 & 0
\\
0 & 0 & 0 & -3 & -3 & 0
\\
-3 & -2 & 0 & 0 & 0 & 0
\\
0 & 0 & 0 & 0 & 0 & -6
\\
0 & 0 & 0 & -21\,\Xaux & -24\,\Xaux & 0
\end{array}
\!\right)\,
\left(\!
\begin{array}{c}
\Kaux_2^1
\\
\Kaux_2^2
\\
\Kaux_2^3
\\
\Kaux_2^4
\\
\Kaux_2^5
\\
\Kaux_2^6
\end{array}
\!\right),
\]
and has the unique solution:
\[
\aligned
\Kaux_2^1
\,=\,
0,
\ \ \ \ \ \ \ \ \ \
\Kaux_2^2
\,=\,
1,
\ \ \ \ \ \ \ \ \ \
\Kaux_2^3
\,=\,
0,
\ \ \ \ \ \ \ \ \ \
\Kaux_2^4
&
\,=\,
0,
\ \ \ \ \ \ \ \ \ \
\Kaux_2^5
&
\,=\,
0,
\ \ \ \ \ \ \ \ \ \
\Kaux_2^6
\,=\,
0.
\endaligned
\]

For the three non-phantom invariants:
\[
\Iaux_{5,0}
\,=:\,
\Xaux,
\ \ \ \ \ \ \ \ \ \ \ \ \ \ \ \ \ \ \ \
\Iaux_{7,0}
\,=:\,
\Yaux,
\ \ \ \ \ \ \ \ \ \ \ \ \ \ \ \ \ \ \ \
\Iaux_{8,0},
\]
the recurrence formul{ae} are:
\[
\left(
\begin{array}{c}
\mathcal{D}_1\Xaux
\\
\mathcal{D}_1\Yaux
\\
\mathcal{D}_1\Iaux_{8,0}
\end{array}
\!\right)
\,\,=\,\,
\left(\!
\begin{array}{c}
0
\\
\Iaux_{8,0}
\\
\Iaux_{9,0}
\end{array}
\!\right)
+
\left(\!
\begin{array}{cccccc}
-6\Xaux & -\Xaux & 0 & 0 & 0 & 0
\\
-8\Yaux & -\Yaux & 0 & 0 & 0 & -105\Xaux
\\
-9\Iaux_{8,0} & -\Iaux_{8,0} & 0 & -36\Yaux & -48\Yaux & 0
\end{array}
\!\right)\,
\left(\!
\begin{array}{c}
\Kaux_1^1
\\
\Kaux_1^2
\\
\Kaux_1^3
\\
\Kaux_1^4
\\
\Kaux_1^5
\\
\Kaux_1^6
\end{array}
\!\right),
\]
and:
\[
\left(
\begin{array}{c}
\mathcal{D}_2\Xaux
\\
\mathcal{D}_2\Yaux
\\
\mathcal{D}_2\Iaux_{8,0}
\end{array}
\!\right)
\,\,=\,\,
\left(\!
\begin{array}{c}
4\Xaux
\\
6\Yaux
\\
7\Iaux_{8,0}
\end{array}
\!\right)
+
\left(\!
\begin{array}{cccccc}
-6\Xaux & -\Xaux & 0 & 0 & 0 & 0
\\
-8\Yaux & -\Yaux & 0 & 0 & 0 & -105\Xaux
\\
-9\Iaux_{8,0} & -\Iaux_{8,0} & 0 & -36\Yaux & -48\Yaux & 0
\end{array}
\!\right)\,
\left(\!
\begin{array}{c}
\Kaux_2^1
\\
\Kaux_2^2
\\
\Kaux_2^3
\\
\Kaux_2^4
\\
\Kaux_2^5
\\
\Kaux_2^6
\end{array}
\!\right).
\]
So we obtain:
\[
\aligned
\mathcal{D}_1\Xaux
&
\,=\,
0,
&
\ \ \ \ \ \ \ \ \ \ \ \ \ \ \ \ \ \ \ \
\mathcal{D}_2\Xaux
&
\,=\,
3\,\Xaux,
\\
\mathcal{D}_1\Yaux
&
\,=\,
\Iaux_{8,0}
-
\frac{35}{2}\,
\Xaux^2,
&
\ \ \ \ \ \ \ \ \ \ \ \ \ \ \ \ \ \ \ \
\mathcal{D}_2\Yaux
&
\,=\,
5\,\Yaux,
\\
\mathcal{D}_1\Iaux_{8,0}
&
\,=\,
\Iaux_{9,0}
-
4\,\frac{\Yaux^2}{\Xaux}\,
,
&
\ \ \ \ \ \ \ \ \ \ \ \ \ \ \ \ \ \ \ \
\mathcal{D}_2\Iaux_{8,0}
&
\,=\,
6\,\Iaux_{8,0}.
\endaligned
\]
Looking at these equations, we see that we can solve the 
$8$\textsuperscript{th} order invariant $\Iaux_{8,0}$
in terms of $\Xaux$, $\Yaux$ and their invariant derivatives.
An elementary induction yields:

\begin{Proposition}
Within the branch $\Saux \neq 0$, $\Waux \equiv 0$,
the algebra of differential invariants 
is generated by  $\Xaux$, $\Yaux$ and their invariant derivatives.

Moreover, $\Yaux$ cannot be obtained from $\Xaux$ and its
invariant derivatives.
\end{Proposition}

\proof
Indeed, both
derivatives $\mathcal{D}_1 \Xaux = 0$ and $\mathcal{D}_2 \Xaux = 5\,
\Xaux$ do not raise the jet order.
\endproof

\Subsection{Commutators of invariant differentials} Besides taking invariant derivatives $\mathcal{D}_1$, $\mathcal{D}_2$, as in~{\cite{Olver-2007}},
there is another way to get syzygies among differential invariants: by means of the commutator
$\big[ \mathcal{D}_1, \mathcal{D}_2 \big]$.
For our group $\SA_3(\R)$, Olver in~{\cite{Olver-2007}} 
obtained the following formulas:
\[
\aligned
\mathcal{D}_3
\,:=\,
\big[
\mathcal{D}_1,\,
\mathcal{D}_2
\big]
\,:=\,
&\,
\mathcal{D}_1
\circ
\mathcal{D}_2
-
\mathcal{D}_2
\circ
\mathcal{D}_1
\\
\,=\,
&\,
\Zaux_1\,\mathcal{D}_1
+
\Zaux_2\,\mathcal{D}_2,
\endaligned
\]
with the two differential invariants:
\[
\aligned
\Zaux_1
\,:=\,
\sum_{1\leqslant\sigma\leqslant 6}\,
\Big(
\frac{\partial\xi_\sigma}{\partial x}
(0,0,0)\,
\Kaux_2^\sigma
-
\frac{\partial\xi_\sigma}{\partial y}
(0,0,0)\,
\Kaux_1^\sigma
\Big)
\,\,=\,\,
\Kaux_2^1
-
\Kaux_1^3,
\\
\Zaux_2
\,:=\,
\sum_{1\leqslant\sigma\leqslant 6}\,
\Big(
\frac{\partial\eta_\sigma}{\partial x}
(0,0,0)\,
\Kaux_2^\sigma
-
\frac{\partial\eta_\sigma}{\partial y}
(0,0,0)\,
\Kaux_1^\sigma
\Big)
\,\,=\,\,
\Kaux_2^5
-
\Kaux_1^2.
\endaligned
\]

Within the branch $\Saux \neq 0$, $\Waux \neq 0$, we have:
\[
\aligned
\Zaux_1
&
\,=\,
0
-
1
\,=\,
-\,1,
\\
\Zaux_2
&
\,=\,
\frac{4}{3}\,
\Waux
-
\Waux
\,=\,
\frac{1}{3}\,
\Waux.
\endaligned
\]
Hence $\big[ \mathcal{D}_1, \mathcal{D}_2 \big] = -\, \mathcal{D}_1 +
\frac{1}{3}\, \Waux\, \mathcal{D}_2$, an operator which can raise the
order by at most $1$.  The commutator does not generate anything other
than the invariant derivations $\mathcal{D}_1$, $\mathcal{D}_2$ would
do.  One can double-check this formula by calculating $\big[
  \mathcal{D}_1, \mathcal{D}_2 \big]\, \Waux = \frac{4}{3}\, \Waux^2 =
-\, \mathcal{D}_1 \Waux + \frac{1}{3}\, \Waux\, \mathcal{D}_2 \Waux$.

Note that when $\mathcal{D}_2 \Maux\neq0$, from
\[
\mathcal{D}_3\Maux=-\mathcal{D}_1\Maux+\frac{1}{3}\Waux\,\mathcal{D}_2\Maux,
\]
the invariant $\Waux$ is solved in terms of $\Maux$ and its invariant
differentials, so are all the other differential invariants; indeed,
we check that $\mathcal{D}_2 \Maux \not\equiv 0$ in general, with
numerator having $107$ differential monomials.

Within the branch $\Saux \neq 0$, $\Waux \equiv 0$, we have:
\[
\aligned
\Zaux_1
&
\,=\,
0-1
\,=\,
-\,1,
\\
\Zaux_2
&
\,=\,
0-0
\,=\,
0.
\endaligned
\]
Hence $\big[ \mathcal{D}_1, \mathcal{D}_2 \big] = -\, \mathcal{D}_1$,
an operator which can raise the order by at most $1$. The commutator
does not generate anything other than the invariant derivations
$\mathcal{D}_1$, $\mathcal{D}_2$ would do.  One can double-check this
formula by calculating $\big[ \mathcal{D}_1, \mathcal{D}_2 \big]\,
\Xaux = 0 = -\, \mathcal{D}_1\Xaux$.


\Section{\bf Explicit Invariant Differentiation Operators 
$\mathcal{D}_1$ and $\mathcal{D}_2$}
\label{explicit-D1-D2}
\HEAD{{\ref{explicit-D1-D2}}.~{\sf Explicit Invariant 
Differentiation Operators $\mathcal{D}_1$ and $\mathcal{D}_2$}
}{
Zhangchi {\sc Chen}, Jo\"el~{\sc Merker}}

In order to double-check the overall theoretical coherency of our
recurrence formulas satisfied by infinitely many differential
invariants, let us raise

\begin{Question}
{\sl 
How to make {\em explicit} the two invariant differentiation
operators $\mathcal{D}_1$ and $\mathcal{D}_2$?}
\end{Question}

To fix ideas, we will content ourselves to examine the
branch $\Waux \neq 0$. One strategy would be to follow the
general theory presented in the previous sections,  
but we already saw in 
Section~{\ref{search-resolved-cross-section-parabolic}}
that it seems impossible to directly compute in terms
of a cross-section, which we were unable
to make explicit.

Therefore, we will trace another route. 
At first, 
applying the power series method 
of~Sections~{\ref{parabolic-pseudostablization}}
and~{\ref{relative-invariant-S-first-invariant-W}},
we compute explicitly the differential invariants
$\Waux \equiv \Iaux_{3,1}$, $\Maux \equiv \Iaux_{5,0}$
(already shown)
and also $\Iaux_{6,0}$, $\Iaux_{5,1}$.
The numerator of $\Iaux_{6,0}$ has $225$ monomials,
that of $\Iaux_{5,1}$ (only) $69$.

In advance, on a computer, we declare some dependent parabolic jets:
\[
u_{0,2}
\,:=\,
\frac{u_{1,1}^2}{u_{2,0}},
\ \ \ \ \ \ \ \ \ \ \ \ \ \ \ \ \ \ \ \ \ \ \ \ \ \
u_{1,2}
\,:=\,
\frac{2\,u_{2,0}u_{1,1}u_{2,1}
-u_{1,1}^2u_{3,0}}{u_{2,0}^2},
\]
and so on, up to $u_{5,2}$ (those which will be useful for 
${\sf D}_y$ below), and we declare the
two total differentiation (not invariant) operators:
\[
\footnotesize
\aligned
{\sf D}_x
&
\,:=\,
\frac{\partial}{\partial x}
\\
&
\ \ \ \ \
+
u_{1,0}\,\frac{\partial}{\partial u}
\\
&
\ \ \ \ \
+
u_{2,0}\,\frac{\partial}{\partial u_{1,0}}
+
u_{1,1}\,\frac{\partial}{\partial u_{0,1}}
\\
&
\ \ \ \ \
+
u_{3,0}\,\frac{\partial}{\partial u_{2,0}}
+
u_{2,1}\,\frac{\partial}{\partial u_{1,1}}
\\
&
\ \ \ \ \
+
u_{4,0}\,\frac{\partial}{\partial u_{3,0}}
+
u_{3,1}\,\frac{\partial}{\partial u_{2,1}}
\\
&
\ \ \ \ \
+
u_{5,0}\,\frac{\partial}{\partial u_{4,0}}
+
u_{4,1}\,\frac{\partial}{\partial u_{3,1}}
\\
&
\ \ \ \ \
+
u_{6,0}\,\frac{\partial}{\partial u_{5,0}}
+
u_{5,1}\,\frac{\partial}{\partial u_{4,1}}
\\
&
\ \ \ \ \
+
u_{7,0}\,\frac{\partial}{\partial u_{6,0}}
+
u_{6,1}\,\frac{\partial}{\partial u_{5,1}}
\endaligned
\ \ \ \ \ \ \ \ \ \ \ \ \ \ \ \ \ \ \ \
\text{and}
\ \ \ \ \ \ \ \ \ \ \ \ \ \ \ \ \ \ \ \
\aligned
{\sf D}_y
&
\,:=\,
\frac{\partial}{\partial y}
\\
&
\ \ \ \ \
+
u_{0,1}\,\frac{\partial}{\partial u}
\\
&
\ \ \ \ \
+
u_{1,1}\,\frac{\partial}{\partial u_{1,0}}
+
u_{0,2}\,\frac{\partial}{\partial u_{0,1}}
\\
&
\ \ \ \ \
+
u_{2,1}\,\frac{\partial}{\partial u_{2,0}}
+
u_{1,2}\,\frac{\partial}{\partial u_{1,1}}
\\
&
\ \ \ \ \
+
u_{3,1}\,\frac{\partial}{\partial u_{3,0}}
+
u_{2,2}\,\frac{\partial}{\partial u_{2,1}}
\\
&
\ \ \ \ \
+
u_{4,1}\,\frac{\partial}{\partial u_{4,0}}
+
u_{3,2}\,\frac{\partial}{\partial u_{3,1}}
\\
&
\ \ \ \ \
+
u_{5,1}\,\frac{\partial}{\partial u_{5,0}}
+
u_{4,2}\,\frac{\partial}{\partial u_{4,1}}
\\
&
\ \ \ \ \
+
u_{6,1}\,\frac{\partial}{\partial u_{6,0}}
+
u_{5,2}\,\frac{\partial}{\partial u_{5,1}}.
\endaligned
\]

According to the general theory, there are coefficients
$\alpha$, $\beta$, $\gamma$, $\delta$ so that:
\[
\mathcal{D}_1
\,=\,
\alpha\,{\sf D}_x
+
\beta\,{\sf D}_y
\ \ \ \ \ \ \ \ \ \ \ \ \ \ \ \ \ \ \ \
\text{and}
\ \ \ \ \ \ \ \ \ \ \ \ \ \ \ \ \ \ \ \
\mathcal{D}_2
\,=\,
\gamma\,{\sf D}_x
+
\delta\,{\sf D}_y.
\]
Consequently, two appropriate pairs of recurrence relations seen in 
Section~{\ref{recurrence-relations-parabolic-surfaces}} read
as two linear systems satisfied by the two pairs
$\{\alpha, \beta\}$ and $\{\gamma, \delta\}$:
\[
\left[
\aligned
\alpha\,{\sf D}_x\Waux
+
\beta\,{\sf D}_y\Waux
&
\,=\,
-\,\tfrac{2}{3}\,\Waux^2,
\\
\alpha\,{\sf D}_x\Maux
+
\beta\,{\sf D}_y\Maux
&
\,=\,
\Iaux_{6,0}
-
14\,\Waux\Maux
+
\tfrac{10}{3}\,
\Waux\Iaux_{5,1},
\endaligned\right.
\]
and:
\[
\left[
\aligned
\gamma\,{\sf D}_x\Waux
+
\delta\,{\sf D}_y\Waux
&
\,=\,
2\,\Waux,
\\
\gamma\,{\sf D}_x\Maux
+
\delta\,{\sf D}_y\Maux
&
\,=\,
\Iaux_{5,1}
-
\Maux
+
\tfrac{80}{9}\,
\Waux^3.
\endaligned\right.
\]

Contrary to what could be expected/hoped,
the common determinant:
\[
\footnotesize
\aligned
\Delta
\,:=\,
\left\vert\!
\begin{array}{cc}
{\sf D}_x\Waux & {\sf D}_y\Waux
\\
{\sf D}_x\Maux & {\sf D}_y\Maux
\end{array}
\!\right\vert
\,=\,
-\,\frac{1}{54}\,
\frac{{\sf complicated\,\,numerator}}{
u_{2,0}^8
(u_{2,0}u_{2,1}-u_{1,1}u_{3,0})^{8/3}
\big(
u_{1,1}u_{2,0}u_{4,0}
-
u_{2,0}^2u_{3,1}
+
2u_{2,0}u_{2,1}u_{3,0}
-
2u_{1,1}u_{3,0}^2
\big)},
\endaligned
\]
is not simple, as its numerator is a non-factorizable 
homogeneous polynomial of degree $15$ in the jet variables
up to order $6$ having {\bf 431} monomials.

For some time, we were afraid
that this unpleasant numerator indicated there
existed a mistake in our recurrence formulas.

But if one really applies the Cramér formulas to the first
linear system:
\[
\alpha
\,:=\,
\frac{
\left\vert\!
\begin{array}{cc}
-\frac{2}{3}\Waux^2 & {\sf D}_y\Waux
\\
\Iaux_{6,0}-14\Waux\Maux+\frac{10}{3}\Waux\Iaux_{5,1} & {\sf D}_y\Maux
\end{array}
\!\right\vert
}{
\left\vert\!
\begin{array}{cc}
{\sf D}_x\Waux & {\sf D}_y\Waux
\\
{\sf D}_x\Maux & {\sf D}_y\Maux
\end{array}
\!\right\vert},
\ \ \ \ \ \ \ \ \ \ \ \ \ \ \ \ \ \ \ \
\beta
\,:=\,
\frac{
\left\vert\!
\begin{array}{cc}
{\sf D}_x\Waux & -\frac{2}{3}\Waux^2
\\
{\sf D}_x\Maux & \Iaux_{6,0}-14\Waux\Maux+\frac{10}{3}\Waux\Iaux_{5,1}
\end{array}
\!\right\vert
}{
\left\vert\!
\begin{array}{cc}
{\sf D}_x\Waux & {\sf D}_y\Waux
\\
{\sf D}_x\Maux & {\sf D}_y\Maux
\end{array}
\!\right\vert},
\]
one realizes that this large complicated
numerator is in fact an {\sl extraneous factor},
namely it cancels as {\em it appears both in numerator
and in denominator places for both $\alpha$ and $\beta$}.

After clearing this factor and cleaning, we receive:
\[
\aligned
\alpha
&
\,:=\,
\frac{1}{6}\,
\,
\frac{1}{
u_{2,0}\big(u_{2,0}u_{2,1}-u_{1,1}u_{3,0}\big)^{2/3}
\big(
u_{1,1}u_{2,0}u_{4,0}
-
u_{2,0}^2u_{3,1}
+
2u_{2,0}u_{2,1}u_{3,0}
-
2u_{1,1}u_{3,0}^2
\big)}
\Big\{
\\
&
\ \ \ \ \ \ \ \ \ \
\Big\{
12\,u_{3,0}u_{2,1}^2u_{2,0}^2
-
6\,u_{3,1}u_{2,1}u_{2,0}^3
-
44\,u_{3,0}^2u_{1,1}u_{2,1}u_{2,0}
+
16\,u_{3,0}u_{3,1}u_{1,1}u_{2,0}^2
+
15\,u_{4,0}u_{1,1}u_{2,1}u_{2,0}^2
\\
&
\ \ \ \ \ \ \ \ \ \
-\,3\,u_{1,1}u_{4,1}u_{2,0}^3
+
32\,u_{3,0}^3u_{1,1}^2
-
25\,u_{4,0}u_{1,1}^2u_{3,0}u_{2,0}
+
3\,u_{5,0}u_{1,1}^2u_{2,0}^2
\Big\},
\endaligned
\]
and:
\[
\aligned
\beta
&
\,:=\,
\frac{1}{6}\,
\,
\frac{1}{
\big(u_{2,0}u_{2,1}-u_{1,1}u_{3,0}\big)^{2/3}
\big(
u_{1,1}u_{2,0}u_{4,0}
-
u_{2,0}^2u_{3,1}
+
2u_{2,0}u_{2,1}u_{3,0}
-
2u_{1,1}u_{3,0}^2
\big)}
\Big\{
\\
&
\ \ \ \ \ \ \ \ \ \
\Big\{
20\,u_{2,0}u_{2,1}u_{3,0}^2
-
10\,u_{3,0}u_{2,0}^2u_{3,1}
-
9\,u_{2,0}^2u_{2,1}u_{4,0}
+
3\,u_{4,1}u_{2,0}^3
-
20\,u_{1,1}u_{3,0}^3
\\
&
\ \ \ \ \ \ \ \ \ \ \ \
+
19\,u_{3,0}u_{1,1}u_{2,0}u_{4,0}
-
3\,u_{1,1}u_{2,0}^2u_{5,0}
\Big\}.
\endaligned
\]
The same extraneous factor from $\Delta$ also disappears
from Cramér's formulas in the second linear system,
and we receive rather neat and simpler expressions for:
\[
\aligned
\gamma
\,:=\,
-\,\frac{u_{2,0}u_{1,1}}{
u_{2,0}u_{2,1}-u_{1,1}u_{3,0}}
\endaligned
\]
and for:
\[
\aligned
\delta
\,:=\,
\frac{u_{2,0}^2}{
u_{2,0}u_{2,1}-u_{1,1}u_{3,0}}.
\endaligned
\]

A computation of the determinant:
\[
\left\vert\!
\begin{array}{cc}
\alpha & \beta
\\
\gamma & \delta
\end{array}
\!\right\vert
\,=\,
\frac{u_{2,0}}{
\big(u_{2,0}u_{2,1}-u_{1,1}u_{3,0}
\big)^{2/3}}.
\] 
confirms a property stated by the general theory:
\[
\Span\,
\big\{
\mathcal{D}_1,\mathcal{D}_2
\big\}
\,=\,
\Span\,
\big\{
{\sf D}_x,{\sf D}_y
\big\}.
\]

Lastly, with these explicit formulas for $\mathcal{D}_1$ and
$\mathcal{D}_2$, we verify on a computer, after
computing independently also, say,
$\Iaux_{7,0}$, $\Iaux_{6,1}$ (and so on), that the next few
recurrence formulas hold identically.

In conclusion, this confirms {\em explicit coherency} 
of all our formulas.

\Section{\bf Relation with the classification of developable surfaces: 
\\
${\sf SA}_3(\mathbb{R})$-invariant PDEs for cylinders and cones}
\label{relation-developable-surfaces}
\HEAD{{\ref{relation-developable-surfaces}}.~{\sf 
Relation with the classification of developable surfaces}
}{
Zhangchi {\sc Chen}, Jo\"el~{\sc Merker}}

A surface in $\mathbb{R}^3$ is called a {\em ruled} surface if it can be parametrized by a family of lines (rulers):
\[
\vec{x}(t,v)=\vec{\alpha}(t)+v\vec{w}(t)
\ \ \ \ \
\quad{\scriptstyle{(t\in(0,1),\,v\in\mathbb{R})}},
\]
where $\vec{\alpha}(t),\vec{w}(t)$ are $C^1$-smooth. It is called {\em developable} if
\[
(\vec{w},\vec{w'},\vec{\alpha'})\equiv0.
\]
Properties of ruled surfaces and developable surfaces can be found in \cite{Guggenheimer-1977, Do-Carmo-2016}. For instance, a $C^2$-smooth surface is developable if and only if its Gaussian curvature is identically 0. In this section, we will assume as before that geometric objects are analytic.

Parabolic surfaces are developable since their Gaussian curvature is constantly zero. A developable surface is called:
\begin{itemize}
\item {\sl cylindrical} if all rules are parallel;
\item {\sl conical} if all rulers pass through the same point;
\item {\sl tangential} if all rulers are tangent to a certain curve in $\mathbb{R}^3$.
\end{itemize}
Near a $C^2$-smooth point, a developable surface is locally cylindrical or conical or tangential~\cite[p.~197]{Do-Carmo-2016}.

\begin{figure}[ht] 
  \begin{minipage}[b]{0.5\linewidth}
    \centering
    \includegraphics[width=.9\linewidth]{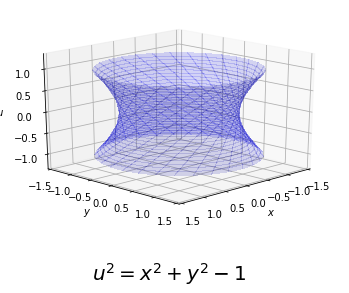} 
    \caption{A non-developable ruled surface} 
  \end{minipage}
  \begin{minipage}[b]{0.5\linewidth}
    \centering
    \includegraphics[width=.9\linewidth]{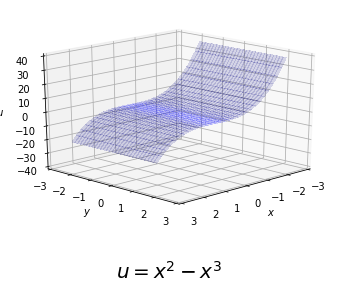} 
    \caption{A cylindrical surface} 
  \end{minipage} 
  \begin{minipage}[b]{0.5\linewidth}
    \centering
    \includegraphics[width=.9\linewidth]{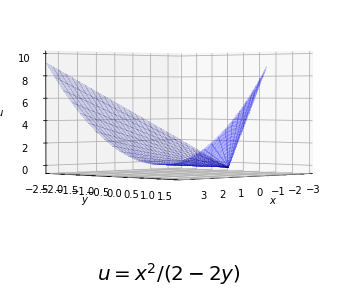} 
    \caption{A conical surface}
  \end{minipage}
  \begin{minipage}[b]{0.5\linewidth}
    \centering
    \includegraphics[width=.9\linewidth]{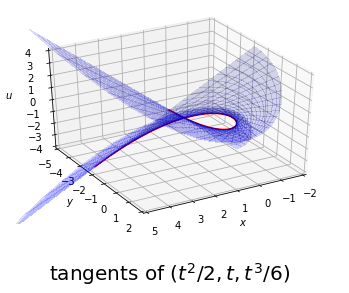} 
    \caption{A tangential surface} 
  \end{minipage} 
\end{figure}

Being locally cylindrical (or conical, or tangential) is clearly a ${\sf SA}_3(\mathbb{R})$-invariant local property. It is expected that such properties can be characterized by some differential ${\sf SA}_3(\mathbb{R})$-invariants. Indeed
\begin{Theorem} A parabolic surface is
\begin{itemize}
\item a cylinder if and only if $\Saux \equiv0$,
\item a cone if and only if $\Saux \neq0$ and $\Waux \equiv 0$,
\item a tangential surface if and only if $\Saux\neq0$ and $\Waux \neq0$.
\end{itemize}
\end{Theorem}

In~{\cite[p.~295]{Guggenheimer-1977}}, the degenerate branches of
cylinders and cones are characterized by the vanishing of some
coefficients in the Cartan's structure equations.

\proof It suffices to prove that at any smooth point of
\begin{enumerate}
\item[{\bf (1)}] any cylinder, $\Saux\equiv 0$,
\item[{\bf (2)}] any cone, $\Saux\neq0$ and $\Waux\equiv0$,
\item[{\bf (3)}] any tangential surface, $\Saux\neq0$ and $\Waux\neq0$.
\end{enumerate}

{\bf (1)} Any cylinder passing by a point $\vec{p}\in\mathbb{R}^3$, after some ${\sf SA}_3(\mathbb{R})$ action, can be viewed as $\mathbb{R}_y$ times a curve in $\mathbb{R}^2_{x,u}$, while $\vec{p}$ is mapped to the origin with the tangent plane $\mathbb{R}^2_{x,y}$. The cylinder is then a graph
\[
u=F(x)=\frac{F_{2,0}}{2}x^2+\frac{F_{3,0}}{6}x^3+\cdots.
\]
We calculate $\Saux=\frac{F_{xx}F_{xxy}-F_{xy}F_{xxx}}{F_{xx}^2}=0$ at the origin. Since $\vec{p}$ is arbitrarily chosen, $\Saux\equiv0$.

\smallskip

{\bf (2)} A cone can be parameterized by
\[
\vec{x}(t,v)=(1-v)\,\vec{\alpha}(t)+v\,\vec{p},
\]
where $v<1,t\in(-1,1)$, $\vec{\alpha}(t)$ parametrizes a smooth directrix and $A:=\vec{p}$ is the apex.

For any smooth marked point $B$ on the cone, we apply the following three steps of ${\sf SA}_3(\mathbb{R})$-actions.

First, we translate $B$ to the origin.

Second, we fix the origin and rotate the cone so that its tangent plane at $B$ is spanned by $\vec{e}_x:=(1,0,0)$ and $\vec{e}_y:=(0,1,0)$, while the generatrix $\vec{BA}$ is parallel to $\vec{e_y}$.

Finally, we apply a dilation of the type $x'=\lambda x$, $y'=\lambda^{-1} y$, $u'=u$ to make sure $|BA|=1$, i.e. $A=(0,1,0)$.

The new cone has the apex $(0,1,0)$ and the marked point $(0,0,0)$. By intersecting this new cone with $\mathbb{R}^2_{x,u} = \{y=0\}$, we get another directrix $\big(t,0,c(t)\big)$ passing the origin. Thus, our cone is ${\sf SA}_3(\mathbb{R})$-equivalent to
\[
\aligned
x(t,v)&=(1-v)\,t,\\
y(t,v)&=v, \\
u(t,v)&=(1-v)\,c(t),
\endaligned
\]
where $c(t)$ is smooth (analytic) and $c(0)=c'(0)=0$. Since the marked point $B$ is arbitrarily chosen, it suffices to check that for the new cone above, the invariant $\Saux\neq0$ and $\Waux=0$ at the origin.

\begin{figure}[ht] 
    \centering
    \includegraphics[width=.5\linewidth]{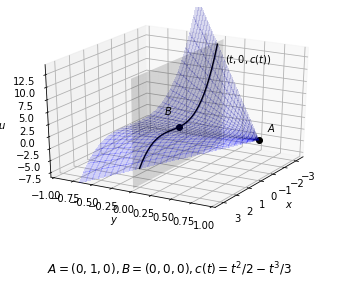} 
    \caption{A ${\sf SA}_3(\mathbb{R})$-normalized cone} 
\end{figure}

Suppose such a conical surface is a graph of $u=F(x,y)$ around the origin. Assuming $c(t)$, $F(x,y)$ analytic, expand
\[
\aligned
c(t)&=\frac{c_{2}}{2}t^2 +\frac{c_{3}}{6}t^3+\cdots,\\
F(x,y)&=F_{1,0} x+ F_{0,1} y +\frac{F_{2,0}}{2}x^2+F_{1,1}x\,y+\frac{F_{0,2}}{2}y^2+\cdots.
\endaligned
\]
In the fundamental identity holding for $t,v$ near $0$:
\[
u(t,v)\equiv F\big(x(t,v),y(t,v)\big),
\]
the Taylor coefficients of all monomials $t^jv^k$ should be the same. Identifying and solving, we get
\[
\aligned
F_{1,0}&=0, & F_{0,1}&=0, & & & & & &\\
F_{2,0}&=c_2, & F_{1,1}&=0, & F_{0,2}&=0, & & & &\\
F_{3,0}&=c_3, & F_{2,1}&=c_2, & F_{1,2}&=0, & F_{0,3}&=0, & &\\
F_{4,0}&=c_4, & F_{3,1}&=2\,c_3, & F_{2,2}&=2\,c_2, & F_{1,3}&=0, & F_{0,4}&=0.
\endaligned
\]
We may then compute
\[
\aligned
\Saux
&
:=
\frac{F_{2,0}\,F_{2,1}-F_{1,1}\,F_{3,0}}{F_{2,0}^2},
\\
\Waux
&
:=
\frac{2F_{1,1}\,F_{3,0}^2+F_{2,0}^2\,F_{3,1}-2F_{2,0}\,F_{2,1}\,F_{3,0}-F_{2,0}\,F_{1,1}\,F_{4,0}}{F_{2,0}^2\,(F_{2,0}\,F_{2,1}-F_{1,1}\,F_{3,0})^{2/3}},
\endaligned
\]
to get $\Saux=1$ and $\Waux=0$.

\smallskip\noindent
{\em Remark.} In the branch $\Waux\equiv0$, there are $2$ generators of ${\sf SA}_3(\mathbb{R})$-invariants: $\Xaux$ of order $5$ and $\Yaux$ of order $7$. Indeed at the origin of the normalized cone
\[
\aligned
\Xaux=\frac{40 c_3^3-45 c_2\, c_3\, c_4+9 c_2^2 \,c_5}{9c_2^4}
\endaligned
\]
is the Monge invariant. The model $u=\frac{x^2}{2(1-y)}$ is a cone with apex $(0,1,0)$ and directrix $(t,0,t^2/2)$. One can verify that $\Waux=0$ and $\Xaux=0$ at the origin.

\smallskip

{\bf (3)} A tangential surface can be parametrized by
\[
\vec{x}(t,v)=\alpha(t)+v\,\alpha'(t),
\]
where $v\in\mathbb{R},t\in(-1,1)$ and $\vec{\alpha}(t)$ parametrizes a smooth (analytic) directrix.

For any smooth point $B=\vec{x}(t_0,v_0)$ on the surface, let $A=\vec{x}(t_0,0)$ be the corresponding point on the directrix. There is a ${\sf SA}_3(\mathbb{R})$-action sending $A$ to $(0,-1,0)$ and $B$ to $(0,0,0)$. The original surface is sent to tangents of a curve passing by $(0,-1,0)$ with tangent direction $\vec{e_y}:=(0,1,0)$. The curve can be locally reparametrized as $\big(a(t),-1+t,c(t)\big)$. Thus our tangential surface is ${\sf SA}_3(\mathbb{R})$-equivalent to
\[
\aligned
x(t,v)&=a(t)+v\,a'(t),\\
y(t,v)&=-1+t+v, \\
u(t,v)&=c(t)+v\,c'(t),\\
\endaligned
\]
where $v\in\mathbb{R}$, $t\in(-1,1)$, $a(t)$ and $c(t)$ are analytic with $a(0)=c(0)=a'(0)=c'(0)=0$. Note that $\vec{x}(0,1)=(0,0,0)$ is the marked point $B$.

\begin{figure}[ht] 
    \centering
    \includegraphics[width=.5\linewidth]{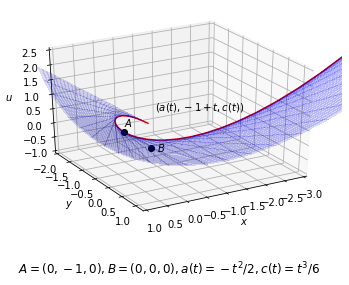} 
    \caption{A ${\sf SA}_3(\mathbb{R})$-normalized tangential surface} 
\end{figure}

We may rotate the $(x,u)$-space while fixing the $y$-axis to make sure the surface is not vertical at the origin. Then the surface is graphed as $u=F(x,y)$. Assuming $a(t)$, $c(t)$, $F(x,y)$ analytic in $t$, we have
\[
\aligned
a(t)&=\frac{a_{2}}{2}t^2 +\frac{a_{3}}{6}t^3+\cdots,\\
c(t)&=\frac{c_{2}}{2}t^2 +\frac{c_{3}}{6}t^3+\cdots,\\
F(x,y)&=F_{1,0} x+ F_{0,1} y +\frac{F_{2,0}}{2}x^2+F_{1,1}x\,y+\frac{F_{0,2}}{2}y^2+\cdots.
\endaligned
\]
By expanding the fundamental identity holding for $(t,v)$ near $(0,1)$
\[
u(t,v)\equiv F\big(x(t,v),y(t,v)\big),
\]
by identifying Taylor coefficients of 
monomials $t^j(v-1)^k$, and by solving, we obtain
\[
\aligned
F_{1,0}&=\frac{c_2}{a_2}, & F_{0,1}&=0, & & \\
F_{2,0}&=\frac{-a_3\,c_2+a_2\,c_3}{a_2^3}, & F_{1,1}&=0, & F_{0,2}&=0,
\endaligned
\]
\[
\aligned
F_{3,0}&=\frac{-a_2\,a_3\,c_2+3a_3^2\,c_2-a_2\,a_4\,c_2+a_2^2\,c_3-3a_2\,a_3\,c_3+a_2^2\,c_4}{a_2^5},\\
F_{2,1}&=\frac{a_3\,c_2-a_2\,c_3}{a_2^3},\ \ \ \ F_{1,2}=0, \ \ \ \ F_{0,3}=0,
\endaligned
\]
\[
\aligned
F_{4,0}&=\frac{1}{a_2^7}\big(-3a_2^2\,a_3\,c_2+10a_2\,a_3^2\,c_2-15a_3^3\,c_2-3a_2^2\,a_4\,c_2+10a_2\,a_3\,a_4\,c_2+a_2^2\,a_5\,c_2+3a_2^3\,c_3\\
& \ \ \ \ \ \ \ \ -10a_2^2\,a_3\,c_3+15a_2\,a_3^2\,c_3-4a_2^2\,a_4\,c_3+3a_2^3\,c_4-6a_2^2\,a_3\,c_4+a_2^3\,c_5\big),\\
F_{3,1}&=\frac{1}{a_2^5}\big(3\,a_2\,a_3\,c_2-6a_3^2\,c_2+2a_2\,a_4\,c_2-3a_2^2\,c_3+6a_2\,a_3\,c_4+a_2^3\,c_5\big),\\
F_{2,2}&=\frac{2(-a_3\,c_2+a_2\,c_3)}{a_2^3},\ \ \ \ F_{1,3}=0,\ \ \ \ F_{0,4}=0.
\endaligned
\]
A computation conducts to the compact expression
\[
\Waux=(a_3\,c_2-a_2\,c_3)^{-1/3}.
\]

On the other hand, observe that 
for the curve $\vec{\alpha}(t)=\big(a(t),-1+t,c(t)\big)$, the torsion is
\[
\tau(t)=\frac{1+a'(t)^2+c'(t)^2}{\big(a''(t)^2+c''(t)^2\big)^2}\big(c''(t)\,a'''(t)-a''(t)\,c'''(t)\big).
\]
We shall exclude the case $\tau(t) \equiv 0$,
since in this degenerate case, the curve is locally planar,
whence the tangential surface is locally flat, contradicting our
assumption $F_{xx} \neq 0$.
If $t=0$ is an isolated zero of $\tau(t)$, then
the tangential surface has a cuspidal edge along $\vec{x}(0,v)$ \cite{Cleave-1979}, contradicting our overall
assumption that the surface is smooth (and analytic) at $\vec{x}(0,1)$.
Thus $\tau(0)\neq 0$ necessarily, and this guarantees that
$c_2\,a_3-c_3\,a_2\neq 0$. In conclusion, $\Waux\neq 0$.
\endproof

\begin{figure}[ht] 
    \centering
    \includegraphics[width=.5\linewidth]{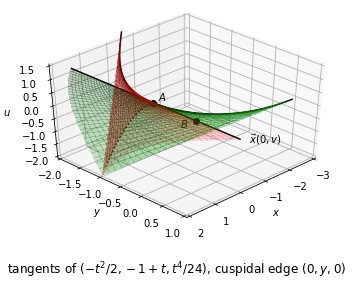} 
    \caption{Tangential surface of a curve with an isolated zero torsion point at $A=(0,-1,0)$.} 
\end{figure}

\Section{\bf Special Affinely Homogeneous Models}
\label{SA-3-homogeneous-models}
\HEAD{{\ref{SA-3-homogeneous-models}}.~{\sf Special Affinely 
Homogeneous Models}
}{
Zhangchi {\sc Chen}, Jo\"el~{\sc Merker}}

The complete classification of $A_3(\R)$-homogeneous surfaces $S^2
\subset \R^3$ was terminated by
Doubrov-Komrakov-Rabinovich~{\cite{Doubrov-Komrakov-Rabinovich-1996}},
who solved a long-standing open problem. Also,
Abdalla-Dillen-Vrancken~{\cite{Abdalla-Dillen-Vrancken-1997}} finished
the delicate classification of affinely homogeneous surfaces in $\R^3$
having {\em vanishing Pick invariant}. This full classification was done again later by Eastwood-Ezhov~{\cite{Eastwood-Ezhov-1999}},
who employed the power series method, without considering algebras of
differential invariants. In the literature, we did not find an answer
to the simple

\begin{Question}
{\sl
What are the {\rm special} affinely homogeneous surfaces
$S^2 \subset \R^3$?}
\end{Question}

Because any differential invariant $\Iaux$ satisfies, according
to Sections~{\ref{moving-frame-method}}
and~{\ref{meaning-differential-invariant}}:
\[
\Iaux\big(g\cdot z^{(n)}\big)
\,=\,
\Iaux\big(z^{(n)}\big),
\]
when the horizontal group action, namely the projection
of $g \cdot z^{(n)}$:
\[
\big(x,{\tt u}(x)\big)
\,\,\longmapsto\,\,
\Big(
\varphi\big(g,x,{\tt u}(x)\big),\,
\psi\big(g,x,{\tt u}(x)\big)
\Big),
\]
is transitive, this forces $\Iaux$ to be constant.

\begin{Observation}
\label{Obs-I-constant-homogeneous}
All differential invariants are constant for
geometric objects that are homogeneous.\qed
\end{Observation}

Consequently, {\em all} invariant derivatives of any order
$\geqslant 1$ are trivially zero:
\[
\mathcal{D}_1^{j_1}\cdots\mathcal{D}_p^{j_p}
\Iaux
\,\equiv\,
0.
\]
This yields interesting simplications in the recurrence relations
of Theorem~{\ref{Thm-recurrence-formulae}}, that are
better rewritten as:
\[
\Iaux_{K,j}^\alpha
\,=\,
\zero{\mathcal{D}_j
\Iaux_K^\alpha}
-
\sum_{1\leqslant\sigma\leqslant r}\,
{\varphi_\sigma}_K^\alpha
\big(
\Iaux^{(K)}
\big)
\cdot
\Kaux_j^\sigma
\big(
\Iaux^{(n_\GG+1)}
\big).
\]

Generally, the power series of the graphing functions
$u^\alpha = F^\alpha(x^1, \dots, x^p)$, $\alpha = 1, \dots, q$,
can be fully determined from these recurrence relations,
when the group action is transitive.
We now illustrate this general idea in our elementary context.

\[
\xymatrix{
&&
&&
\Paux\,\equiv\,0
&&
\Caux\,\equiv\,0
&&
\\
&&
\ar[urr]
\Saux\,\equiv\,0
\ar[rr]
&&
\ar[urr]
\Paux\,\neq\,0
\ar[rr]
&&
\Caux\,\neq\,0
&&
\\
\ar[urr]
\ar[rr]
\boxed{\root}
&&
\ar[rr]
\Saux\,\neq\,0
\ar[ddrr]
&&
\ar[rr]
\Waux\,\equiv\,0
\ar[drr]
&&
\Xaux\,\equiv\,0
&&
\\
&&
&&
&&
\ar[rr]
\Xaux\,\neq\,0
\ar[drr]
&&
\Yaux\,\equiv\,0
\\
&&
&&
\ar[rr]
\Waux\,\neq\,0
\ar[drr]
&&
\Maux\,\equiv\,0
&&
\Yaux\,\neq\,0
\\
&&
&&
&&
\Maux\,\neq\,0
&&
}
\]

Coming back to the complete branching diagram of
Section~{\ref{presentation-results}} copied above, 
we remember that the branch
$\Saux \equiv 0$ corresponds to {\sl cylindrical surfaces}, namely
surfaces that are the product of a curve $C^1 \subset \R_{x,u}$ with
$\R_y^1$. We also observe that ${\sf SA}_3(\R)$-equivalence classes of
such cylindrical surfaces are in one-to-one correspondence with ${\sf
A}_2(\R)$-equivalence classes of the corresponding curves, because
one can always use appropriate dilations along the dumb axes $\R_y^1$
to adjust preservation of the volume, so that the volume-preserving
condition in $\R^3$ is {\em not} transmitted as an aera-preserving
condition in $\R^2$.

Thus, let us discuss how ${\sf A}_2(\R)$-homogeneous curves can be
determined by employing the recurrence relations for differential
invariants shown in Section~{\ref{affine-invariants-curves-R-2}}.

\smallskip\noindent$\bullet$\,
Clearly, the relative invariant $\Iaux_2 := F_{xx}$
vanishes identically if and only if the curve is affinely equivalent
to the straight line $\{ u = 0\}$.

\smallskip\noindent$\bullet$\,
On the branch $\Iaux_2 \neq 0$, there is the {\em relative} 
differential invariant $\Iaux_4 := \frac{1}{3}\,
\frac{-5F_{xxx}^2+3F_{xx}F_{xxxx}}{F_{xx}^2}$.
Then Lemma~{\ref{LM-vanishing-H-M}} {\bf (1)}
already showed that $\Iaux_4 \equiv 0$ if and only if
the curve is ${\sf A}_2(\R)$-equivalent to the parabola
$\{u = x^2\}$.

\begin{center}
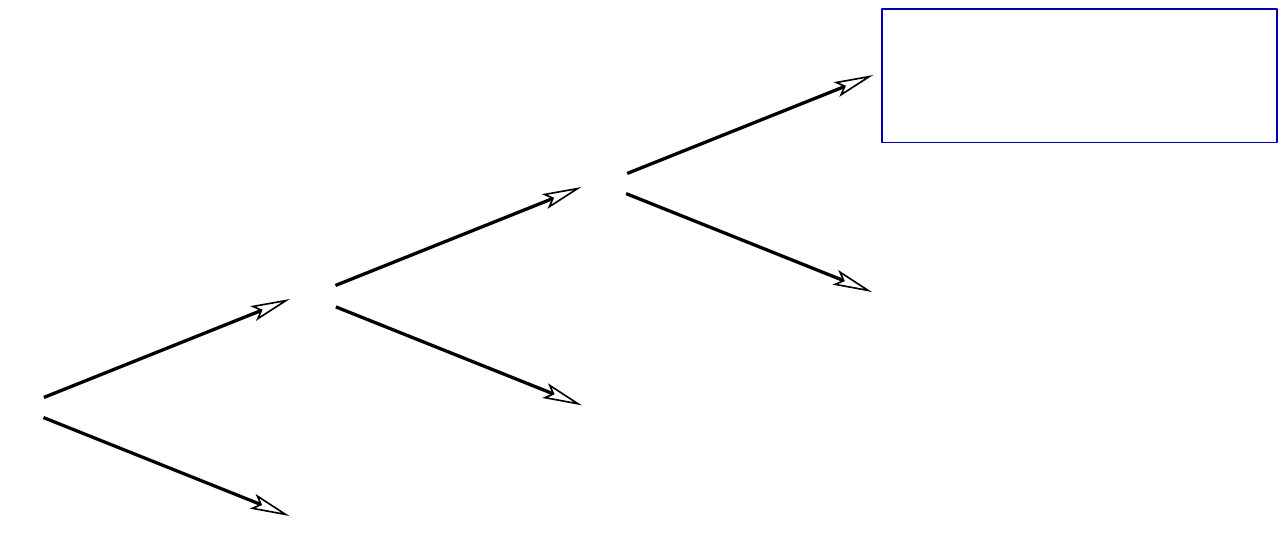
\end{center}

\smallskip\noindent$\bullet$\,
On the branch $\Iaux_2 \neq 0 \neq \Iaux_4$, 
with $\varepsilon := \pm 1$ denoting the {\em sign}
of $\Iaux_4$, 
there comes the
first {\em absolute} ({\em i.e.} not relative)
differential invariant:
\[
I_5
\,:=\,
\frac{1}{\sqrt{3}}\,
\frac{9\,F_{xx}^2
F_{xxxxx}
-
45\,F_{xx}\,
F_{xxx}\,F_{xxxx}
+
40\,F_{xxx}^3}{
(\varepsilon\,
[3\,F_{xx}\,F_{xxxx}
-
5\,F_{xxx}^2])^{3/2}}.
\]
Then Lemma~{\ref{LM-vanishing-H-M}} {\bf (2)} already
showed that $\{ u = F(x)\}$ is contained in a nondegenerate
conic (hence not a parabola) if and only if $\Iaux_5 \equiv 0$.
One can verify that only two normalized equations occur:
\[
u
\,=\,
\varepsilon\,
\sqrt{1+\varepsilon\,x^2}
-
\varepsilon,
\]
a circle for $\varepsilon = 1$, and a hyperbola for $\varepsilon = -
1$. Both are well known to be affinely homogeneous.

\smallskip\noindent$\bullet$\,
Lastly, it remains to consider the branch $\Iaux_5 \neq 0$.
Then according to Observation~{\ref{Obs-I-constant-homogeneous}}, 
for homogeneity to hold,
$\Iaux_5 \equiv \aaux \in \R \backslash \{0\}$ must
necessarily be constant. Furthermore, 
the recurrence relations written at the end
of Section~{\ref{affine-invariants-curves-R-2}}
become:
\[
\aligned
\Iaux_6
&
\,=\,
\zero{\mathcal{D}_x(\Iaux_5)}
\pm
{\textstyle{\frac{3}{2}}}
\,
\Iaux_5^2
+
5,
&
&
\,=\,
\pm\tfrac{3}{2}\,\aaux^2
+
5,
\\
\Iaux_7
&
\,=\,
\zero{\mathcal{D}_x(\Iaux_6)}
\pm
2\,\Iaux_5\Iaux_6
\pm
7\,\Iaux_5
&
&
\,=\,
3\,\aaux^3
\pm
17\,\aaux.
\endaligned
\]
Beyond, we have for every $k \geqslant 5$:
\[
\Iaux_{k+1}
\,=\,
\zero{\mathcal{D}_x\big(\Iaux_k\big)}
-
\sum_{1\leqslant\kappa\leqslant 4}\,
\inv\big(
\Phi_\kappa^k
\big)\,
\Raux^\kappa,
\]
hence knowing the values of the $\Raux^\kappa$, computing the values of
$\inv\big( \Phi_\kappa^k \big)$, we realize that {\em all} $\Iaux_{k
\geqslant 5}$ are uniquely determined polynomials in $\aaux$.
Thus, the power series reads:
\begin{equation}
\label{hom-curve-a-param}
u
\,=\,
\frac{x^2}{2!}
\pm
\frac{x^4}{4!}
+
\aaux\,
\frac{x^5}{5!}
+
\Big(
5
\pm
\frac{3}{2}\,\aaux^2
\Big)\,
\frac{x^6}{6!}
+
\Big(
3\,\aaux^3
\pm
17\,\aaux
\Big)\,
\frac{x^7}{7!}
+\cdots.
\end{equation}
By employing the Cauchy majorant method, one could prove
that this power series has a radius of convergence
$> 0$ for any parameter $\aaux$, but this would be tedious.

Another more elementary and straightforward method
is to test whether an affine infinitesimal transformation
with $6$ unknowns:
\[
L
\,:=\,
\big(
\AA\,x+\BB\,u+\CC
\big)\,
\tfrac{\partial}{\partial x}
+
\big(
\EE\,x+\FF\,u+\GG
\big)\,
\tfrac{\partial}{\partial u},
\]
is tangent to the above graph, and to realize by examining
Taylor coefficients only up to order $5$, that one 
comes to the single solution\,\,---\,\,up to dilation\,\,---:
\[
L
\,:=\,
\big(
\pm 1
-
\tfrac{1}{2}\,\aaux\,x
-
\tfrac{1}{3}\,u
\big)\,
\tfrac{\partial}{\partial x}
+
\big(
\pm x
-
\aaux\,u
\big)\,
\tfrac{\partial}{\partial u}.
\]
Furthermore, one can verify that this vector field is indeed tangent
to the graph~{\eqref{hom-curve-a-param}} truncated to any order, and
that the recurrence relations are {\em equivalent} to such a tangency
condition.

Therefore, the homogeneous curve can simply be taken as the {\em
orbit} of the origin $0 \in \R^2$ by the flow of the vector field
$L$, and since $L$ is analytic with $L(0) \neq 0$, its flow and its
local orbits are analytic too.

\begin{figure}[ht] 
  \begin{minipage}[b]{0.5\linewidth}
    \centering
    \includegraphics[width=.9\linewidth]{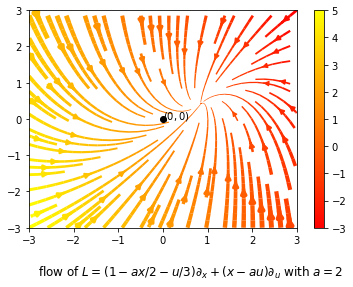} 
  \end{minipage}
  \begin{minipage}[b]{0.5\linewidth}
    \centering
    \includegraphics[width=.9\linewidth]{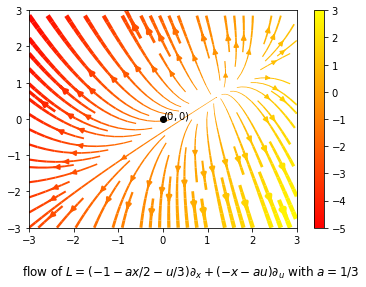} 
  \end{minipage} 
  \caption{Flow of $L$}
\end{figure}

By construction, different values of the parameter $\aaux \in \R
\backslash \{0\}$ conduct to affinely {\em inequivalent} homogeneous
curves.  This is the main interest of the use of differential
invariants.

Alternative classifications provide {\em closed, explicit expressions}
for homogeneous graphs $\{ u = F(x) \}$, but this then requires
discussions ({\em see} {\em e.g.}~{\cite{Merker-Nurowski-2020}}) about
ranges of incoming parameters in order to determine in which precise
(invariant) branches do sit the corresponding homogeneous models.
We will not touch these aspects here.

\medskip

Now, we come to the {\em non-cylindrical} surfaces, those having
$\Saux \neq 0$. From Theorem~{\ref{Thm-flat-cone}}, 
we already know that when $\Waux \equiv 0 \equiv \Xaux$, there
is, up to ${\sf SA}_3(\R)$, the single graph:
\[
u
\,=\,
\frac{1}{2}\,
\frac{x^2}{1-y},
\]
which is equivalent to a smooth part of the standard straight cone
$\{ x_1^2 + x_2^2 - x_3^2 = 0\}$.
Furthermore, one can verify that 
this model {\em is} special affinely
homogeneous, with its algebra of infinitesimal transformations
contained in $\mathfrak{sa}_3(\R)$ 
being generated by:
\[
\aligned
e_1
&
\,:=\,
-\,u\,\partial_x
+
x\,\partial_y,
\\
e_2
&
\,:=\,
(1-y)\,\partial_x
+
x\,\partial_u,
\\
e_3
&
\,:=\,
(1-y)\,\partial_y
+
u\,\partial_u,
\endaligned
\ \ \ \ \ \ \ \ \ \ \ \ \ \ \ \ \ \ \ \
\text{with Lie structure}
\ \ \ \ \ \ \ \ \ \ \ \ \ \ \ \ \ \ \ \
\aligned
{}
[e_1,e_2]
&
\,=\,
-e_3,
\\
[e_1,e_3]
&
\,=\,
-\,e_1,
\\
[e_2,e_3]
&
\,=\,
e_2,
\endaligned
\]
isomorphic to $\mathfrak{sl}_2(\R)$.

\begin{Observation}
Excepting the straight cone $\big\{ u = \frac{1}{2}\, 
\frac{x^2}{1-y} \big\}$, there are no non-cylindrical
special affinely homogeneous parabolic surfaces
$S^2 \subset \R^3$.
\end{Observation}

\proof
Coming back to the complete branching diagram copied above, 
it remains to examine the two branches
$\Waux \equiv 0 \neq \Xaux$ and $\Waux \neq 0$, 
with no need of going to subbranches.

\smallskip\noindent$\bullet$\,
When $\Waux \equiv 0 \neq \Xaux$, by homogeneity 
we have $\Xaux =: \aaux \in \R \backslash \{0\}$ constant,
but one of the recurrence relations shown at the end of
Subsection~{\ref{Subsection-W-zero}} immediately brings
a contradiction:
\[
0
\,=\,
\mathcal{D}_2\Xaux
\,=\,
3\,\Xaux.
\]

\smallskip\noindent$\bullet$\,
When $\Waux \neq 0$, 
we also have $\Waux =: \baux \in \R \backslash \{0\}$ constant,
but one of the recurrence relations shown at the end of
Subsection~{\ref{Subsection-W-nonzero}} again brings
an immediate contradiction:
\[
0
\,=\,
\mathcal{D}_2\Waux
\,=\,
2\,\Waux.
\qedhere
\]
\endproof


\vfill\end{document}

%% file: macros.tex
\newtheorem{Theorem}[equation]{Theorem}

\newtheorem{Proposition}[equation]{Proposition}

\newtheorem{Lemma}[equation]{Lemma}

\newtheorem{Assertion}[equation]{Assertion}

\newtheorem{Observation}[equation]{Observation}

\theoremstyle{definition}

\newtheorem{Hypothesis}[equation]{Hypothesis}
\newtheorem{Definition}[equation]{Definition}

\newtheorem{Terminology}[equation]{Terminology}
\newtheorem{Convention}[equation]{Convention}

\newtheorem{Example}[equation]{Example}

\newtheorem{Question}[equation]{Question}

\newtheorem{Problem}[equation]{Problem}

\newcommand{\C}{\mathbb{C}}

\newcommand{\N}{\mathbb{N}}

\newcommand{\R}{\mathbb{R}}

\renewcommand{\AA}{\text{\sc a}}
\newcommand{\BB}{\text{\sc b}}
\newcommand{\CC}{\text{\sc c}}

\newcommand{\EE}{\text{\sc e}}
\newcommand{\FF}{\text{\sc f}}
\newcommand{\GG}{\text{\sc g}}

\newcommand{\NN}{\text{\sc n}}

\newcommand{\aaux}{{\text{\usefont{T1}{qcs}{m}{sl}a}}}
\newcommand{\baux}{{\text{\usefont{T1}{qcs}{m}{sl}b}}}

\newcommand{\Caux}{{\text{\usefont{T1}{qcs}{m}{sl}C}}}
\newcommand{\Daux}{{\text{\usefont{T1}{qcs}{m}{sl}D}}}

\newcommand{\Faux}{{\text{\usefont{T1}{qcs}{m}{sl}F}}}

\newcommand{\Haux}{{\text{\usefont{T1}{qcs}{m}{sl}H}}}
\newcommand{\Iaux}{{\text{\usefont{T1}{qcs}{m}{sl}I}}}
\newcommand{\Jaux}{{\text{\usefont{T1}{qcs}{m}{sl}J}}}
\newcommand{\Kaux}{{\text{\usefont{T1}{qcs}{m}{sl}K}}}

\newcommand{\Maux}{{\text{\usefont{T1}{qcs}{m}{sl}M}}}

\newcommand{\Paux}{{\text{\usefont{T1}{qcs}{m}{sl}P}}}
\newcommand{\Qaux}{{\text{\usefont{T1}{qcs}{m}{sl}Q}}}
\newcommand{\Raux}{{\text{\usefont{T1}{qcs}{m}{sl}R}}}
\newcommand{\Saux}{{\text{\usefont{T1}{qcs}{m}{sl}S}}}

\newcommand{\Waux}{{\text{\usefont{T1}{qcs}{m}{sl}W}}}
\newcommand{\Xaux}{{\text{\usefont{T1}{qcs}{m}{sl}X}}}
\newcommand{\Yaux}{{\text{\usefont{T1}{qcs}{m}{sl}Y}}}
\newcommand{\Zaux}{{\text{\usefont{T1}{qcs}{m}{sl}Z}}}

\definecolor{blue}{cmyk}{1.,1.,0.,0.63}
\definecolor{red}{cmyk}{0.,1.,1.,0.63}
\definecolor{green}{cmyk}{1.,0.,1.,0.63}
\definecolor{black}{cmyk}{1.,1.,1.,1.}

\newcommand{\blue}{\textcolor{blue}}
\newcommand{\green}{\textcolor{green}}
\newcommand{\red}{\textcolor{red}}

\makeatletter
\renewcommand{\@fnsymbol}[1]
{\ensuremath{\ifcase#1\or $*$\or $**$\or $***$\or $****$\or $*****$
\else\@ctrerr\fi}}
\makeatother

\newcommand{\HEAD}[2]{%
\pagestyle{fancy}
\fancyhead[RO]{\tiny\sf\thepage}
\fancyhead[CO]{{\tiny\sf #1}}
\fancyhead[LE]{\tiny\sf\thepage}
\fancyhead[CE]{{\tiny\sf #2}}
\fancyfoot{}}

\numberwithin{equation}{section}

\newcommand{\Section}[1]{
\renewcommand{\thesection}{\bf\arabic{section}}
\section{#1}
\renewcommand{\thesection}{\arabic{section}}}

\newcommand{\Subsection}[1]{
\refstepcounter{equation}
\medskip\noindent{\bf\arabic{section}.\arabic{equation}.~#1.}}

\newcommand{\style}[1]{\text{\footnotesize{\sf #1}}}

\newcommand{\stylesmall}[1]{{\sf #1}}

\newcommand{\Card}{\style{Card}}

\renewcommand{\cos}{\style{cos}}

\renewcommand{\det}{\style{det}}

\renewcommand{\dim}{\style{dim}}

\renewcommand{\exp}{\style{exp}}

\newcommand{\formula}{\style{formula}}

\newcommand{\GL}{\style{GL}}

\newcommand{\Hessian}{\style{Hessian}}

\newcommand{\inv}{\style{inv}}

\newcommand{\Iso}{\style{Iso}}

\newcommand{\Jac}{\style{Jac}}

\newcommand{\Lie}{\style{Lie}}

\renewcommand{\lim}{\style{lim}}

\renewcommand{\max}{\style{max}}

\newcommand{\nonzero}{\style{nonzero}}

\newcommand{\rank}{\style{rank}}

\renewcommand{\root}{\style{root}}

\newcommand{\SA}{\style{SA}}

\newcommand{\SE}{\style{SE}}

\renewcommand{\sin}{\style{sin}}

\newcommand{\SL}{\style{SL}}

\newcommand{\SO}{\style{SO}}

\newcommand{\Span}{\style{Span}}

\newcommand{\stabsmall}{\stylesmall{stab}}

\newcommand{\medcup}{\mathbin{\scalebox{1.5}{\ensuremath{\cup}}}}

\newcommand{\oorder}{{\scriptstyle{\mathcal{O}}}}

\newcommand{\smallsum}[1]{
\underset{#1}{\raisebox{1pt}{$\sum$\,}}
}

\newcommand{\vf}{\vfill\end{document}}

\newcommand{\zero}[1]{\underline{#1}_{\red{\circ}}}


%% file: n-order-parabolic-jets.pdf_t
\begin{picture}(0,0)%
\includegraphics{n-order-parabolic-jets.pdf}%
\end{picture}%
\setlength{\unitlength}{4144sp}%
\begingroup\makeatletter\ifx\SetFigFont\undefined%
\gdef\SetFigFont#1#2#3#4#5{%
  \reset@font\fontsize{#1}{#2pt}%
  \fontfamily{#3}\fontseries{#4}\fontshape{#5}%
  \selectfont}%
\fi\endgroup%
\begin{picture}(2656,1605)(742,-1610)
\put(1001,-1538){\makebox(0,0)[lb]{\smash{{\SetFigFont{6}{7.2}{\familydefault}{\mddefault}{\updefault}{\color[rgb]{0,0,0}$(1,\!0)$}%
}}}}
\put(1238,-1540){\makebox(0,0)[lb]{\smash{{\SetFigFont{6}{7.2}{\familydefault}{\mddefault}{\updefault}{\color[rgb]{0,0,0}$(2,\!0)$}%
}}}}
\put(769,-1537){\makebox(0,0)[lb]{\smash{{\SetFigFont{6}{7.2}{\familydefault}{\mddefault}{\updefault}{\color[rgb]{0,0,0}$(0,\!0)$}%
}}}}
\put(757,-1311){\makebox(0,0)[lb]{\smash{{\SetFigFont{6}{7.2}{\familydefault}{\mddefault}{\updefault}{\color[rgb]{0,0,0}$(0,\!1)$}%
}}}}
\put(991,-1313){\makebox(0,0)[lb]{\smash{{\SetFigFont{6}{7.2}{\familydefault}{\mddefault}{\updefault}{\color[rgb]{0,0,0}$(1,\!1)$}%
}}}}
\put(2566,-1321){\makebox(0,0)[lb]{\smash{{\SetFigFont{6}{7.2}{\familydefault}{\mddefault}{\updefault}{\color[rgb]{0,0,0}$(n\!-\!1,\!1)$}%
}}}}
\put(2791,-1546){\makebox(0,0)[lb]{\smash{{\SetFigFont{6}{7.2}{\familydefault}{\mddefault}{\updefault}{\color[rgb]{0,0,0}$(n,\!0)$}%
}}}}
\put(1846,-196){\makebox(0,0)[lb]{\smash{{\SetFigFont{9}{10.8}{\familydefault}{\mddefault}{\updefault}{\color[rgb]{0,0,0}\green{$j\!+\!k\!=\!n$}}%
}}}}
\put(936,-152){\makebox(0,0)[lb]{\smash{{\SetFigFont{8}{9.6}{\familydefault}{\mddefault}{\updefault}{\color[rgb]{0,0,0}\blue{$k$}}%
}}}}
\put(3276,-1358){\makebox(0,0)[lb]{\smash{{\SetFigFont{8}{9.6}{\familydefault}{\mddefault}{\updefault}{\color[rgb]{0,0,0}\blue{$j$}}%
}}}}
\end{picture}%

%% file: straightening-G-orbits.pdf_t
\begin{picture}(0,0)%
\includegraphics{straightening-G-orbits.pdf}%
\end{picture}%
\setlength{\unitlength}{4144sp}%
\begingroup\makeatletter\ifx\SetFigFont\undefined%
\gdef\SetFigFont#1#2#3#4#5{%
  \reset@font\fontsize{#1}{#2pt}%
  \fontfamily{#3}\fontseries{#4}\fontshape{#5}%
  \selectfont}%
\fi\endgroup%
\begin{picture}(5894,2071)(429,-1883)
\put(4982,-915){\makebox(0,0)[lb]{\smash{{\SetFigFont{9}{10.8}{\familydefault}{\mddefault}{\updefault}{\color[rgb]{0,0,0}\red{$0$}}%
}}}}
\put(4995, 78){\makebox(0,0)[lb]{\smash{{\SetFigFont{9}{10.8}{\familydefault}{\mddefault}{\updefault}{\color[rgb]{0,0,0}\red{$\R^{\NN-s}$}}%
}}}}
\put(6199,-923){\makebox(0,0)[lb]{\smash{{\SetFigFont{9}{10.8}{\familydefault}{\mddefault}{\updefault}{\color[rgb]{0,0,0}$\R^s$}%
}}}}
\put(1065, 22){\makebox(0,0)[lb]{\smash{{\SetFigFont{9}{10.8}{\familydefault}{\mddefault}{\updefault}{\color[rgb]{0,0,0}\red{$T^n$}}%
}}}}
\put(2071,-477){\makebox(0,0)[lb]{\smash{{\SetFigFont{9}{10.8}{\familydefault}{\mddefault}{\updefault}{\color[rgb]{0,0,0}$z^{(n)}$}%
}}}}
\put(3174,-735){\makebox(0,0)[lb]{\smash{{\SetFigFont{9}{10.8}{\familydefault}{\mddefault}{\updefault}{\color[rgb]{0,0,0}\blue{${\sf straightening}$}}%
}}}}
\put(1058,-835){\makebox(0,0)[lb]{\smash{{\SetFigFont{9}{10.8}{\familydefault}{\mddefault}{\updefault}{\color[rgb]{0,0,0}\red{$z_0^{(n)}$}}%
}}}}
\end{picture}%

%% file: T-cross-section.pdf_t
\begin{picture}(0,0)%
\includegraphics{T-cross-section.pdf}%
\end{picture}%
\setlength{\unitlength}{4144sp}%
\begingroup\makeatletter\ifx\SetFigFont\undefined%
\gdef\SetFigFont#1#2#3#4#5{%
  \reset@font\fontsize{#1}{#2pt}%
  \fontfamily{#3}\fontseries{#4}\fontshape{#5}%
  \selectfont}%
\fi\endgroup%
\begin{picture}(4374,2828)(1104,-2693)
\put(4244,-650){\makebox(0,0)[lb]{\smash{{\SetFigFont{9}{10.8}{\familydefault}{\mddefault}{\updefault}{\color[rgb]{0,0,0}$z^{(n)}$}%
}}}}
\put(3248,-2087){\makebox(0,0)[lb]{\smash{{\SetFigFont{8}{9.6}{\familydefault}{\mddefault}{\updefault}{\color[rgb]{0,0,0}\green{$\rho(z^{(n)})\in G$}}%
}}}}
\put(2403, 16){\makebox(0,0)[lb]{\smash{{\SetFigFont{9}{10.8}{\familydefault}{\mddefault}{\updefault}{\color[rgb]{0,0,0}\red{$T^n$}}%
}}}}
\put(2061,-1511){\makebox(0,0)[lb]{\smash{{\SetFigFont{7}{8.4}{\familydefault}{\mddefault}{\updefault}{\color[rgb]{0,0,0}\red{$z_0^{(n)}$}}%
}}}}
\put(2334,-1088){\makebox(0,0)[lb]{\smash{{\SetFigFont{9}{10.8}{\familydefault}{\mddefault}{\updefault}{\color[rgb]{0,0,0}\red{$\rho(z^{(n)})\cdot z^{(n)}$}}%
}}}}
\end{picture}%

%% file: T-target.pdf_t
\begin{picture}(0,0)%
\includegraphics{T-target.pdf}%
\end{picture}%
\setlength{\unitlength}{4144sp}%
\begingroup\makeatletter\ifx\SetFigFont\undefined%
\gdef\SetFigFont#1#2#3#4#5{%
  \reset@font\fontsize{#1}{#2pt}%
  \fontfamily{#3}\fontseries{#4}\fontshape{#5}%
  \selectfont}%
\fi\endgroup%
\begin{picture}(5938,2831)(900,-2684)
\put(6682,-1507){\makebox(0,0)[lb]{\smash{{\SetFigFont{9}{10.8}{\familydefault}{\mddefault}{\updefault}{\color[rgb]{0,0,0}\blue{$w_{\nu_1}^{(n_1)}$}}%
}}}}
\put(2095,-1403){\makebox(0,0)[lb]{\smash{{\SetFigFont{9}{10.8}{\familydefault}{\mddefault}{\updefault}{\color[rgb]{0,0,0}\red{$z_0^{(n)}$}}%
}}}}
\put(2345,-276){\makebox(0,0)[lb]{\smash{{\SetFigFont{9}{10.8}{\familydefault}{\mddefault}{\updefault}{\color[rgb]{0,0,0}\red{$T^n$}}%
}}}}
\put(3707,-856){\makebox(0,0)[lb]{\smash{{\SetFigFont{9}{10.8}{\familydefault}{\mddefault}{\updefault}{\color[rgb]{0,0,0}\blue{$w^{(n)}\big(g,z^{(n)}\big)$}}%
}}}}
\put(5140, 15){\makebox(0,0)[lb]{\smash{{\SetFigFont{9}{10.8}{\familydefault}{\mddefault}{\updefault}{\color[rgb]{0,0,0}\blue{$w_\nu^{(n)}$}}%
}}}}
\put(6069,-22){\makebox(0,0)[lb]{\smash{{\SetFigFont{9}{10.8}{\familydefault}{\mddefault}{\updefault}{\color[rgb]{0,0,0}\red{$T_{\sf target}^n$}}%
}}}}
\put(5682,-619){\makebox(0,0)[lb]{\smash{{\SetFigFont{9}{10.8}{\familydefault}{\mddefault}{\updefault}{\color[rgb]{0,0,0}\blue{$w_{\nu_r}^{(n_r)}$}}%
}}}}
\put(5682,-1653){\makebox(0,0)[lb]{\smash{{\SetFigFont{9}{10.8}{\familydefault}{\mddefault}{\updefault}{\color[rgb]{0,0,0}\red{$c_1$}}%
}}}}
\put(5282,-1196){\makebox(0,0)[lb]{\smash{{\SetFigFont{9}{10.8}{\familydefault}{\mddefault}{\updefault}{\color[rgb]{0,0,0}\red{$c_r$}}%
}}}}
\end{picture}%

%% file: g-differential-invariant.pdf_t
\begin{picture}(0,0)%
\includegraphics{g-differential-invariant.pdf}%
\end{picture}%
\setlength{\unitlength}{4144sp}%
\begingroup\makeatletter\ifx\SetFigFont\undefined%
\gdef\SetFigFont#1#2#3#4#5{%
  \reset@font\fontsize{#1}{#2pt}%
  \fontfamily{#3}\fontseries{#4}\fontshape{#5}%
  \selectfont}%
\fi\endgroup%
\begin{picture}(6137,2247)(654,-2736)
\put(4016,-2667){\makebox(0,0)[lb]{\smash{{\SetFigFont{9}{10.8}{\familydefault}{\mddefault}{\updefault}{\color[rgb]{0,0,0}\blue{${\tt y}(g,x)$}}%
}}}}
\put(4061,-1905){\makebox(0,0)[lb]{\smash{{\SetFigFont{9}{10.8}{\familydefault}{\mddefault}{\updefault}{\color[rgb]{0,0,0}\blue{$\big({\tt y}(g,x),{\tt v}\big(g,{\tt y}(g,x)\big)\big)$}}%
}}}}
\put(1122,-2642){\makebox(0,0)[lb]{\smash{{\SetFigFont{9}{10.8}{\familydefault}{\mddefault}{\updefault}{\color[rgb]{0,0,0}\blue{$x$}}%
}}}}
\put(2165,-1352){\makebox(0,0)[lb]{\smash{{\SetFigFont{9}{10.8}{\familydefault}{\mddefault}{\updefault}{\color[rgb]{0,0,0}\blue{$J_{x,u}^n$}}%
}}}}
\put(2932,-914){\makebox(0,0)[lb]{\smash{{\SetFigFont{9}{10.8}{\familydefault}{\mddefault}{\updefault}{\color[rgb]{0,0,0}$g\cdot$}%
}}}}
\put(2931,-2130){\makebox(0,0)[lb]{\smash{{\SetFigFont{9}{10.8}{\familydefault}{\mddefault}{\updefault}{\color[rgb]{0,0,0}$g\cdot$}%
}}}}
\put(1082,-2222){\makebox(0,0)[lb]{\smash{{\SetFigFont{9}{10.8}{\familydefault}{\mddefault}{\updefault}{\color[rgb]{0,0,0}\blue{$\big(x,{\tt u}(x)\big)$}}%
}}}}
\put(6776,-2539){\makebox(0,0)[lb]{\smash{{\SetFigFont{9}{10.8}{\familydefault}{\mddefault}{\updefault}{\color[rgb]{0,0,0}\blue{$\R_y^p$}}%
}}}}
\put(2503,-2545){\makebox(0,0)[lb]{\smash{{\SetFigFont{9}{10.8}{\familydefault}{\mddefault}{\updefault}{\color[rgb]{0,0,0}\blue{$\R_x^p$}}%
}}}}
\put(732,-873){\makebox(0,0)[lb]{\smash{{\SetFigFont{9}{10.8}{\familydefault}{\mddefault}{\updefault}{\color[rgb]{0,0,0}\blue{$\big(x,{\tt u}^\alpha(x),{\tt u}_{x^J}^\beta(x)\big)$}}%
}}}}
\put(3665,-881){\makebox(0,0)[lb]{\smash{{\SetFigFont{9}{10.8}{\familydefault}{\mddefault}{\updefault}{\color[rgb]{0,0,0}\blue{$\big({\tt y}(g,x),{\tt v}^\gamma\big(g,{\tt y}(g,x)\big), {\tt v}_{y^K}^\delta\big(g,{\tt y}(g,x)\big)\big)$}}%
}}}}
\put(6448,-1350){\makebox(0,0)[lb]{\smash{{\SetFigFont{9}{10.8}{\familydefault}{\mddefault}{\updefault}{\color[rgb]{0,0,0}\blue{$J_{y,v}^n$}}%
}}}}
\end{picture}%

%% file: O-order-parabolic-jets.pdf_t
\begin{picture}(0,0)%
\includegraphics{O-order-parabolic-jets.pdf}%
\end{picture}%
\setlength{\unitlength}{4144sp}%
\begingroup\makeatletter\ifx\SetFigFont\undefined%
\gdef\SetFigFont#1#2#3#4#5{%
  \reset@font\fontsize{#1}{#2pt}%
  \fontfamily{#3}\fontseries{#4}\fontshape{#5}%
  \selectfont}%
\fi\endgroup%
\begin{picture}(3781,1599)(742,-1604)
\put(1001,-1538){\makebox(0,0)[lb]{\smash{{\SetFigFont{6}{7.2}{\familydefault}{\mddefault}{\updefault}{\color[rgb]{0,0,0}$(1,\!0)$}%
}}}}
\put(1238,-1540){\makebox(0,0)[lb]{\smash{{\SetFigFont{6}{7.2}{\familydefault}{\mddefault}{\updefault}{\color[rgb]{0,0,0}$(2,\!0)$}%
}}}}
\put(769,-1537){\makebox(0,0)[lb]{\smash{{\SetFigFont{6}{7.2}{\familydefault}{\mddefault}{\updefault}{\color[rgb]{0,0,0}$(0,\!0)$}%
}}}}
\put(757,-1311){\makebox(0,0)[lb]{\smash{{\SetFigFont{6}{7.2}{\familydefault}{\mddefault}{\updefault}{\color[rgb]{0,0,0}$(0,\!1)$}%
}}}}
\put(991,-1313){\makebox(0,0)[lb]{\smash{{\SetFigFont{6}{7.2}{\familydefault}{\mddefault}{\updefault}{\color[rgb]{0,0,0}$(1,\!1)$}%
}}}}
\put(3436,-1318){\makebox(0,0)[lb]{\smash{{\SetFigFont{6}{7.2}{\familydefault}{\mddefault}{\updefault}{\color[rgb]{0,0,0}$(\oorder\!-\!1,\!1)$}%
}}}}
\put(3710,-1539){\makebox(0,0)[lb]{\smash{{\SetFigFont{6}{7.2}{\familydefault}{\mddefault}{\updefault}{\color[rgb]{0,0,0}$(\oorder,\!0)$}%
}}}}
\put(2711,-186){\makebox(0,0)[lb]{\smash{{\SetFigFont{9}{10.8}{\familydefault}{\mddefault}{\updefault}{\color[rgb]{0,0,0}\green{$k\!+\!l\!=\!\oorder$}}%
}}}}
\put(4436,-1379){\makebox(0,0)[lb]{\smash{{\SetFigFont{8}{9.6}{\familydefault}{\mddefault}{\updefault}{\color[rgb]{0,0,0}\blue{$k$}}%
}}}}
\put(936,-152){\makebox(0,0)[lb]{\smash{{\SetFigFont{8}{9.6}{\familydefault}{\mddefault}{\updefault}{\color[rgb]{0,0,0}\blue{$l$}}%
}}}}
\end{picture}%

%% file: u-F-x-v-G-y.pdf_t
\begin{picture}(0,0)%
\includegraphics{u-F-x-v-G-y.pdf}%
\end{picture}%
\setlength{\unitlength}{4144sp}%
\begingroup\makeatletter\ifx\SetFigFont\undefined%
\gdef\SetFigFont#1#2#3#4#5{%
  \reset@font\fontsize{#1}{#2pt}%
  \fontfamily{#3}\fontseries{#4}\fontshape{#5}%
  \selectfont}%
\fi\endgroup%
\begin{picture}(6569,1888)(519,-1463)
\put(937,312){\makebox(0,0)[lb]{\smash{{\SetFigFont{9}{10.8}{\familydefault}{\mddefault}{\updefault}{\color[rgb]{0,0,0}$u$}%
}}}}
\put(602,306){\makebox(0,0)[lb]{\smash{{\SetFigFont{9}{10.8}{\familydefault}{\mddefault}{\updefault}{\color[rgb]{0,0,0}$\R^2$}%
}}}}
\put(3053,-1371){\makebox(0,0)[lb]{\smash{{\SetFigFont{9}{10.8}{\familydefault}{\mddefault}{\updefault}{\color[rgb]{0,0,0}$x$}%
}}}}
\put(4494,307){\makebox(0,0)[lb]{\smash{{\SetFigFont{9}{10.8}{\familydefault}{\mddefault}{\updefault}{\color[rgb]{0,0,0}$\R^2$}%
}}}}
\put(4812,331){\makebox(0,0)[lb]{\smash{{\SetFigFont{9}{10.8}{\familydefault}{\mddefault}{\updefault}{\color[rgb]{0,0,0}$v$}%
}}}}
\put(6940,-1355){\makebox(0,0)[lb]{\smash{{\SetFigFont{9}{10.8}{\familydefault}{\mddefault}{\updefault}{\color[rgb]{0,0,0}$y$}%
}}}}
\put(5864,-513){\makebox(0,0)[lb]{\smash{{\SetFigFont{9}{10.8}{\familydefault}{\mddefault}{\updefault}{\color[rgb]{0,0,0}\blue{$\big\{v=G(y)\big\}$}}%
}}}}
\put(1971,-500){\makebox(0,0)[lb]{\smash{{\SetFigFont{9}{10.8}{\familydefault}{\mddefault}{\updefault}{\color[rgb]{0,0,0}\blue{$\big\{u=F(x)\big\}$}}%
}}}}
\put(3396,-415){\makebox(0,0)[lb]{\smash{{\SetFigFont{9}{10.8}{\familydefault}{\mddefault}{\updefault}{\color[rgb]{0,0,0}${\sf special}\,\,{\sf affine}$}%
}}}}
\put(3626,-850){\makebox(0,0)[lb]{\smash{{\SetFigFont{9}{10.8}{\familydefault}{\mddefault}{\updefault}{\color[rgb]{0,0,0}\red{${\sf inverse}$}}%
}}}}
\end{picture}%

%% file: 4-proj-parabolic.pdf_t
\begin{picture}(0,0)%
\includegraphics{4-proj-parabolic.pdf}%
\end{picture}%
\setlength{\unitlength}{4144sp}%
\begingroup\makeatletter\ifx\SetFigFont\undefined%
\gdef\SetFigFont#1#2#3#4#5{%
  \reset@font\fontsize{#1}{#2pt}%
  \fontfamily{#3}\fontseries{#4}\fontshape{#5}%
  \selectfont}%
\fi\endgroup%
\begin{picture}(1801,1369)(742,-1605)
\put(757,-1311){\makebox(0,0)[lb]{\smash{{\SetFigFont{6}{7.2}{\familydefault}{\mddefault}{\updefault}{\color[rgb]{0,0,0}$(0,\!1)$}%
}}}}
\put(991,-1313){\makebox(0,0)[lb]{\smash{{\SetFigFont{6}{7.2}{\familydefault}{\mddefault}{\updefault}{\color[rgb]{0,0,0}$(1,\!1)$}%
}}}}
\put(1230,-1315){\makebox(0,0)[lb]{\smash{{\SetFigFont{6}{7.2}{\familydefault}{\mddefault}{\updefault}{\color[rgb]{0,0,0}$(2,\!1)$}%
}}}}
\put(1467,-1314){\makebox(0,0)[lb]{\smash{{\SetFigFont{6}{7.2}{\familydefault}{\mddefault}{\updefault}{\color[rgb]{0,0,0}$(3,\!1)$}%
}}}}
\put(1001,-1538){\makebox(0,0)[lb]{\smash{{\SetFigFont{6}{7.2}{\familydefault}{\mddefault}{\updefault}{\color[rgb]{0,0,0}$(1,\!0)$}%
}}}}
\put(1238,-1540){\makebox(0,0)[lb]{\smash{{\SetFigFont{6}{7.2}{\familydefault}{\mddefault}{\updefault}{\color[rgb]{0,0,0}$(2,\!0)$}%
}}}}
\put(769,-1537){\makebox(0,0)[lb]{\smash{{\SetFigFont{6}{7.2}{\familydefault}{\mddefault}{\updefault}{\color[rgb]{0,0,0}$(0,\!0)$}%
}}}}
\put(1466,-1541){\makebox(0,0)[lb]{\smash{{\SetFigFont{6}{7.2}{\familydefault}{\mddefault}{\updefault}{\color[rgb]{0,0,0}$(3,\!0)$}%
}}}}
\put(1698,-1541){\makebox(0,0)[lb]{\smash{{\SetFigFont{6}{7.2}{\familydefault}{\mddefault}{\updefault}{\color[rgb]{0,0,0}$(4,\!0)$}%
}}}}
\put(1143,-637){\makebox(0,0)[lb]{\smash{{\SetFigFont{6}{7.2}{\familydefault}{\mddefault}{\updefault}{\color[rgb]{0,0,0}\green{$k\!+\!l\!=\! 4$}}%
}}}}
\put(2386,-1366){\makebox(0,0)[lb]{\smash{{\SetFigFont{8}{9.6}{\familydefault}{\mddefault}{\updefault}{\color[rgb]{0,0,0}\blue{$k$}}%
}}}}
\put(931,-383){\makebox(0,0)[lb]{\smash{{\SetFigFont{8}{9.6}{\familydefault}{\mddefault}{\updefault}{\color[rgb]{0,0,0}\blue{$l$}}%
}}}}
\end{picture}%

%% file: 3-1-S-parabolic.pdf_t
\begin{picture}(0,0)%
\includegraphics{3-1-S-parabolic.pdf}%
\end{picture}%
\setlength{\unitlength}{4144sp}%
\begingroup\makeatletter\ifx\SetFigFont\undefined%
\gdef\SetFigFont#1#2#3#4#5{%
  \reset@font\fontsize{#1}{#2pt}%
  \fontfamily{#3}\fontseries{#4}\fontshape{#5}%
  \selectfont}%
\fi\endgroup%
\begin{picture}(1801,1369)(742,-1605)
\put(991,-1313){\makebox(0,0)[lb]{\smash{{\SetFigFont{6}{7.2}{\familydefault}{\mddefault}{\updefault}{\color[rgb]{0,0,0}$(1,\!1)$}%
}}}}
\put(1001,-1538){\makebox(0,0)[lb]{\smash{{\SetFigFont{6}{7.2}{\familydefault}{\mddefault}{\updefault}{\color[rgb]{0,0,0}$(1,\!0)$}%
}}}}
\put(1238,-1540){\makebox(0,0)[lb]{\smash{{\SetFigFont{6}{7.2}{\familydefault}{\mddefault}{\updefault}{\color[rgb]{0,0,0}$(2,\!0)$}%
}}}}
\put(769,-1537){\makebox(0,0)[lb]{\smash{{\SetFigFont{6}{7.2}{\familydefault}{\mddefault}{\updefault}{\color[rgb]{0,0,0}$(0,\!0)$}%
}}}}
\put(1466,-1541){\makebox(0,0)[lb]{\smash{{\SetFigFont{6}{7.2}{\familydefault}{\mddefault}{\updefault}{\color[rgb]{0,0,0}$(3,\!0)$}%
}}}}
\put(1698,-1541){\makebox(0,0)[lb]{\smash{{\SetFigFont{6}{7.2}{\familydefault}{\mddefault}{\updefault}{\color[rgb]{0,0,0}$(4,\!0)$}%
}}}}
\put(757,-1311){\makebox(0,0)[lb]{\smash{{\SetFigFont{6}{7.2}{\familydefault}{\mddefault}{\updefault}{\color[rgb]{0,0,0}$(0,\!1)$}%
}}}}
\put(2386,-1366){\makebox(0,0)[lb]{\smash{{\SetFigFont{8}{9.6}{\familydefault}{\mddefault}{\updefault}{\color[rgb]{0,0,0}\blue{$k$}}%
}}}}
\put(931,-383){\makebox(0,0)[lb]{\smash{{\SetFigFont{8}{9.6}{\familydefault}{\mddefault}{\updefault}{\color[rgb]{0,0,0}\blue{$l$}}%
}}}}
\end{picture}%

%% file: 4-1-W-parabolic.pdf_t
\begin{picture}(0,0)%
\includegraphics{4-1-W-parabolic.pdf}%
\end{picture}%
\setlength{\unitlength}{4144sp}%
\begingroup\makeatletter\ifx\SetFigFont\undefined%
\gdef\SetFigFont#1#2#3#4#5{%
  \reset@font\fontsize{#1}{#2pt}%
  \fontfamily{#3}\fontseries{#4}\fontshape{#5}%
  \selectfont}%
\fi\endgroup%
\begin{picture}(1801,1369)(742,-1605)
\put(1001,-1538){\makebox(0,0)[lb]{\smash{{\SetFigFont{6}{7.2}{\familydefault}{\mddefault}{\updefault}{\color[rgb]{0,0,0}$(1,\!0)$}%
}}}}
\put(1238,-1540){\makebox(0,0)[lb]{\smash{{\SetFigFont{6}{7.2}{\familydefault}{\mddefault}{\updefault}{\color[rgb]{0,0,0}$(2,\!0)$}%
}}}}
\put(769,-1537){\makebox(0,0)[lb]{\smash{{\SetFigFont{6}{7.2}{\familydefault}{\mddefault}{\updefault}{\color[rgb]{0,0,0}$(0,\!0)$}%
}}}}
\put(1466,-1541){\makebox(0,0)[lb]{\smash{{\SetFigFont{6}{7.2}{\familydefault}{\mddefault}{\updefault}{\color[rgb]{0,0,0}$(3,\!0)$}%
}}}}
\put(1698,-1541){\makebox(0,0)[lb]{\smash{{\SetFigFont{6}{7.2}{\familydefault}{\mddefault}{\updefault}{\color[rgb]{0,0,0}$(4,\!0)$}%
}}}}
\put(991,-1313){\makebox(0,0)[lb]{\smash{{\SetFigFont{6}{7.2}{\familydefault}{\mddefault}{\updefault}{\color[rgb]{0,0,0}$(1,\!1)$}%
}}}}
\put(757,-1311){\makebox(0,0)[lb]{\smash{{\SetFigFont{6}{7.2}{\familydefault}{\mddefault}{\updefault}{\color[rgb]{0,0,0}$(0,\!1)$}%
}}}}
\put(1218,-1314){\makebox(0,0)[lb]{\smash{{\SetFigFont{6}{7.2}{\familydefault}{\mddefault}{\updefault}{\color[rgb]{0,0,0}$(2,\!1)$}%
}}}}
\put(2386,-1366){\makebox(0,0)[lb]{\smash{{\SetFigFont{8}{9.6}{\familydefault}{\mddefault}{\updefault}{\color[rgb]{0,0,0}\blue{$k$}}%
}}}}
\put(931,-383){\makebox(0,0)[lb]{\smash{{\SetFigFont{8}{9.6}{\familydefault}{\mddefault}{\updefault}{\color[rgb]{0,0,0}\blue{$l$}}%
}}}}
\end{picture}%

%% file: destroy.pdf_t
\begin{picture}(0,0)%
\includegraphics{destroy.pdf}%
\end{picture}%
\setlength{\unitlength}{4144sp}%
\begingroup\makeatletter\ifx\SetFigFont\undefined%
\gdef\SetFigFont#1#2#3#4#5{%
  \reset@font\fontsize{#1}{#2pt}%
  \fontfamily{#3}\fontseries{#4}\fontshape{#5}%
  \selectfont}%
\fi\endgroup%
\begin{picture}(4844,1468)(669,-1468)
\put(5404,-1379){\makebox(0,0)[lb]{\smash{{\SetFigFont{8}{9.6}{\familydefault}{\mddefault}{\updefault}{\color[rgb]{0,0,0}$x$}%
}}}}
\put(756,-147){\makebox(0,0)[lb]{\smash{{\SetFigFont{8}{9.6}{\familydefault}{\mddefault}{\updefault}{\color[rgb]{0,0,0}$y$}%
}}}}
\put(4176,-147){\makebox(0,0)[lb]{\smash{{\SetFigFont{8}{9.6}{\familydefault}{\mddefault}{\updefault}{\color[rgb]{0,0,0}$y$}%
}}}}
\end{picture}%

%% file: curves-branches.pdf_t
\begin{picture}(0,0)%
\includegraphics{curves-branches.pdf}%
\end{picture}%
\setlength{\unitlength}{4144sp}%
\begingroup\makeatletter\ifx\SetFigFont\undefined%
\gdef\SetFigFont#1#2#3#4#5{%
  \reset@font\fontsize{#1}{#2pt}%
  \fontfamily{#3}\fontseries{#4}\fontshape{#5}%
  \selectfont}%
\fi\endgroup%
\begin{picture}(5850,2442)(932,-1829)
\put(995,-1267){\makebox(0,0)[lb]{\smash{{\SetFigFont{9}{10.8}{\familydefault}{\mddefault}{\updefault}{\color[rgb]{0,0,0}$I_2$}%
}}}}
\put(1395,-992){\makebox(0,0)[lb]{\smash{{\SetFigFont{9}{10.8}{\familydefault}{\mddefault}{\updefault}{\color[rgb]{0,0,0}$\neq0$}%
}}}}
\put(1344,-1587){\makebox(0,0)[lb]{\smash{{\SetFigFont{9}{10.8}{\familydefault}{\mddefault}{\updefault}{\color[rgb]{0,0,0}$\equiv0$}%
}}}}
\put(2282,-1774){\makebox(0,0)[lb]{\smash{{\SetFigFont{9}{10.8}{\familydefault}{\mddefault}{\updefault}{\color[rgb]{0,0,0}$u=0$}%
}}}}
\put(2311,-762){\makebox(0,0)[lb]{\smash{{\SetFigFont{9}{10.8}{\familydefault}{\mddefault}{\updefault}{\color[rgb]{0,0,0}$I_4$}%
}}}}
\put(3614,-1256){\makebox(0,0)[lb]{\smash{{\SetFigFont{9}{10.8}{\familydefault}{\mddefault}{\updefault}{\color[rgb]{0,0,0}$u=x^2$}%
}}}}
\put(2817,-448){\makebox(0,0)[lb]{\smash{{\SetFigFont{9}{10.8}{\familydefault}{\mddefault}{\updefault}{\color[rgb]{0,0,0}$\neq0$}%
}}}}
\put(2717,-1094){\makebox(0,0)[lb]{\smash{{\SetFigFont{9}{10.8}{\familydefault}{\mddefault}{\updefault}{\color[rgb]{0,0,0}$\equiv0$}%
}}}}
\put(3634,-245){\makebox(0,0)[lb]{\smash{{\SetFigFont{9}{10.8}{\familydefault}{\mddefault}{\updefault}{\color[rgb]{0,0,0}$I_5$}%
}}}}
\put(4132, 61){\makebox(0,0)[lb]{\smash{{\SetFigFont{9}{10.8}{\familydefault}{\mddefault}{\updefault}{\color[rgb]{0,0,0}$\neq0$}%
}}}}
\put(4051,-577){\makebox(0,0)[lb]{\smash{{\SetFigFont{9}{10.8}{\familydefault}{\mddefault}{\updefault}{\color[rgb]{0,0,0}$\equiv0$}%
}}}}
\put(5005,439){\makebox(0,0)[lb]{\smash{{\SetFigFont{9}{10.8}{\familydefault}{\mddefault}{\updefault}{\color[rgb]{0,0,0}Two $1$-parameter families}%
}}}}
\put(5006,246){\makebox(0,0)[lb]{\smash{{\SetFigFont{9}{10.8}{\familydefault}{\mddefault}{\updefault}{\color[rgb]{0,0,0}of homogeneous models $\aaux\!\neq\!0$}%
}}}}
\put(5005, 44){\makebox(0,0)[lb]{\smash{{\SetFigFont{9}{10.8}{\familydefault}{\mddefault}{\updefault}{\color[rgb]{0,0,0}$u\!=\!\frac{x^2}{2}+\varepsilon\frac{x^4}{4!}+\aaux x^5+\cdots$}%
}}}}
\put(4953,-741){\makebox(0,0)[lb]{\smash{{\SetFigFont{9}{10.8}{\familydefault}{\mddefault}{\updefault}{\color[rgb]{0,0,0}$u=\varepsilon\sqrt{1+\varepsilon x^2}\!-\!\varepsilon$}%
}}}}
\put(947,466){\makebox(0,0)[lb]{\smash{{\SetFigFont{9}{10.8}{\familydefault}{\mddefault}{\updefault}{\color[rgb]{0,0,0}$\varepsilon\!=\,\pm 1$}%
}}}}
\end{picture}%